\documentclass[10pt,journal,compsoc]{IEEEtran}

\usepackage{times}
\usepackage{amsmath}
\usepackage{amssymb}
\usepackage{url}
\usepackage{amsthm} 
\usepackage{algorithm}
\usepackage{diagbox} 
\usepackage{bbm}
\usepackage{bm}
\usepackage{diagbox}
\usepackage{graphicx}
\usepackage{subcaption}
\usepackage{epstopdf}
\usepackage{booktabs}
\usepackage{pdfpages}

\usepackage{courier}
\usepackage{color}
\definecolor{colorpink}{RGB}{251,53,155}
\definecolor{colorblue}{RGB}{0,148,200}
\definecolor{colorgreen}{RGB}{0,150,0}

\def\beq{\begin{eqnarray}}
\def\eeq{\end{eqnarray}}
\def\noi{\noindent}
\def\nn{\nonumber}
\def\la{\langle}
\def\ra{\rangle}

\newtheorem{lemma}{Lemma}
\newtheorem{theorem}{Theorem}

\newcommand{\bbb}[1]{\boldsymbol{\mathbf{#1}}}

\usepackage{multirow}
\usepackage{rotating}

\newcommand{\tabsizeall}{\fontsize{9.25}{12.5}\selectfont}

\def\J{\mathcal{J}}
\def\Lag{\mathcal{L}}

\def\beq{\begin{eqnarray}}
\def\eeq{\end{eqnarray}}
\def\noi{\noindent}
\def\nn{\nonumber}
\def\la{\langle}
\def\ra{\rangle}

\def\ghs{\hspace{0.05cm}} 
\def\tabsize{\scriptsize}
\def\J{\mathcal{J}}

\def\fourfigwid{0.235\textwidth}

\def\imgheibisect{0.135\textheight}
\def\imgheisubgraph{0.135\textheight}
\def\imgheiclustering{0.135\textheight}
\def\imgheimrf{0.135\textheight}
\def\imgheiconv{0.135\textheight}

\def\imgwidseg{0.185\textwidth}
\def\imgheiseg{0.086\textheight}
\def\ghs{\hspace{0.05cm}} 
\def\ghst{\hspace{0.1cm}} 

\usepackage{alphalph}

\usepackage[font=small, justification=centering]{caption}
\DeclareCaptionFont{tiny}{\tiny}

\newenvironment{psmallmatrix}
  {\left(\begin{smallmatrix}}
  {\end{smallmatrix}\right)}
\ifCLASSOPTIONcompsoc
  \usepackage[nocompress]{cite}
\else
  \usepackage{cite}
\fi

\begin{document}

\title{Binary Optimization via Mathematical Programming with Equilibrium Constraints}

\author{Ganzhao~Yuan,~Bernard~Ghanem
\IEEEcompsocitemizethanks{\IEEEcompsocthanksitem Ganzhao~Yuan (yuanganzhao@gmail.com) is with Sun Yat-sen University (SYSU), China. Bernard Ghanem (bernard.ghanem@kaust.edu.sa) is with King Abdullah University of Science and Technology (KAUST), Saudi Arabia.\protect\\
}
\thanks{Manuscript received April 19, 2005; revised August 26, 2015.}}

\markboth{}
{Shell \MakeLowercase{\textit{et al.}}: Bare Demo of IEEEtran.cls for Computer Society Journals}

\IEEEtitleabstractindextext{%
\begin{abstract}
Binary optimization is a central problem in mathematical optimization and its applications are abundant. To solve this problem, we propose a new class of continuous optimization techniques which is based on Mathematical Programming with Equilibrium Constraints (MPECs). We first reformulate the binary program as an equivalent augmented biconvex optimization problem with a bilinear equality constraint, then we propose two penalization/regularization methods (exact penalty and alternating direction) to solve it. The resulting algorithms seek desirable solutions to the original problem via solving a sequence of linear programming convex relaxation subproblems. In addition, we prove that both the penalty function and augmented Lagrangian function, induced by adding the complementarity constraint to the objectives, are exact, i.e., they have the same local and global minima with those of the original binary program when the penalty parameter is over some threshold. The convergence of both algorithms can be guaranteed, since they essentially reduce to block coordinate descent in the literature. Finally, we demonstrate the effectiveness and versatility of our methods on several important problems, including graph bisection, constrained image segmentation, dense subgraph discovery, modularity clustering and Markov random fields. Extensive experiments show that our methods outperform existing popular techniques, such as iterative hard thresholding, linear programming relaxation and semidefinite programming relaxation.

\end{abstract}

}

\maketitle
\IEEEdisplaynontitleabstractindextext
\IEEEpeerreviewmaketitle

\renewcommand{\arraystretch}{1}

\begin{table*}[!ht]
\tabsizeall
\center
\caption{Existing continuous methods for binary optimization.}
\vspace{-5pt}
\begin{center}
\scalebox{1.15}{\begin{tabular}{|p{0.18cm}|p{5.9cm}|p{8.3cm}|}
\hline
\multicolumn{1}{|c|}{} & Method and Reference  & Description  \\
\hline
\multirow{3}{*}{\begin{sideways} {\tabsize Relaxed Approximation~~~} \end{sideways}}
&  { spectral relaxation \cite{cour2007solving,ShiM00}} & {$ \{-1,+1\}^n \approx \{ \bbb{x}~|~\|\bbb{x}\|_2^2= n\}$}   \\
\cline{2-3}
&  { linear programming relaxation \cite{komodakis2007approximate,hsieh2015pu}} & {$\{-1,+1\}^n \approx \{ \bbb{x}~|~-\bbb{1}\leq \bbb{x} \leq \bbb{1} \}$}   \\
\cline{2-3}
& \multirow{2}{*}{SDP relaxation \cite{BQP2015Wang,Wen2010,keuchel2003binary}} & {$ \{0,+1\}^n \approx \{ \bbb{x}~|~\bbb{X} \succeq \bbb{x}\bbb{x}^T,diag(\bbb{X})=\bbb{x}\}$}   \\
\cline{3-3}
&   & {$ \{-1,+1\}^n \approx \{ \bbb{x}~|~\bbb{X}\succeq\bbb{x}\bbb{x}^T,~diag(\bbb{X})=\bbb{1}\}$}   \\
\cline{2-3}
&  {  doubly positive relaxation \cite{HuangCG14,Wen2010}} & {$ \{0,+1\}^n \approx \{ \bbb{x}~|~\bbb{X}\succeq\bbb{x}\bbb{x}^T,diag(\bbb{X})=\bbb{x},~\bbb{x}\geq \bbb{0},~\bbb{X}\geq \bbb{0}\}$}   \\
\cline{2-3}
&  {  completely positive relaxation   \cite{burer2010optimizing,burer2009copositive}} & {$  \{0,+1\}^n \approx \{ \bbb{x}~|~\bbb{X}\succeq\bbb{x}\bbb{x}^T,diag(\bbb{X})=\bbb{x},~\bbb{x}\geq \bbb{0},~\bbb{X}$~is CP$\} $}   \\
\cline{2-3}
&  {  SOCP relaxation \cite{KumarKT09,GhaddarVA11}} &  {$ \{-1,+1\}^n \approx \{ \bbb{x}~|~\la \bbb{X} - \bbb{x}\bbb{x}^T, \bbb{LL}^T\ra \geq 0,~diag(\bbb{X})=\bbb{1}\},~\forall~\bbb{L}$}   \\
\hline
\multirow{3}{*}{\begin{sideways} {\tabsize Equivalent Optimization~~~~~~~~~~}\end{sideways}}
&  {  iterative hard thresholding \cite{Yuan013,BeckE13} } & {$\min_{\bbb{x}}~\|\bbb{x}-\bbb{x}'\|_2^2,~s.t.~\bbb{x} \in \{-1,+1\}^n$}   \\
\cline{2-3}
&  {  $\ell_0$ norm reformulation \cite{LuZ13,yuan2016l0mpec}} & {$ \{-1,+1\}^n \Leftrightarrow \{ \bbb{x}~|~\|\bbb{x} + \bbb{1}\|_0 + \|\bbb{x} - \bbb{1}\|_0 \leq  n  \}$}   \\
\cline{2-3}
&  {  piecewise separable reformulation \cite{zhang2007binary}} & {$ \{-1,+1\}^n \Leftrightarrow \{ \bbb{x}~| ~ (\bbb{1}+\bbb{x}) \odot (\bbb{1}-\bbb{x}) = \bbb{0}\}$}   \\
\cline{2-3}
&  {  $\ell_2$ box non-separable reformulation \cite{Raghavachari1969,MurrayN10}} & {$\{-1,+1\}^n \Leftrightarrow \{ \bbb{x}~|~-\bbb{1}\leq \bbb{x} \leq   \bbb{1},~\|\bbb{x}\|_2^2= n  \}$}   \\
\cline{2-3}
&  {  $\ell_p$ box non-separable reformulation \cite{WuG16a}} & {$\{-1,+1\}^n \Leftrightarrow \{ \bbb{x}~|~-\bbb{1}\leq \bbb{x} \leq \bbb{1},~\|\bbb{x}\|_p^p= n,~0<p<\infty  \}$}   \\
\cline{2-3}
&  {  $\ell_{\infty}$ box separable MPEC \hspace*{\fill}[This paper]}  & {$ \{-1,+1\}^n \Leftrightarrow \{ \bbb{x}~|~\bbb{x} \odot \bbb{v}=\bbb{1},~\|\bbb{x}\|_{\infty}\leq 1,~\|\bbb{v}\|_{\infty}\leq 1,~\forall \bbb{v} \}$}\\
\cline{2-3}
&  {  $\ell_2$ box separable MPEC \hspace*{\fill}[This paper]} & {$ \{-1,+1\}^n \Leftrightarrow \{ \bbb{x}~|~ \bbb{x} \odot \bbb{v}=1,~\|\bbb{x}\|_{\infty}\leq 1,~\|\bbb{v}\|_2^2\leq n,~\forall \bbb{v} \}$} \\
\cline{2-3}
&  {  $\ell_{\infty}$ box non-separable MPEC \hspace*{\fill}[This paper]} & {$ \{-1,+1\}^n \Leftrightarrow \{ \bbb{x}~|~\la \bbb{x},\bbb{v}\ra=n,~\|\bbb{x}\|_{\infty}\leq 1,~\|\bbb{v}\|_{\infty}\leq 1,~\forall \bbb{v} \}$}\\
\cline{2-3}
&  {  $\ell_2$ box non-separable MPEC \hspace*{\fill} [This paper]} & {$ \{-1,+1\}^n \Leftrightarrow \{ \bbb{x}~|~\la \bbb{x},\bbb{v}\ra=n,~\|\bbb{x}\|_{\infty}\leq 1,~\|\bbb{v}\|_2^2\leq n,~\forall \bbb{v} \}$} \\
\hline
\end{tabular}}
\end{center}
\label{tab:binary:method}
\vspace{-3pt}
\end{table*}

\section{Introduction}
In this paper, we mainly focus on the following binary optimization problem:
\beq \label{eq:main}
\min_{\bbb{x}}~f(\bbb{x}),~s.t.~\bbb{x}\in\{-1,1\}^{n},~\bbb{x}\in \Omega
\eeq
\noi where the objective function $f:\mathbb{R}^{n}\rightarrow\mathbb{R}$ is convex (but not necessarily smooth) on some convex set $\Omega$, and the non-convexity of (\ref{eq:main}) is only caused by the binary constraints. In addition, we assume $\{-1,1\}^{n} \cap \Omega \neq \emptyset$.

The optimization in (\ref{eq:main}) describes many applications of interest in both computer vision and machine learning, including graph bisection \cite{Goemans95,keuchel2003binary}, image (co-)segmentation \cite{ShiM00,keuchel2003binary,JoulinBP10}, Markov random fields \cite{Boykov2001}, permutation problem \cite{FogelJBd15}, graph matching \cite{cour2007balanced,yan2015discrete,toshev2007image}, binary matrix completion \cite{davenport20141,hsieh2015pu}, hashing coding \cite{LiuWKC11,WangLKC16}, image registration \cite{BQP2015Wang}, multimodal feature learning \cite{ShrivastavaRSCD15}, multi-target tracking \cite{ShiLXH13}, visual alignment \cite{shokrollahi2015unsupervised}, and social network analysis (e.g. subgraphs discovery \cite{Yuan013,Ames2015}, biclustering \cite{Ames2014}, planted k-disjoint-clique discover \cite{Ames2014b}, planted clique and biclique discovery \cite{Ames2011}, community discovery \cite{hejoint,ChanY11}), etc.

The binary optimization problem is difficult to solve, since it is NP-hard. One type of method to solve this problem is continuous in nature. The simple way is to relax the binary constraint with Linear Programming (LP) relaxation constraints $-\bbb{1} \leq \bbb{x} \leq \bbb{1}$ and round the entries of the resulting continuous solution to the nearest integer at the end. However, not only may this solution not be optimal, it may not even be feasible and violate some constraint. Another type of optimization focuses on the cutting-plane and branch-and-cut method. The cutting plane method solves the LP relaxation and then adds linear constraints that drive the solution towards integers. The branch-and-cut method partially develops a binary tree and iteratively cuts out the nodes having a lower bound that is worse than the current upper bound, while the lower bound can be found using convex relaxation, Lagrangian duality, or Lipschitz continuity. However, this class of method ends up solving all $2^n$ convex subproblems in the worst case. Our algorithm aligns with the first research direction. It solves a convex LP relaxation subproblem iteratively, but it provably terminates in polynomial iterations.

In non-convex optimization, good initialization is very important to the quality of the solution. Motivated by this, several papers design smart initialization strategies and establish optimality qualification of the solutions for non-convex problems. For example, the work of \cite{zhang2010analysis} considers a multi-stage convex optimization algorithm to refine the global solution by the initial convex method; the work of \cite{CandesLS15} starts with a careful initialization obtained by a spectral method and improves this estimate by gradient descent; the work of \cite{jain2013low} uses the top-$k$ singular vectors of the matrix as initialization and provides theoretical guarantees through biconvex alternating minimization. The proposed method also uses a similar initialization strategy, as it reduces to convex LP relaxation in the first iteration.

The contributions of this paper are three-fold. \textbf{(a)} We reformulate the binary program as an equivalent augmented optimization problem with a bilinear equality constraint via a variational characterization of the binary constraint. Then, we propose two penalization/regularization methods (exact penalty and alternating direction) to solve it. The resulting algorithms seek desirable solutions to the original binary program. \textbf{(b)} We prove that both the penalty function and augmented Lagrangian function, induced by adding the complementarity constraint to the objectives, are exact, i.e., the set of their globally optimal solutions coincide with that of (\ref{eq:main}) when the penalty parameter is over some threshold. Thus, the convergence of both algorithms can be guaranteed since they reduce to block coordinate descent in the literature \cite{tseng2001convergence,Bolte2014}. To our best knowledge, this is the first attempt to solve general non-smooth binary program with guaranteed convergence. \textbf{(c)} We provide numerical comparisons with state-of-the-art techniques, such as iterative hard thresholding \cite{Yuan013}, linear programming relaxation \cite{komodakis2007approximate,KumarKT09} and semidefinite programming relaxation \cite{BQP2015Wang} on a variety of concrete computer vision and machine learning problems. Extensive experiments have demonstrated the effectiveness of our proposed methods. A preliminary version of this paper appeared in AAAI 2017 \cite{YuanG17}.



This paper is organized as follows. Section \ref{sect:related} provides a brief description of the related work. Section \ref{sect:mpecopt} presents our MPEC-based optimization framework. Section \ref{sect:discussion} discusses some features of our methods. Section \ref{sect:exp} summarizes the experimental results. Finally, Section \ref{sect:conc} concludes this paper. Throughout this paper, we use lowercase and uppercase boldfaced letters to denote real vectors and matrices respectively. The Euclidean inner product between $\bbb{x}$ and $\bbb{y}$ is denoted by $\la \bbb{x},\bbb{y}\ra$ or $\bbb{x}^T\bbb{y}$. We use $\bbb{I}_n$ to denote an identity matrix of size $n$, where sometimes the subscript is dropped when $n$ is known from the context. $\bbb{X}\succeq \bbb{0}$ means that matrix $\bbb{X}$ is positive semi-definite. Finally, ${\rm sign}$ is a signum function with ${\rm sign}(0)=\pm 1$.

\renewcommand{\arraystretch}{1}
\section{Related Work} \label{sect:related}
This paper proposes a new continuous method for binary optimization. We briefly review existing related work in this research direction in the literature (see Table \ref{tab:binary:method}).

There are generally two types of methods in the literature. One is the relaxed approximation method. Spectral relaxation \cite{cour2007solving,olsson2007solving,ShiM00,lin2013general} replaces the binary constraint with a spherical one and solves the problem using eigen decomposition. Despite its computational merits, it is difficult to generalize to handle linear or nonlinear constraints. Linear programming relaxation \cite{komodakis2007approximate,KumarKT09} transforms the NP-hard optimization problem into a convex box-constrained optimization problem, which can be solved by well-established optimization methods and software. Semi-Definite Programming (SDP) relaxation \cite{HuangCG14} uses a lifting technique $\bbb{X}=\bbb{x}\bbb{x}^T$ and relaxes to a convex conic $\bbb{X}\succeq \bbb{x}\bbb{x}^T$ \footnote{Using Schur complement lemma, one can rewrite $\bbb{X}\succeq\bbb{x}\bbb{x}^T$ as $\begin{psmallmatrix} \bbb{X} & \bbb{x} \\\bbb{x}^T  & 1\end{psmallmatrix} \succeq \bbb{0}$.} to handle the binary constraint. Combining this with a unit-ball randomized rounding algorithm, the work of \cite{Goemans95} proves that at least a factor of 87.8\% to the global optimal solution can be achieved for the graph bisection problem. Since the original paper of \cite{Goemans95}, SDP has been applied to develop numerous approximation algorithms for NP-hard problems. As more constraints lead to tighter bounds for the objective, doubly positive relaxation considers constraining both the eigenvalues and the elements of the SDP solution to be nonnegative, leading to better solutions than canonical SDP methods. In addition, Completely Positive (CP) relaxation \cite{burer2010optimizing,burer2009copositive} further constrains the entries of the factorization of the solution $\bbb{X}=\bbb{L}\bbb{L}^T$ to be nonnegative $\bbb{L}\geq\bbb{0}$. It can be solved by tackling its associated dual co-positive program, which is related to the study of indefinite optimization and sum-of-squares optimization in the literature. Second-Order Cone Programming (SOCP) relaxes the SDP conic into the nonnegative orthant \cite{KumarKT09} using the fact that $\la \bbb{X} - \bbb{x}\bbb{x}^T, \bbb{LL}^T\ra \geq 0,~\forall~\bbb{L}$, resulting in tighter bound than the LP method, but looser than that of the SDP method. Therefore it can be viewed as a balance between efficiency and efficacy.

Another type of methods for binary optimization relates to equivalent optimization. The iterative hard thresholding method directly handles the non-convex constraint via projection and it has been widely used due to its simplicity and efficiency \cite{Yuan013}. However, this method is often observed to obtain sub-optimal accuracy and it is not directly applicable, when the objective is non-smooth. A piecewise separable reformulation has been considered in \cite{zhang2007binary}, which can exploit existing smooth optimization techniques. Binary optimization can be reformulated as an $\ell_0$ norm semi-continuous optimization problem\footnote{One can rewrite the $\ell_0$ norm constraint $\|\bbb{x} + \bbb{1}\|_0 + \|\bbb{x} - \bbb{1}\|_0 \leq  n$ in a compact matrix form $\|\bbb{Ax}-\bbb{b}\|_0\leq n$ with $\bbb{A} = [\bbb{I}_n~|~\bbb{I}_n]^T \in \mathbb{R}^{2n \times n},~\bbb{b}=[\bbb{1}^T|-\bbb{1}^T]^T \in \mathbb{R}^{2n}$.}. Thus, existing $\ell_0$ norm sparsity constrained optimization techniques such as quadratic penalty decomposition method \cite{LuZ13} and multi-stage convex optimization method \cite{zhang2010analysis,yuan2016l0mpec} can be applied. A continuous $\ell_2$ box non-separable reformulation\footnote{They replace $\bbb{x} \in \{0,1\}^{n}$ with $\bbb{0} \leq \bbb{x} \leq \bbb{1},~\bbb{x}^T(\bbb{1}-\bbb{x}) = 0$. We extend this strategy to replace $\{-1,+1\}^n$ with $-\bbb{1} \leq \bbb{x} \leq \bbb{1},~(\bbb{1}+\bbb{x})^T(\bbb{1}-\bbb{x}) = 0$ which reduces to $\|\bbb{x}\|_{\infty} \leq {1},~\|\bbb{x}\|_2^2 = n$. Appendix \ref{app:equ} provides a simple proof.} has been used in the literature \cite{Raghavachari1969,Kalantari1982,Pardalos1987}. A second-order interior point method \cite{MurrayN10,de2012continuous} has been developed to solve the continuous reformulation optimization problem. A continuous $\ell_p$ box non-separable reformulation has recently been used in \cite{WuG16a}, where an interesting geometric illustration of $\ell_p$-box intersection has been shown\footnote{We adapt their formulation to our $\{-1,+1\}$ formulation.}. In addition, they infuse this equivalence into the optimization framework of Alternating Direction Method of Multipliers (ADMM). However, their guarantee of convergence is weak. In this paper, to tackle the binary optimization problem, we propose a new framework that is based on Mathematical Programming with Equilibrium Constraints (MPECs) (refer to the proposed MPEC reformulations in Table \ref{tab:binary:method}). Our resulting algorithms are theoretically convergent and empirically effective.

Mathematical programs with equilibrium constraints\footnote{In fact, we focus on mathematical programs with complementarity constraints (MPCC), a subclass of MPECs where the original discrete optimization problem can be formulated as a complementarity problem and therefore as a nonlinear program. Here we use the term MPECs for the purpose of generality and historic conventions.} are optimization problems where the constraints include complementarities or variational inequalities. They are difficult to deal with because their feasible region may not necessarily be convex or even connected. Motivated by recent development of MPECs for non-convex optimization \cite{YuanG2015l0tv,yuan2016proximal,yuan2016l0mpec,bi2014exact,luo1996mathematical}, we consider continuous MPEC reformulations of the binary optimization problem. Since our preliminary experimental results show that $\ell_2$ box non-separable MPEC in Table \ref{tab:binary:method} often presents the best performance in terms of accuracy, we only focus on this specific formulation in our forthcoming algorithm design and numerical experiments. For future references, other possible separable/non-separable\footnote{Note that separable MPEC has $n$ complementarity constraints and non-separable MPEC has one complementarity constraint. The terms separable and non-separable are related to whether the constraints can be decomposed to independent components.} MPEC reformulations are presented here. Their mathematical derivation and mathematical treatment are similar.

\section{Proposed Optimization Algorithms}\label{sect:mpecopt}

This section presents our MPEC-based optimization algorithms. We first propose an equivalent reformulation of binary optimization program, and then we consider two algorithms (exact penalty method and alternating direction method) to solve it.

\subsection{Equivalent Reformulation}

First of all, we present a variational reformulation of the binary constraint. The proofs for other reformulations of binary constraints can be found in Appendix A.

\begin{lemma} \label{lemma:binary:mpec}
\bbb{$\ell_2$ box non-separable MPEC.} We define $\Theta \triangleq \{ (\bbb{x},\bbb{v})~|~ \bbb{x}^T\bbb{v}=n,~\|\bbb{x}\|_{\infty}\leq 1,~\|\bbb{v}\|_{2}^2\leq n \}$. Assume that $(\bbb{x},\bbb{v}) \in \Theta$, then we have $\bbb{x}\in \{-1,+1\}^n$, $\bbb{v}\in \{-1,+1\}^n$, and $\bbb{x}=\bbb{v}$.

\begin{proof}
(i) Firstly, we prove that $\bbb{x} \in \{-1,+1\}^n$. Using the definition of $\Theta$ and the Cauchy-Schwarz Inequality, we have: $n =  \bbb{x}^T\bbb{v} \leq  \|\bbb{x}\|_{2}  \|\bbb{v}\|_2 \leq \|\bbb{x}\|_{2}  \sqrt{n} = \sqrt{n \bbb{x}^T \bbb{x} }   \leq \sqrt{n\|\bbb{x}\|_1 \|\bbb{x}\|_{\infty} } \leq \sqrt{n\|\bbb{x}\|_1 }$. Thus, we obtain $\|\bbb{x}\|_1\geq n$. We define $\bbb{z}=|\bbb{x}|$. Combining $\|\bbb{x}\|_{\infty}\leq 1$, we have the following constraint sets for $\bbb{z}$: $\sum_{i}~ \bbb{z}_i \geq n,~\bbb{0}\leq \bbb{z}\leq \bbb{1}$. Therefore, we have $\bbb{z}=\bbb{1}$ and it holds that $\bbb{x} \in \{-1,+1\}^n$.
(ii) Secondly, we prove that $\bbb{v} \in \{-1,+1\}^n$. We have:
\beq \label{eq:hahalast}
\textstyle n =  \bbb{x}^T\bbb{v}\leq \|\bbb{x}\|_{\infty} \|\bbb{v}\|_{1} \leq  \|\bbb{v}\|_{1} =  |\bbb{v}|^T\bbb{1}\leq \|\bbb{v}\|_2 \|\bbb{1}\|_2
\eeq
\noi Thus, we obtain $\|\bbb{v}\|_2\geq \sqrt{n}$. Combining $\|\bbb{v}\|_2^2\leq{n}$, we have $\|\bbb{v}\|_2=\sqrt{n}$ and $\|\bbb{v}\|_2 \|\bbb{1}\|_2=n$. By the Squeeze Theorem, all the equalities in (\ref{eq:hahalast}) hold automatically. Using the equality condition for Cauchy-Schwarz Inequality, we have $|\bbb{v}| = \bbb{1}$ and it holds that $\bbb{v} \in \{-1,+1\}^n$.

(iii) Finally, since $\bbb{x}\in \{-1,+1\}^n$, $\bbb{v}\in \{-1,+1\}^n$, and $\la \bbb{x},\bbb{v}\ra=n$, we obtain $\bbb{x}=\bbb{v}$.

\end{proof}
\end{lemma}

Using Lemma \ref{lemma:binary:mpec}, we can rewrite (\ref{eq:main}) in an equivalent form as follows.
\beq \label{eq:mpec}
\min_{-\bbb{1}\leq \bbb{x}\leq \bbb{1},~ \|\bbb{v}\|_2^2\leq n}~f(\bbb{x}),~s.t.~\bbb{x}^T\bbb{v}=n,~\bbb{x}\in \Omega
\eeq
\noi We remark that $\bbb{x}^T\bbb{v}=n$ is called complementarity (or equilibrium) constraint in the literature \cite{luo1996mathematical,ralph2004some,steffensen2010new} and it always holds that $\bbb{x}^T\bbb{v}\leq \|\bbb{x}\|_{\infty}\|\bbb{v}\|_{1}\leq \sqrt{n}\|\bbb{v}\|_{2}\leq n$ for any feasible $\bbb{x}$ and $\bbb{v}$.

\subsection{Exact Penalty Method} \label{sect:epm}

We now present our exact penalty method for solving the optimization problem in (\ref{eq:mpec}). It is worthwhile to point out that there are many studies on exact penalty for MPECs (refer to \cite{luo1996mathematical,hu2004convergence,ralph2004some,yuan2016l0mpec} for examples), but they do not afford the exactness of our penalty problem. In an exact penalty method, we penalize the complementary error directly by a penalty function. The resulting objective $\mathcal{J}: \mathbb{R}^{n}\times\mathbb{R}^{m}\to\mathbb{R}$ is defined in (\ref{eq:mpec:pen}), where $\rho$ is the penalty parameter that is iteratively increased to enforce the bilinear constraint.
\beq \label{eq:mpec:pen}
\begin{split}
&&\J_{\rho}(\bbb{x},\bbb{v}) = f(\bbb{x}) + \rho (n-\bbb{x}^T\bbb{v})~~~~\\
&&s.t.~-\bbb{1}\leq \bbb{x}\leq \bbb{1},~ \|\bbb{v}\|_2^2\leq n,~\bbb{x}\in \Omega
\end{split}
\eeq
\noi  In each iteration, we minimize over $\bbb{x}$ and $\bbb{v}$ alternatingly \cite{tseng2001convergence,Bolte2014}, while fixing the parameter $\rho$. We summarize our exact penalty method in Algorithm \ref{alg:epm}. The parameter $T$ is the number of inner iterations for solving the biconvex problem and the parameter $L$ is the Lipschitz constant of the objective function $f(\cdot)$. We make the following observations about the algorithm.

\begin{algorithm} [!h]
\caption{\label{alg:epm} {MPEC-EPM: An Exact Penalty Method for Solving MPEC Problem (\ref{eq:mpec})}}
(S.0) Set $t=0,~\bbb{x}^0= \bbb{v}^0 = \bbb{0}$, ~$\rho >0$,~ $\sigma>1$.

(S.1) Solve the following $\bbb{x}$-subproblem [primal step]:
     \beq \label{eq:subprob:epm:x}
     \bbb{x}^{t+1}= \mathop{\arg\min}_{\bbb{x}}~\J(\bbb{x},\bbb{v}^{t}),~s.t.~\bbb{x} \in [-1,+1]^n\cap \Omega
     \eeq

(S.2) Solve the following $\bbb{v}$-subproblem [dual step]:
     \beq \label{eq:subprob:epm:v}
      \bbb{v}^{t+1} =\mathop{\arg\min}_{ \bbb{v}}~\J(\bbb{x}^{t+1},\bbb{v}),~s.t.~\|\bbb{v}\|_2^2\leq n
     \eeq

(S.3) Update the penalty in every $T$ iterations:
                \beq \label{eq:update:epm:rho}
                \rho \Leftarrow \min(2L,\rho \times \sigma )
                \eeq
(S.4) Set $t:=t+1$ and then go to Step (S.1)
 \end{algorithm}

\vspace{1pt}\noi \bbb{(a) Initialization}. We initialize $\bbb{v}^0$ to $\bbb{0}$. This is for the sake of finding a reasonable local minimum in the first iteration, as it reduces to LP convex relaxation \cite{komodakis2007approximate} for the binary optimization problem.

\vspace{1pt}\noi \bbb{(b) Exact property}. One remarkable feature of our method is the boundedness of the penalty parameter $\rho$ (see Theorem \ref{theorem:1}). Therefore, we terminate the optimization when the threshold is reached (see (\ref{eq:update:epm:rho}) in Algorithm \ref{alg:epm}). This distinguishes it from the quadratic penalty method \cite{LuZ13}, where the penalty may become arbitrarily large for non-convex problems.

\vspace{1pt}\noi \bbb{(c) $\bbb{v}$-Subproblem}. Variable $\bbb{v}$ in (\ref{eq:subprob:epm:v}) is updated by solving the following convex problem:
\beq \label{eq:u:subp}
\begin{split}
\textstyle \bbb{v}^{t+1} =\mathop{\arg\min} ~   \la \bbb{v},-\bbb{x}^{t+1} \ra,~s.t.~\|\bbb{v}\|_2^2\leq n
\end{split}
\eeq
\noi When $\bbb{x}^{t+1} = 0$, any feasible solution is also an optimal solution; when $\bbb{x}^{t+1} \neq 0$, the optimal solution will be achieved in the boundary with $\|\bbb{v}\|_2^2= n$ and (\ref{eq:u:subp}) is equivalent to solving: $\mathop{\min}_{\|\bbb{v}\|_2^2=n}~ \frac{1}{2}\|\bbb{v}\|_2^2- \la \bbb{v},\bbb{x}^{t+1} \ra$. Thus, we have the following optimal solution for $\bbb{v}$:
\beq
\bbb{v}^{t+1} =
 \left\{
   \begin{array}{ll}
     \sqrt{n} \cdot {\bbb{x}^{t+1}}/{\|\bbb{x}^{t+1}\|_2}, & \hbox{$\bbb{x}^{t+1} \neq 0$;} \\
     \text{any}~\bbb{v}~\text{with}~\|\bbb{v}\|_2^2\leq n, & \hbox{otherwise.}
   \end{array}
 \right.
 \eeq

\noi \bbb{(d) $\bbb{x}$-Subproblem}. Variable $\bbb{x}$ in (\ref{eq:subprob:epm:x}) is updated by solving box constrained convex problem which has no closed-form solution. However, it can be solved using Nesterov's proximal gradient method \cite{Nesterov03} or classical/linearized ADM \cite{HeY12}.

\textbf{Theoretical Analysis.} In the following, we present some theoretical analysis of our exact penalty method. The following lemma is very crucial and useful in our proofs.

\begin{lemma}\label{eq:sign:bound}
We define:
\beq \label{projection:lowbound}
\begin{split}
h^* \triangleq \left( \min_{\bbb{x}\in \mathbb{R}^n,~-\bbb{1}\leq \bbb{x}\leq \bbb{1},~\text{sign}(\bbb{x}) \neq \bbb{x}}~h(\bbb{x})\right)   \\
\text{with}~~h(\bbb{x})\triangleq \frac{n-\sqrt{n}\|\bbb{x}\|_2}{\|\text{sign}(\bbb{x}) - \bbb{x}\|_2} ~~~~~~~~
\end{split}
\eeq
\noi It holds that $h^* ~\geq~ 1/2$.
\begin{proof}
We define $B$ and $N$ as the index set of the binary and non-binary variable for any $\bbb{x} \in \mathbb{R}^n$, respectively. Therefore, we have $|\bbb{x}_i|=1,~\forall i\in B$ and $|\bbb{x}_j|\neq1,~\forall j\in N$. Moreover, we define $s \triangleq |N|$ and we have $|B|=n-s$. We can rewrite $h(\bbb{x})$ as
\beq
h(\bbb{x}) &=& \frac{n-\sqrt{n}\sqrt{ \|\bbb{x}_B\|_2^2 + \|\bbb{x}_N\|_2^2 }}{\sqrt{\|\text{sign}(\bbb{x}_B) - \bbb{x}_B\|^2_2 + \|\text{sign}(\bbb{x}_N) - \bbb{x}_N\|_2^2 } } \nn \\
&=& \frac{n-\sqrt{n}\sqrt{ |B| + \|\bbb{x}_N\|_2^2 }}{\sqrt{0+\|\text{sign}(\bbb{x}_N) - \bbb{x}_N\|_2^2 } } \nn
\eeq
\noi We have the following problem which is equivalent to (\ref{projection:lowbound})
\beq \label{projection:lowbound2}
\begin{split}
\min_{\bbb{z} \in \mathbb{R}^{s}}~h'(\bbb{z})=\frac{n-\sqrt{n}\sqrt{ |B| + \|\bbb{z}\|_2^2 }}{\sqrt{\|\text{sign}(\bbb{z}) -\bbb{z}\|_2^2 } }\\
s.t.~{-\bbb{1}\leq \bbb{z}\leq \bbb{1}, ~|\bbb{z}|_i\neq 1, \forall i}
\end{split}
\eeq

\noi We notice that the function $h'(\bbb{z})$ has the property that $h'(\bbb{z})=h'(\bbb{-z})=h'(\bbb{Pz})$ for any perturbation matrix $\bbb{P}$ with $\bbb{P}\in\{0,1\}^{s\times s}, ~\bbb{P1}=\bbb{1},~\bbb{1P}=\bbb{1}$. Therefore, the optimal solution $\bbb{z}^*$ for problem (\ref{projection:lowbound2}) will be achieved with $|\bbb{z}_1^*|=|\bbb{z}_2^*|=,...,=|\bbb{z}_{s}^*|=\delta$ for some unknown constant $\delta$ with $0<\delta<1$. Therefore, we have the following optimization problem which is equivalent to (\ref{projection:lowbound}) and (\ref{projection:lowbound2})
\beq \label{projection:lowbound3}
\min_{s,\delta}~ J(s,\delta) \triangleq \frac{n-\sqrt{n}\sqrt{ n-s + s \delta^2}}{\sqrt{ s (1-\delta)^2 } } ,~s.t.~s \in \{1,2.,,,.n\}
\eeq
\noi Since $J(s,\delta)$ is an increasing function for any $\delta\in(0,1)$, the optimal solution $s^*$ for (\ref{projection:lowbound3}) with be achieved with $s^*=1$. Therefore, we derive the following inequalities:
\beq \label{eq:lowerbound}
h^* \geq\frac{n - \sqrt{n} \sqrt{(n-1) + \delta^2}} {\sqrt{(1-\delta)^2}} \triangleq p(\delta)
 \eeq

\noi In what follows, we prove that $p(\delta)>1/2$ on $0<\delta<1$ for any $n$. First of all, it is not hard to validate that $p(\delta)$ is a decreasing function on $0<\delta<1$ for any $n$. We now consider the two limit points for $p(\delta)$ in (\ref{eq:lowerbound}). (i) We have $\lim_{\delta \rightarrow 0^{+}} p(\delta) = n - \sqrt{n^2-n} > 1/2$ due to the following inequalities: $1/4>0 \Rightarrow  n^2-n+{1}/{4}>n^2-n \Rightarrow (n-1/2)^2>{n^2-n} \Rightarrow (n-1/2)>\sqrt{n^2-n} \Rightarrow n-\sqrt{n^2-n}  >1/2$. (ii) Moreover, we have $\lim_{\delta \rightarrow 1^{-}} p(\delta)= \lim_{\delta \rightarrow 1^{-}} \frac{n - \sqrt{n} \sqrt{(n-1) + \delta^2}} {\sqrt{(1-\delta)^2}} = \lim_{\delta \rightarrow 1^{-}} \frac{n -  \sqrt{n^2-n + \delta^2 n}} {1-\delta} > 1/2$. Hereby we finish the proof of this lemma.
\end{proof}
\end{lemma}

The following lemma is useful in establishing the exactness property of the penalty function in Algorithm \ref{alg:epm}.

\begin{lemma} \label{lemma:eq:hold}
Consider the following optimization problem:
\beq \label{eq:opt:pen}
(\bbb{x}_{\rho}^*,\bbb{v}_{\rho}^*) = \arg\min_{-\bbb{1}\leq \bbb{x}\leq \bbb{1},~\|\bbb{v}\|_2^2\leq n,~\bbb{x}\in \Omega}~\J_{\rho}(\bbb{x},\bbb{v}).
\eeq
\noi Assume that $f(\cdot)$ is a $L$-Lipschitz continuous convex function on $-\bbb{1} \leq \bbb{x}\leq \bbb{1}$. When $\rho>2L$, $\la\bbb{x}_{\rho}^*,\bbb{v}_{\rho}^*\ra =n$ will be achieved for any local optimal solution of (\ref{eq:opt:pen}).

\begin{proof}
To finish the proof of this lemma, we consider two cases: $\bbb{x}=\text{sign}(\bbb{x})$ and $\bbb{x} \neq \text{sign}(\bbb{x})$. For the first case, $\bbb{x}=\text{sign}(\bbb{x})$ implies that the solution is binary, any $\rho$ guarantees that the complimentary constraint is satisfied. The conclusion of this lemma clearly holds. Now we consider the case that $\bbb{x} \neq \text{sign}(\bbb{x})$.

First of all, we focus on the \bbb{v}-subproblem in the optimization problem in (\ref{eq:opt:pen}).
\beq
\bbb{v}_{\rho}^* = \arg \min_{\bbb{v}}~-\bbb{x}^T\bbb{v},~s.t.~\|\bbb{v}\|_2^2 \leq n
\eeq
\noi Assume that $\bbb{x}_{\rho}^*\neq\bbb{0}$, we have $\bbb{v}_{\rho}^*=\sqrt{n}\cdot {\bbb{x}_{\rho}^*}/{\|\bbb{x}_{\rho}^*\|_2}$. Then the biconvex optimization problem reduces to the following:
\beq \label{eq:subproblem:x}
\bbb{x}_{\rho}^* = \arg \min_{-\bbb{1}\leq \bbb{x}\leq \bbb{1}} p(\bbb{x}) \triangleq  f(\bbb{x}) + \rho (n-\sqrt{n}\|\bbb{x}\|_2)
\eeq

\noi For any $\bbb{x}_{\rho}^* \in \Omega$, we derive the following inequalities:
\beq
&&0.5\rho\| \text{sign}(\bbb{x}_{\rho}^*) - \bbb{x}_{\rho}^*  \|_2 \nn\\
&\leq& \rho(n-\sqrt{n} \|\bbb{x}_{\rho}^*\|_2) \nn\\
&=&[\rho(n-\sqrt{n} \|\bbb{x}_{\rho}^*\|_2) + f(\bbb{x}_{\rho}^*)] - f(\bbb{x}_{\rho}^*)  \nn \\
&\leq&[\rho(n-\sqrt{n} \|  \text{sign}(\bbb{\bbb{x}_{\rho}^*}) \|_2) + f(\text{sign}(\bbb{\bbb{x}_{\rho}^*}))] f(\bbb{x}_{\rho}^*)\nn\\
&=&  f(\text{sign}(\bbb{\bbb{x}_{\rho}^*})) - f(\bbb{x}_{\rho}^*)\nn\\
&=&L \|\text{sign}(\bbb{x}_{\rho}^*)- \bbb{x}_{\rho}^*\|_2 \label{eq:bound:xsignx}
\eeq
\noi where the first step uses Lemma \ref{eq:sign:bound} that $\|\text{sign}(\bbb{x}) - \bbb{x}\|_2 \leq 2 (n-\sqrt{n}\|\bbb{x}\|_2)$ for any $\bbb{x}$ in $\|\bbb{x}\|_{\infty} \leq {1}$; the third steps uses the optimality of $\bbb{x}_{\rho}^*$ in (\ref{eq:subproblem:x}) that $p(\bbb{x}_{\rho}^*) \leq p(\bbb{y})$ for any $\bbb{y}$ with $-\bbb{1}\leq \bbb{y}\leq \bbb{1},~\bbb{y}\in \Omega$; the fourth step uses the fact that $\text{sign}(\bbb{x}_{\rho}) \in \{-1,+1\}^n$ and $\sqrt{n}\|\text{sign}(\bbb{x}_{\rho})\|_2=n$; the last step uses the Lipschitz continuity of $f(\cdot)$.

From (\ref{eq:bound:xsignx}), we have $\|\bbb{x}_{\rho}^* - \text{sign}(\bbb{x}_{\rho}^*)\|_2\cdot (\rho - 2L) \leq 0$. Since $\rho - 2L>0$, we conclude that it always holds that $\|\bbb{x}_{\rho}^* - \text{sign}(\bbb{x}_{\rho}^*)\|_2=0$. Thus, $\bbb{x}_{\rho}^*\in\{-1,+1\}^n$. Finally, we have $\bbb{x}_{\rho}^* =\sqrt{n}\cdot {\bbb{x}_{\rho}^*}/{\|\bbb{x}_{\rho}^*\|_2} = \bbb{v}_{\rho}^*$ and $\la \bbb{x}_{\rho}^*, \bbb{v}_{\rho}^*\ra = n$.

\end{proof}
\end{lemma}


The following theorem shows that when the penalty parameter $\rho$ is larger than some threshold, the biconvex objective function in (\ref{eq:mpec:pen}) is equivalent to the original constrained MPEC problem in (\ref{eq:mpec}). This essentially implies the theoretical convergence of the algorithm since it reduces to block coordinate descent in the literature\footnote{Specifically, using Tseng's convergence results of block coordinate descent for non-differentiable minimization \cite{tseng2001convergence}, one can guarantee that every clustering point of Algorithm \ref{alg:adm} is also a stationary point. In addition, stronger convergence results \cite{Bolte2014,yuan2016l0mpec} can be obtained by combining a proximal strategy and Kurdyka-\L ojasiewicz inequality assumption on $\mathcal{J}(\cdot)$.}.

\begin{theorem} \label{theorem:1}
\textbf{Exactness of the Penalty Function.} Assume that $f(\cdot)$ is a $L$-Lipschitz continuous convex function on $-\bbb{1} \leq \bbb{x}\leq \bbb{1}$. When $\rho>2L$, the biconvex optimization in (\ref{eq:mpec:pen}) has the same local and global minima with the original problem in (\ref{eq:mpec}).

\begin{proof}
We let $\bbb{x}^*$ be any global minimizer of (\ref{eq:mpec}) and $(\bbb{x}^*_\rho,\bbb{v}^*_\rho)$ be any global minimizer of (\ref{eq:mpec:pen}) for some $\rho>2L$.

(i) We now prove that $\bbb{x}^*$ is also a global minimizer of (\ref{eq:mpec:pen}). For any feasible $\bbb{x}$ and $\bbb{v}$ that $\|\bbb{x}\|_{\infty}\leq 1,~ \|\bbb{v}\|_2^2\leq n$, we derive the following inequalities:
\beq
&&\J(\bbb{x},\bbb{v},\rho) \nn\\
&\geq& \min_{\|\bbb{x}\|_{\infty}\leq 1,~ \|\bbb{v}\|_2^2\leq n,~\bbb{x}\in\Omega}~f(\bbb{x}) + \rho (n-\bbb{x}^T\bbb{v})\nn\\
&=& \min_{\|\bbb{x}\|_{\infty}\leq 1,~ \|\bbb{v}\|_2^2\leq n,~\bbb{x}\in\Omega}~f(\bbb{x}),~s.t.~ \bbb{x}^T\bbb{v} = n \nn\\
&=& f(\bbb{x}^*) + \rho (n- \bbb{x}^{*T} \bbb{v}^* ) \nn\\
&=& \J(\bbb{x}^*,\bbb{v}^*,\rho) \nn
\eeq
\noi where the first equality holds due to the fact that the constraint $\bbb{x}^T\bbb{v} = n$ is satisfied at the local optimal solution when $\rho>2L$ (see Lemma \ref{lemma:eq:hold}). Therefore, we conclude that any optimal solution of (\ref{eq:mpec}) is also an optimal solution of (\ref{eq:mpec:pen}).

(ii) We now prove that $\bbb{x}_{\rho}^*$ is also a global minimizer of (\ref{eq:mpec}). For any feasible $\bbb{x}$ and $\bbb{v}$ that $\|\bbb{x}\|_{\infty}\leq 1,~ \|\bbb{v}\|_2^2\leq n,~\bbb{x}^T\bbb{v}=n,~\bbb{x}\in\Omega$, we naturally have the following inequalities:
\beq
&&f(\bbb{x}_\rho^*) - f(\bbb{x})\nn\\
&=& f(\bbb{x}_\rho^*) +  \rho (n- \bbb{x}_{\rho}^{*T} \bbb{v}_{\rho}^* )  - f(\bbb{x}) -  \rho (n- \bbb{x}^{T} \bbb{v} ) \nn\\
&=& \J_\rho(\bbb{x}_\rho^*,\bbb{v}_\rho^*) - \J_\rho(\bbb{x},\bbb{v}) \nn\\
&\leq& 0\nn
\eeq
\noi where the first equality uses Lemma \ref{lemma:eq:hold}. Therefore, we conclude that any optimal solution of (\ref{eq:mpec:pen}) is also an optimal solution of (\ref{eq:mpec}).

Finally, we conclude that when $\rho>2L$, the biconvex optimization in (\ref{eq:mpec:pen}) has the same local and global minima with the original problem in (\ref{eq:mpec}).

%

\end{proof}
\end{theorem}

The following the theorem characterizes the convergence rate and asymptotic monotone property the algorithm.

\begin {theorem}\label{theorem:2}
\textbf{Convergence Rate and Asymptotic Monotone Property of Algorithm \ref{alg:epm}.} Assume that $f(\cdot)$ is a $L$-Lipschitz continuous convex function on $-\bbb{1} \leq \bbb{x}\leq \bbb{1}$. Algorithm \ref{alg:epm} will converge to the first-order KKT point in at most $\lceil  (\ln(L\sqrt{2n}) - \ln(\epsilon \rho^0)) / \ln \sigma \rceil$ outer iterations\footnote{Every time we increase $\rho$, we call it one outer iteration.} with the accuracy at least $ n - \bbb{x}^T\bbb{v}\leq \epsilon$. Moreover, after $\la \bbb{x},\bbb{v} \ra=n$ is obtained, the sequence of $\{f(\bbb{x}^t)\}$ generated by Algorithm \ref{alg:epm} is monotonically non-increasing.

\begin{proof}
We denote $s$ and $t$ as the outer iterations counter and inner iteration counter in Algorithm \ref{alg:epm}, respectively.

(i) we now prove the convergence rate of Algorithm \ref{alg:epm}. Assume that Algorithm \ref{alg:epm} takes $s$ outer iterations to converge. We denote $f'(\bbb{x})$ as the sub-gradient of $f(\cdot)$ in $\bbb{x}$. According the the $\bbb{x}$-subproblem in (\ref{eq:subproblem:x}), if $\bbb{x}^*$ solves (\ref{eq:subproblem:x}), then we have the following variational inequality \cite{HeY12}:
\beq
\forall \bbb{x} \in [-1,+1]^n\cap\Omega,~ \la \bbb{x}-\bbb{x}^*,f'(\bbb{x}^*) \ra +\nn\\
\rho(n-\sqrt{n}\|\bbb{x}\|_2) - \rho(n-\sqrt{n}\|\bbb{x}^*\|_2)  \geq 0\nn
\eeq
\noi Letting $\bbb{x}$ be any feasible solution that $\bbb{x}\in\{-1,+1\}^n \cap\Omega$, we have the following inequalities:
\beq \label{eq:conv:rate}
&&(n-\sqrt{n}\|\bbb{x}^*\|_2) \nn\\
&\leq& (n-\sqrt{n}\|\bbb{x}\|_2) + \tfrac{1}{\rho}\la\bbb{x}-\bbb{x}^*, f'(\bbb{x}^*)\ra\nn\\
&\leq&\tfrac{1}{\rho}\|\bbb{x}-\bbb{x}^*\|_2 \cdot \|f'(\bbb{x}^*)\|_2\nn\\
&\leq&L\sqrt{2n}/{\rho}
\eeq
\noi where the second inequality is due to the Cauchy-Schwarz Inequality, the third inequality is due to the fact that $\|\bbb{x}-\bbb{y}\|_2\leq \sqrt{2n},~\forall -\bbb{1} \leq \bbb{x}~,\bbb{y}\leq \bbb{1}$ and the Lipschitz continuity of $f(\cdot)$ that $\|f'(\bbb{x}^*)\|_2\leq L$.

The inequality in (\ref{eq:conv:rate}) implies that when $\rho^s\geq L\sqrt{2n}/\epsilon$, Algorithm \ref{alg:epm} achieves accuracy at least $n - \sqrt{n}\|\bbb{x}\|_2 \leq \epsilon$. Noticing that $\rho^s=\sigma^s\rho^0$, we have that $\epsilon$ accuracy will be achieved when
\beq
&&\sigma^s \rho^0\geq \frac{L\sqrt{2n}}{\epsilon}\nn\\
&\Rightarrow&\sigma^s\geq\frac{L\sqrt{2n}}{\epsilon \rho^0}\nn\\
&\Rightarrow&  s \geq (\ln(L\sqrt{2n}) - \ln(\epsilon \rho^0)) / \ln \sigma \nn
\eeq

(ii) we now prove the asymptotic monotone property of Algorithm \ref{alg:epm}. We naturally derive the following inequalities:
\beq
&&f(\bbb{x}^{t+1}) - f(\bbb{x}^{t})\nn\\
&\leq&  \rho(n - \la \bbb{x}^{t},\bbb{v}^t\ra) - \rho(n-\la \bbb{x}^{t+1},\bbb{v}^t\ra)\nn\\
&=& \rho\left( \la \bbb{x}^{t+1},\bbb{v}^t\ra - \la \bbb{x}^{t},\bbb{v}^t\ra  \right)\nn\\
&\leq& \rho\left( \la \bbb{x}^{t+1},\bbb{v}^{t+1}\ra - \la \bbb{x}^{t},\bbb{v}^t\ra  \right)\nn\\
&=&0\nn
\eeq
\noi where the first step uses the fact that $f(\bbb{x}^{t+1}) + \rho (n-\la \bbb{x}^{t+1},\bbb{v}^t \ra )\leq f(\bbb{x}^{t}) + \rho (n-\la \bbb{x}^{t},\bbb{v}^t \ra )$ holds due to $\bbb{x}^{t+1}$ is the optimal solution of (\ref{eq:subprob:epm:x}); the third step uses the fact $-\la \bbb{x}^{t+1},\bbb{v}^{t+1}\ra \leq -\la \bbb{x}^{t+1},\bbb{v}^t\ra$ holds due to $\bbb{v}^{t+1}$ is the optimal solution of (\ref{eq:subprob:epm:v}); the last step uses $\la \bbb{x},\bbb{v} \ra=n$. Note that the equality $\la \bbb{x},\bbb{v} \ra=n$ together with the feasible set $-\bbb{1}\leq \bbb{x}\leq \bbb{1},~\|\bbb{v}\|_2^2\leq n$ also implies that $\bbb{x}\in \{-1,+1\}^n$.
\end{proof}
\end{theorem}

\begin{algorithm}  [!t]
\caption{\label{alg:adm} { MPEC-ADM: An Alternating Direction Method for Solving MPEC Problem (\ref{eq:mpec})}}
(S.0) Set $t=0,~\bbb{x}^0=\bbb{v}^0 = \bbb{0},~\rho^{0}= 0,~\alpha>0,~\sigma>1$.

(S.1) Solve the following $\bbb{x}$-subproblem [primal step]:
     \beq \label{eq:subprob:adm:x}
     \bbb{x}^{t+1}= \mathop{\arg\min}_{\bbb{x}}~\mathcal{L}(\bbb{x},\bbb{v}^{t},\rho^t),~s.t.~\bbb{x} \in [-1,+1]^n\cap \Omega
     \eeq

(S.2) Solve the following $\bbb{v}$-subproblem [dual step]:
     \beq \label{eq:subprob:adm:v}
      \bbb{v}^{t+1} =\mathop{\arg\min}_{ \bbb{v}}~\mathcal{L} (\bbb{x}^{t+1},\bbb{v},\rho^t)~s.t.~\|\bbb{v}\|_2^2\leq n
     \eeq

(S.3) Update the Lagrange multiplier:
                \beq \label{eq:update:adm:pi}
                \rho^{t+1}=\rho^{t} + \alpha  (  n - \la \bbb{x}^{t+1},\bbb{v}^{t+1} \ra   )
                \eeq

(S.4) Update the penalty in every $T$ iterations (if necessary):
                \beq \label{eq:update:adm:alpha}
                \alpha  \Leftarrow \alpha \times \sigma
                \eeq
(S.5) Set $t:=t+1$ and then go to Step (S.1).
 \end{algorithm}

We have a few remarks on the theorems above. We assume that the objective function is $L$-Lipschitz continuous. However, such hypothesis is not strict. Because the solution $\bbb{x}$ is defined on the compact set, the Lipschitz constant can always be computed for any continuous objective (e.g. norm function, min/max envelop function). In fact, it is equivalent to say that the (sub-) gradient of the objective is bounded by $L$\footnote{For example, for the quadratic function $f(\bbb{x}) = 0.5\bbb{x}^T\bbb{Ax}+\bbb{x}^T\bbb{b}$ with $\bbb{A}\in\mathbb{R}^{n\times n}$ and $\bbb{b}\in\mathbb{R}^{n}$, the Lipschitz constant is bounded by $L \leq  \|\bbb{Ax}+\bbb{b}\|\leq \|\bbb{A}\|\|\bbb{x}\| + \|\bbb{b}\| \leq \|\bbb{A}\|\sqrt{n}+\|\bbb{b}\|$; for the $\ell_1$ regression function $f(\bbb{x})=\|\bbb{Ax}-\bbb{b}\|_1$ with $\bbb{A}\in\mathbb{R}^{m\times n}$ and $\bbb{b}\in\mathbb{R}^{m}$, the Lipschitz constant is bounded by $L \leq \|\bbb{A}^T \partial |\bbb{Ax}-\bbb{b}|\|\leq\|\bbb{A}^T\| \sqrt{m}$.}. Although exact penalty method has been study in the literature \cite{han1979exact,di1989exact,di1994exact}, their results cannot directly apply here. The theoretical bound $2 L$ (on the penalty parameter $\rho$) heavily depends on the specific structure of the optimization problem. Moreover, we also establish the convergence rate and asymptotic monotone property of our algorithm.

\subsection{ Alternating Direction Method} \label{sect:adm}
This section presents an alternating direction method (of multipliers) (ADM or ADMM) for solving (\ref{eq:mpec}). This is mainly motivated by the recent popularity of ADM in the optimization literature \cite{HeY12,Wen2010,YuanG2015l0tv,yuan2016proximal}.

We first form the augmented Lagrangian $\mathcal{L}: \mathbb{R}^{n}\times\mathbb{R}^{m}\times\mathbb{R}^m\to\mathbb{R}$ in (\ref{eq:mpec:lag}) as follows:
\beq \label{eq:mpec:lag}
\begin{split}
\mathcal{L}(\bbb{x},\bbb{v},\rho) \triangleq f(\bbb{x}) + \rho (n- \bbb{x}^T\bbb{v}) + \frac{\alpha}{2} (n- \bbb{x}^T\bbb{v})^2\\
s.t.~-\bbb{1}\leq \bbb{x}\leq \bbb{1},~ \|\bbb{v}\|_2^2\leq n,~\bbb{x}\in \Omega~~~~~~~~~~~~
\end{split}
\eeq
\noi where $\rho$ is the Lagrange multiplier associated with the complementarity constraint $n-\la \bbb{x},\bbb{v}\ra=0$, and $\alpha>0$ is the penalty parameter. Interestingly, we find that the augmented Lagrangian function can be viewed as adding an elastic net regularization \cite{Hui2005} on the complementarity error. We detail the ADM iteration steps for (\ref{eq:mpec:lag}) in Algorithm \ref{alg:adm}, which has the following properties.

\vspace{1pt}\noi \bbb{(a) Initialization}. We set $\bbb{v}^0= \bbb{0}$ and $\rho^0=0$. This finds a reasonable local minimum in the first iteration, as it reduces to LP relaxation for the $\bbb{x}$-subproblem.

\vspace{1pt}\noi \bbb{(b) Monotone property}. For any feasible solution $\bbb{x}$ and $\bbb{v}$ in (\ref{eq:mpec:lag}), it holds that $ n - \bbb{x}^T\bbb{v} \geq {0}$. Using the fact that $\alpha^t>0$ and due to the $\rho^t$ update rule, $\rho^t$ is monotone increasing.

\vspace{1pt}\noi \bbb{(c) $\bbb{v}$-Subproblem}. Variable $\bbb{v}$ in (\ref{eq:subprob:adm:v}) is updated by solving the following problem:
\beq \label{eq:V:subproblem:1}
\begin{split}
 \bbb{v}^{t+1} =\mathop{\arg\min}_{\bbb{v}}   ~ \frac{1}{2} \bbb{v}^T \bbb{a} \bbb{a}^T \bbb{v} + \la \bbb{v},\bbb{b}\ra,~s.t.~\|\bbb{v}\|_2 \leq \sqrt{n}\nn
\end{split}
\eeq
\noi where $\bbb{a}\triangleq \bbb{x}^{t+1}/\sqrt{\alpha}$ and $\bbb{b}\triangleq -(\rho + \alpha n )\cdot \bbb{x}^{t+1}$. This problem is also known as constrained eigenvalue problem in the literature \cite{gander1989constrained}, which has efficient solution. Since the quadratic matrix is of rank one, this problem has nearly closed-form solution (please refer to Appendix B for detailed discussions).

\noi \bbb{(d) $\bbb{x}$-Subproblem}. Variable $\bbb{x}$ in (\ref{eq:mpec:lag}) is updated by solving a box constrained optimization problem. Similar to (\ref{eq:subprob:epm:x}), it can be solved using Nesterov's proximal gradient method or classical/linearized ADM.

\textbf{Theoretical Analysis.} In the following, we present some theoretical analysis of our alternating direction method. The proofs of theoretical convergence are similar to our previous analysis for the exact penalty method. The following lemma is useful in building the exactness property of the augmented Lagrangian function in Algorithm \ref{alg:adm}.

\begin{lemma} \label{lemma:eq:hold2}

Consider the following optimization problem:
\beq \label{eq:opt:pen2}
(\bbb{x}_{\rho}^*,\bbb{v}_{\rho}^*) = \arg\min_{-\bbb{1}\leq \bbb{x}\leq \bbb{1},\|\bbb{v}\|_2^2\leq n,~\bbb{x}\in\Omega}~\Lag(\bbb{x},\bbb{v},\rho).
\eeq
\noi Assume that $f(\cdot)$ is a $L$-Lipschitz continuous convex function on $-\bbb{1} \leq \bbb{x}\leq \bbb{1}$. When $\rho>2L$, $\la\bbb{x}_{\rho}^*,\bbb{v}_{\rho}^*\ra =n$ will be achieved for any local optimal solution of (\ref{eq:opt:pen2}).

\begin{proof}
The proof of this theorem is based on Lemma \ref{lemma:eq:hold}. We observe that $n-\bbb{x}^T\bbb{v}=0 \Leftrightarrow (n-\bbb{x}^T\bbb{v})^2 = 0$. We define $h(\bbb{x}) \triangleq f(\bbb{x}) + \frac{\alpha}{2}(n-\bbb{x}^T\bbb{v})^2$ and denote $L_h$ as the Lipschtz constant of $h(\cdot)$. We replace $f(\cdot)$ with $h(\cdot)$ in Lemma \ref{lemma:eq:hold} and conclude that when $\rho>2 L_g$, we have $n-\bbb{x}^T\bbb{v}=0$. Thus, the term $\frac{\alpha}{2}(n-\bbb{x}^T\bbb{v})^2$ in $h(\bbb{x})$ reduces to zero and $L_h=L$. We conclude that when $\rho>2L$, $\la\bbb{x}_{\rho}^*,\bbb{v}_{\rho}^*\ra =n$ will be achieved for any local optimal solution.
\end{proof}
\end{lemma}

Our MPEC-ADM has excellent convergence property. The following theorem shows that when the multiplier $\rho$ is larger than some threshold, the biconvex objective function in (\ref{eq:mpec:lag}) is equivalent to the original constrained MPEC problem in (\ref{eq:mpec}).

\begin{theorem} \label{theorem:3}

\textbf{Exactness of the augmented Lagrangian Function.} Assume that $f(\cdot)$ is a $L$-Lipschitz continuous convex function on $-\bbb{1} \leq \bbb{x}\leq \bbb{1}$. When $\rho>2L$, the biconvex optimization problem $\min_{\bbb{x},~\bbb{v}}~\Lag(\bbb{x},\bbb{v},\rho),~s.t.~-\bbb{1}\leq \bbb{x}\leq \bbb{1},~\|\bbb{v}\|_2^2\leq n,~\bbb{x}\in \Omega$ in (\ref{eq:mpec:lag}) has the same local and global minima with the original problem in (\ref{eq:mpec}).

\begin{proof}
The proof of this theorem is similar to Theorem \ref{theorem:1}. We let $\bbb{x}^*$ be any global minimizer of (\ref{eq:mpec}) and $(\bbb{x}^*_\rho,\bbb{v}^*_\rho)$ be any global minimizer of (\ref{eq:mpec:pen}) for some $\rho>2L$.

(i) We now prove that $\bbb{x}^*$ is also a global minimizer of (\ref{eq:mpec:pen}). For any feasible $\bbb{x}$ and $\bbb{v}$ that $\|\bbb{x}\|_{\infty}\leq 1,~ \|\bbb{v}\|_2^2\leq n,~\bbb{x}\in\Omega$, we derive the following inequalities:
\beq
&&\Lag(\bbb{x},\bbb{v},\rho) \nn\\
&\geq& \min_{\|\bbb{x}\|_{\infty}\leq 1,~ \|\bbb{v}\|_2\leq n}~f(\bbb{x}) + \rho (n-\bbb{x}^T\bbb{v}) + \frac{\alpha}{2}(n- \bbb{x}^T\bbb{v})^2\nn\\
&=& \min_{\|\bbb{x}\|_{\infty}\leq 1,~ \|\bbb{v}\|_2\leq n}~f(\bbb{x}),~s.t.~ \bbb{x}^T\bbb{v} = n \nn\\
&=& f(\bbb{x}^*) + \rho (n- \bbb{x}^{*T} \bbb{v}^* ) + \frac{\alpha}{2}(n- \bbb{x}^{*T} \bbb{v}^* )^2 \nn\\
&=& \Lag(\bbb{x}^*,\bbb{v}^*,\rho) \nn
\eeq
\noi where the first equality holds due to the fact that the constraint $\bbb{x}^T\bbb{v} = n$ is satisfied at the local optimal solution when $\rho>2L$ (see Lemma \ref{lemma:eq:hold}). Therefore, we conclude that any optimal solution of (\ref{eq:mpec}) is also an optimal solution of (\ref{eq:mpec:pen}).

(ii) We now prove that $\bbb{x}_{\rho}^*$ is also a global minimizer of (\ref{eq:mpec}). For any feasible $\bbb{x}$ and $\bbb{v}$ that $\|\bbb{x}\|_{\infty}\leq 1,~ \|\bbb{v}\|_2^2\leq n,~\bbb{x}^T\bbb{v}=n,~\bbb{x}\in\Omega$, we naturally have the following inequalities:
\beq
&&f(\bbb{x}_\rho^*) - f(\bbb{x})\nn\\
&=& f(\bbb{x}_\rho^*) +  \rho (n- \bbb{x}_{\rho}^{*T} \bbb{v}_{\rho}^* ) + \frac{\alpha}{2}(n- \bbb{x}_{\rho}^{*T} \bbb{v}_{\rho}^* )^2 \nn\\
&&  - f(\bbb{x}) -  \rho (n- \bbb{x}^{T} \bbb{v} ) - \frac{\alpha}{2}(n- \bbb{x}^{T} \bbb{v} )^2\nn\\
&=& \Lag(\bbb{x}_\rho^*,\bbb{v}_\rho^*,\rho) - \Lag(\bbb{x},\bbb{v},\rho) \nn\\
&\leq& 0\nn
\eeq
\noi Therefore, we conclude that any optimal solution of (\ref{eq:mpec:pen}) is also an optimal solution of (\ref{eq:mpec}).

Finally, we conclude that when $\rho>2L$, the biconvex optimization in (\ref{eq:mpec:pen}) has the same local and global minima with the original problem in (\ref{eq:mpec}).

\end{proof}
\end{theorem}



We have some remarks on the theorem above. We use the same assumption as in the exact penalty method. Although alternating direction method has been used in the literature \cite{lai2014splitting,WuG16a}, their convergence results are much weaker than ours. The main novelties of our method is the self-penalized feature owning to the boundedness and monotone property of multiplier.

\begin{figure*}[!th]
\captionsetup[subfigure]{justification=centering}
    \centering
      \begin{subfigure}{\fourfigwid}\includegraphics[width=\textwidth,height=\imgheibisect]{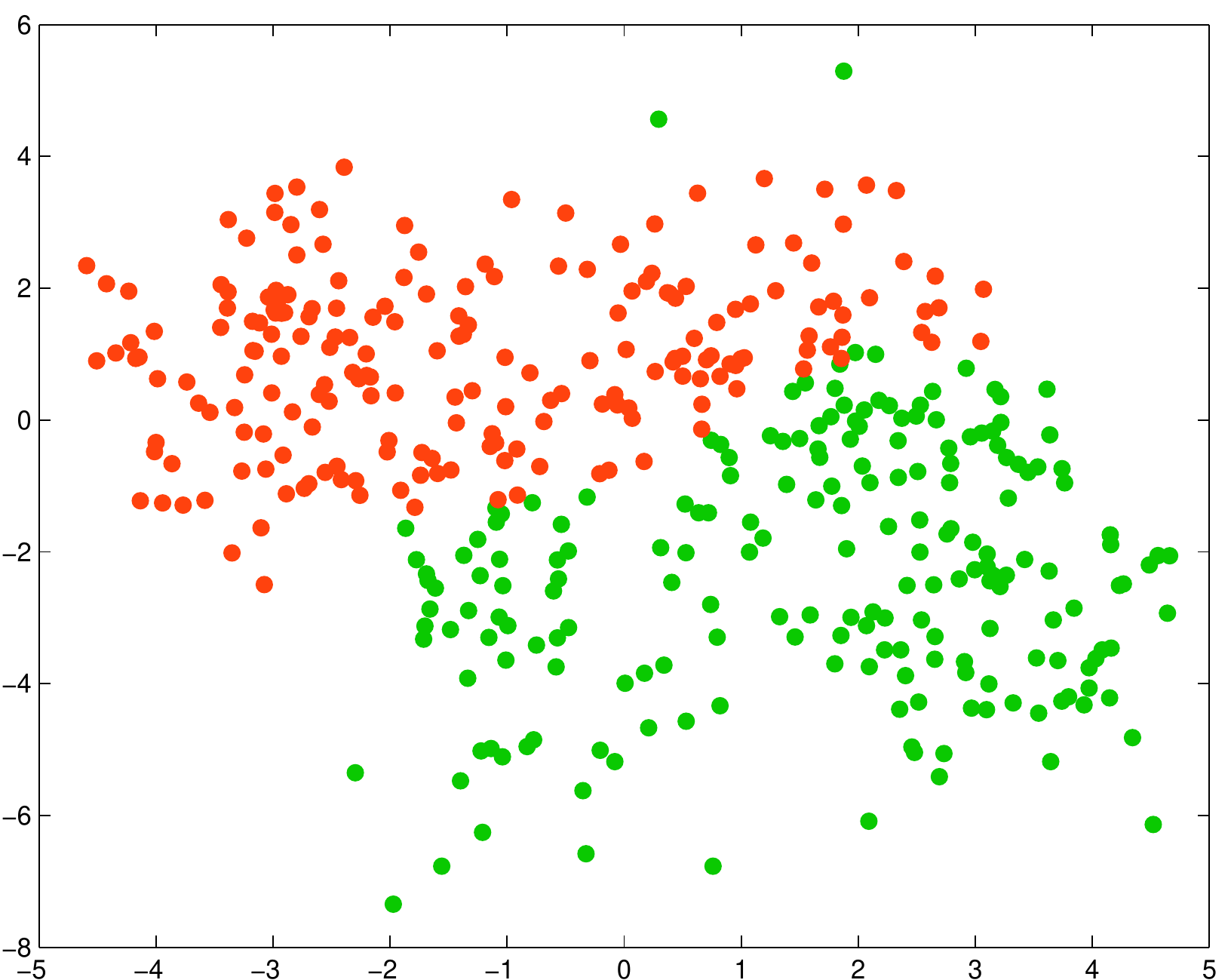}\caption{\text{LP}\\$f = 473.646$}\end{subfigure}\ghs
      \begin{subfigure}{\fourfigwid}\includegraphics[width=\textwidth,height=\imgheibisect]{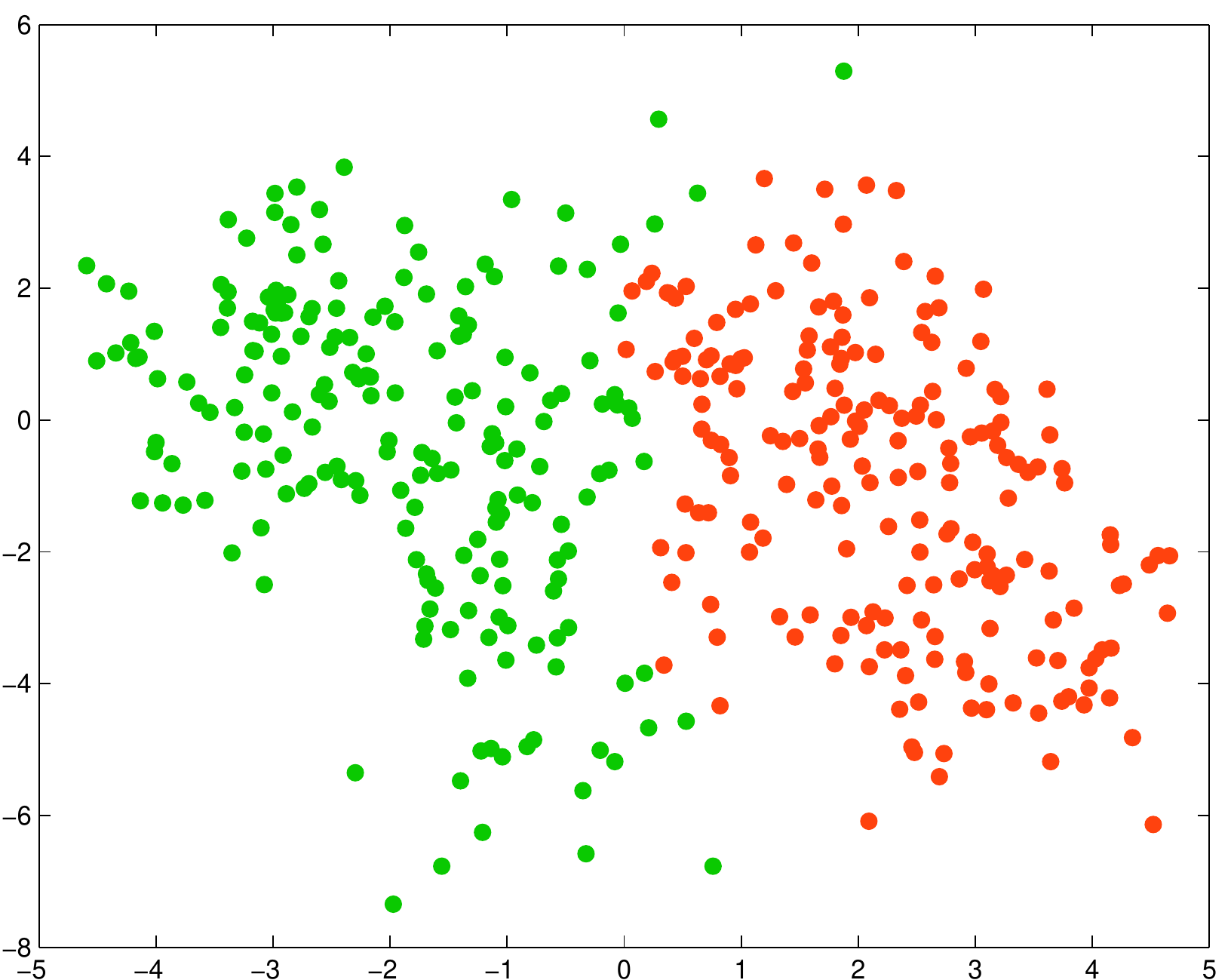}\caption{\text{NCUT}\\$f = 230.049$}\end{subfigure}\ghs
      \begin{subfigure}{\fourfigwid}\includegraphics[width=\textwidth,height=\imgheibisect]{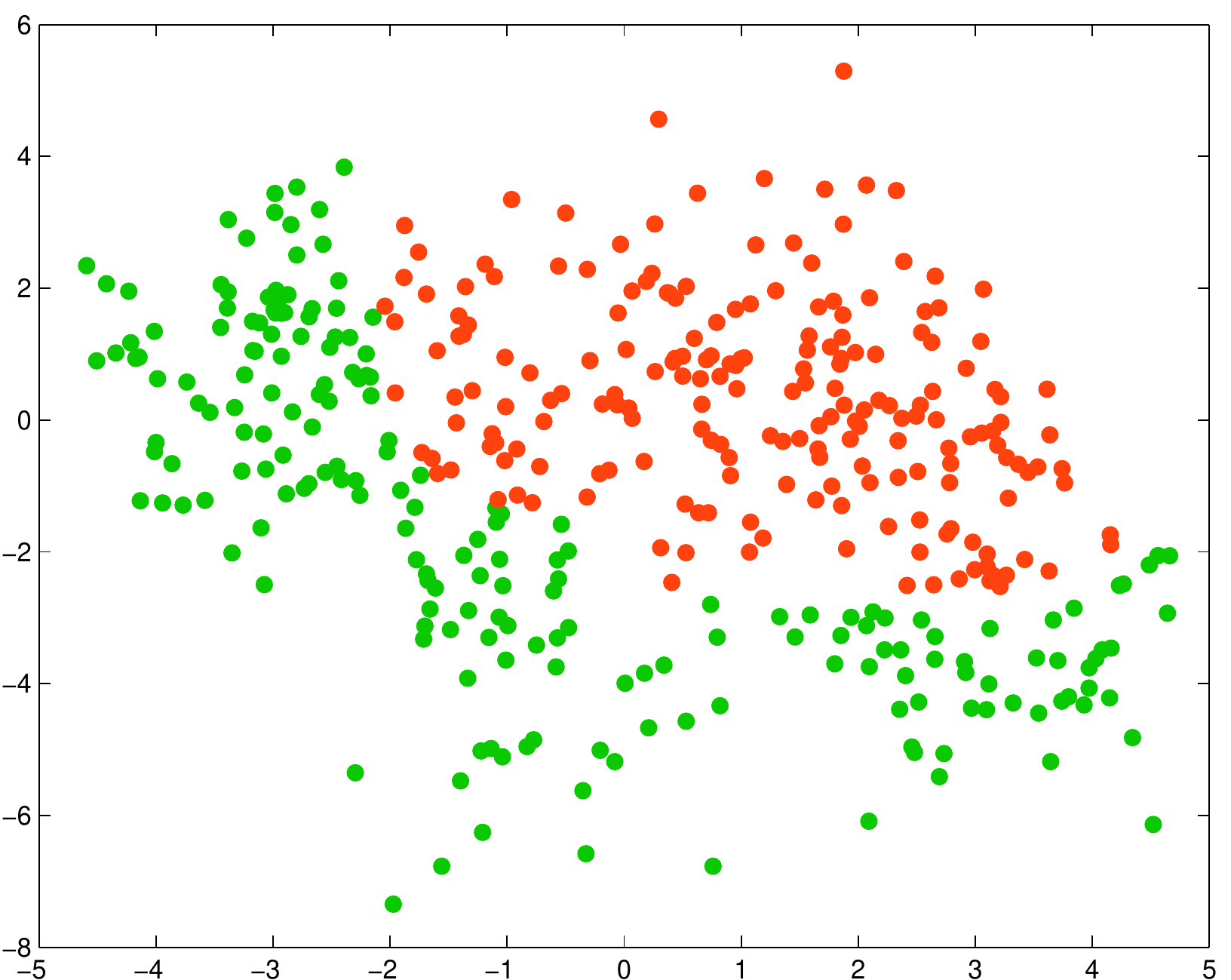}\caption{\text{RCUT}\\$f = 548.964$}\end{subfigure}\ghs
      \begin{subfigure}{\fourfigwid}\includegraphics[width=\textwidth,height=\imgheibisect]{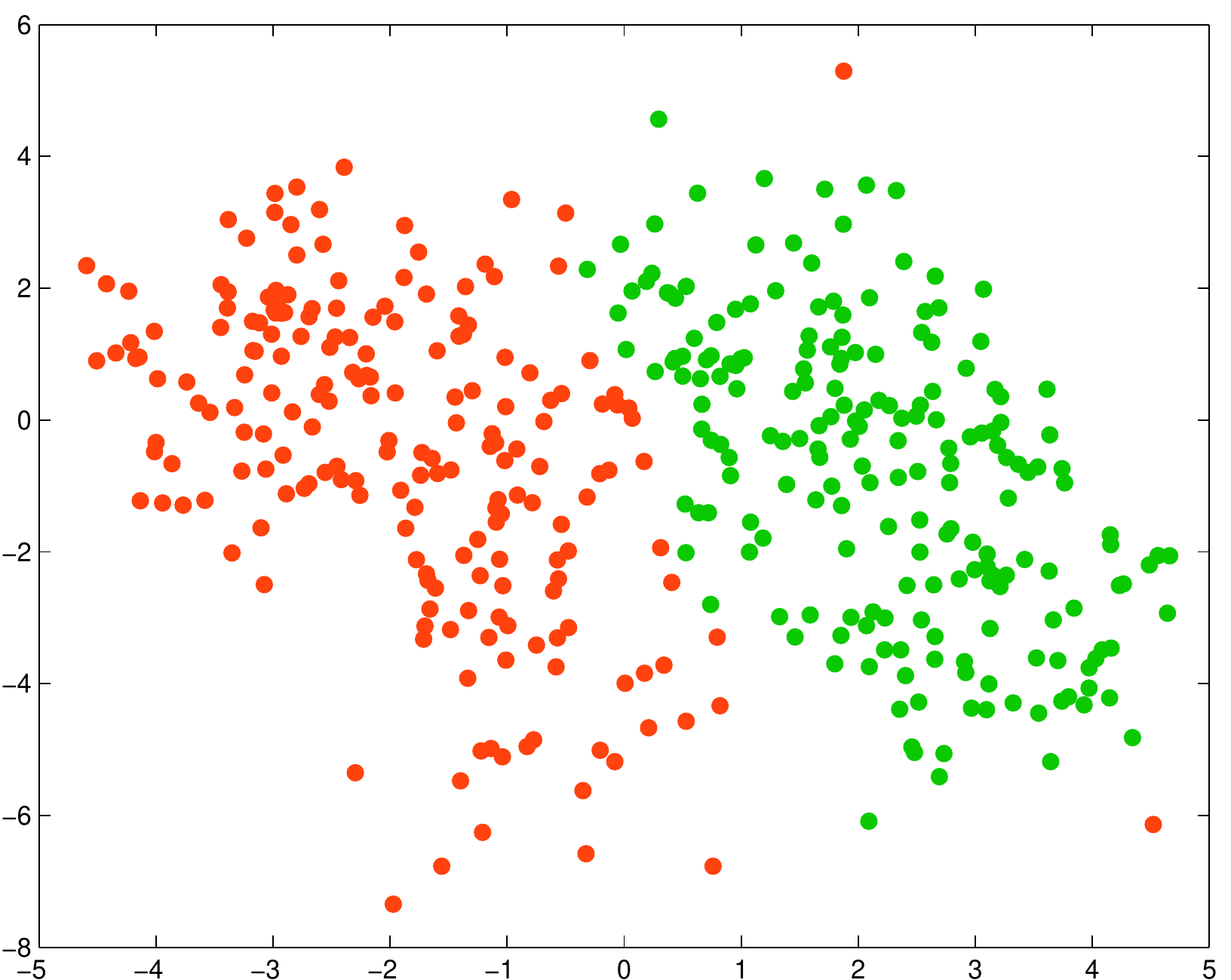}\caption{\text{SDP}\\$f = 194.664$}\end{subfigure}

      \begin{subfigure}{\fourfigwid}\includegraphics[width=\textwidth,height=\imgheibisect]{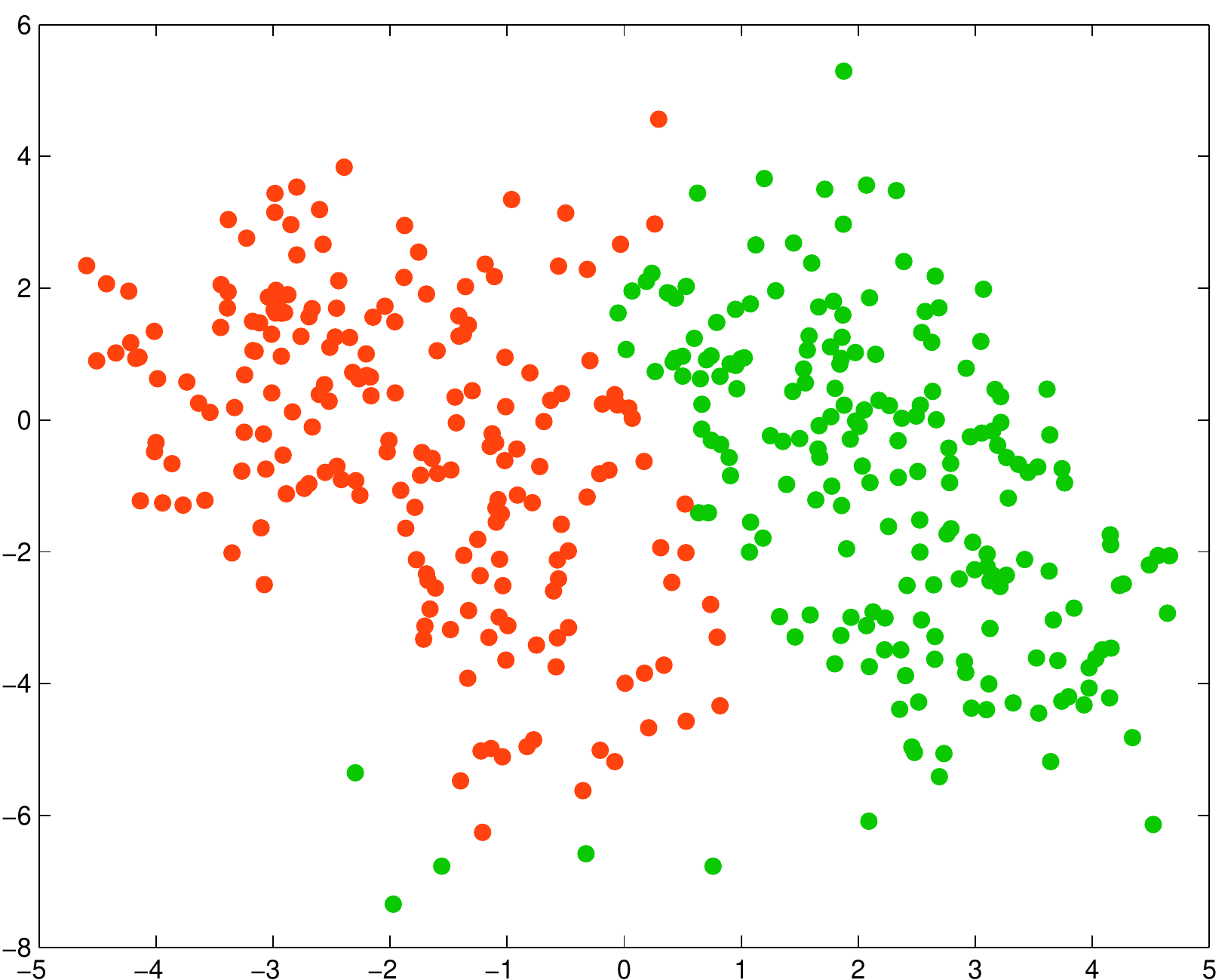}\caption{\text{L2box-ADMM}\\$f = 196.964$}\end{subfigure}\ghs
      \begin{subfigure}{\fourfigwid}\includegraphics[width=\textwidth,height=\imgheibisect]{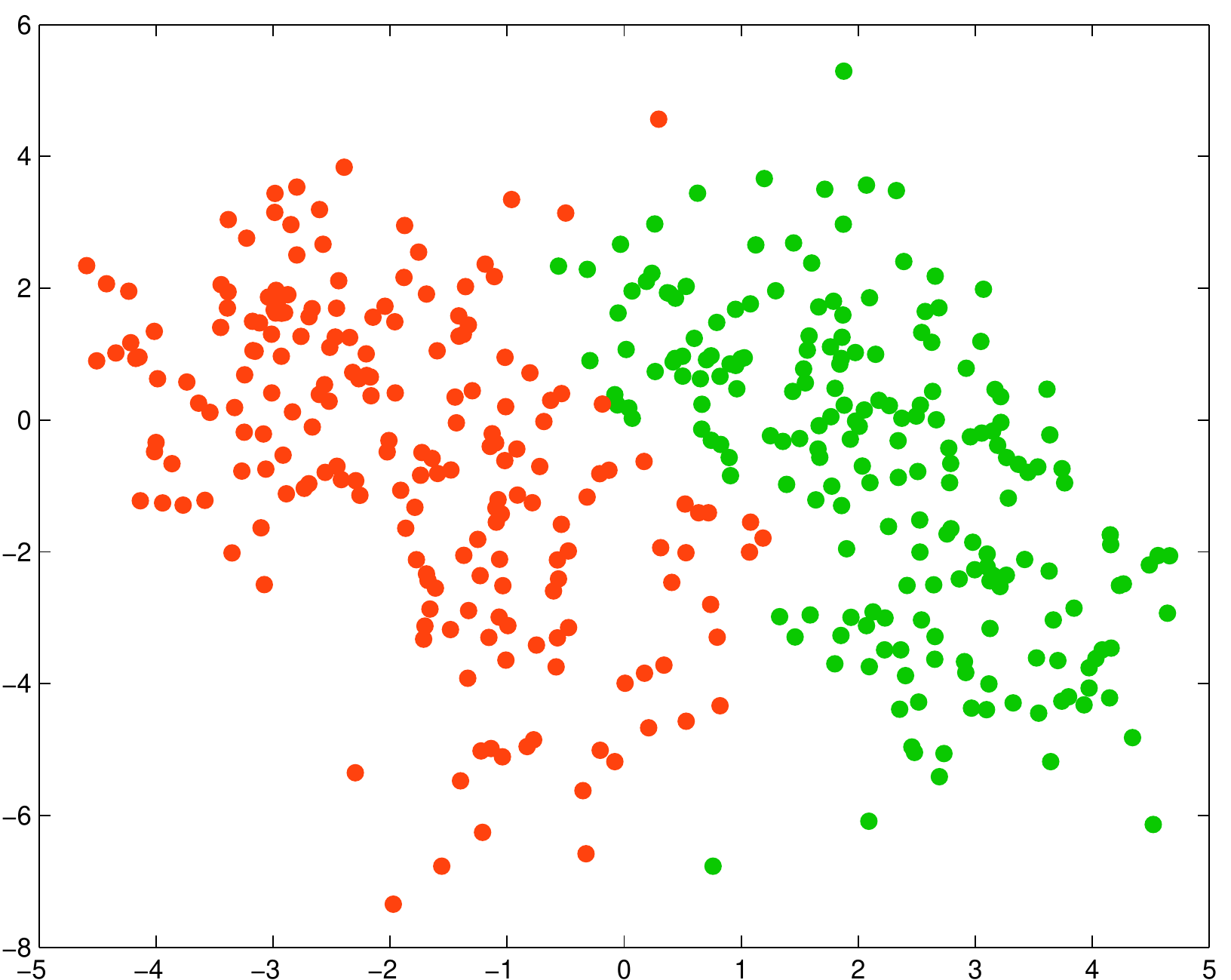}\caption{\text{MPEC-EPM}\\$f = 186.926$}\end{subfigure}\ghs
      \begin{subfigure}{\fourfigwid}\includegraphics[width=\textwidth,height=\imgheibisect]{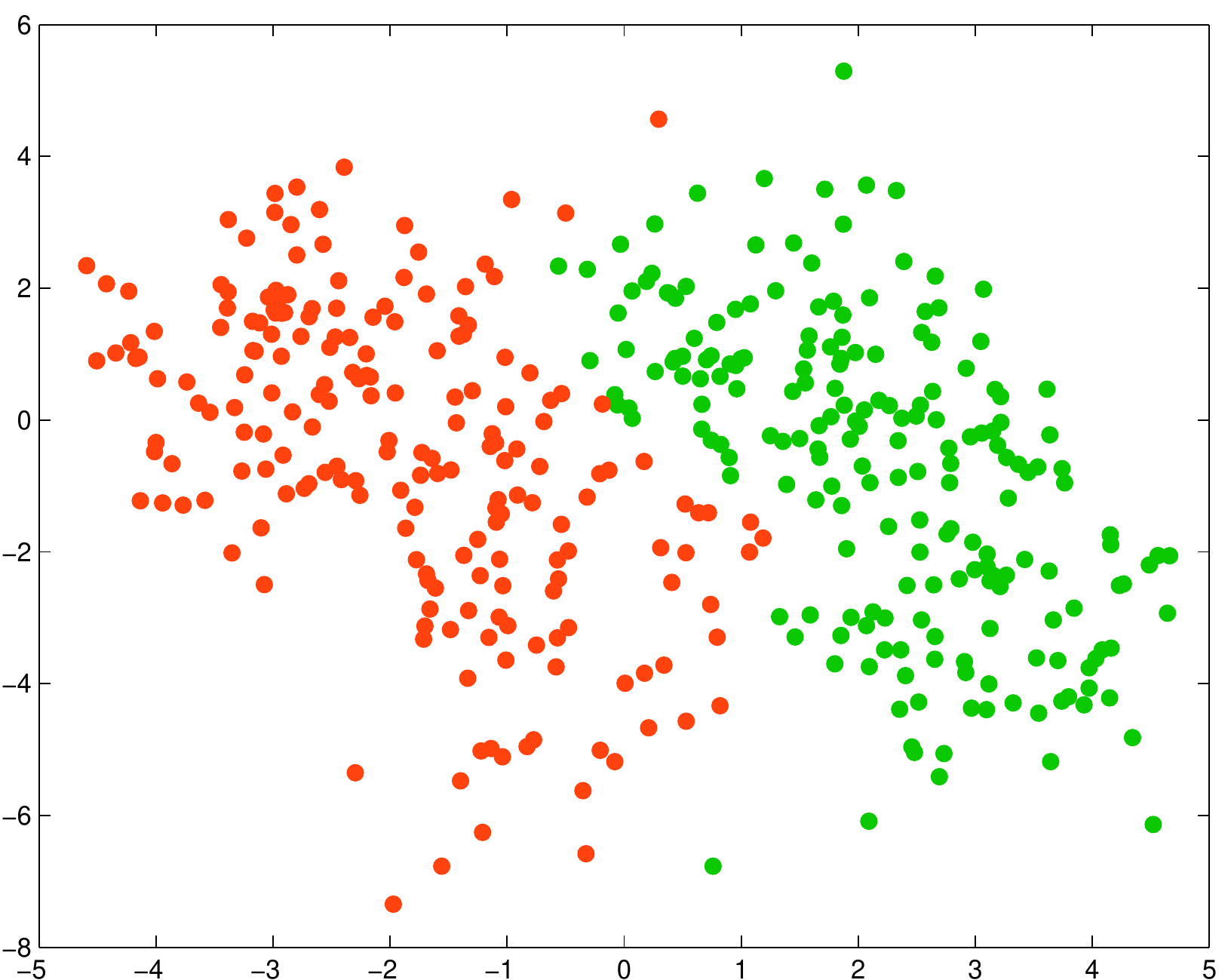}\caption{\text{MPEC-ADM}\\$f = 186.926$}\end{subfigure}
\caption{Graph bisection on the `4gauss' data set.}
\label{fig:bsect}
\end{figure*}

\section{Discussions} \label{sect:discussion}

This section discusses the comparisons of MPEC-EPM and MPEC-ADM, the merits of our methods, and the extensions to zero-one constraint and orthogonality constraint.

\textbf{MPEC-EPM vs. MPEC-ADM.} The proposed MPEC-EPM and MPEC-ADM have their own advantages. {(a)} MPEC-EPM is more simple and elegant and it can directly use existing LP relaxation optimization solver. {(b)} MPEC-EPM may be less adaptive since the penalty parameter $\rho$ is monolithically increased until a threshold is achieved. In comparison, MPEC-ADM is more adaptive, since a constant penalty also guarantees monotonically non-decreasing multipliers and convergence. 

\textbf{Merits of our methods.} There are several merits behind our MPEC-based penalization/regularization methods. {(a)} They exhibit strong convergence guarantees since they essentially reduce to block coordinate descent in the literature \cite{tseng2001convergence,Bolte2014}. {(b)} They seek desirable solutions since the LP convex relaxation methods in the first iteration provide good initializations. {(c)} They are efficient since they are amenable to the use of existing convex methods to solve the sub-problems. {(d)} They have a monotone/greedy property due to the complimentary constraints brought on by MPECs. We penalize the complimentary error and ensure that the error is decreasing in every iteration, leading to binary solutions.


\textbf{Extensions to Zero-One Constraint and Orthogonality Constraint.} For convenience, we define $\Delta = \{0,1\}^n$ and $\mathbb{O}=\{\bbb{X}\in \mathbb{R}^{n\times r}$$~$$|~$$\bbb{X}^T\bbb{X}=\bbb{I}_r~(n\geq r)\}$. Noticing that $\bbb{y}=(\bbb{x}+\bbb{1})/2 \in \{0,1\}^n$, we can extend $\ell_{\infty}$ box non-separable MPEC\footnote{$\Delta\Leftrightarrow \{\bbb{x}~|~\bbb{0}\leq \bbb{x}\leq \bbb{1},~\bbb{0}\leq\bbb{v}\leq \bbb{1},~\la 2\bbb{x}-\bbb{1},2\bbb{v}-\bbb{1} \ra = n,~\forall \bbb{v} \}$} and $\ell_2$ box non-separable MPEC\footnote{$\Delta\Leftrightarrow\{\bbb{x}~|~\bbb{0}\leq \bbb{x}\leq \bbb{1},~\|2\bbb{v}-1\|_2^2\leq n,~\la 2\bbb{x}-\bbb{1},2\bbb{v}-\bbb{1} \ra = n,~\forall \bbb{v}\}$} to handle zero-one binary optimization. Moreover, observing that binary constraint $\bbb{x}\in\{-1,+1\}^n \Leftrightarrow |\bbb{x}|=\bbb{1}$ is analogous to orthogonality constraint since $\bbb{X} \in \mathbb{O} \Leftrightarrow \sigma(\bbb{X})=\bbb{1}$, we can extend our $\ell_{\infty}$ box non-separable MPEC\footnote{$\mathbb{O} \Leftrightarrow \{\bbb{X}~|~\bbb{X}^T\bbb{X} \preceq\bbb{I},~\bbb{V}^T\bbb{V} \preceq\bbb{I},~ \la \bbb{X},\bbb{V}\ra = r,~\forall \bbb{V}\}$} and $\ell_2$ box non-separable MPEC\footnote{$\mathbb{O} \Leftrightarrow \{\bbb{X}~|~\bbb{X}^T\bbb{X} \preceq\bbb{I}_r,~\|\bbb{V}\|_F^2\leq r,~\la \bbb{X},\bbb{V}\ra = r,~\forall \bbb{V}\}$} to the optimization problem with orthogonality constraint \cite{wen2013feasible,Chen2016}. 


\begin{figure*}[!t]
\captionsetup[subfigure]{justification=centering}
    \centering
      \begin{subfigure}{\imgwidseg}\includegraphics[width=\textwidth,height=\imgheiseg]{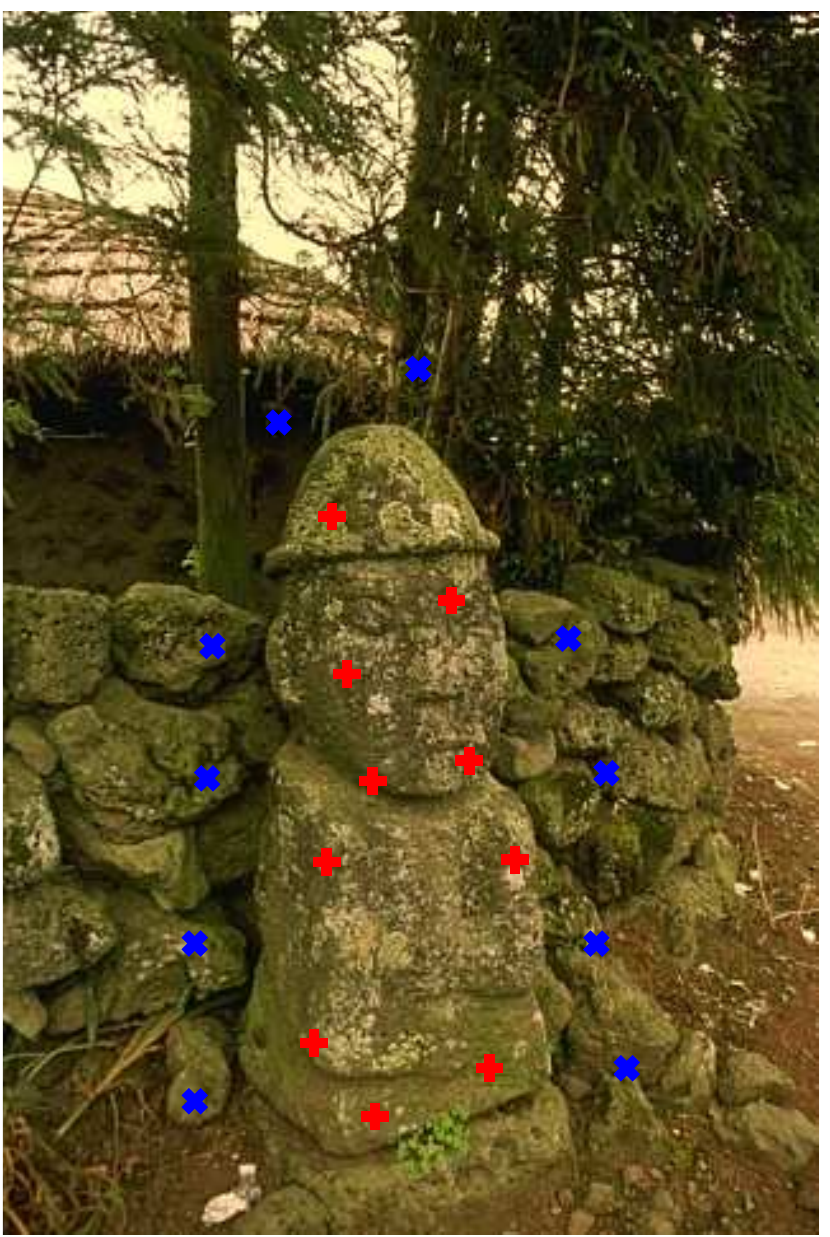}\end{subfigure}\ghst
      \begin{subfigure}{\imgwidseg}\includegraphics[width=\textwidth,height=\imgheiseg]{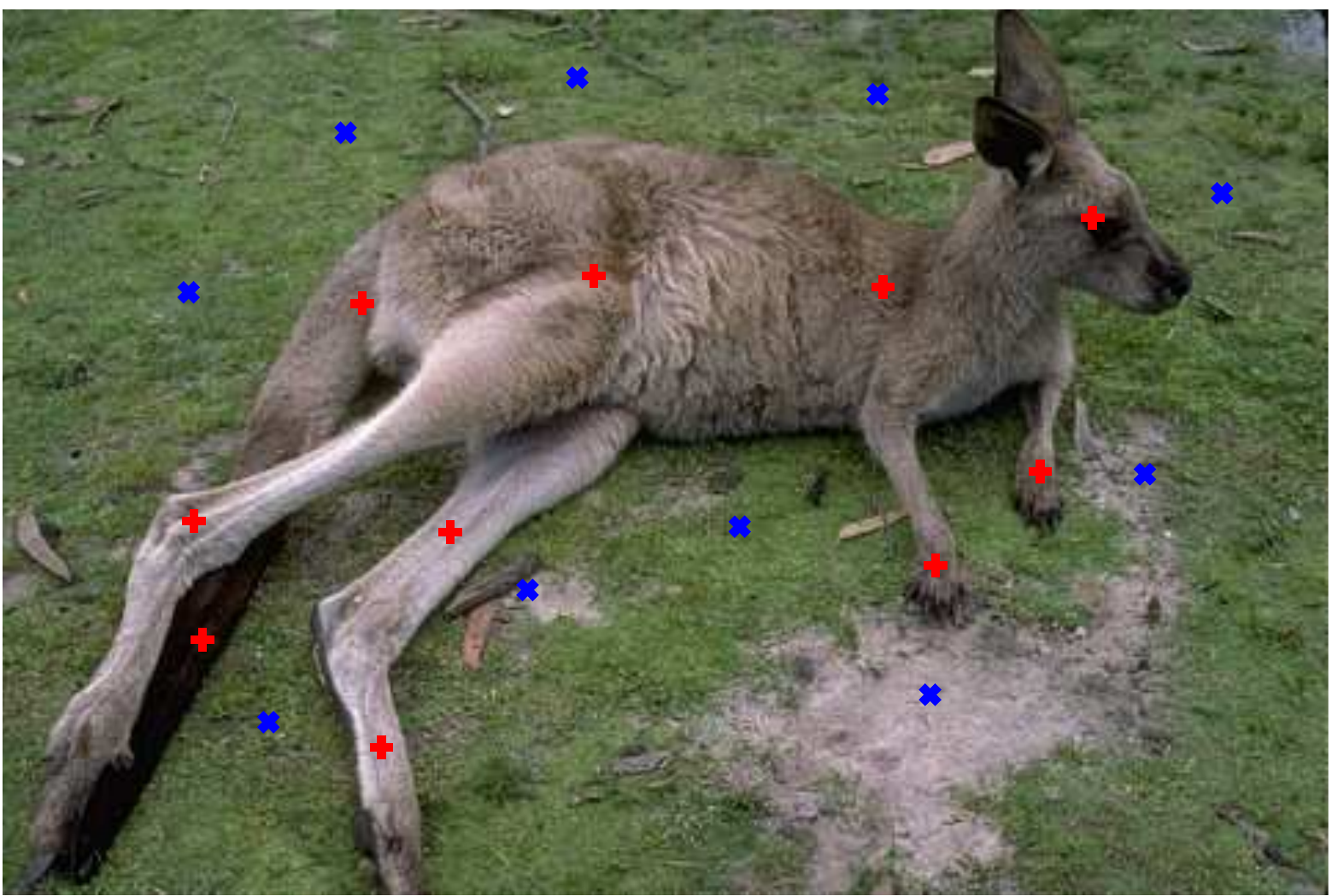}\end{subfigure}\ghst
      \begin{subfigure}{\imgwidseg}\includegraphics[width=\textwidth,height=\imgheiseg]{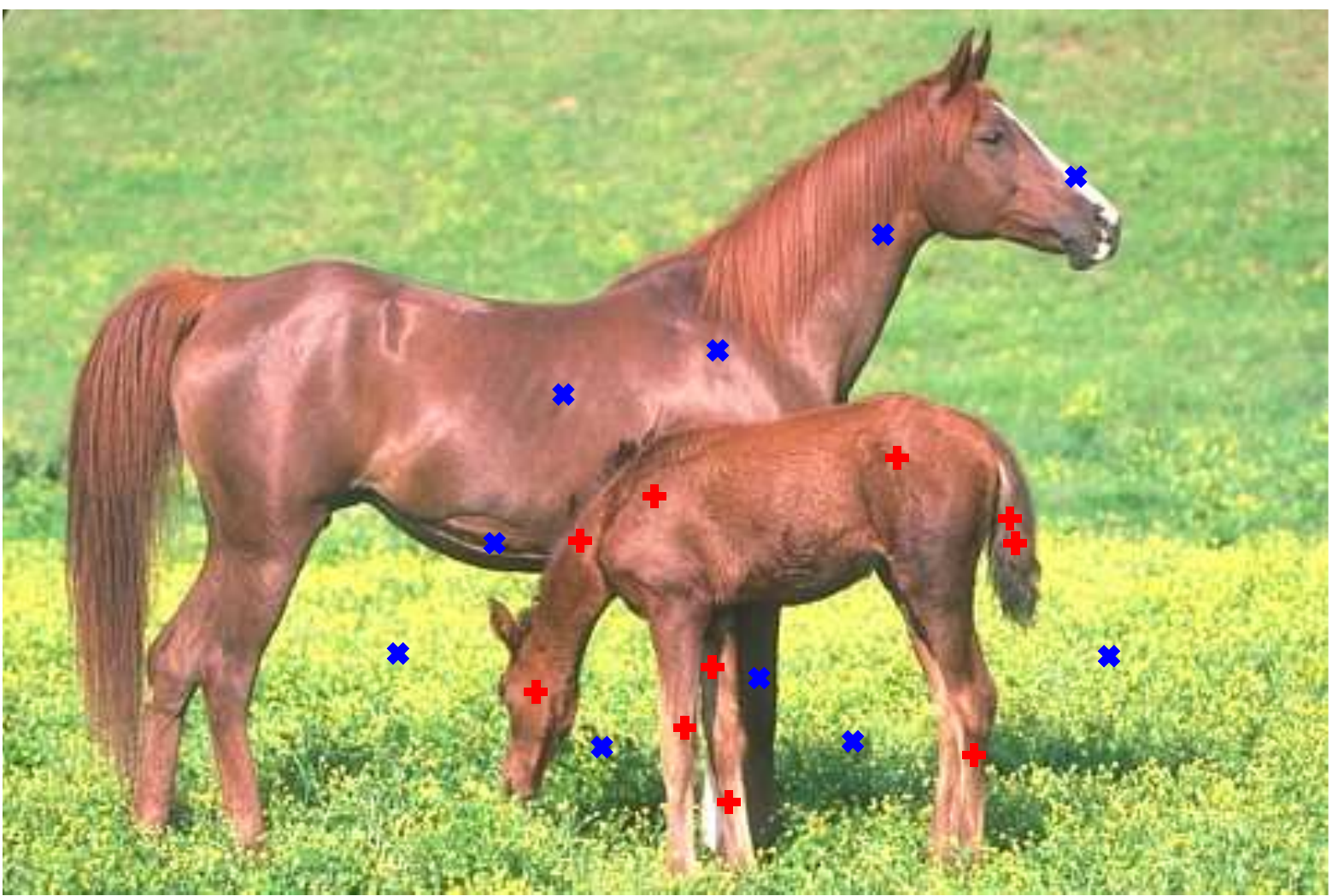}\end{subfigure}\ghst
      \begin{subfigure}{\imgwidseg}\includegraphics[width=\textwidth,height=\imgheiseg]{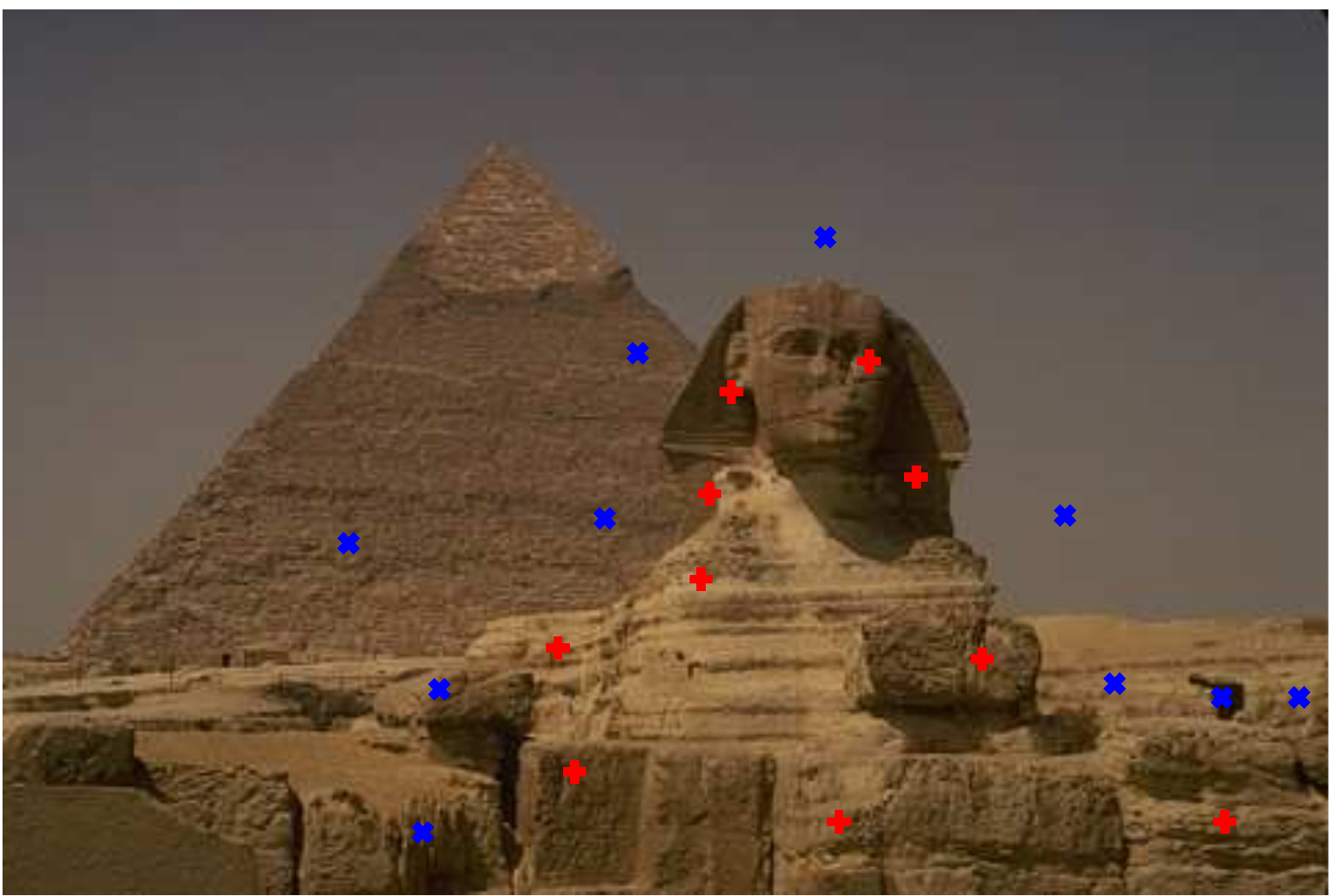}\end{subfigure}\ghst
      \begin{subfigure}{\imgwidseg}\includegraphics[width=\textwidth,height=\imgheiseg]{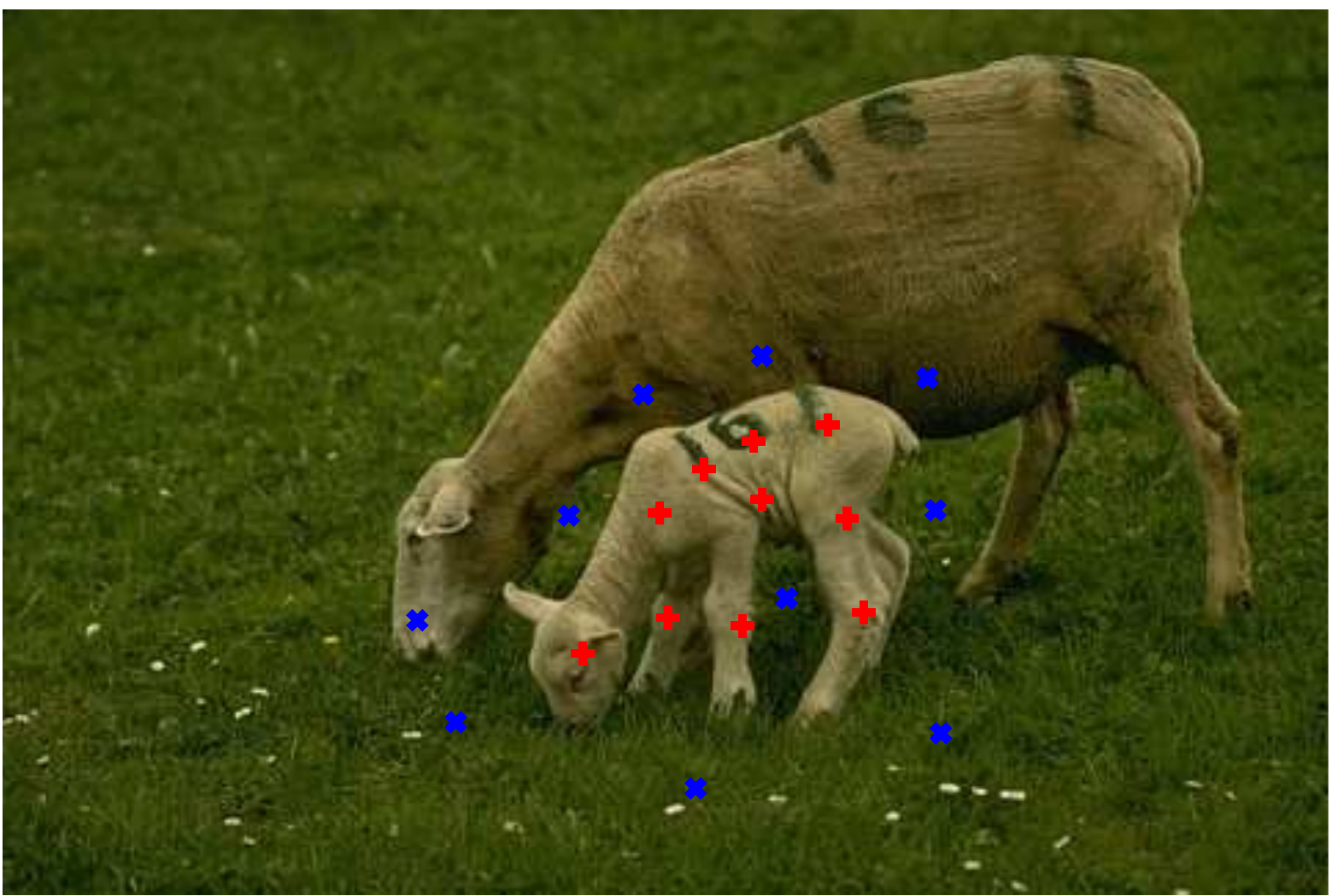}\end{subfigure}
\caption{Images in our constrained image segmentation experiments. 10 foreground pixels and 10 background pixels are annotated by red and blue markers respectively.}
\label{fig:cseg1}
\end{figure*}

\begin{figure*}[!th]
\captionsetup[subfigure]{justification=centering,labelfont=scriptsize}

    \centering
      \begin{subfigure}{\imgwidseg}\includegraphics[width=\textwidth,height=\imgheiseg]{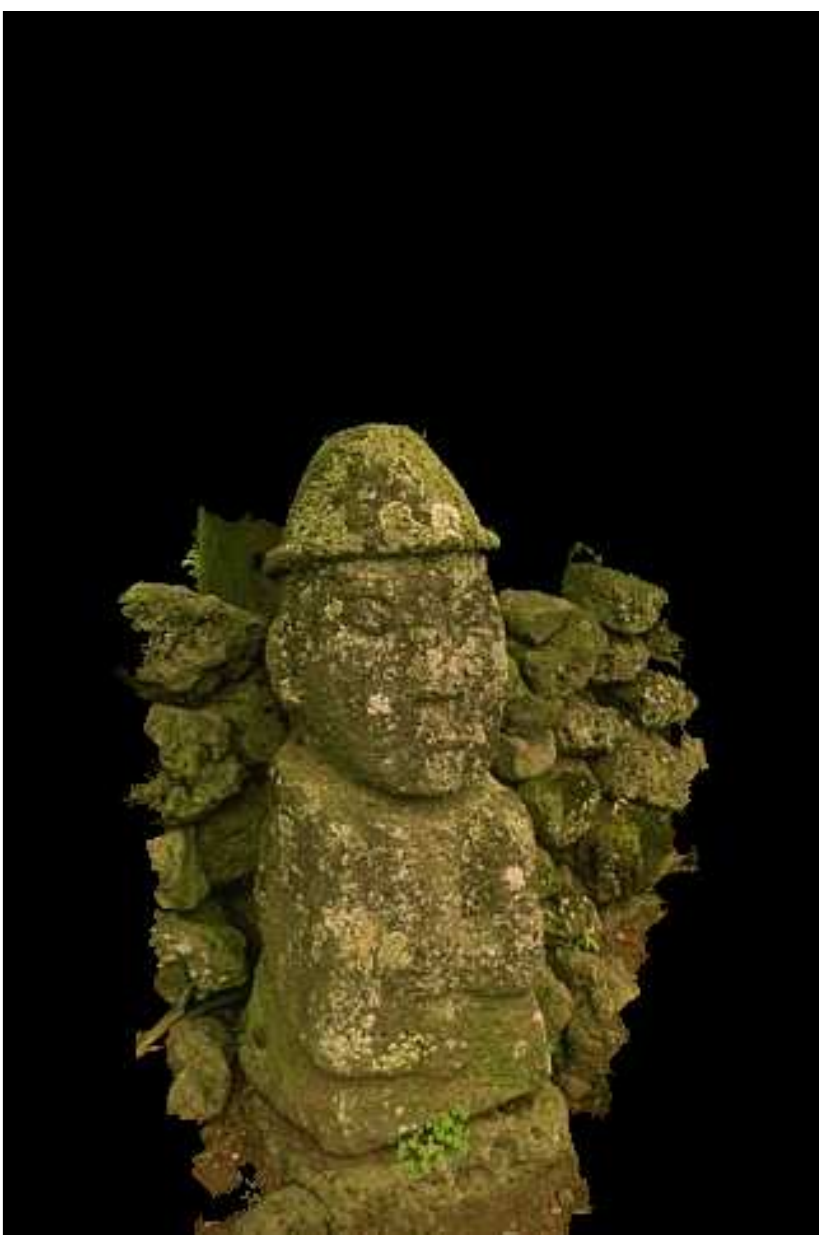}\vspace{-2pt}\caption{{\scriptsize BNCUT,~$f = 89.09$}}\end{subfigure}\ghst
      \begin{subfigure}{\imgwidseg}\includegraphics[width=\textwidth,height=\imgheiseg]{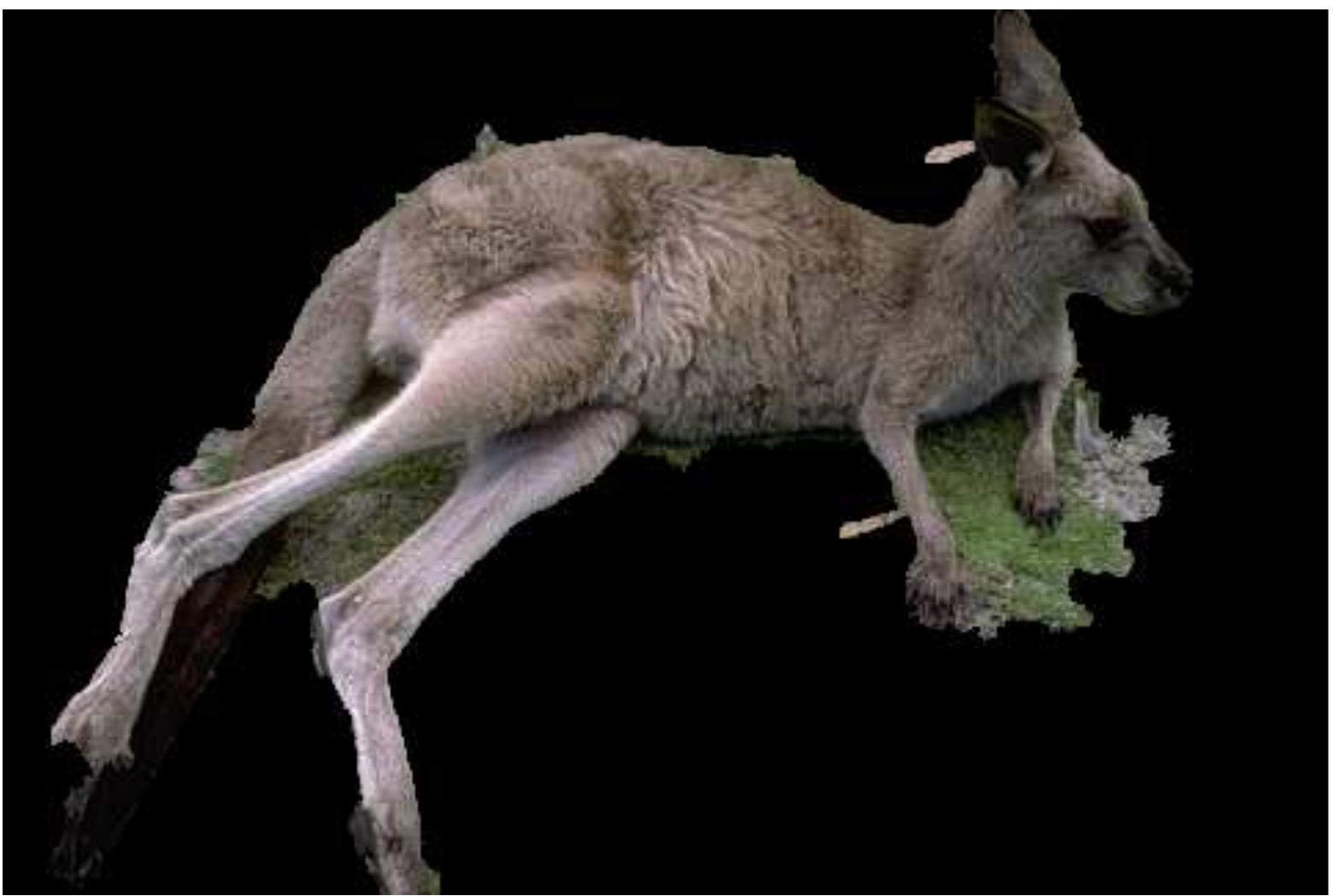}\vspace{-2pt}\caption{{\scriptsize BNCUT,~$f = 23.08$}}\end{subfigure}\ghst
      \begin{subfigure}{\imgwidseg}\includegraphics[width=\textwidth,height=\imgheiseg]{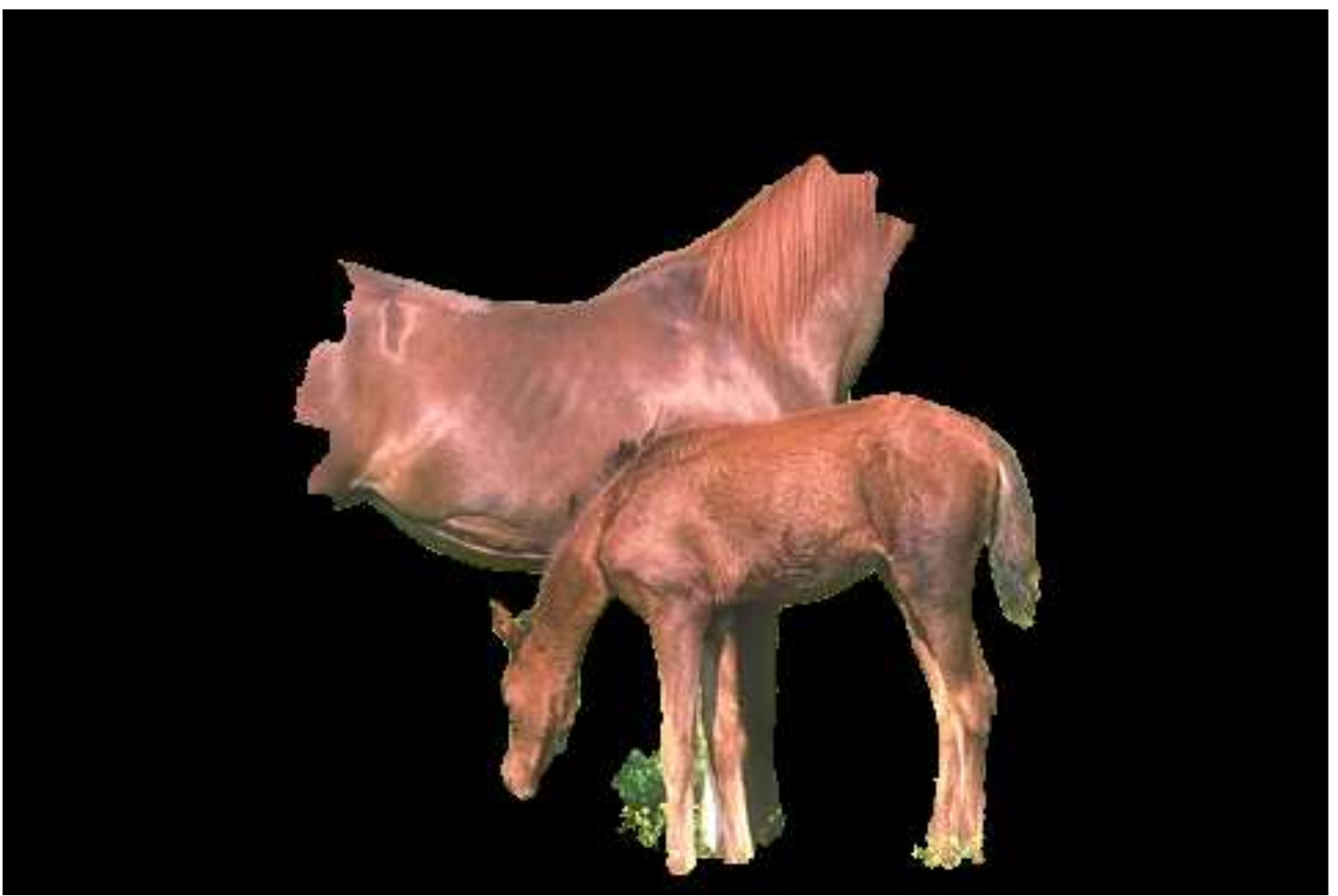}\vspace{-2pt}\caption{{\scriptsize BNCUT,~$f = 26.67 $}}\end{subfigure}\ghst
      \begin{subfigure}{\imgwidseg}\includegraphics[width=\textwidth,height=\imgheiseg]{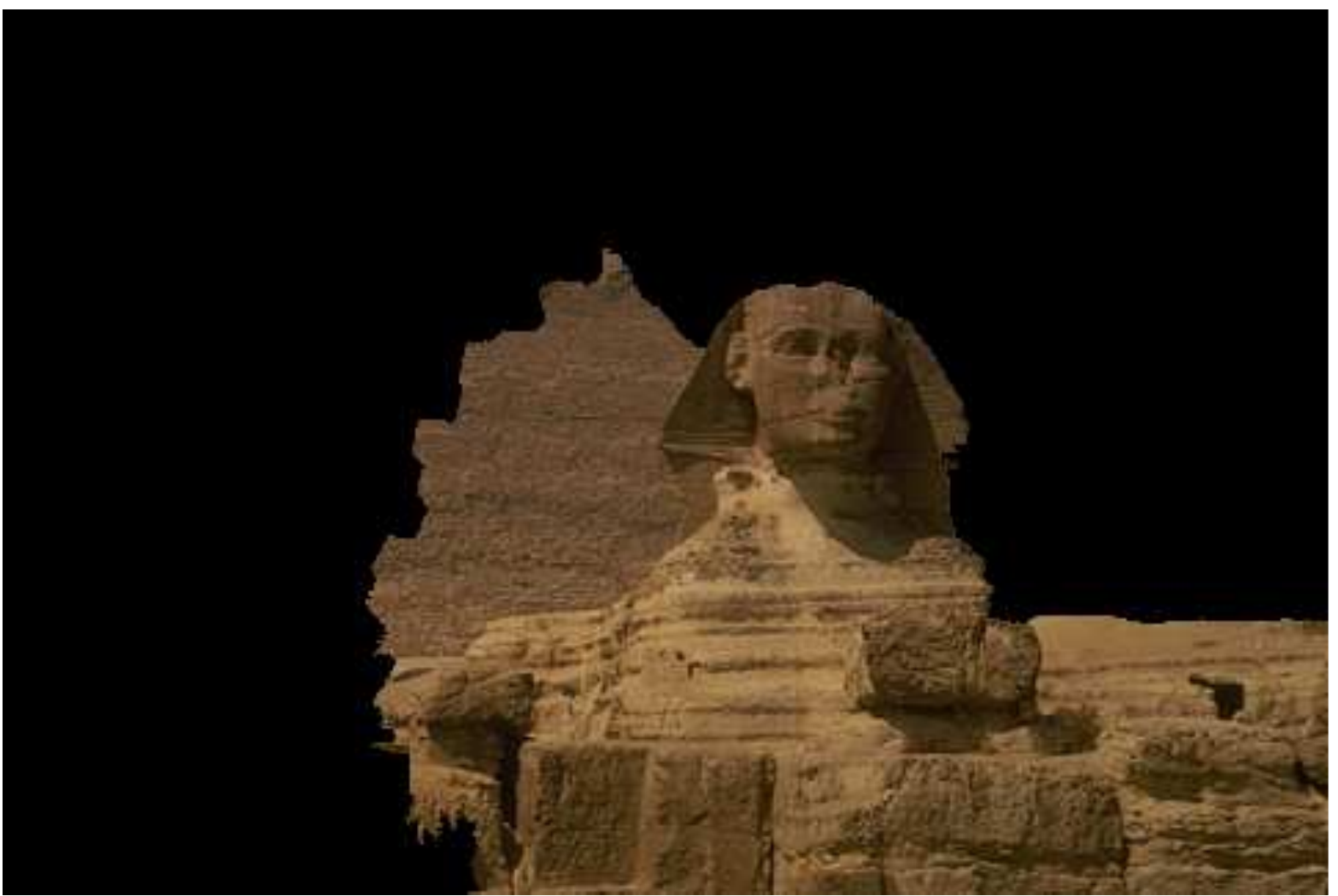}\vspace{-2pt}\caption{{\scriptsize BNCUT,~$f = 56.02$}}\end{subfigure}\ghst
      \begin{subfigure}{\imgwidseg}\includegraphics[width=\textwidth,height=\imgheiseg]{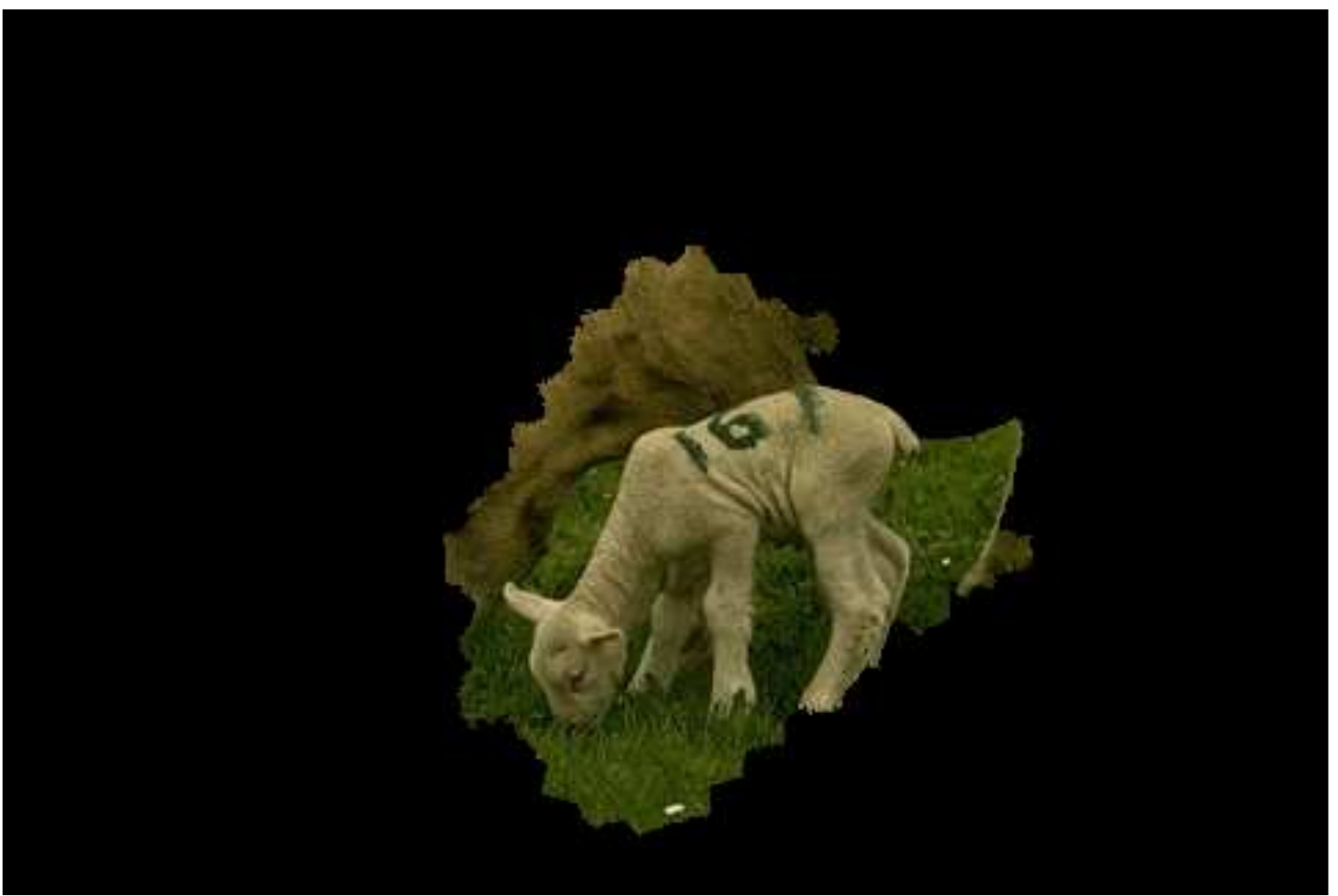}\vspace{-2pt}\caption{{\scriptsize BNCUT,~$f = 75.09 $}}\end{subfigure}

\vspace{5pt}

      \begin{subfigure}{\imgwidseg}\includegraphics[width=\textwidth,height=\imgheiseg]{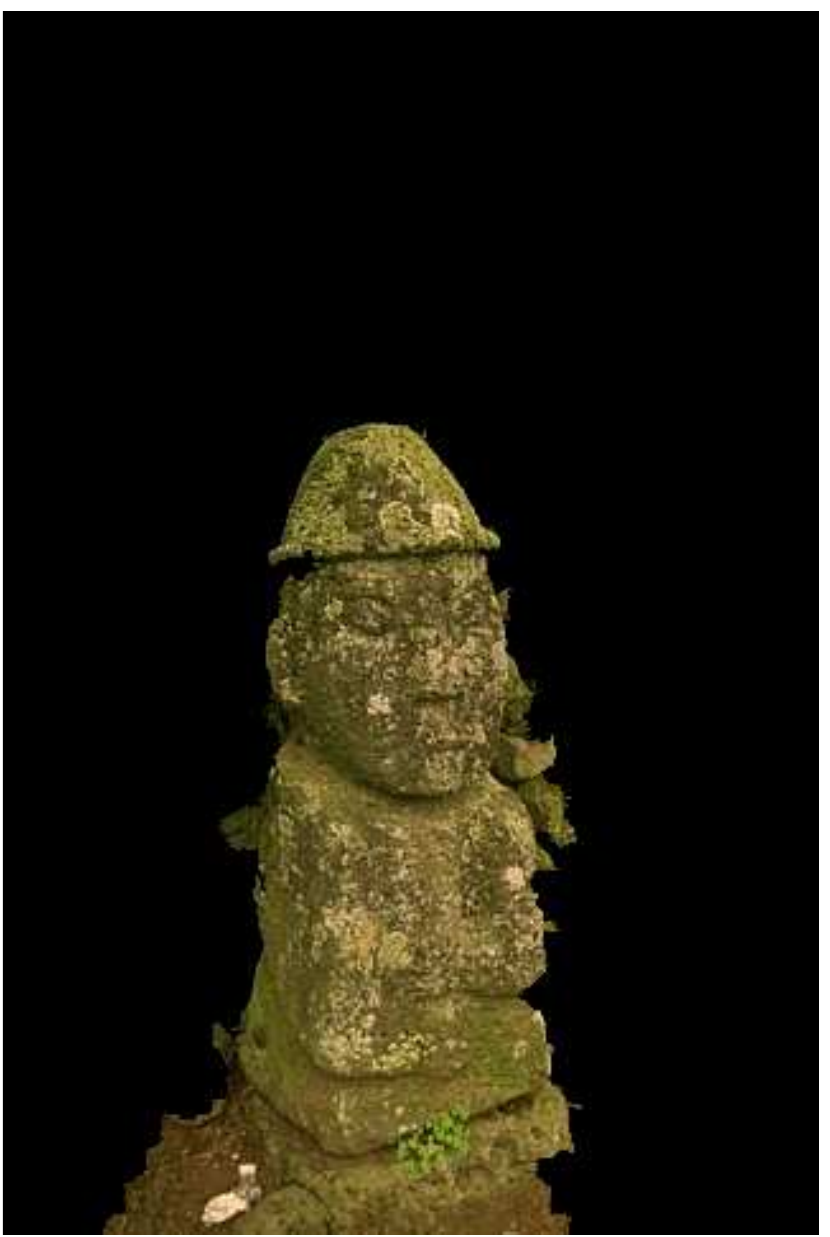}\vspace{-2pt}\caption{{\scriptsize LP,~$f = 65.70$}}\end{subfigure}\ghst
      \begin{subfigure}{\imgwidseg}\includegraphics[width=\textwidth,height=\imgheiseg]{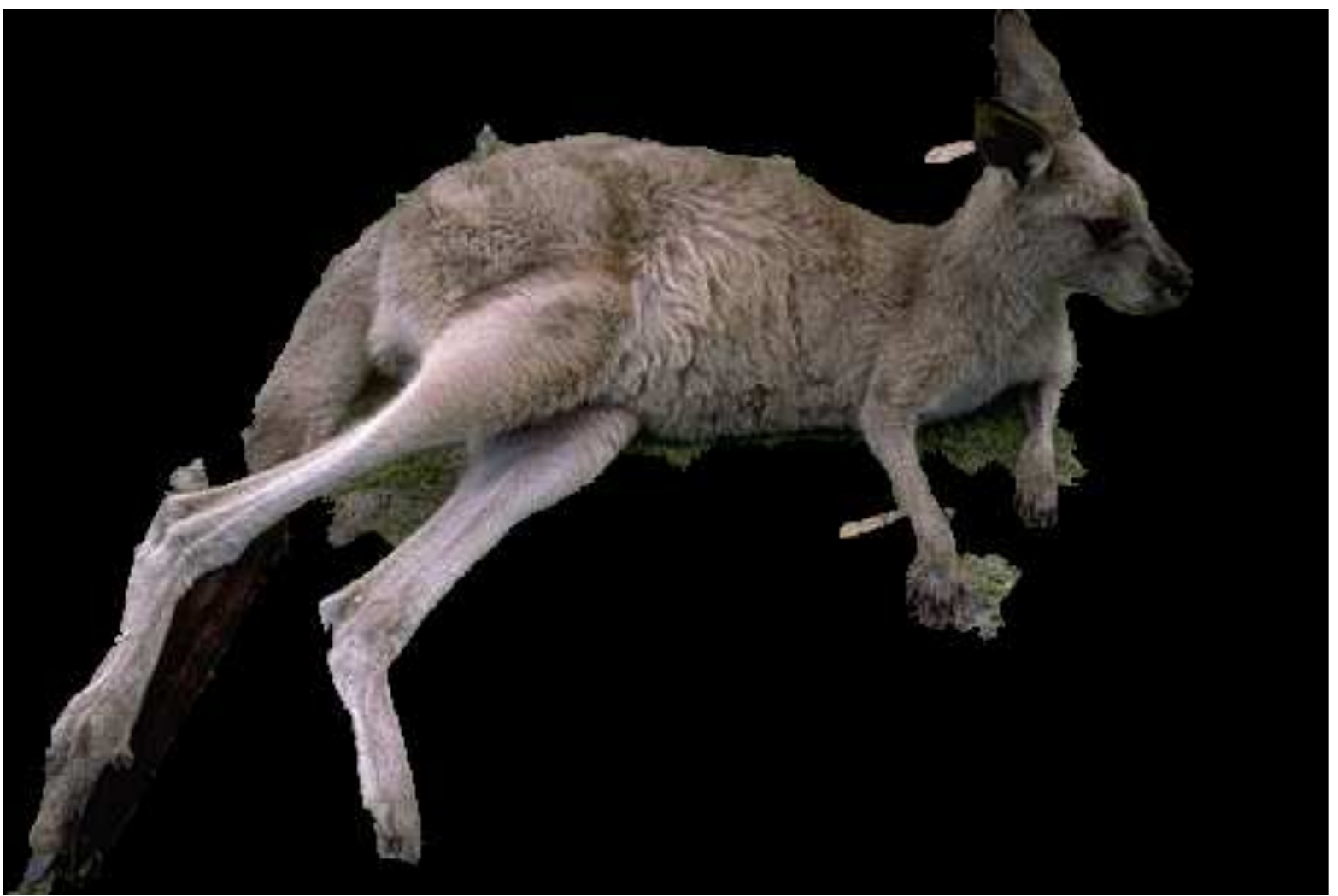}\vspace{-2pt}\caption{{\scriptsize LP,~$f = 27.53$}}\end{subfigure}\ghst
      \begin{subfigure}{\imgwidseg}\includegraphics[width=\textwidth,height=\imgheiseg]{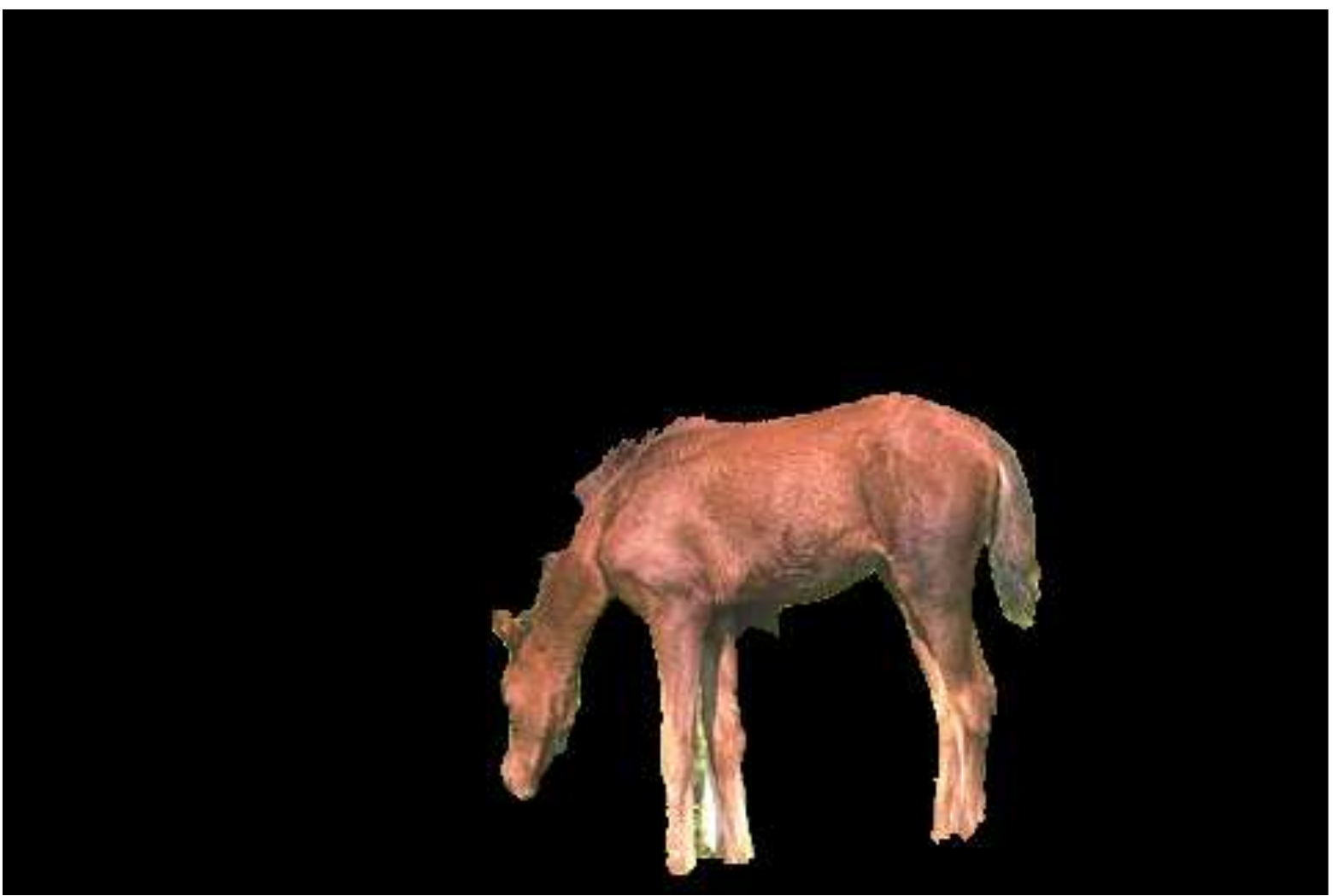}\vspace{-2pt}\caption{{\scriptsize LP,~$f = 14.81$}}\end{subfigure}\ghst
      \begin{subfigure}{\imgwidseg}\includegraphics[width=\textwidth,height=\imgheiseg]{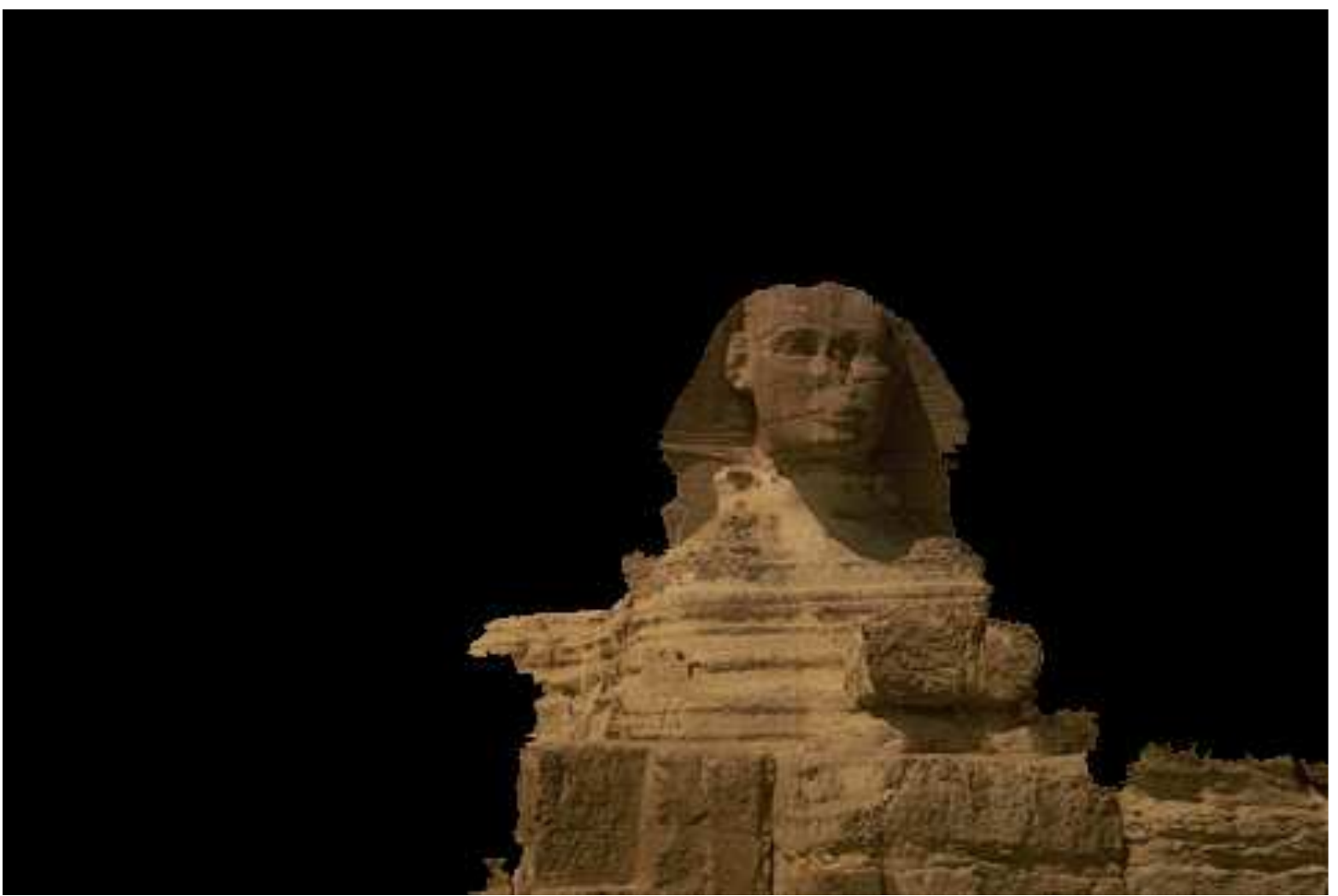}\vspace{-2pt}\caption{{\scriptsize LP,~$f = 27.50$}}\end{subfigure}\ghst
      \begin{subfigure}{\imgwidseg}\includegraphics[width=\textwidth,height=\imgheiseg]{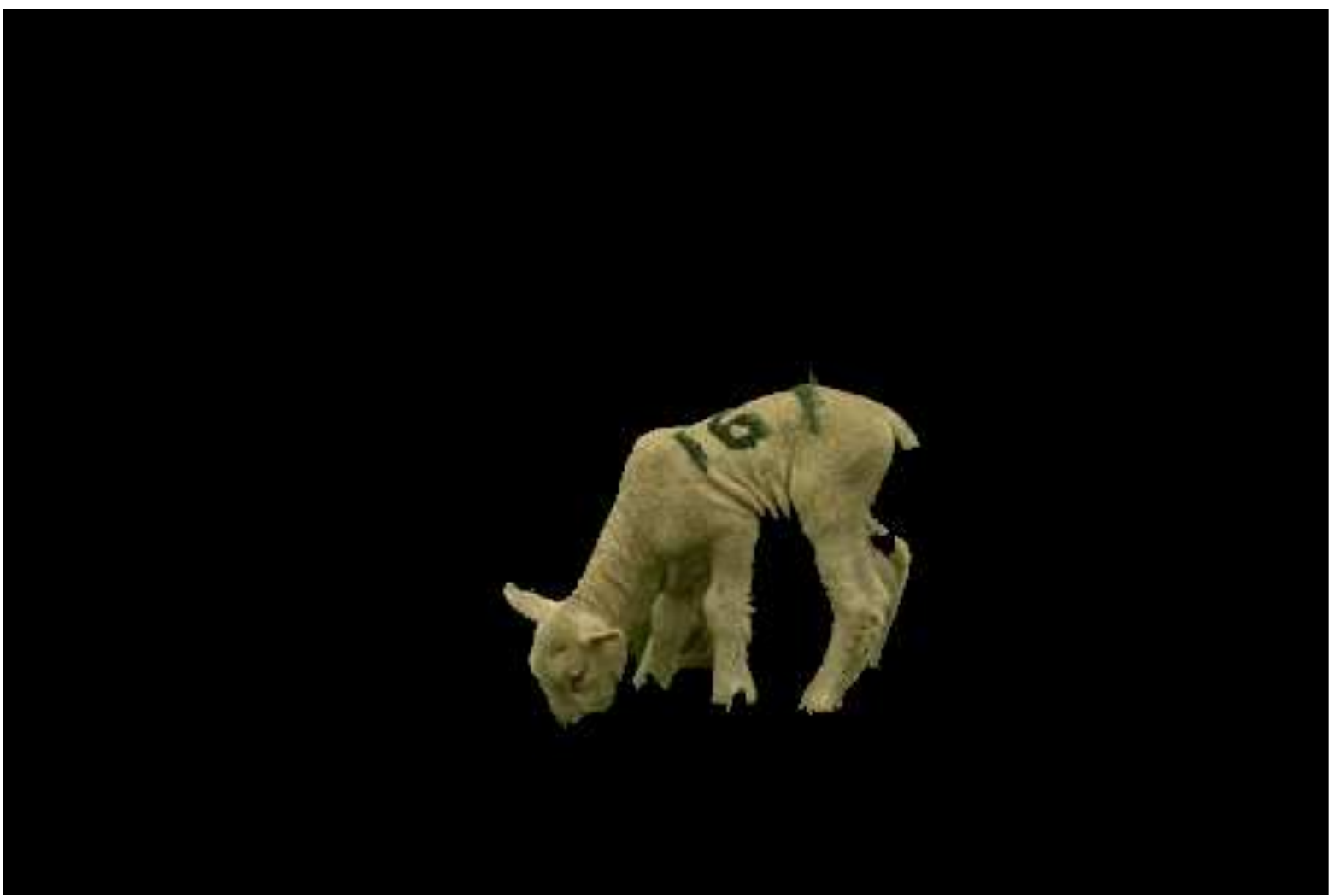}\vspace{-2pt}\caption{{\scriptsize LP,~$f = 6.81$}}\end{subfigure}

\vspace{5pt}

      \begin{subfigure}{\imgwidseg}\includegraphics[width=\textwidth,height=\imgheiseg]{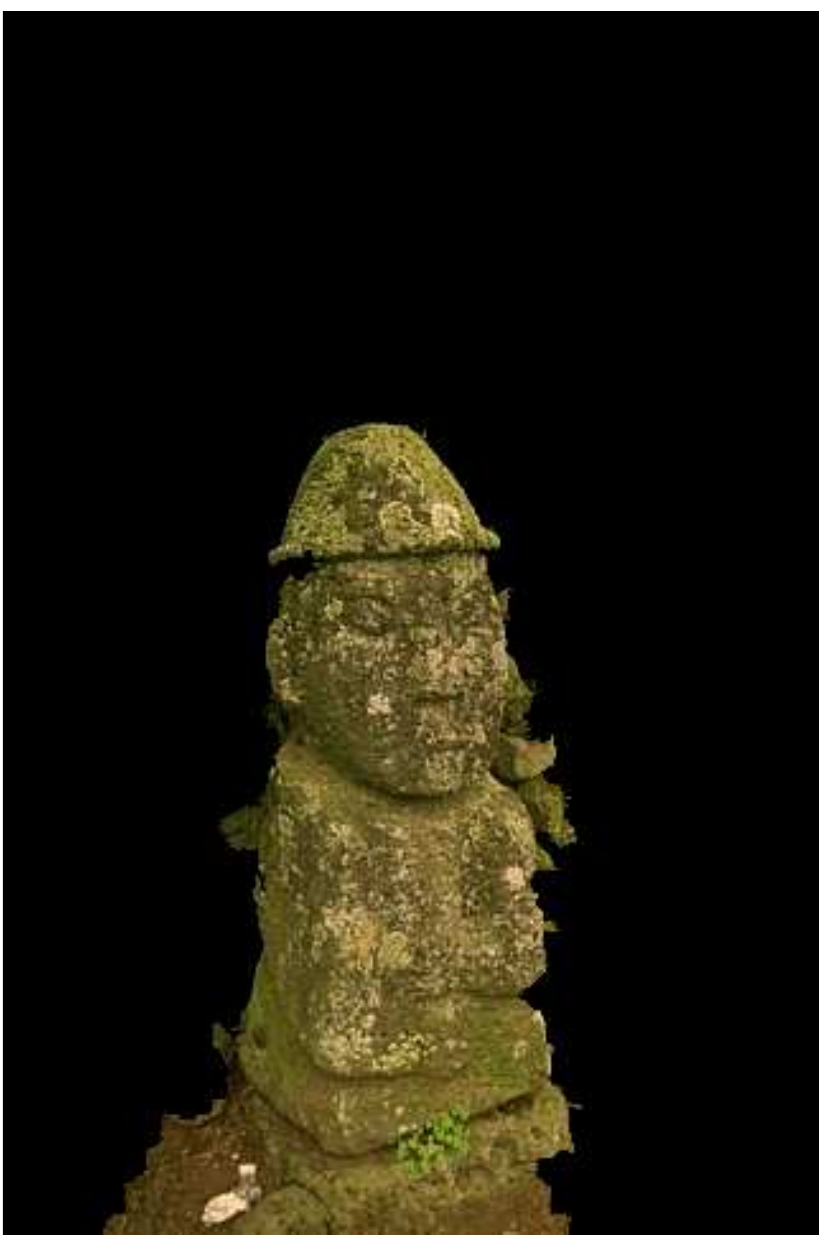}\vspace{-2pt}\caption{{\scriptsize SDP,~$f = 64.40$}}\end{subfigure}\ghst
      \begin{subfigure}{\imgwidseg}\includegraphics[width=\textwidth,height=\imgheiseg]{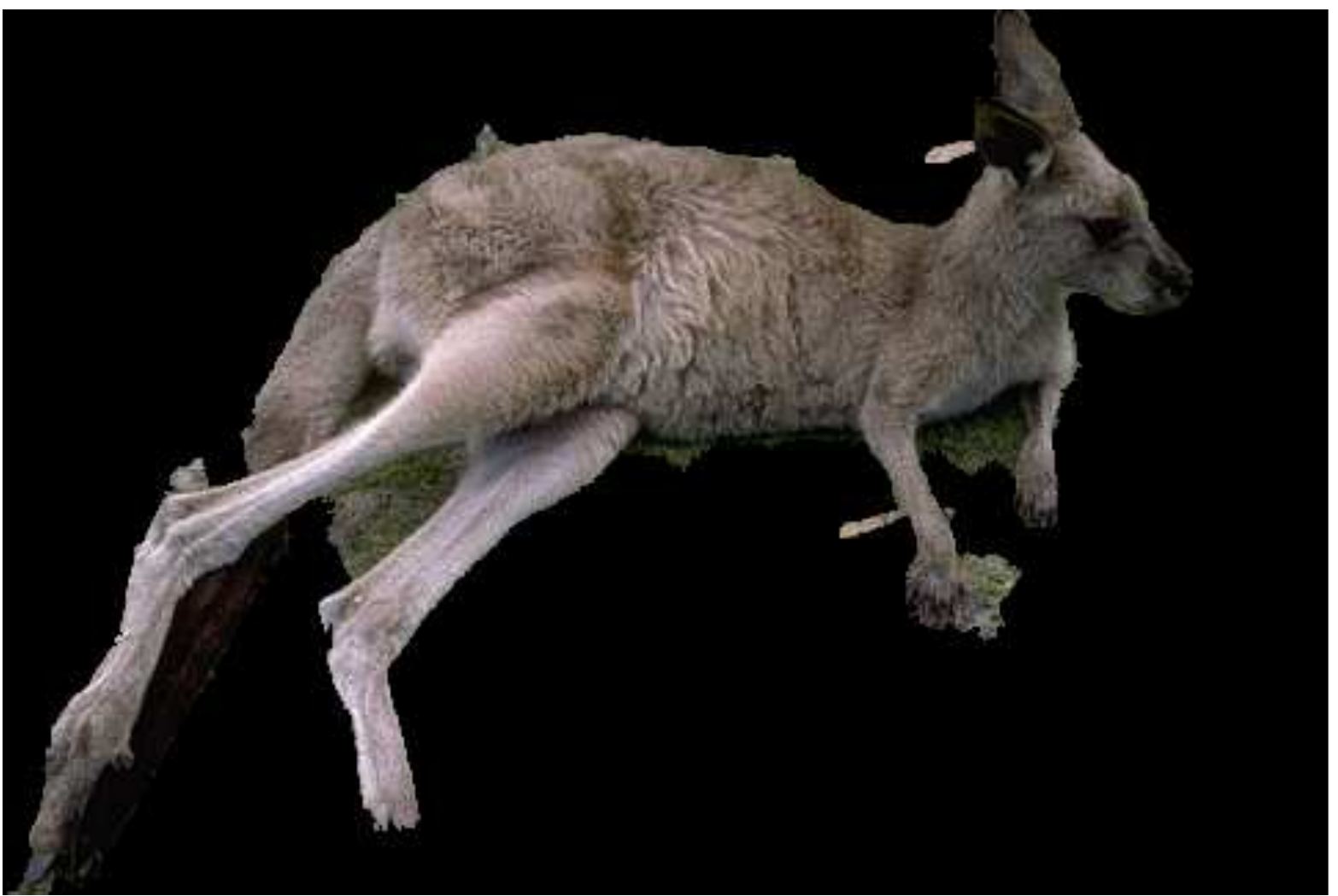}\vspace{-2pt}\caption{{\scriptsize SDP,~$f =19.01 $}}\end{subfigure}\ghst
      \begin{subfigure}{\imgwidseg}\includegraphics[width=\textwidth,height=\imgheiseg]{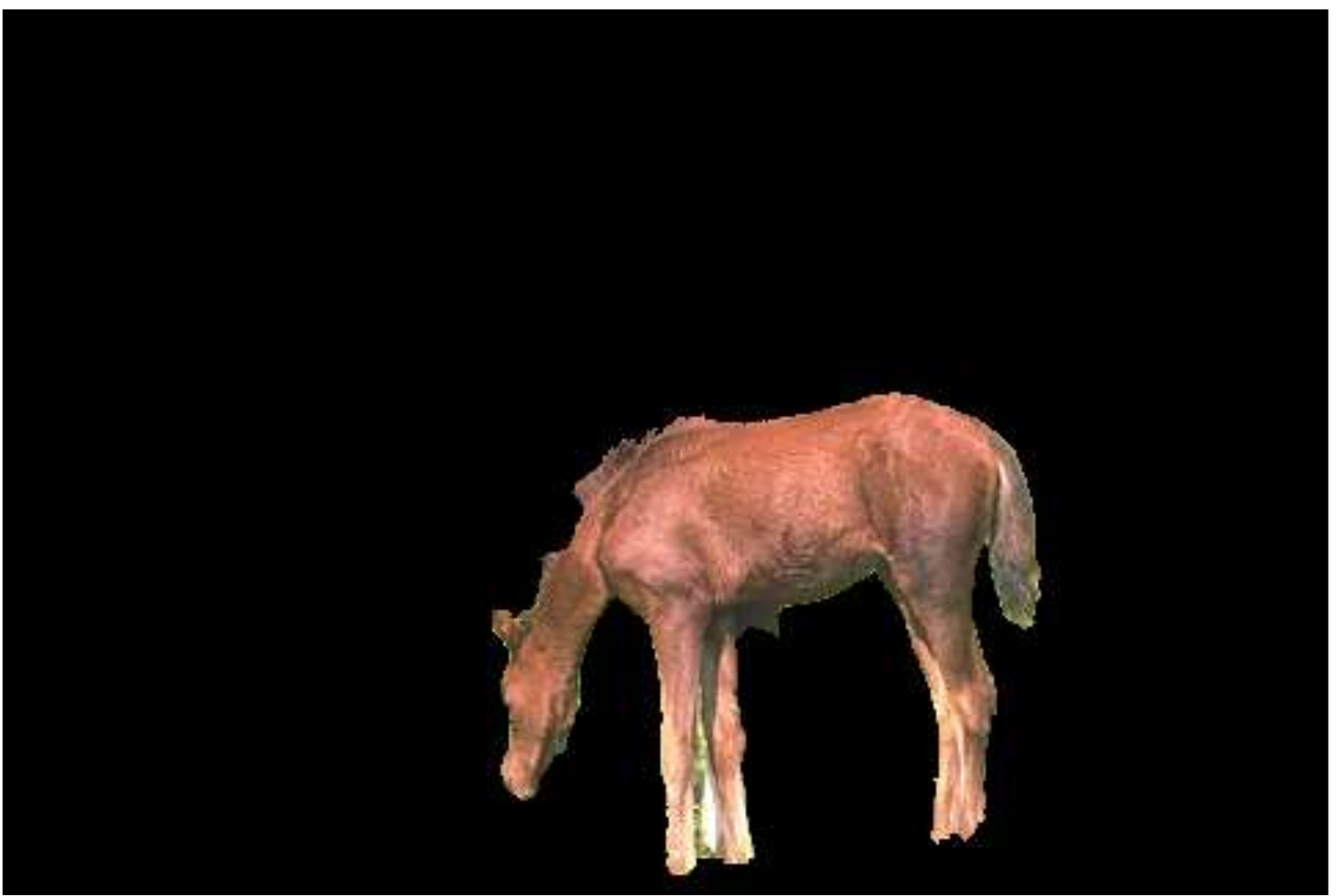}\vspace{-2pt}\caption{{\scriptsize SDP,~$f = 14.81$}}\end{subfigure}\ghst
      \begin{subfigure}{\imgwidseg}\includegraphics[width=\textwidth,height=\imgheiseg]{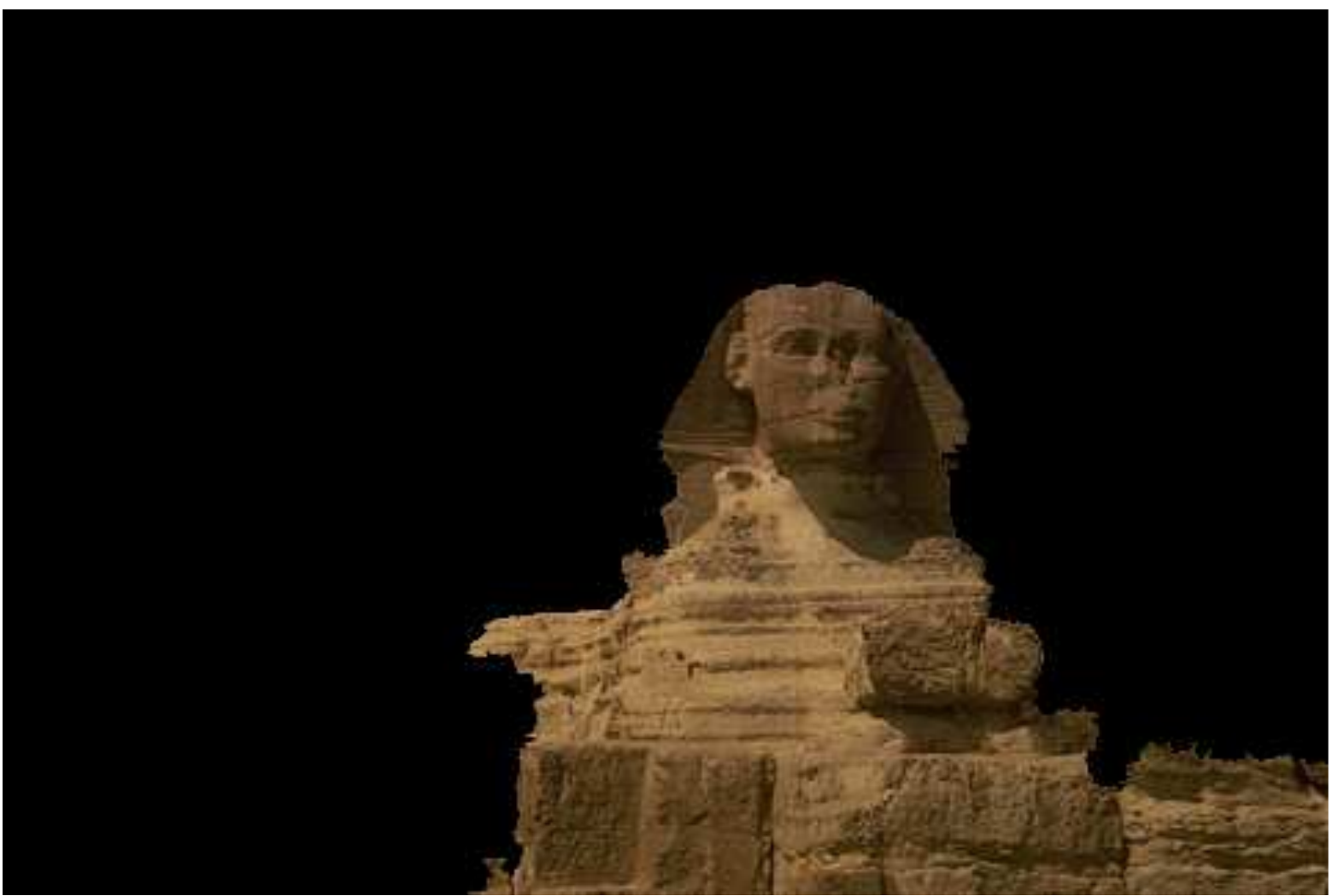}\vspace{-2pt}\caption{{\scriptsize SDP,~$f = 26.97$}}\end{subfigure}\ghst
      \begin{subfigure}{\imgwidseg}\includegraphics[width=\textwidth,height=\imgheiseg]{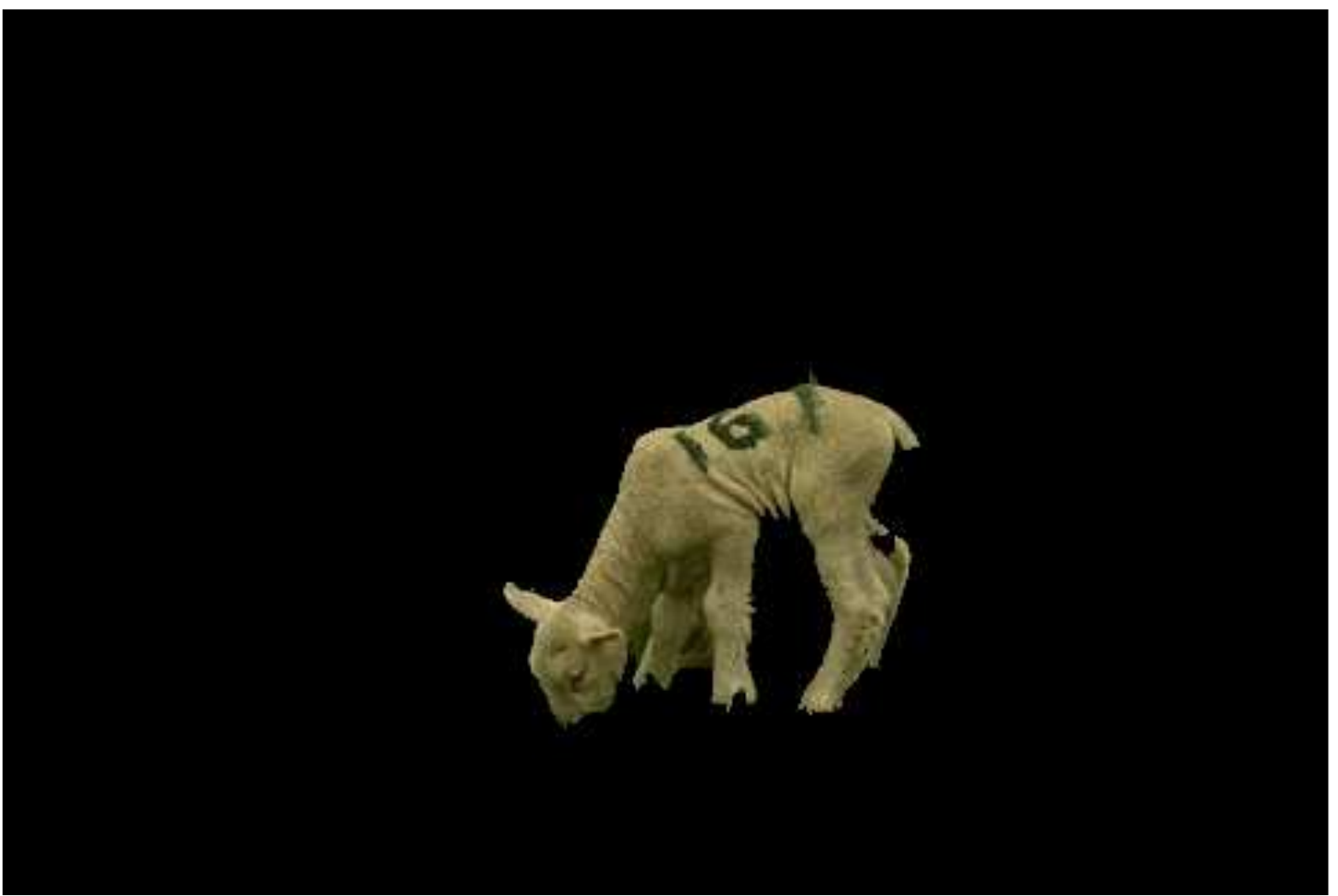}\vspace{-2pt}\caption{{\scriptsize SDP,~$f = 6.81$}}\end{subfigure}

\vspace{5pt}

      \begin{subfigure}{\imgwidseg}\includegraphics[width=\textwidth,height=\imgheiseg]{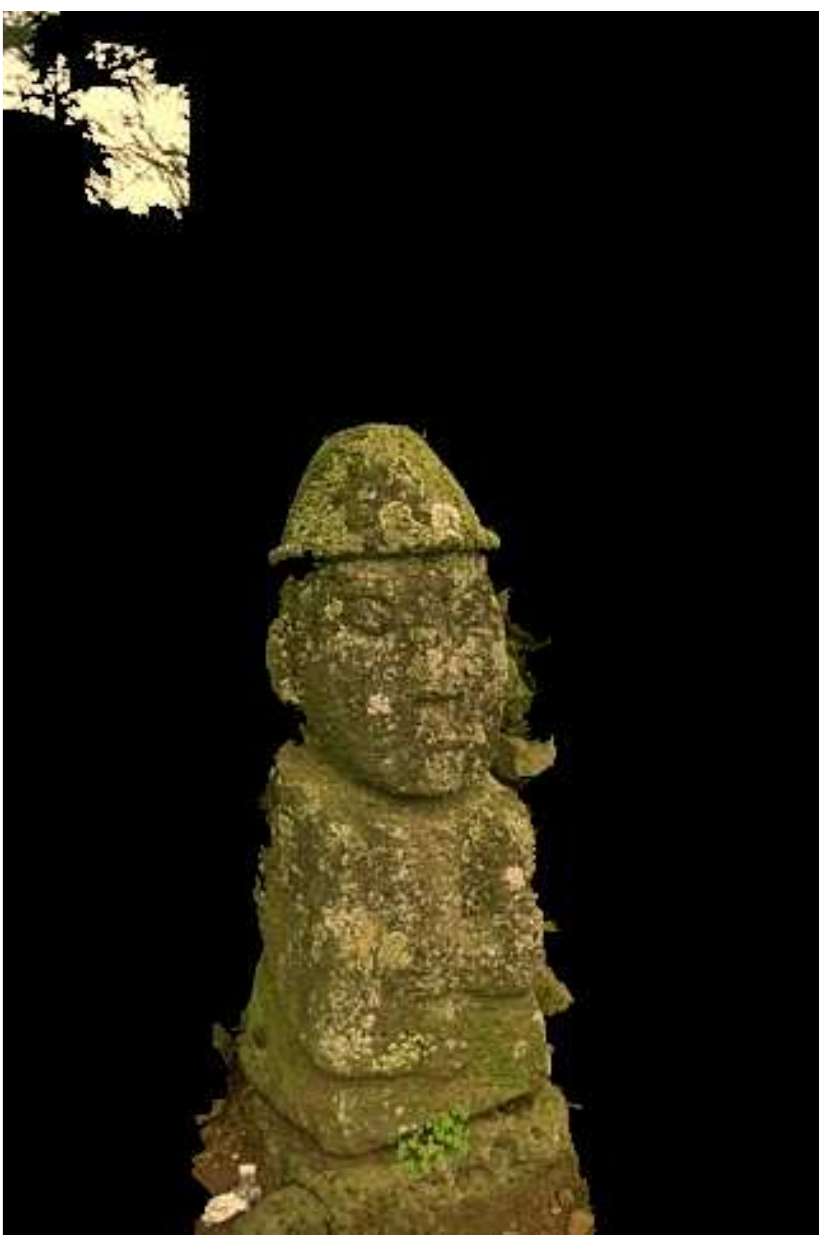}\vspace{-2pt}\caption{{\scriptsize L2box-ADMM,~$f = 56.65$}}\end{subfigure}\ghst
      \begin{subfigure}{\imgwidseg}\includegraphics[width=\textwidth,height=\imgheiseg]{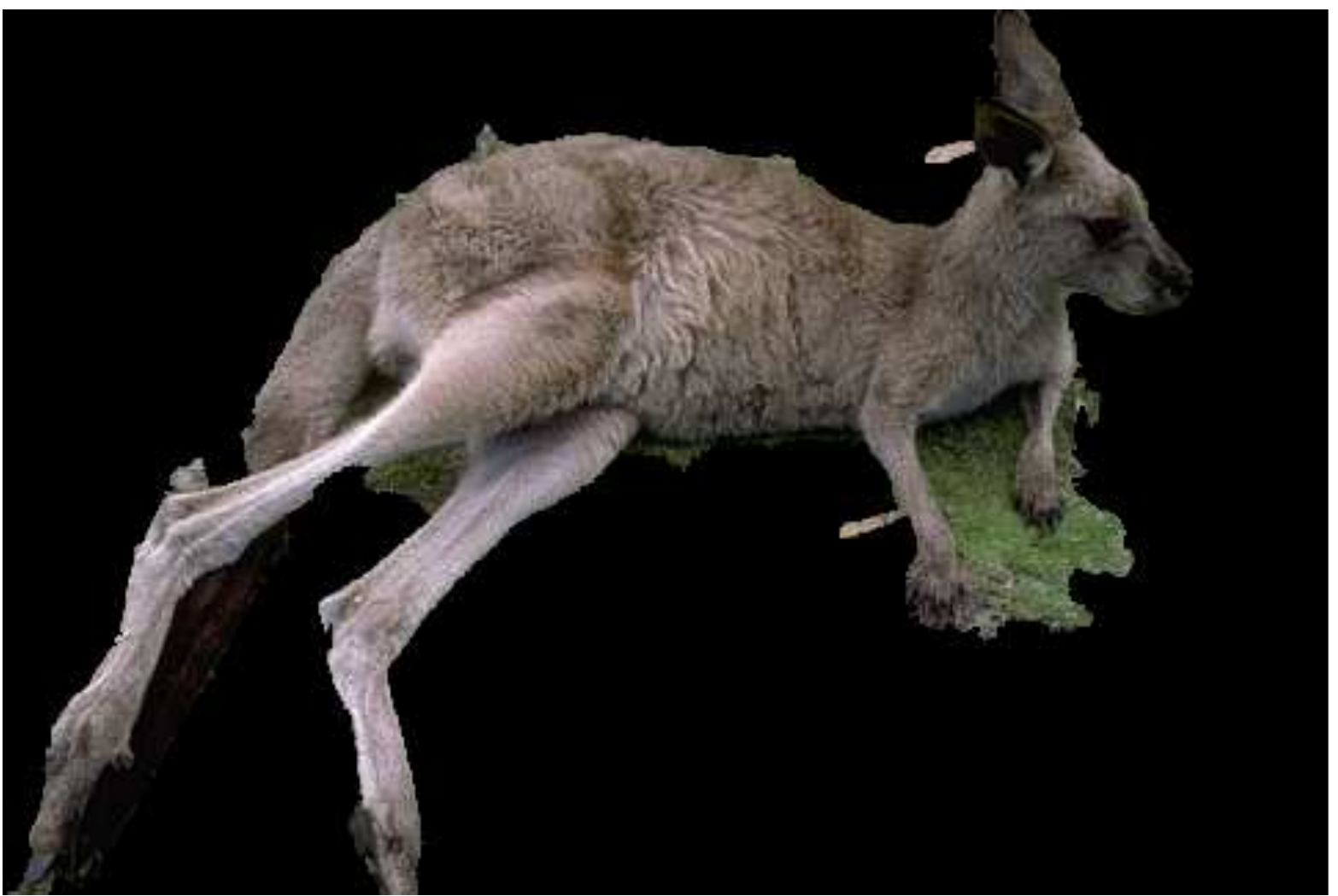}\vspace{-2pt}\caption{{\scriptsize L2box-ADMM,~$f = 19.70$}}\end{subfigure}\ghst
      \begin{subfigure}{\imgwidseg}\includegraphics[width=\textwidth,height=\imgheiseg]{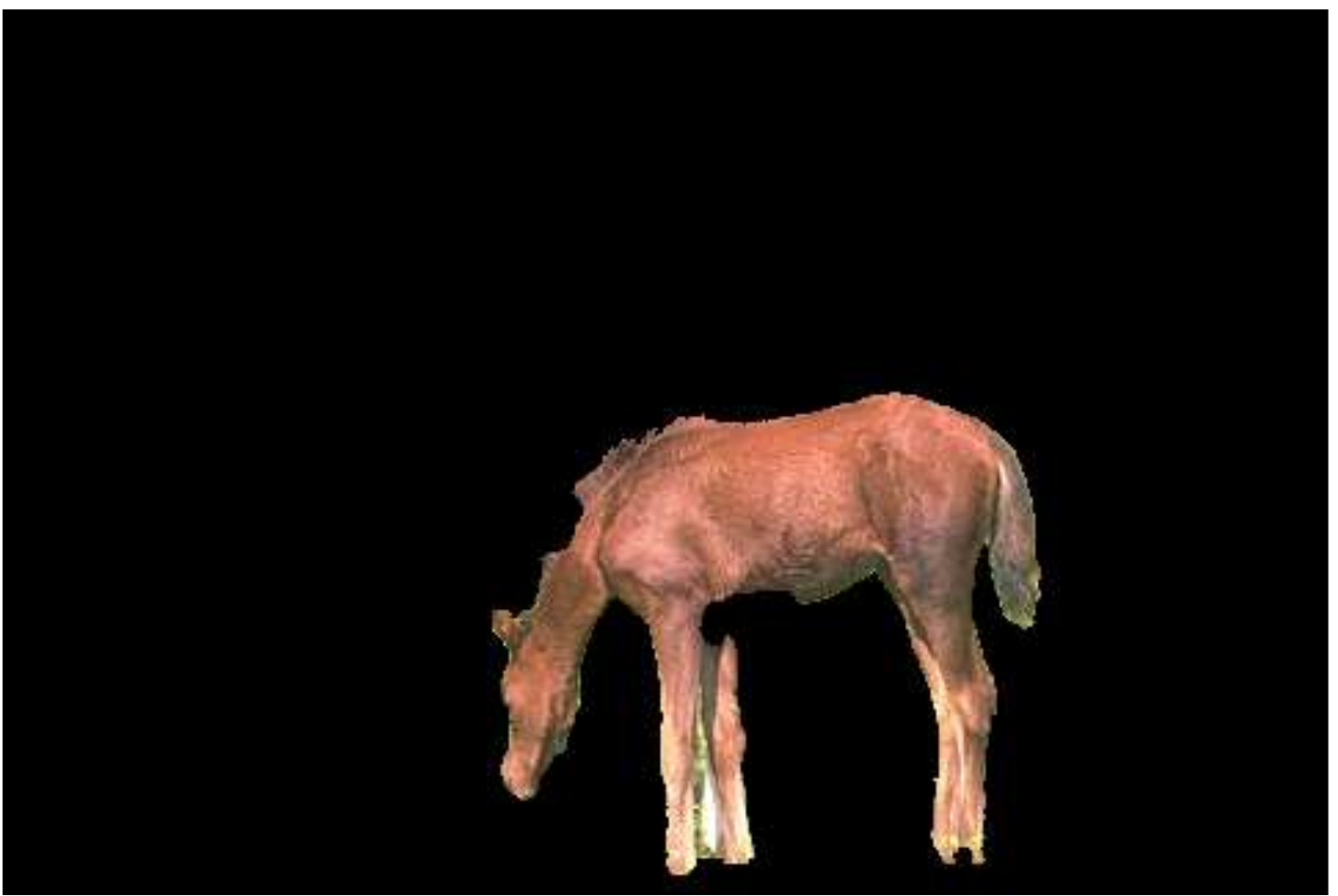}\vspace{-2pt}\caption{{\scriptsize L2box-ADMM,~$f = 14.32$}}\end{subfigure}\ghst
      \begin{subfigure}{\imgwidseg}\includegraphics[width=\textwidth,height=\imgheiseg]{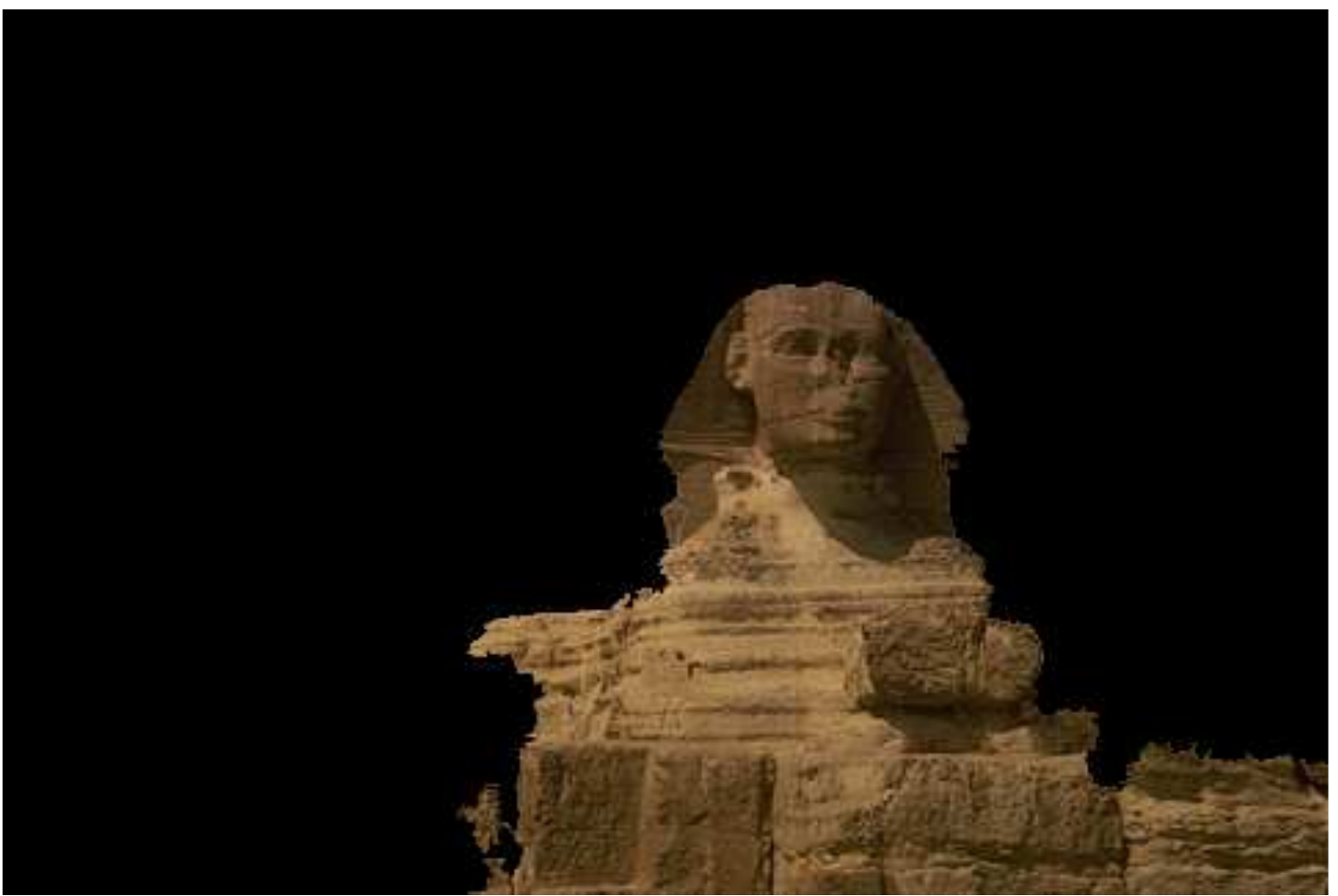}\vspace{-2pt}\caption{{\scriptsize L2box-ADMM,~$f = 25.77$}}\end{subfigure}\ghst
      \begin{subfigure}{\imgwidseg}\includegraphics[width=\textwidth,height=\imgheiseg]{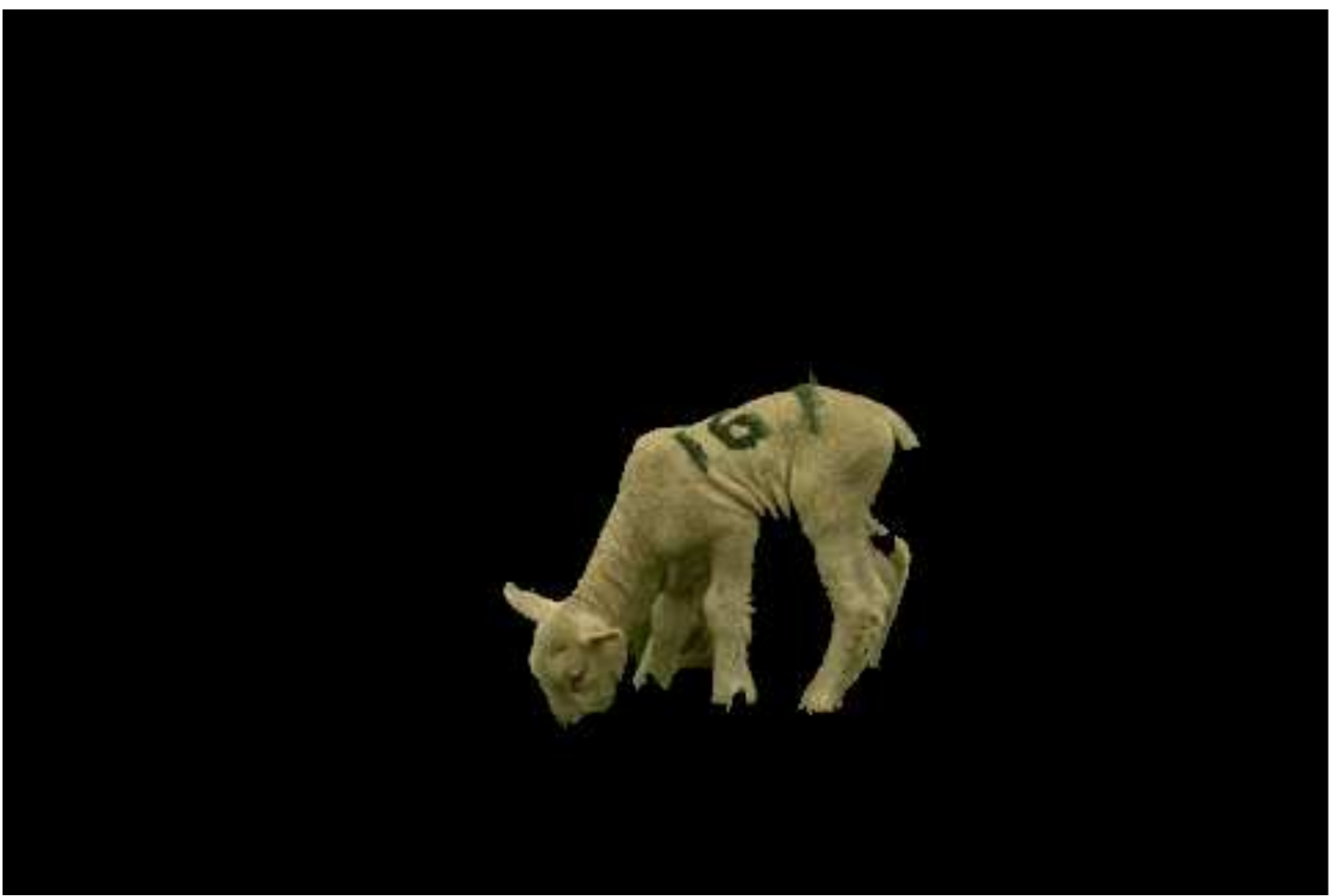}\vspace{-2pt}\caption{{\scriptsize L2box-ADMM,~$f = 6.81$}}\end{subfigure}

\vspace{5pt}

      \begin{subfigure}{\imgwidseg}\includegraphics[width=\textwidth,height=\imgheiseg]{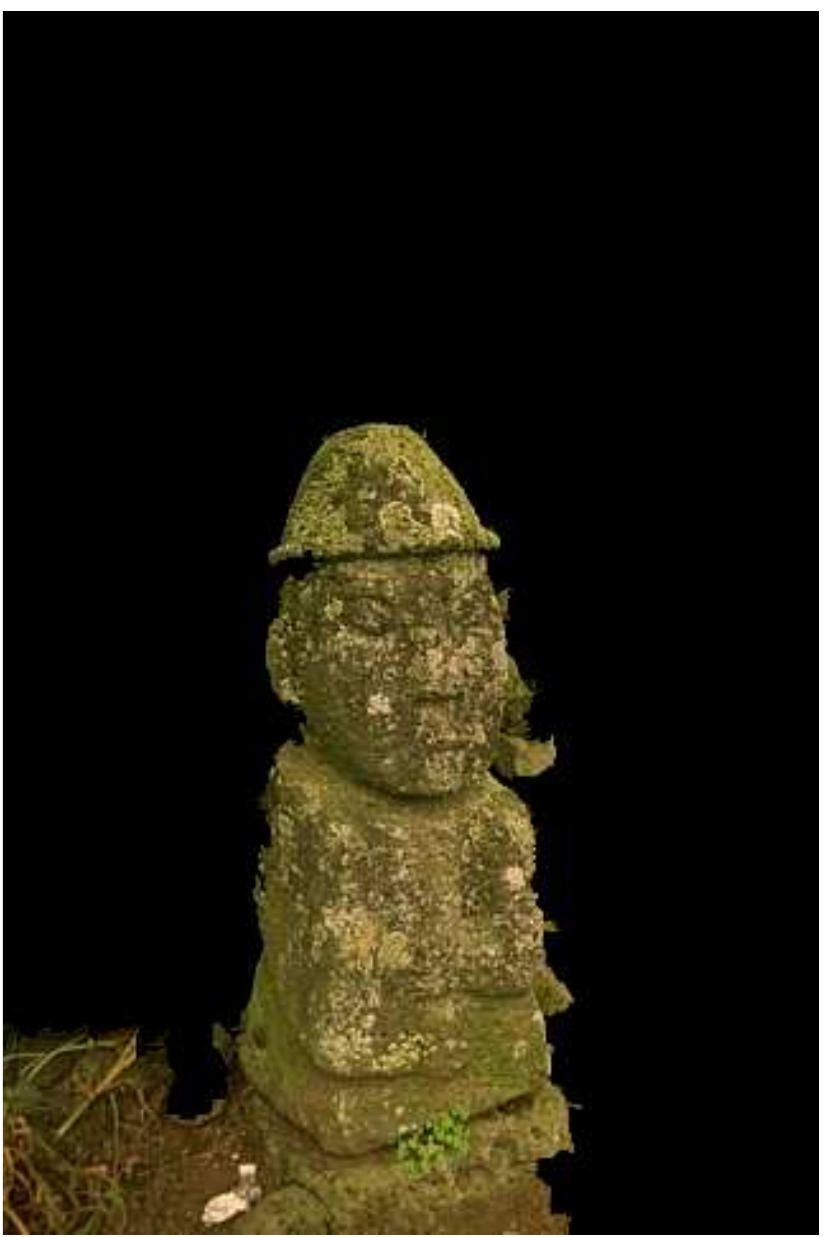}\vspace{-2pt}\caption{{\scriptsize MPEC-EPM,~$f = 49.56$}}\end{subfigure}\ghst
      \begin{subfigure}{\imgwidseg}\includegraphics[width=\textwidth,height=\imgheiseg]{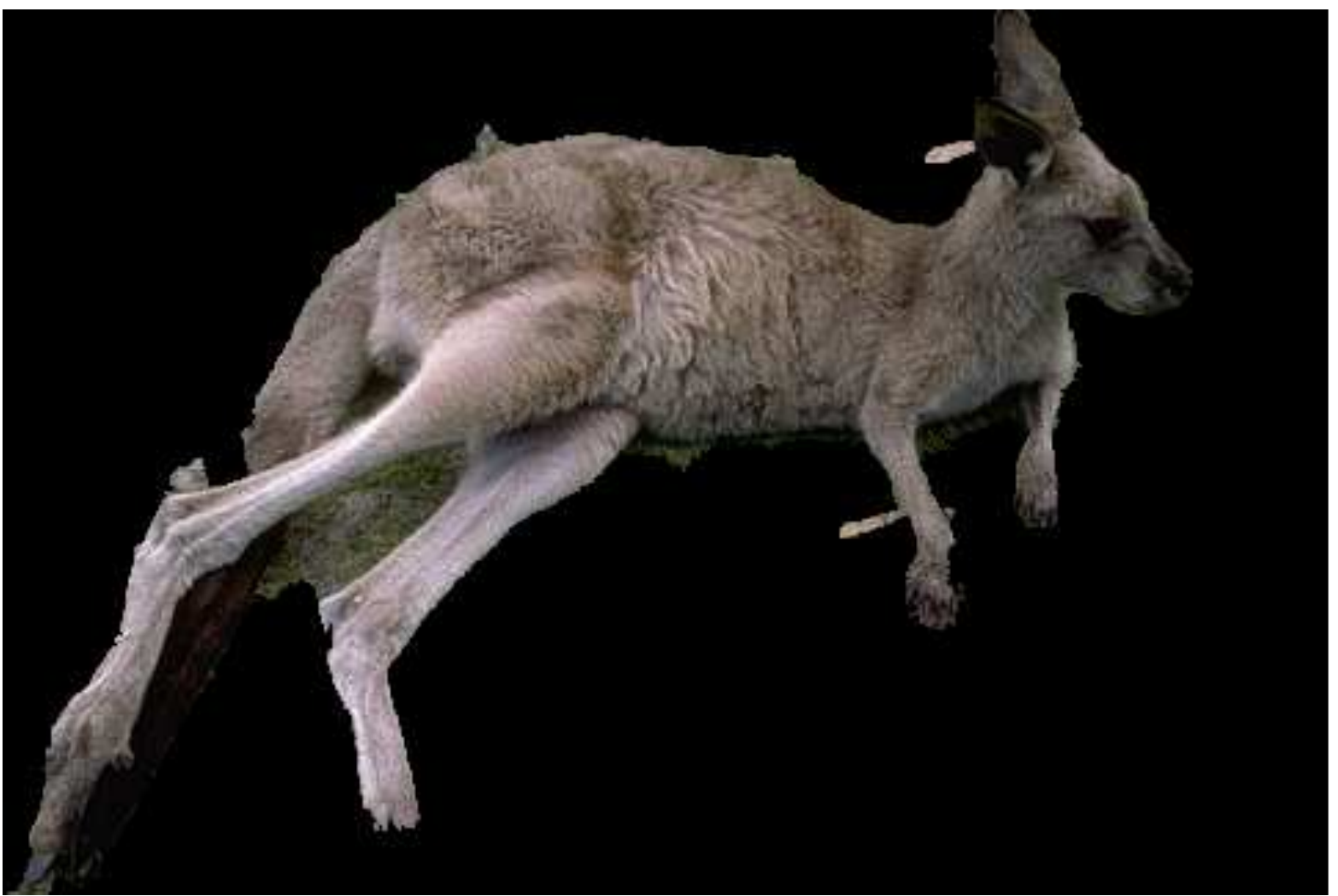}\vspace{-2pt}\caption{{\scriptsize MPEC-EPM,~$f = 16.41$}}\end{subfigure}\ghst
      \begin{subfigure}{\imgwidseg}\includegraphics[width=\textwidth,height=\imgheiseg]{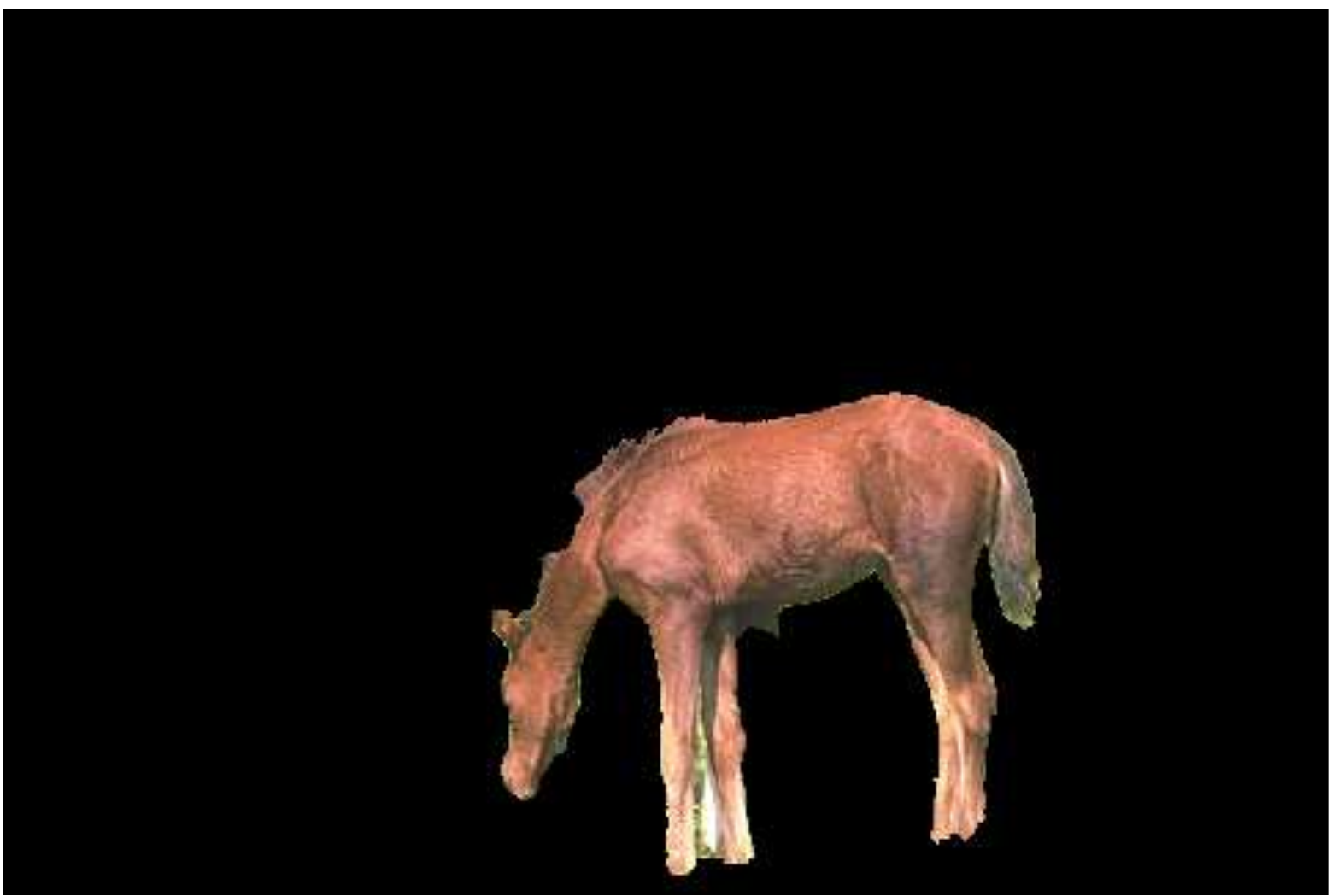}\vspace{-2pt}\caption{{\scriptsize MPEC-EPM,~$f = 14.81$}}\end{subfigure}\ghst
      \begin{subfigure}{\imgwidseg}\includegraphics[width=\textwidth,height=\imgheiseg]{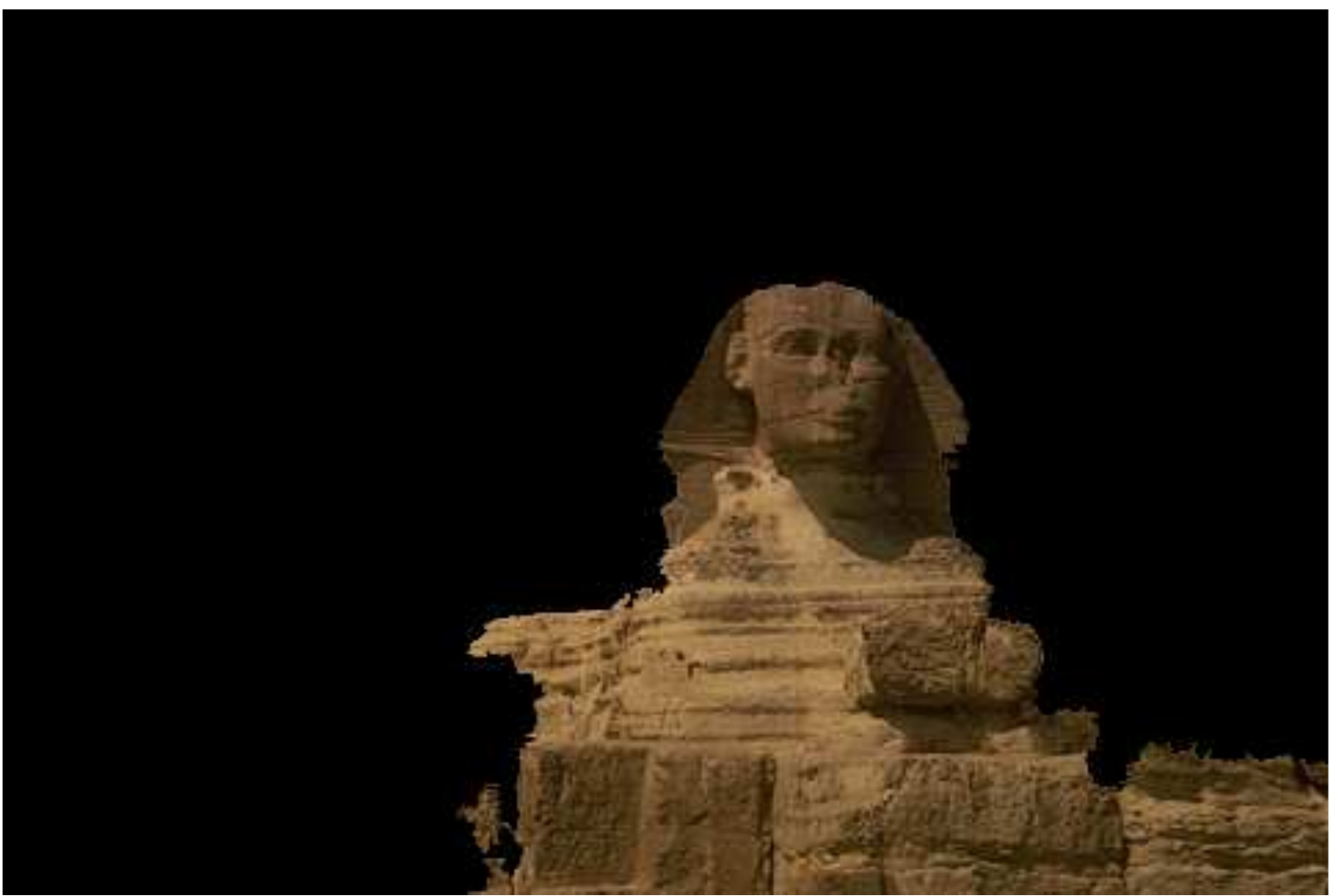}\vspace{-2pt}\caption{{\scriptsize MPEC-EPM,~$f = 25.77$}}\end{subfigure}\ghst
      \begin{subfigure}{\imgwidseg}\includegraphics[width=\textwidth,height=\imgheiseg]{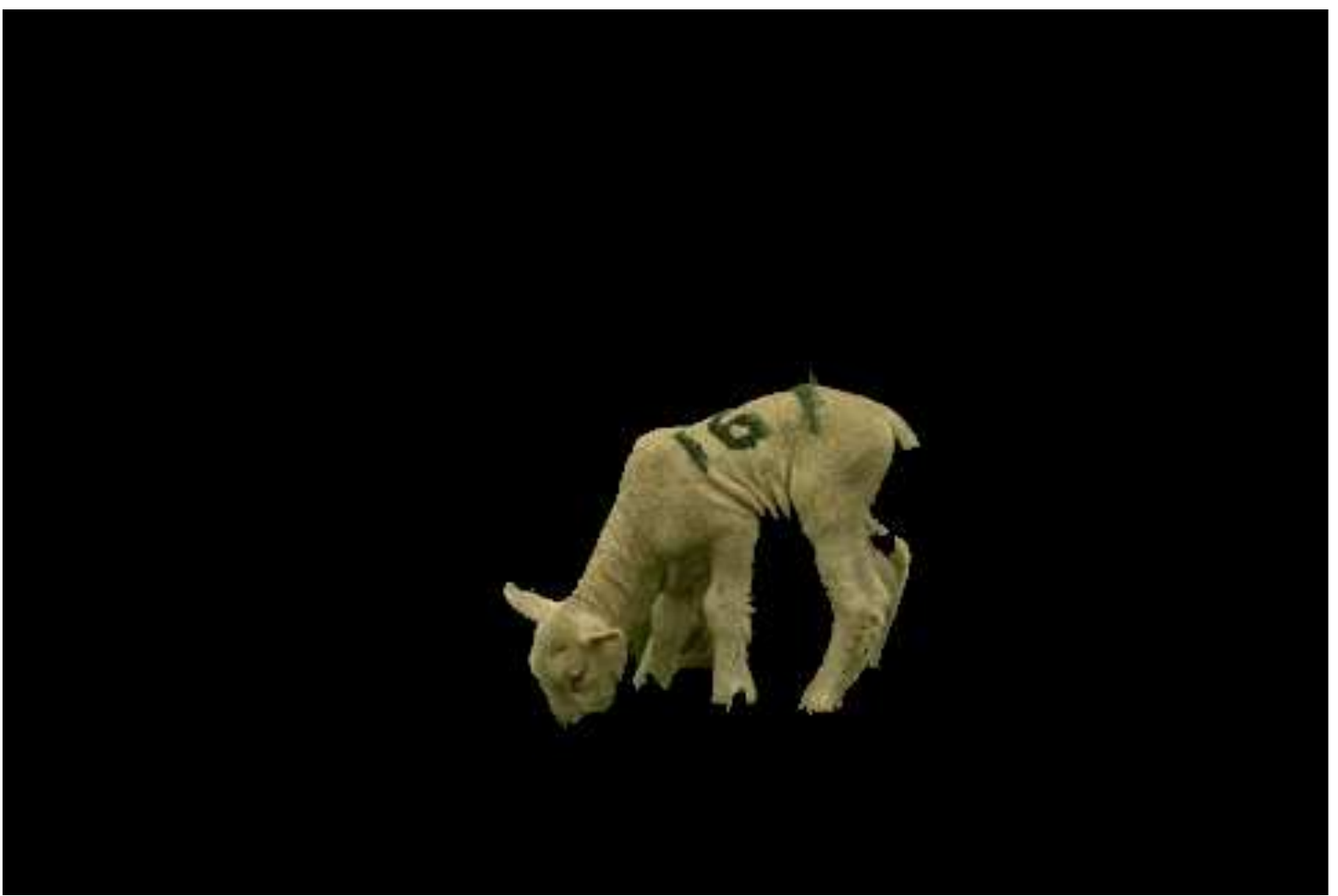}\vspace{-2pt}\caption{{\scriptsize MPEC-EPM,~$f = 6.81$}}\end{subfigure}

\vspace{5pt}

      \begin{subfigure}{\imgwidseg}\includegraphics[width=\textwidth,height=\imgheiseg]{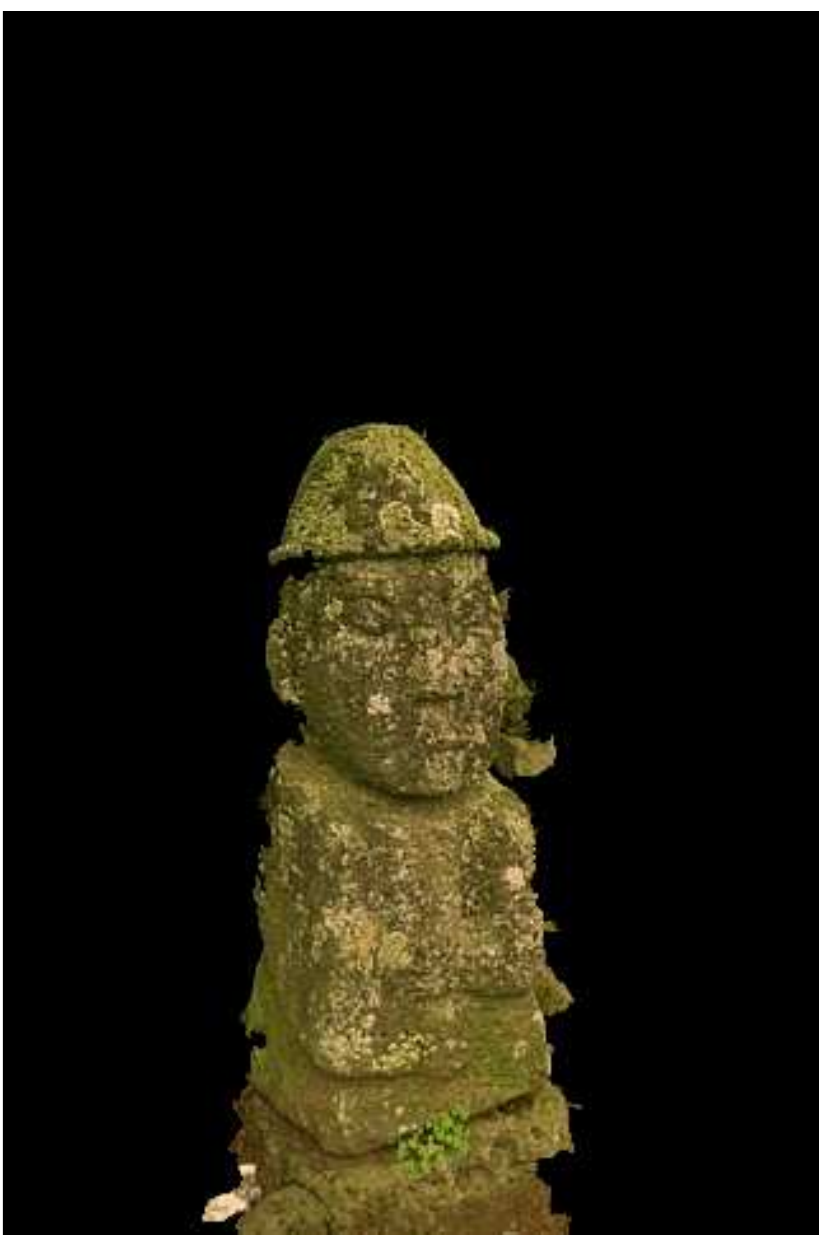}\vspace{-2pt}\caption{{\scriptsize MPEC-ADM,~$f = 50.46$}}\end{subfigure}\ghst
      \begin{subfigure}{\imgwidseg}\includegraphics[width=\textwidth,height=\imgheiseg]{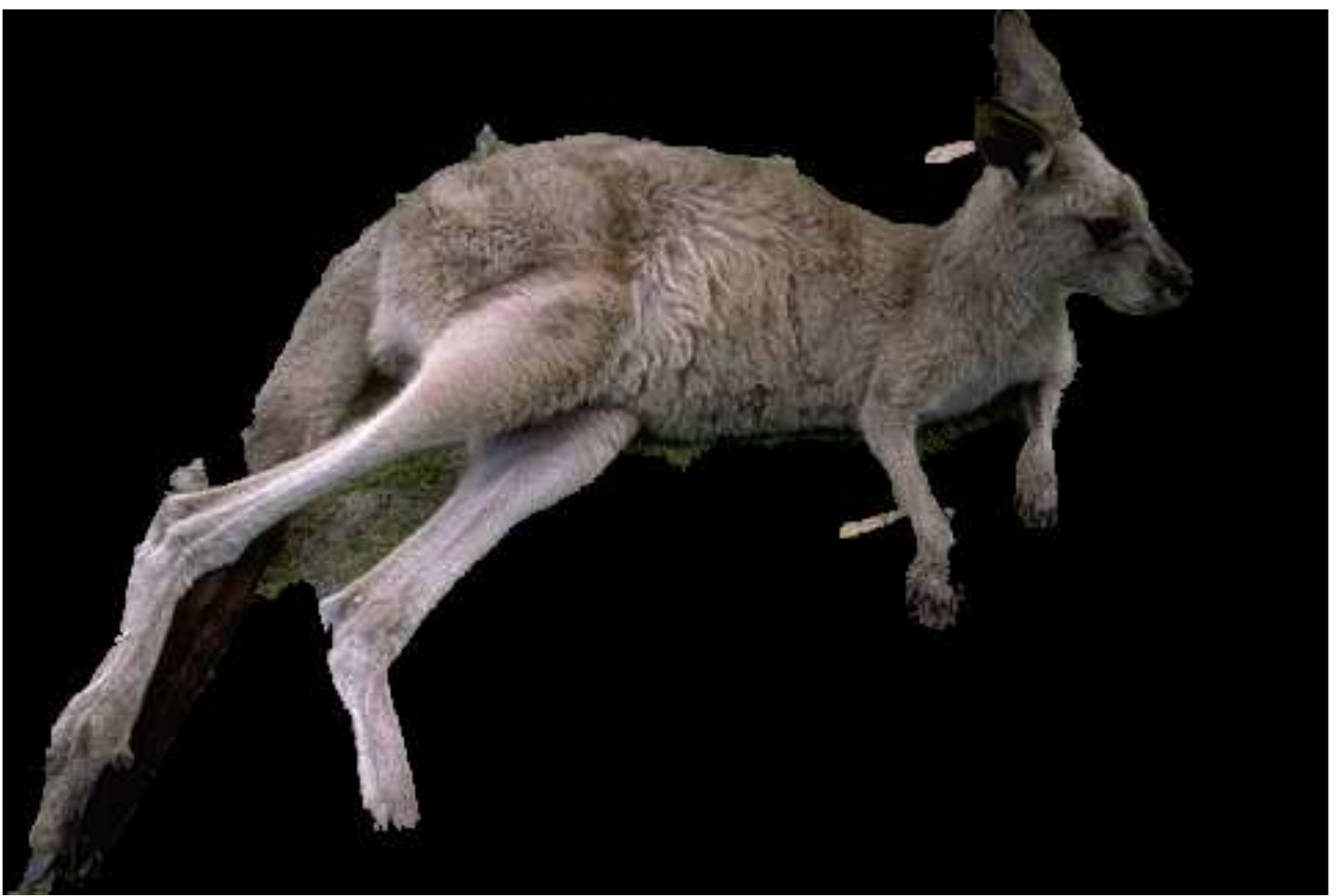}\vspace{-2pt}\caption{{\scriptsize MPEC-ADM,~$f = 16.41$}}\end{subfigure}\ghst
      \begin{subfigure}{\imgwidseg}\includegraphics[width=\textwidth,height=\imgheiseg]{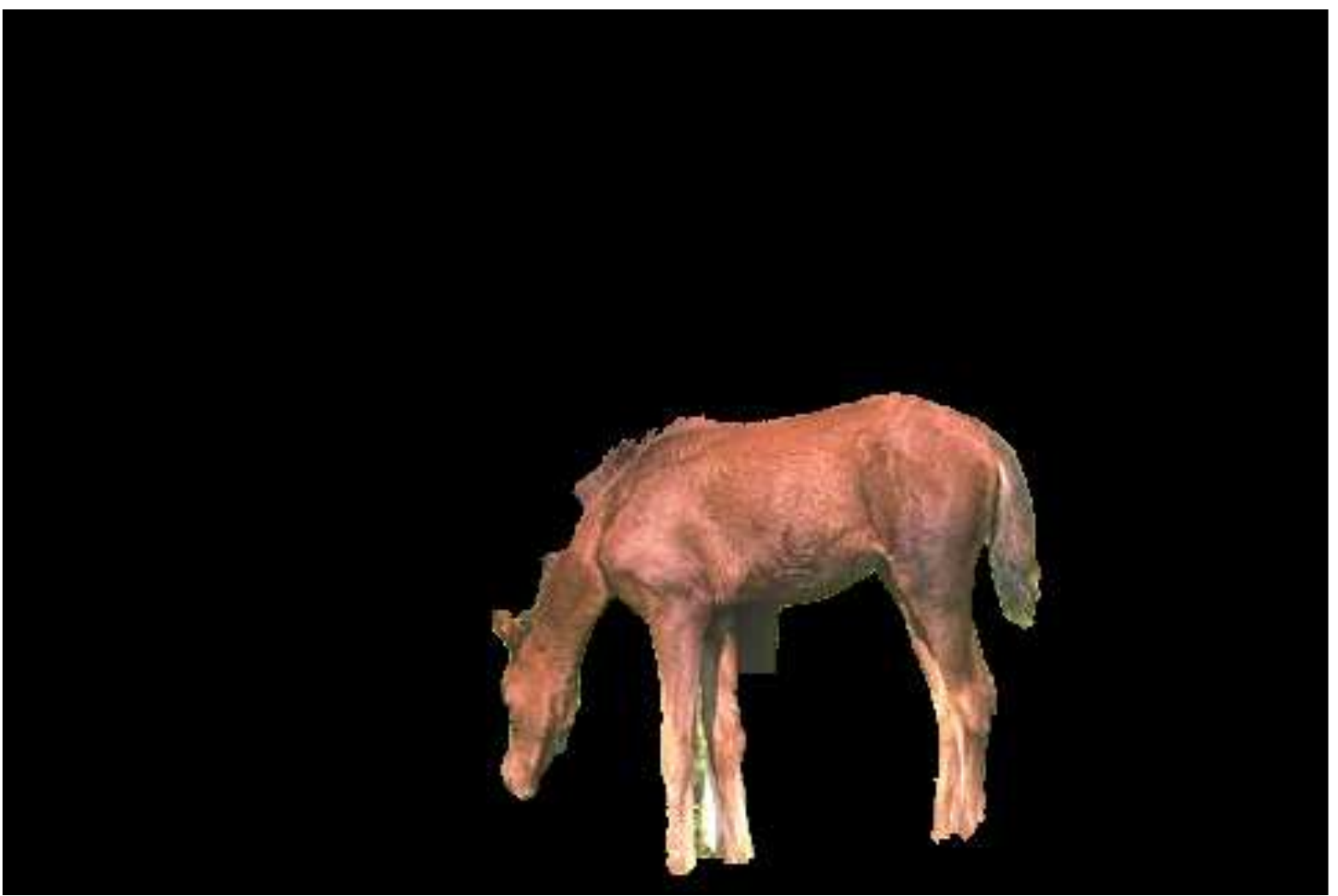}\vspace{-2pt}\caption{{\scriptsize MPEC-ADM,~$f = 13.50$}}\end{subfigure}\ghst
      \begin{subfigure}{\imgwidseg}\includegraphics[width=\textwidth,height=\imgheiseg]{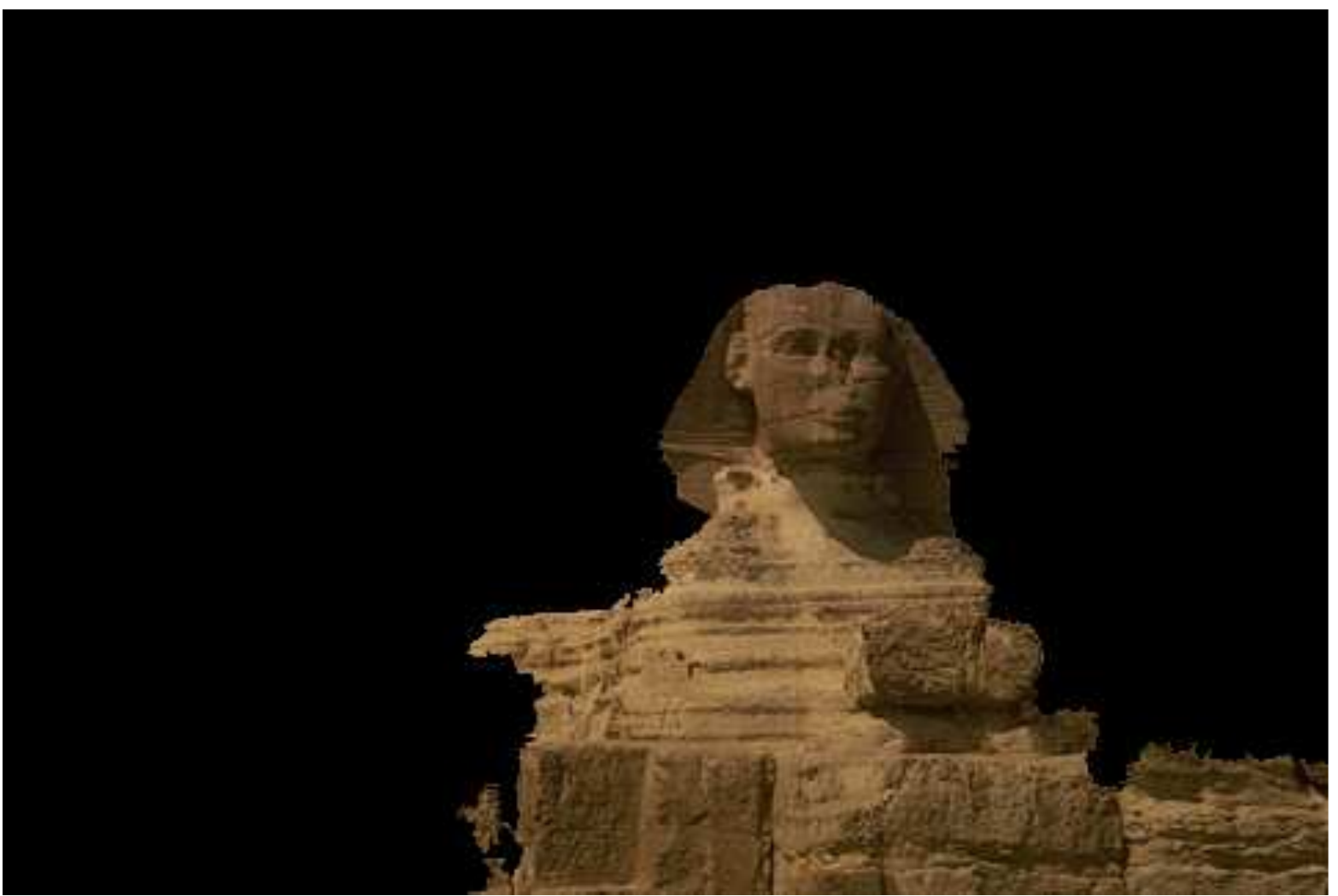}\vspace{-2pt}\caption{{\scriptsize MPEC-ADM,~$f = 24.97$}}\end{subfigure}\ghst
      \begin{subfigure}{\imgwidseg}\includegraphics[width=\textwidth,height=\imgheiseg]{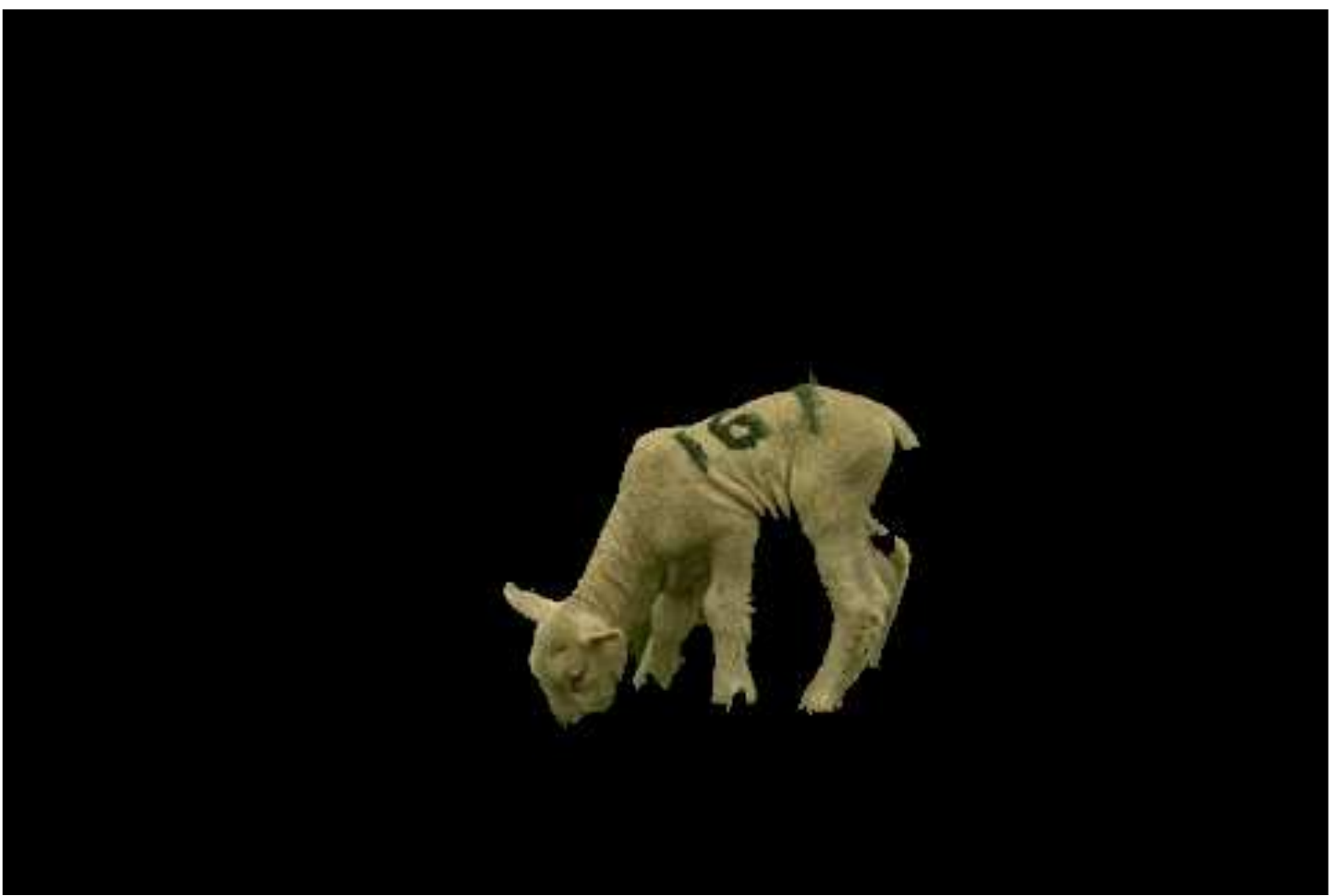}\vspace{-2pt}\caption{{\scriptsize MPEC-ADM,~$f = 6.20$}}\end{subfigure}
\caption{Images used in our constrained image segmentation experiments.}
\label{fig:cseg2}
\end{figure*}

\section{Experimental Validation} \label{sect:exp}

In this section, we demonstrate the effectiveness of our algorithms (MPEC-EPM and MPEC-ADM) on 5 binary optimization tasks, namely graph bisection, constrained image segmentation, dense subgraph discovery, modularity clustering and Markov random fields. All codes are implemented in MATLAB using a 3.20GHz CPU and 8GB RAM. For the purpose of reproducibility, we provide our MATLAB code at: \url{http://isee.sysu.edu.cn/~yuanganzhao/}.

\subsection{Graph Bisection} \label{sect:graph}
Graph bisection aims at separating the vertices of a weighted undirected graph into two disjoint sets with minimal cut edges with equal size. Mathematically, it can be formulated as the following optimization problem \cite{keuchel2003binary,BQP2015Wang}:
\beq
\min_{\bbb{x} \in \{-1,+1\}^n}~\bbb{x}^T\bbb{L}\bbb{x},~s.t.~\bbb{x}^T\bbb{1}=0\nn
\eeq
\noi where $\bbb{L}=\bbb{D}-\bbb{W} \in \mathbb{R}^{n\times n}$ is the Laplacian matrix, $\bbb{W}\in \mathbb{R}^{n\times n}$ is the affinity matrix and $\bbb{D}=diag(\bbb{W1})\in \mathbb{R}^{n\times n}$ is the degree matrix.

\bbb{Compared Methods.} We compare MPEC-EPM and MPEC-ADM against 5 methods on the `4gauss' data set (see Figure \ref{fig:bsect}). (i) LP relaxation simply relaxes the binary constraint to $-\bbb{1}\leq \bbb{x} \leq \bbb{1}$ and solves the following problem:
\beq
\bbb{LP:}~~~~\min_{\bbb{x}}~\bbb{x}^T\bbb{L}\bbb{x},~s.t.~\bbb{x}^T\bbb{1}=0,~-\bbb{1}\leq \bbb{x} \leq \bbb{1}\nn
\eeq
\noi (ii) Ratio Cut (RCUT) and Normalize Cut (NCUT) relax the binary constraint to $\|\bbb{x}\|_2^2 = n$ and solve the following problems \cite{ShiM00,cour2007solving}:
\beq
\bbb{RCut:}~&&\min_{\bbb{x}}~\bbb{x}^T\bbb{L}\bbb{x},~s.t.~\la \bbb{x},\bbb{1}\ra =0,~\|\bbb{x}\|_2^2 = n\nn\\
\bbb{NCut:}~&&\min_{\bbb{x}}~\bbb{x}^T\bar{\bbb{L}}\bbb{x},~s.t.~\la \bbb{x},\bbb{D}^{1/2}\bbb{1}\ra=0,~\|\bbb{x}\|_2^2 = n\nn
\eeq
\noi where $\bar{\bbb{L}}=\bbb{D}^{-1/2}\bbb{L}\bbb{D}^{-1/2}$. The optimal solution of RCut (or NCut) is the second smallest eigenvectors of $\bbb{L}$ (or $\bar{\bbb{L}}$), see e.g. \cite{keuchel2003binary,BQP2015Wang}. (iv) SDP relaxation solves the following convex optimization problem \footnote{\noi For SDP method, we use the randomized rounding strategy in \cite{Goemans95,BQP2015Wang} to get a discrete solution from $\bbb{X}$. Specifically, we sample a random vector $\bbb{x}\in \mathbb{R}^n$ from a Gaussian distribution with mean 0 and covariance $\bbb{X}$, and perform $\bbb{x}^* = \text{sign}(\bbb{x}-\text{median}(\bbb{x}))$. This process is repeated many times and the final solution with the largest objective value is selected.}:
\beq
\bbb{SDP:}~~~~~ \min_{\bbb{X}}~\la \bbb{L},\bbb{X}\ra,~s.t.~diag(\bbb{X})=\bbb{1},~\la \bbb{X},\bbb{1}\bbb{1}^T \ra = 0\nn
\eeq
\noi Finally, we also compare with L2box-ADMM \cite{WuG16a} which applies ADMM directly to the $\ell_2$ box non-separable reformulation:
\beq
\min_{\bbb{x}}~\bbb{x}^T\bbb{L}\bbb{x},~s.t.~\bbb{x}^T\bbb{1}=0,~-\bbb{1}\leq \bbb{x}\leq \bbb{1},~\|\bbb{x}\|_2^2=n.\nn
\eeq
\noi We remark that it is a splitting method that introduces auxiliary variables to separate the two constrained set and then performs block coordinate descend on each variable.

\bbb{Experimental Results.} Several observations can be drawn from Figure \ref{fig:bsect}. (i) The LP, RCUT and NCUT relaxation methods fail to appropriately separate the `4gauss' data set and they result in large objective values. (ii) SDP relaxation provides a good approximation and achieves a lower objective value than LP, RCUT, NCUT and L2box-ADMM. (iii) The proposed methods MPEC-EPM and MPEC-ADM achieve the same lowest objective values among all the compared methods.

\subsection{Constrained Image Segmentation}\label{sect:cons}

In graph-based partition, image is modeled as a weighted undirected graph where nodes corresponds to pixels (or pixel regions) and edges encode feature similarities between the node pairs. Image segmentation can be treated as seeking a partition $\bbb{x}\in \{-1,+1\}^n$ to cut the edges of the graph with minimal weights. Prior label information on the vertices of the graph can be incorporated to improve performance, leading to the following optimization problem:
\beq
\min_{\bbb{x} \in \{-1,+1\}^n}~\bbb{x}^T\bbb{L}\bbb{x},~s.t.~\bbb{x}_F=1,~\bbb{x}_B=-1\nn
\eeq
\noi where $F$ and $B$ denote the index of foreground pixels set and background pixels set, respectively; $\bbb{L}\in \mathbb{R}^{n\times n}$ is the graph Laplacian matrix. Since SDP method can not solve large scale image segmentation problems, we over-segment the images into SLIC pixel regions using the `VLFeat' toolbox \cite{Vedaldi2010}. The affinity matrix $\bbb{W}$ is constructed based on the color similarities and spatial adjacencies between pixel regions.

\bbb{Compared Methods.} We compare MPEC-EPM and MPEC-ADM against 4 methods on the Weizman horses and MSRC datasets \cite{BQP2015Wang} (see Figure \ref{fig:cseg1}). (i) Biased normalized cut (BNCut) \cite{cour2007solving} extends Normalized Cut \cite{ShiM00} to encode the labelled foreground pixels as a quadratic constraint on the solution $\bbb{x}$. The solution of BNCut is a linear combination of the eigenvectors of normalized Laplacian matrix \cite{ShiM00}. (ii) LP relaxation simply replaces the binary constraint with a soft constraint and solves a quadratic programming problem: $\min_{\bbb{x}}~\bbb{x}^T\bbb{L}\bbb{x},~s.t.~-\bbb{1}\leq \bbb{x}\leq \bbb{1},~\bbb{x}_F=1,~\bbb{x}_B=-1$. (iii) SDP relaxation method considers the following optimization problem:
\beq \label{eq:img:seg:sdp}
\min_{\bbb{X} \succeq 0}~\la \bbb{L},\bbb{X}\ra,~s.t.~diag(\bbb{X})=\bbb{1},~\bbb{X}_{I} = 1,~\bbb{X}_J=-1\nn
\eeq
\noi with $\bbb{X}\in\mathbb{R}^{n\times n}$, and ${I}$ and $J$ are the index pairs of similarity and dissimilarity, respectively. Therefore, it contains $n+ \begin{psmallmatrix} 2\\ |B| \end{psmallmatrix}+ \begin{psmallmatrix} 2\\ |F| \end{psmallmatrix}+ |B|\cdot |F|$ linear equality constraints. We use `cvx' optimization software \cite{cvx} to solve this problem. (iv) L2box-ADMM considers the following optimization problem: $\min_{\bbb{x}}~\bbb{x}^T\bbb{L}\bbb{x},~s.t.~-\bbb{1}\leq \bbb{x}\leq \bbb{1},~\bbb{x}_F=1,~\bbb{x}_B=-1,~\|\bbb{x}\|_2^2=n$.

\bbb{Experimental Results.} Several observations can be drawn from Figure \ref{fig:cseg2}. (i) LP relaxation generates better image segmentation results than BNCUT except in the second image. Moreover, we found that a lower objective value does not always necessarily result in better view for image segmentation. We argue that the parameter in constructing the similarity graph is responsible for this result. (ii) SDP method and L2box-ADMM method generate better solutions than BNCUT and LP. (iii) The proposed MPEC-EPM and MPEC-ADM generally obtain lower objective values and outperform all the other compared methods.

\subsection{Dense Subgraph Discovery}

\begin{table}
\small
\caption{The statistics of the web graph data sets used in our dense subgraph discovery experiments.}
\center
\begin{tabular}{|c|c|c|c|}
  \hline
 Graph & \# Nodes & \# Arcs & Avg. Degree \\
  \hline
wordassociation & 10617 & 72172 & 6.80 \\
 enron & 69244 & 276143 & 3.99 \\
 uk-2007-05 & 100000 & 3050615 & 30.51 \\
cnr-2000 & 325557 & 3216152 & 9.88 \\
 dblp-2010 & 326186 & 1615400 & 4.95 \\
 in-2004 & 1382908 & 16917053 & 12.23 \\
 amazon-2008 & 735323 & 5158388 & 7.02 \\
 dblp-2011 & 986324 & 6707236 & 6.80 \\
\hline
\end{tabular}
\label{tab:graph}
\end{table}

\begin{figure*} 
\captionsetup[subfigure]{justification=centering}
    \centering
      \begin{subfigure}{0.7\textwidth}\includegraphics[width=\textwidth,height=13pt]{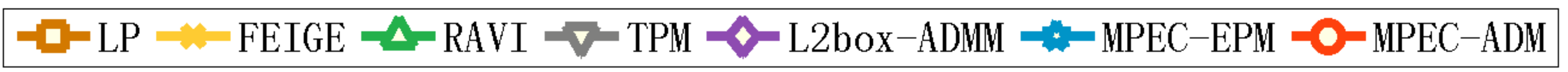}\end{subfigure}

      \vspace{3pt}

      \begin{subfigure}{\fourfigwid}\includegraphics[width=\textwidth,height=\imgheisubgraph]{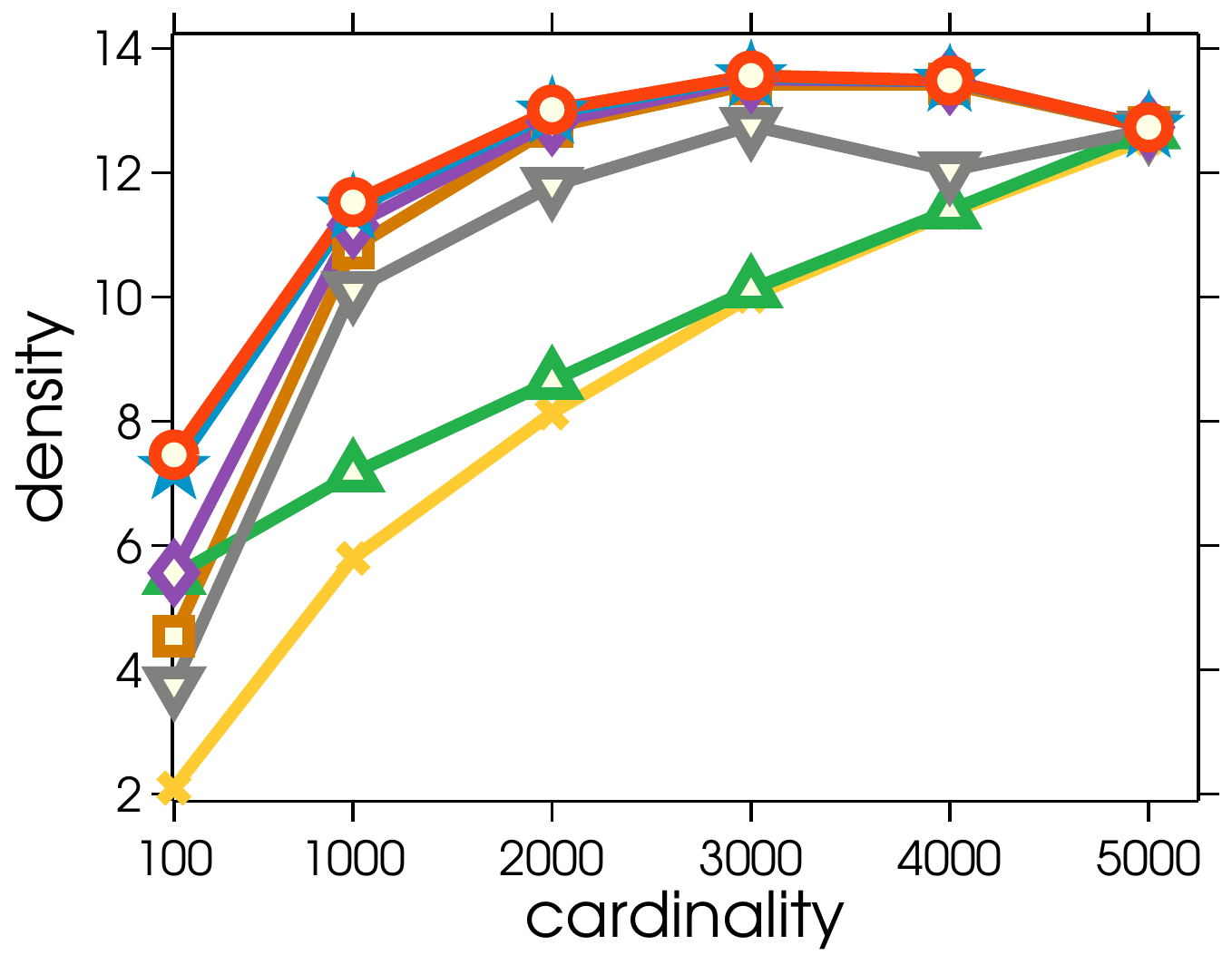}\caption{wordassociation}\end{subfigure}\ghs
      \begin{subfigure}{\fourfigwid}\includegraphics[width=\textwidth,height=\imgheisubgraph]{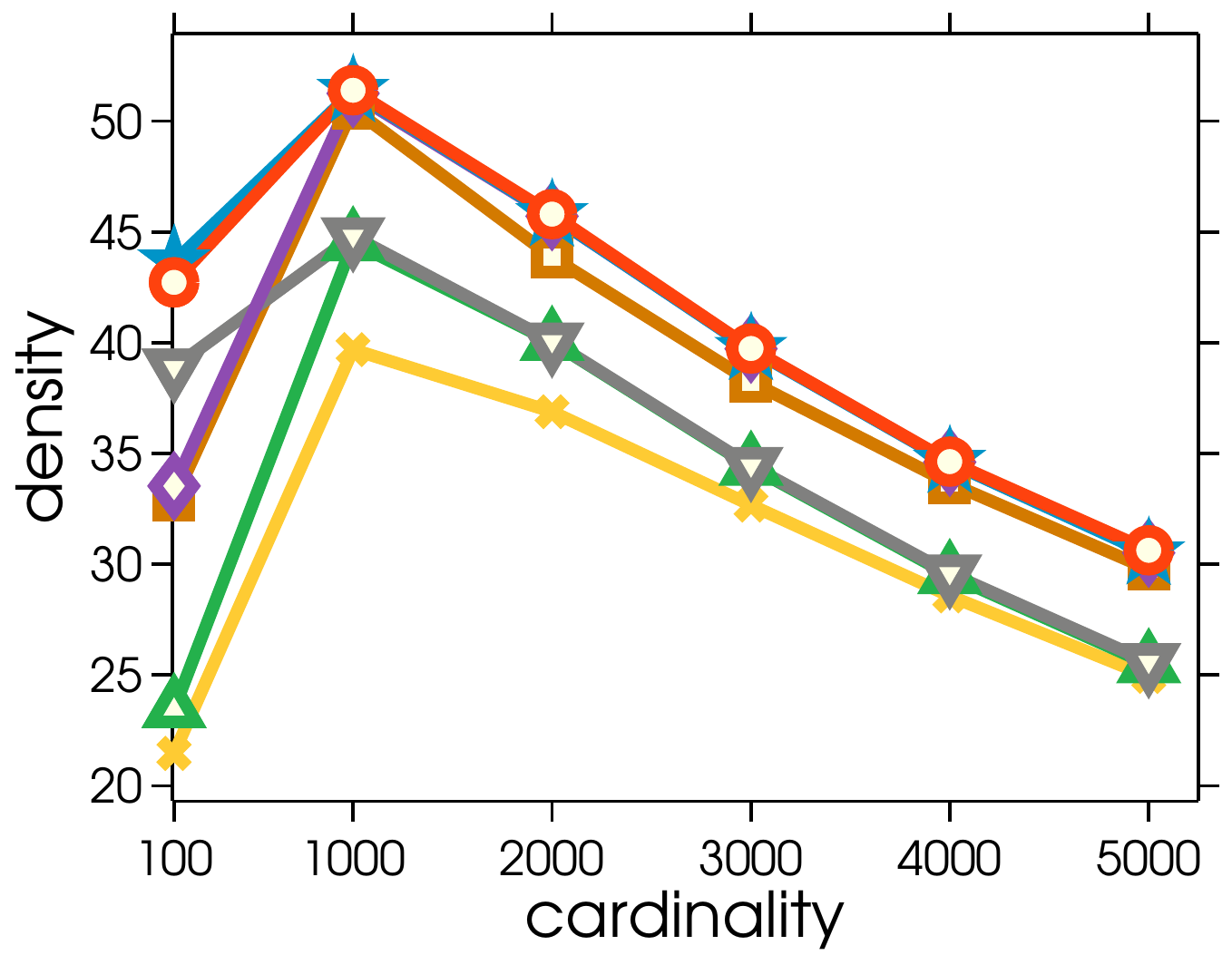}\caption{uk-2007-05}\end{subfigure}\ghs
      \begin{subfigure}{\fourfigwid}\includegraphics[width=\textwidth,height=\imgheisubgraph]{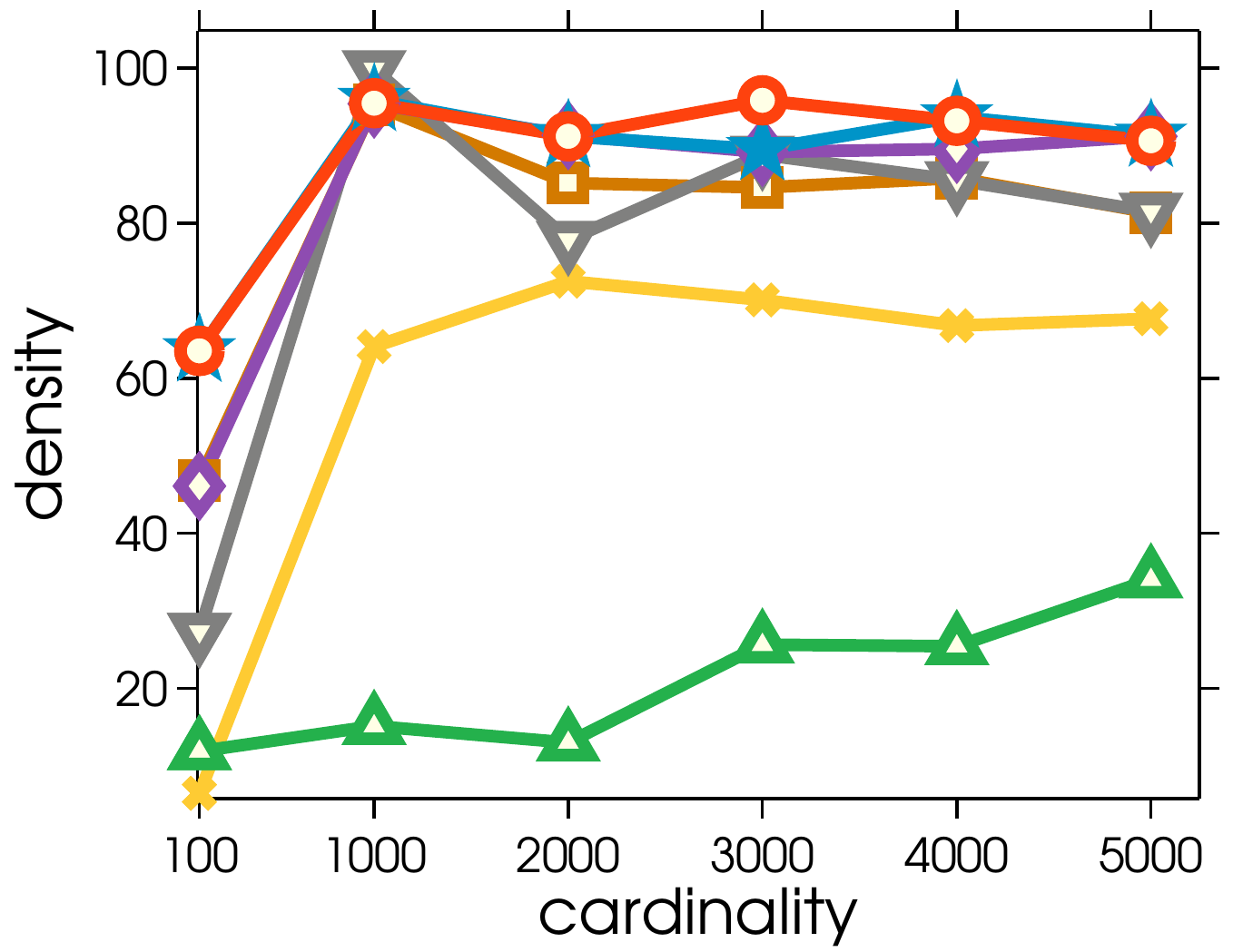}\caption{uk-2014-tpd}\end{subfigure}\ghs
      \begin{subfigure}{\fourfigwid}\includegraphics[width=\textwidth,height=\imgheisubgraph]{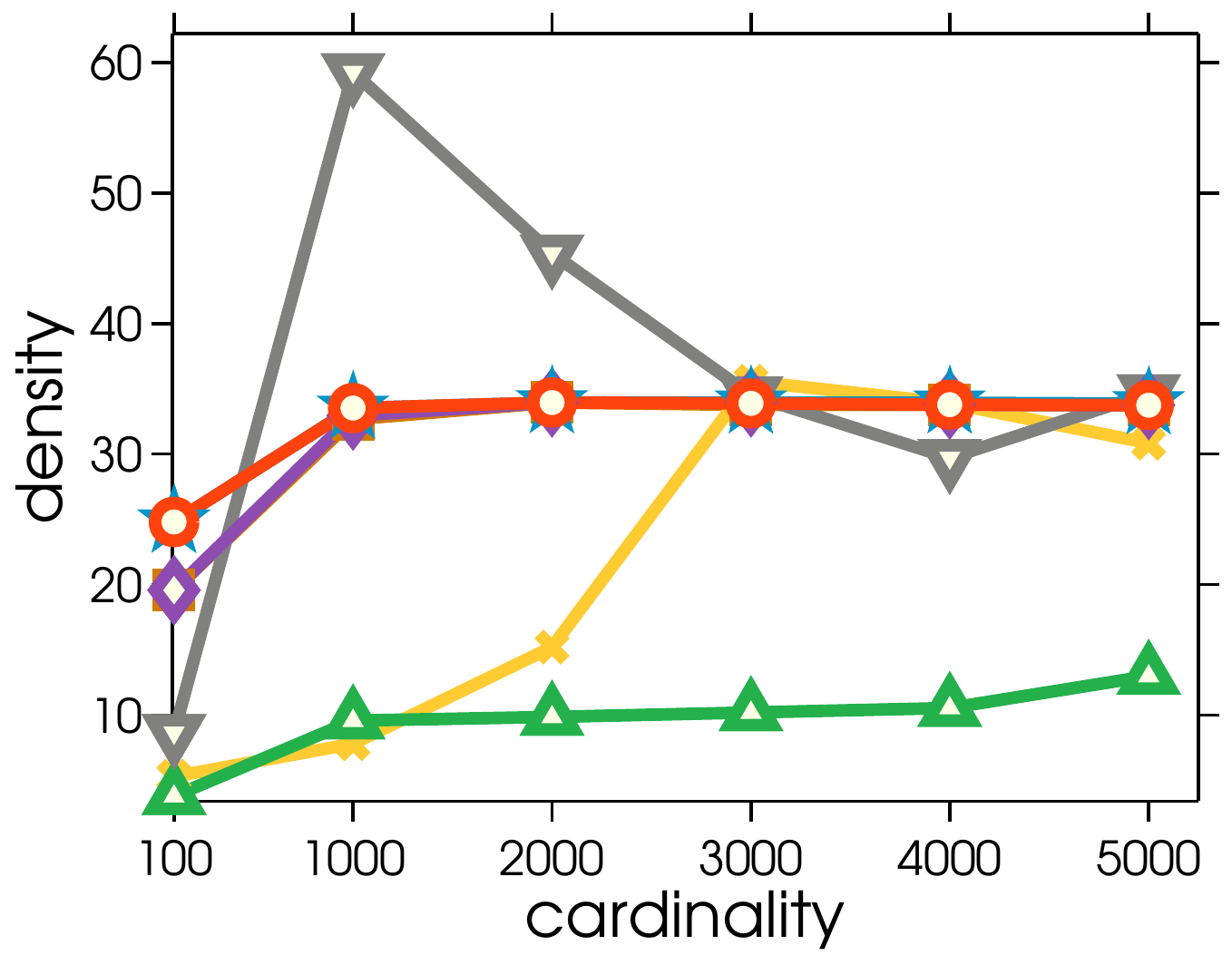}\caption{in-2004}\end{subfigure}

      \begin{subfigure}{\fourfigwid}\includegraphics[width=\textwidth,height=\imgheisubgraph]{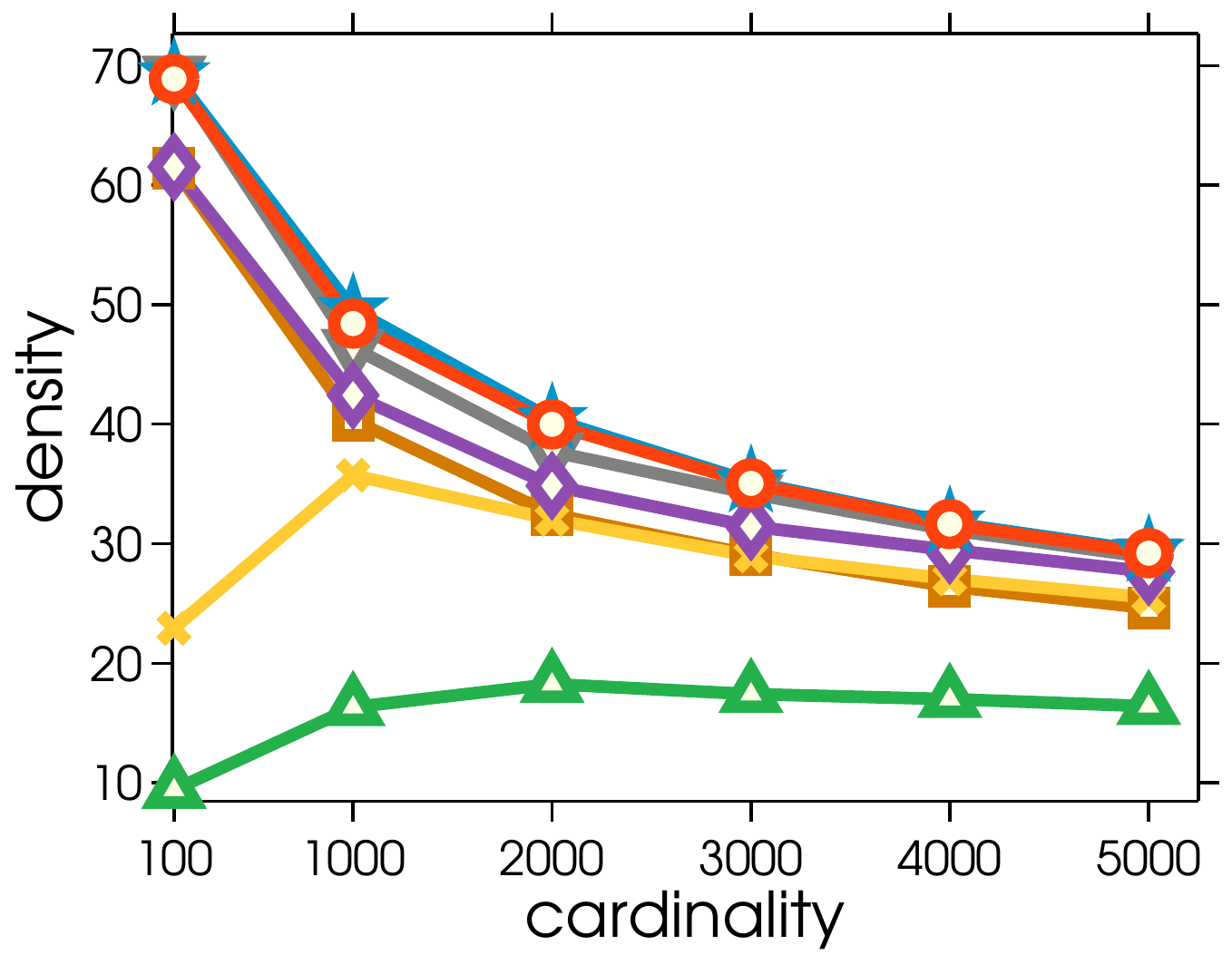}\caption{amazon-2008}\end{subfigure}\ghs
      \begin{subfigure}{\fourfigwid}\includegraphics[width=\textwidth,height=\imgheisubgraph]{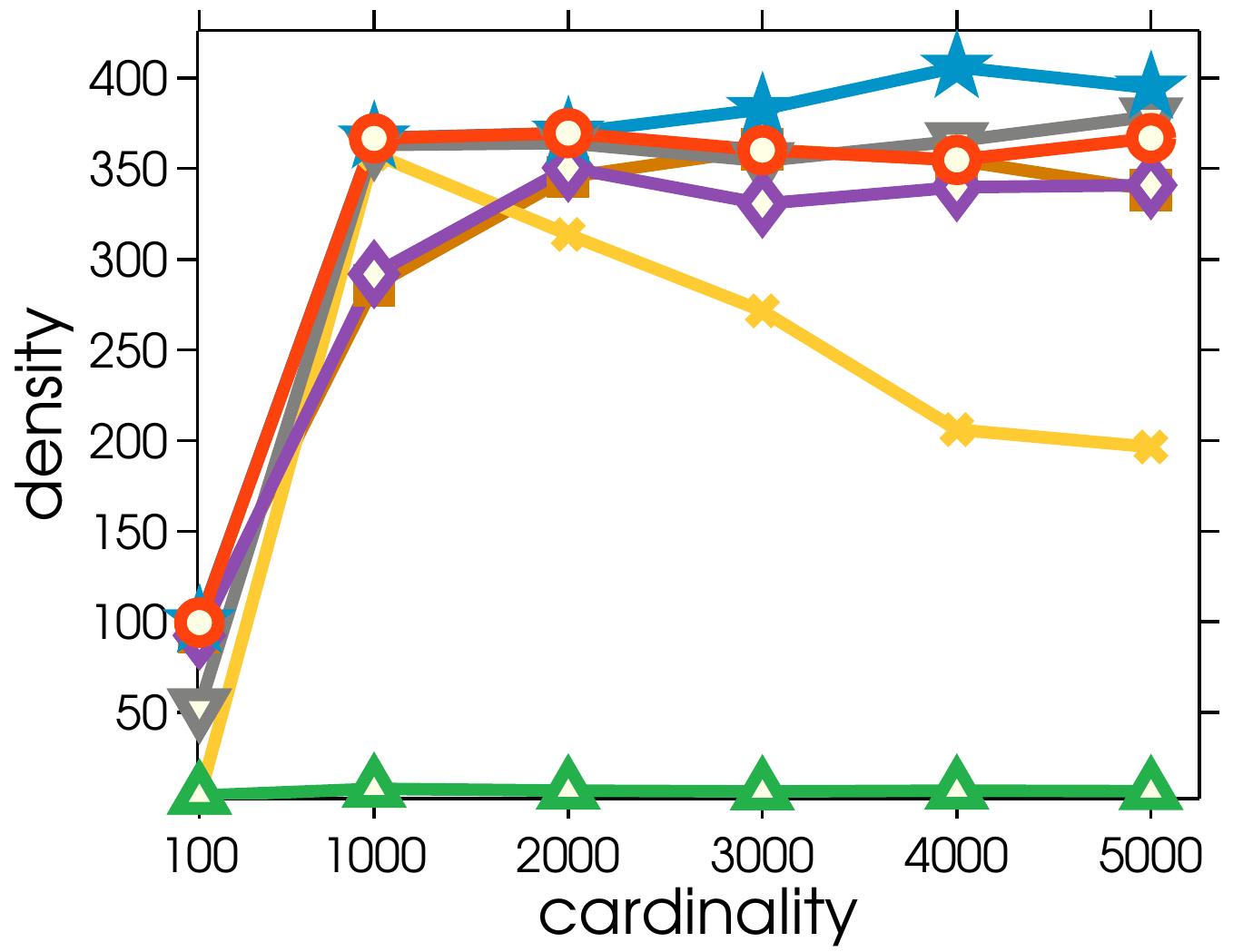}\caption{enron}\end{subfigure}\ghs
      \begin{subfigure}{\fourfigwid}\includegraphics[width=\textwidth,height=\imgheisubgraph]{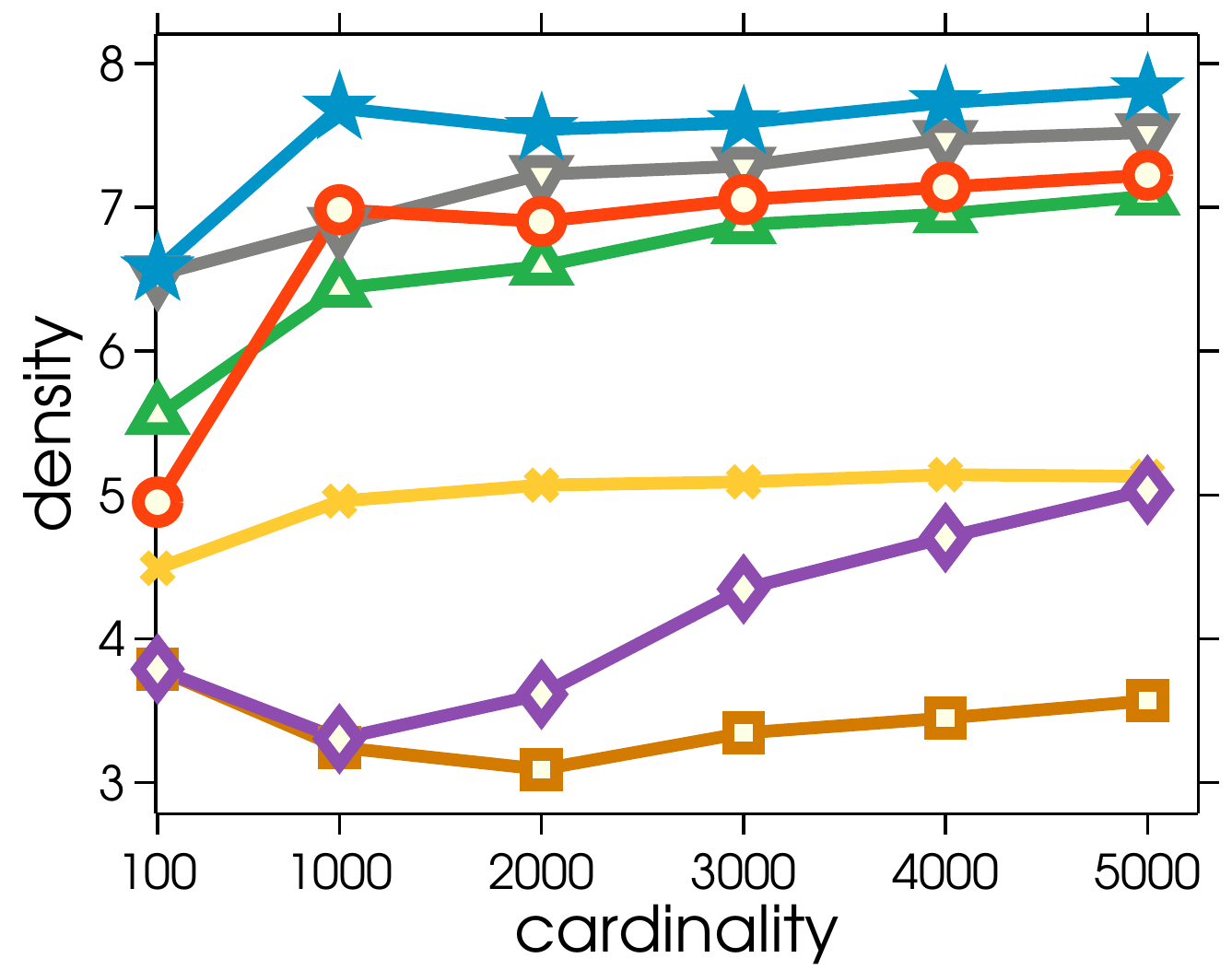}\caption{cnr-2000}\end{subfigure}\ghs
      \begin{subfigure}{\fourfigwid}\includegraphics[width=\textwidth,height=\imgheisubgraph]{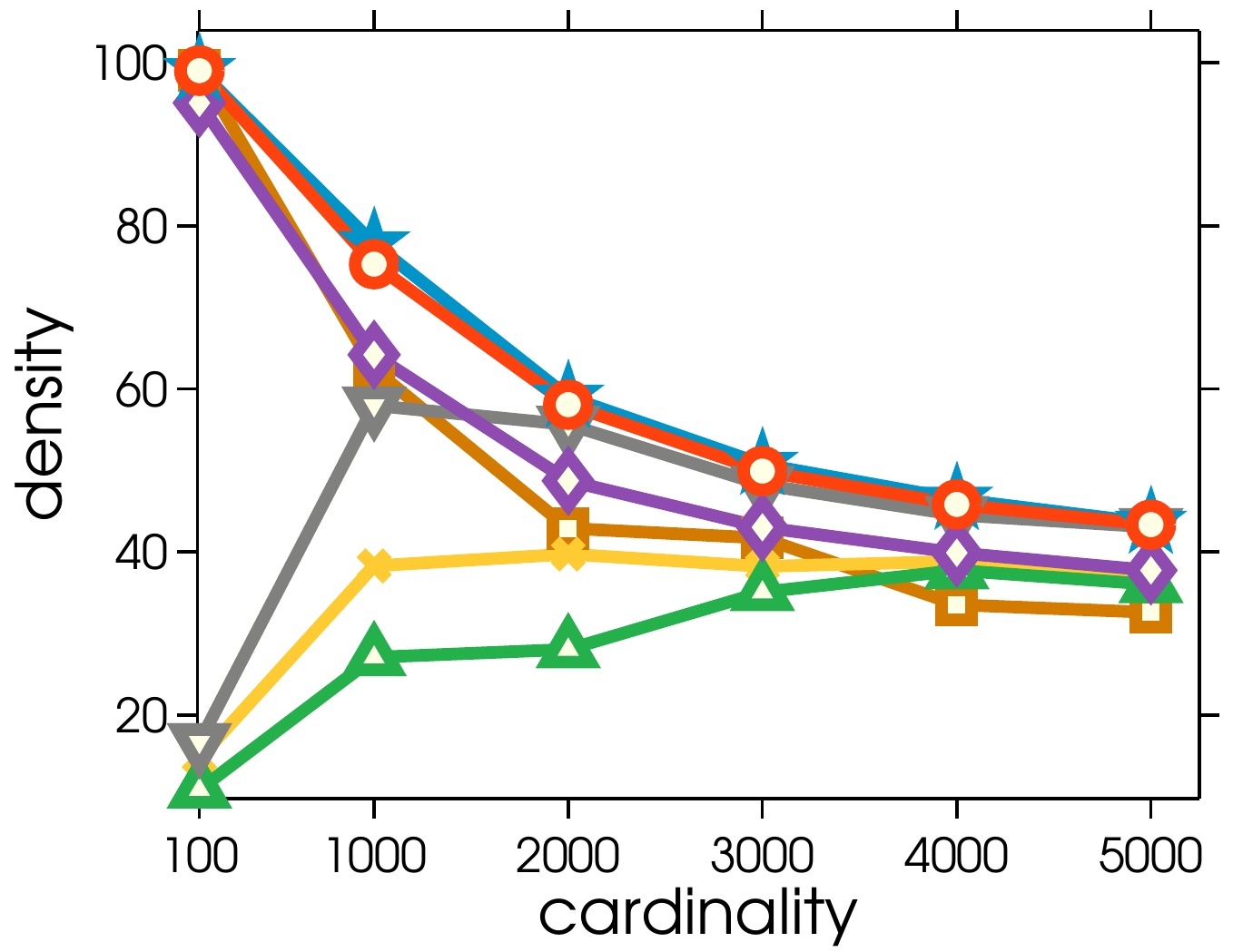}\caption{dblp-2010}\end{subfigure}

\caption{Experimental results for dense subgraph discovery.}
\label{fig:subgraph}
\end{figure*}
Dense subgraphs discovery \cite{ravi1994heuristic,feige2001dense,Yuan013} is a fundamental graph-theoretic problem, as it captures numerous graph mining applications, such as community finding, regulatory motifs detection, and real-time story identification. It aims at finding the maximum density subgraph on $k$ vertices, which can be formulated as the following binary program:
\beq \label{eq:subgraph}
\max_{\bbb{x} \in \{0,1\}^n}~\bbb{x}^T\bbb{W}\bbb{x},~s.t.~\bbb{x}^T\bbb{1} = k
\eeq
\noi where $\bbb{W}\in \mathbb{R}^{n\times n}$ is the adjacency matrix of the graph. Although the objective function in (\ref{eq:subgraph}) may not be convex, one can append an additional term $\lambda \bbb{x}^T\bbb{x}$ to the objective with a sufficiently large $\lambda$ such that $\lambda \bbb{I}-\bbb{W} \succeq 0$. This is equivalent to adding a constant to the objective since $\lambda \bbb{x}^T\bbb{x} = \lambda k$ in the effective domain. Therefore, we have the following optimization problem which is equivalent to (\ref{eq:subgraph}):
\beq \label{eq:subgraph2}
\min_{\bbb{x} \in \{0,1\}^n}~\bbb{x}^T(\lambda \bbb{I}-\bbb{W})\bbb{x},~s.t.~\bbb{x}^T\bbb{1} = k\nn
\eeq
\noi In the experiments, $\lambda$ is set to the largest eigenvalue of $\bbb{W}$.

\bbb{Compared Methods.} We compare MPEC-EPM and MPEC-ADM against 5 methods on 8 datasets\footnote{\url{http://law.di.unimi.it/datasets.php}} (see Table \ref{tab:graph}). (i) Feige's greedy algorithm (GEIGE) \cite{feige2001dense} is included in our comparisons. This method is known to achieve the best approximation ratio for general $k$. (ii) Ravi's greedy algorithm  (RAVI)  \cite{ravi1994heuristic} starts from a heaviest edge and repeatedly adds a vertex to the current subgraph to maximize the weight of the resulting new subgraph. It has asymptotic performance guarantee of $\pi/2$ when the weights satisfy the triangle inequality. (iii) LP relaxation solves a capped simplex problem by standard quadratic programming technique: $\min_{\bbb{x} }~\bbb{x}^T(\lambda \bbb{I}-\bbb{W})\bbb{x},~s.t.~\bbb{0}\leq \bbb{x} \leq \bbb{1},~\bbb{x}^T\bbb{1} = k.$ (iii) L2box-ADMM solves a spherical constraint optimization problem: $\min_{\bbb{x} }~\bbb{x}^T(\lambda \bbb{I}-\bbb{W})\bbb{x},~s.t.~\bbb{0}\leq \bbb{x} \leq \bbb{1},~\bbb{x}^T\bbb{1} = k,~\|2\bbb{x}-1\|_2^2=n.$ (iv) Truncated Power Method (TPM) \cite{Yuan013} considers an iterative procedure that combines power iteration and hard-thresholding truncation. It works by greedily decreasing the objective while maintaining the desired binary property for the intermediate solutions. We use the code\footnote{\url{https://sites.google.com/site/xtyuan1980/publications}} provided by the authors. As suggested in \cite{Yuan013}, the initial solution is set to the indicator vector of the vertices with the top $k$ weighted degrees of the graph in our experiments.

\bbb{Experimental Results.} Several observations can be drawn from Figure \ref{fig:subgraph}. (i) Both FEIGE and RAVI generally fail to solve the dense subgraph discovery problem and they lead to solutions with low density. (ii) LP relaxation gives better performance than state-off-the-art technique TPM in some cases. (iii) L2box-ADMM outperforms LP relaxation for all cases but it generates unsatisfying accuracy in `enron', `cnr-2000' and `dblp-2010'. (iv) Our proposed method MPEC-EPM and MPEC-ADM generally outperforms all the compared methods, while MPEC-EPM seems to present slightly better results than MPEC-ADM in this group of experiments.

\begin{table}[!h]
\small
\caption{The statistics of the web graph data sets used in our modularity clustering experiments.}
\center
\begin{tabular}{|c|c|c|c|}
  \hline
 Graph & \# Nodes & \# Arcs & Avg. Degree \\
  \hline
   karate & 34 & 78 & 4.59 \\
 collab & 235 & 415 & 3.53 \\
 dolphins & 62 & 159 & 5.13 \\
 email & 1133 & 5451 & 9.62 \\
 lesmis & 77 & 820 & 21.30 \\
 polbooks & 105 & 441 & 8.40 \\
 afootball & 115 & 616 & 10.71 \\
 jazz & 198 & 2742 & 27.70 \\
  \hline
\end{tabular}
\label{tab:clustering}
\end{table}

\begin{figure*}
\captionsetup[subfigure]{justification=centering}
    \centering
    \begin{subfigure}{0.7\textwidth}\includegraphics[width=\textwidth,height=13pt]{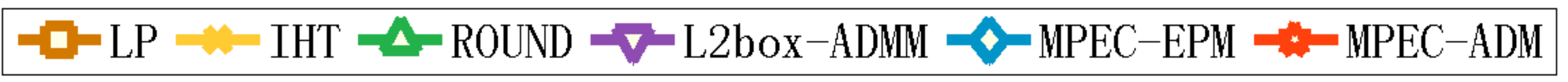}\end{subfigure}

      \vspace{3pt}

      \begin{subfigure}{\fourfigwid}\includegraphics[width=\textwidth,height=\imgheiclustering]{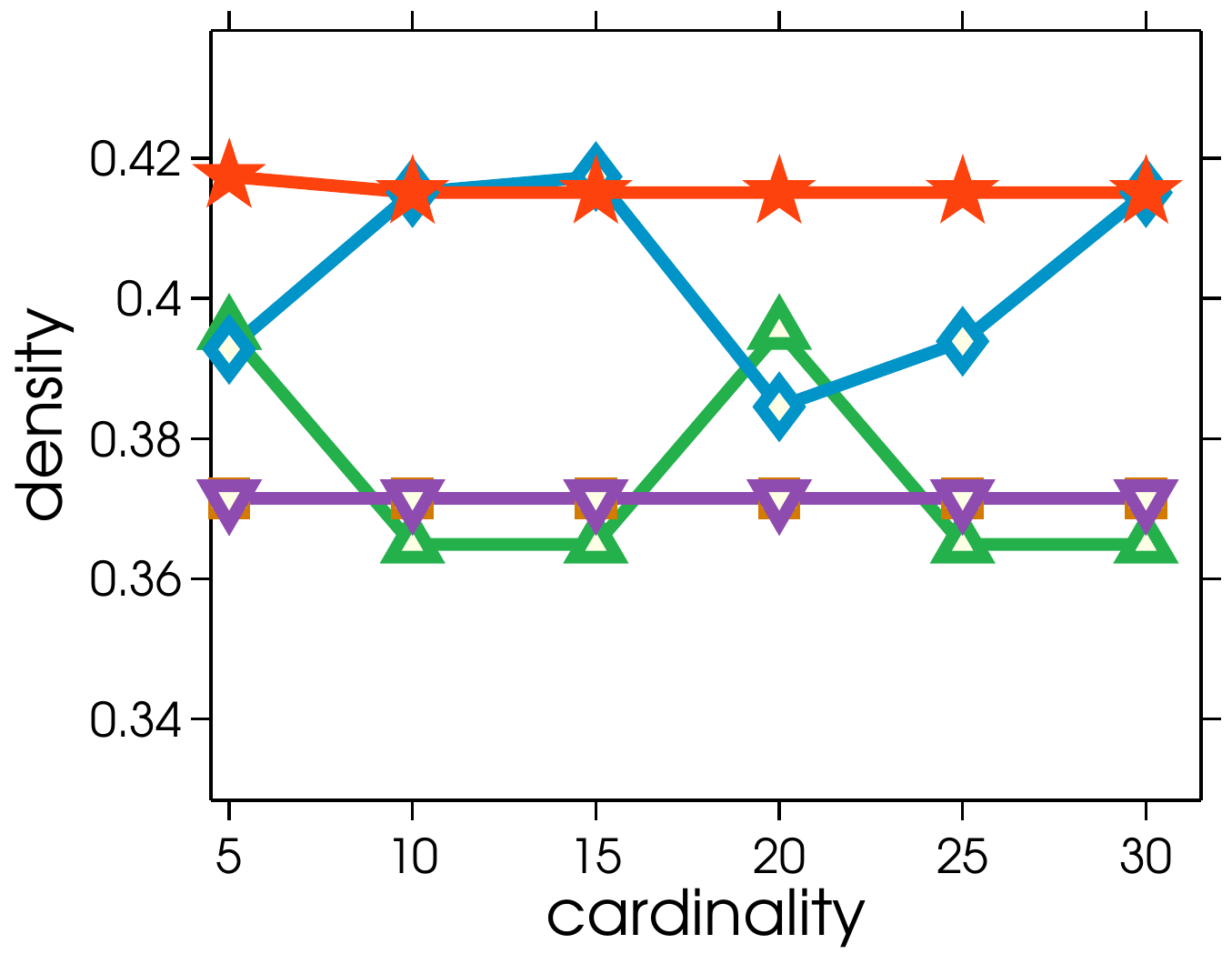}\caption{karate}\end{subfigure}\ghs
      \begin{subfigure}{\fourfigwid}\includegraphics[width=\textwidth,height=\imgheiclustering]{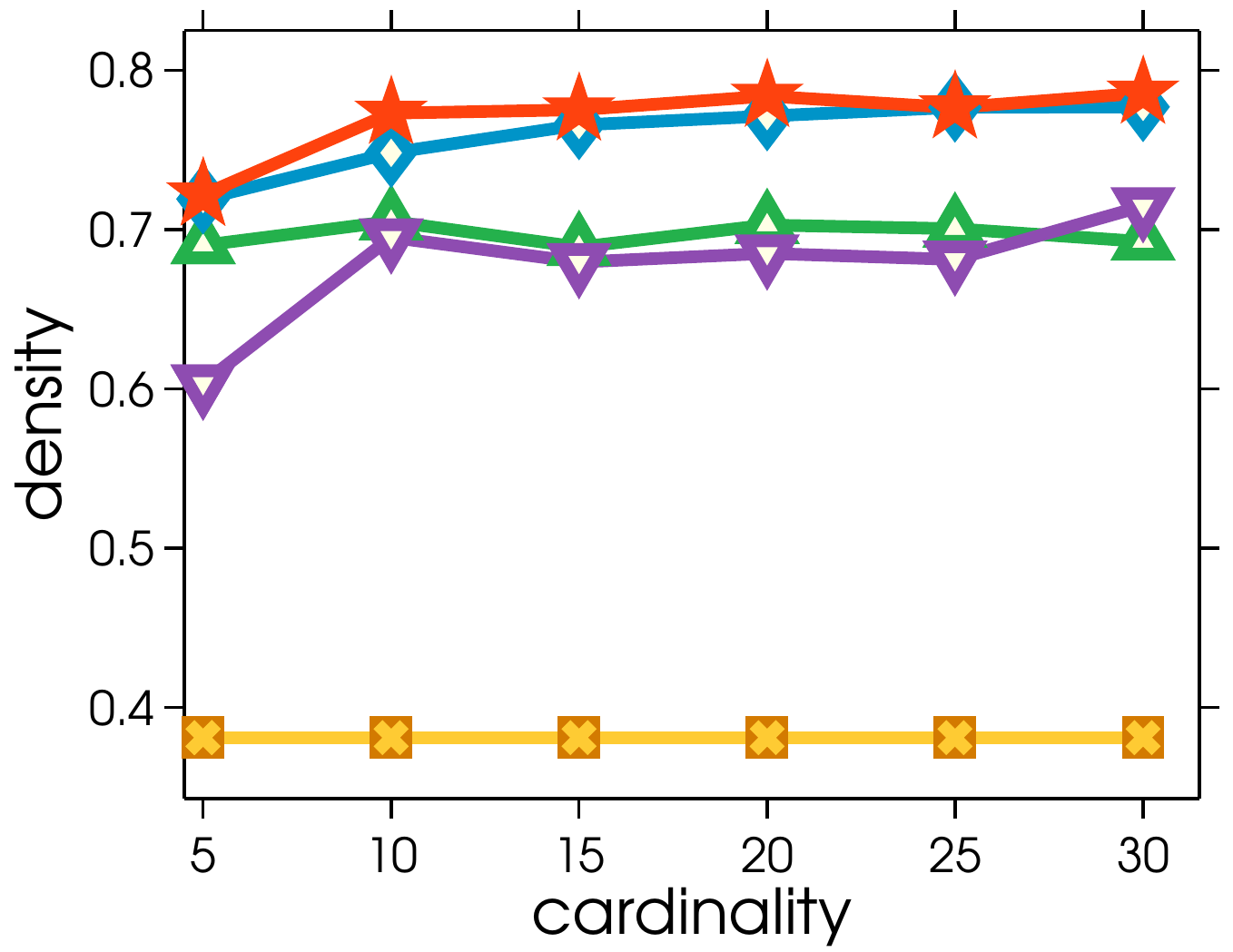}\caption{collab}\end{subfigure}\ghs
      \begin{subfigure}{\fourfigwid}\includegraphics[width=\textwidth,height=\imgheiclustering]{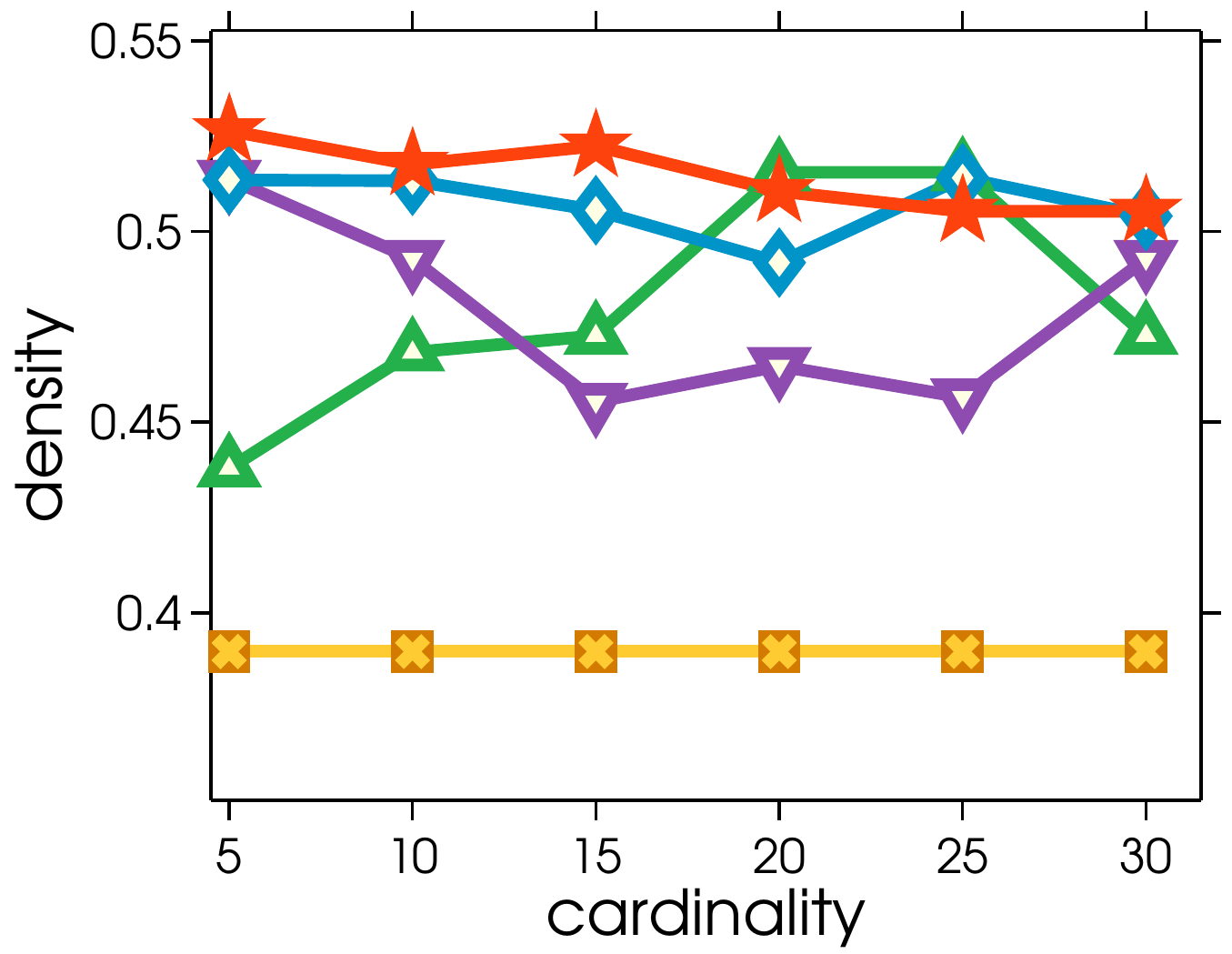}\caption{dolphins}\end{subfigure}\ghs
      \begin{subfigure}{\fourfigwid}\includegraphics[width=\textwidth,height=\imgheiclustering]{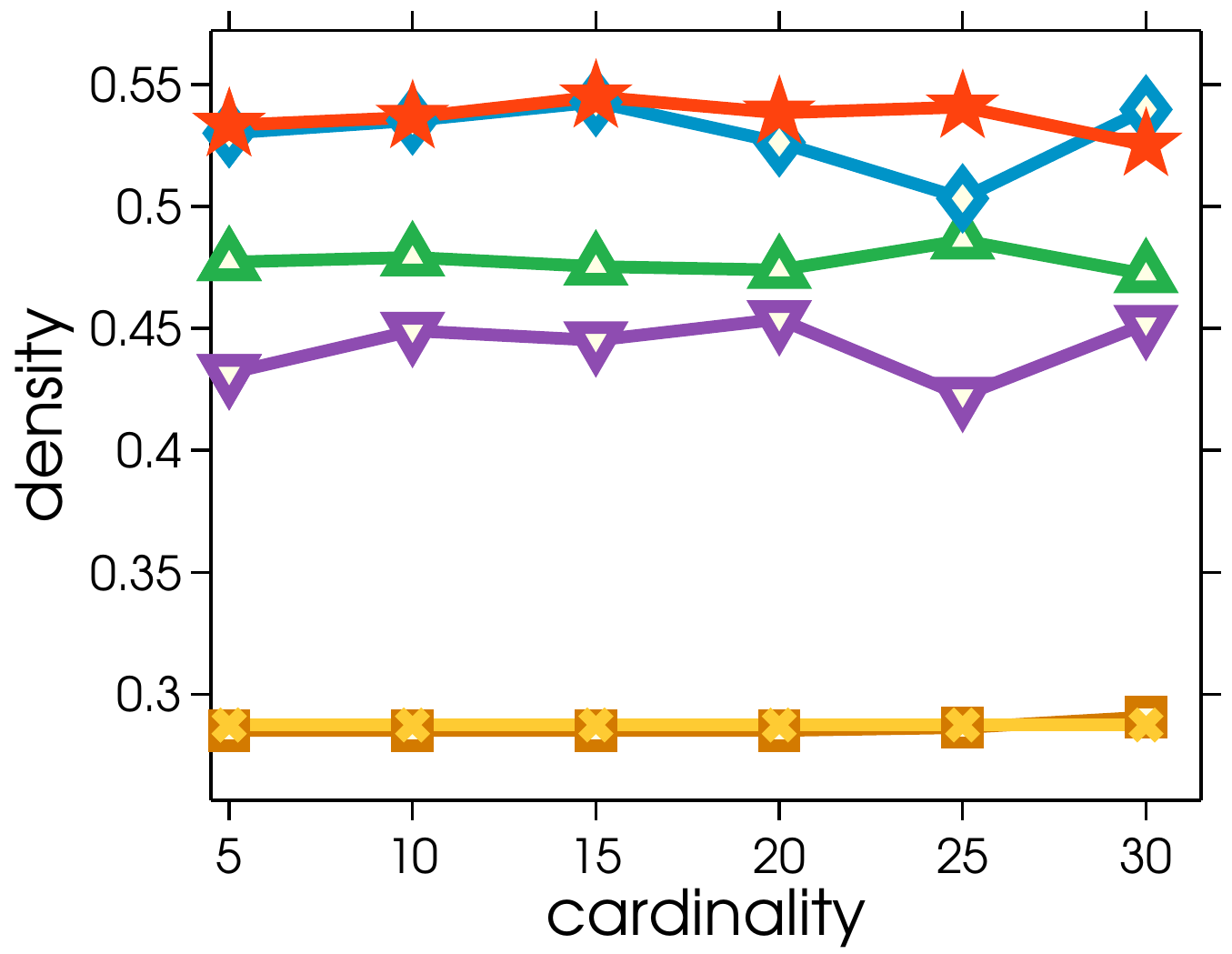}\caption{email}\end{subfigure}

      \begin{subfigure}{\fourfigwid}\includegraphics[width=\textwidth,height=\imgheiclustering]{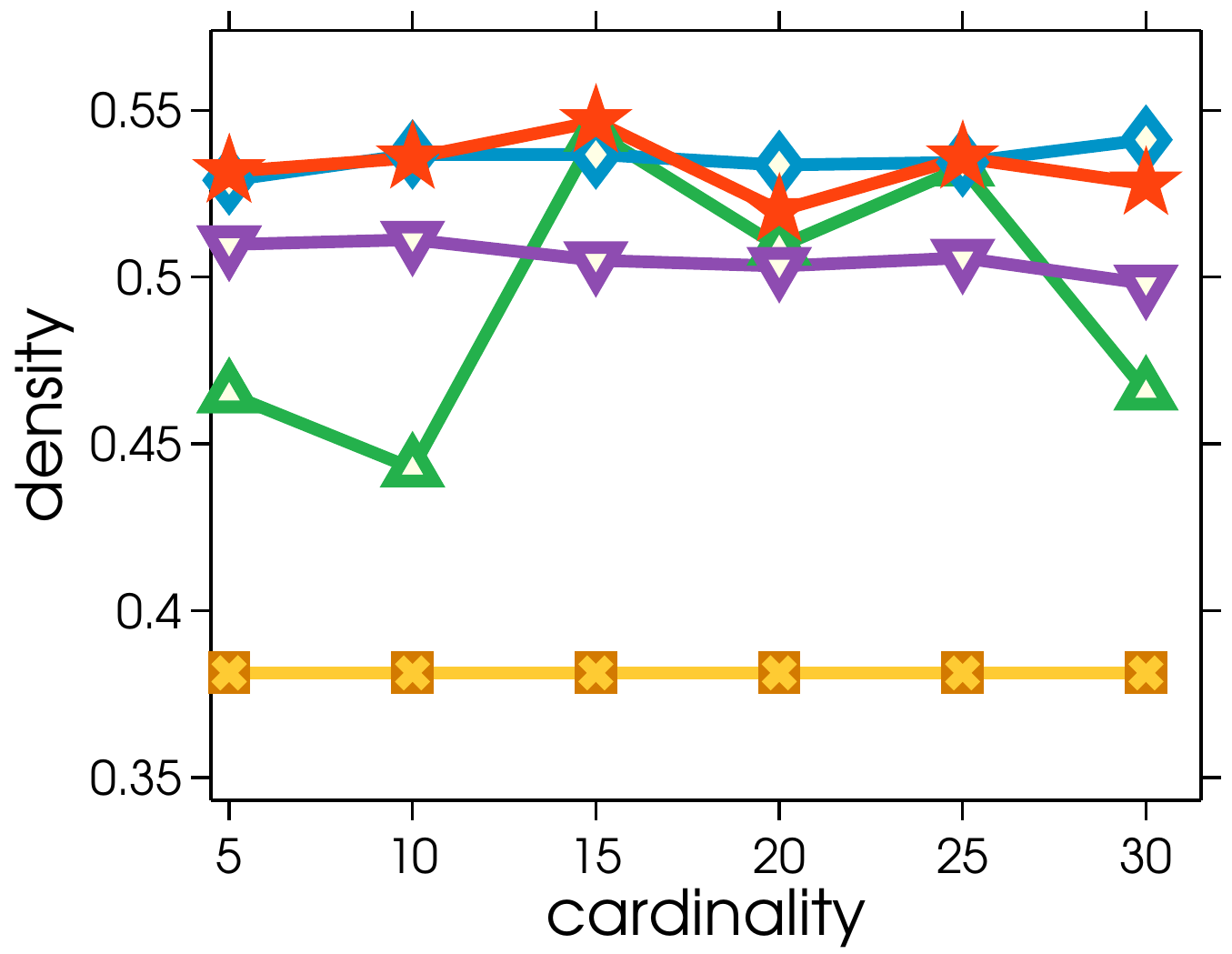}\caption{lesmis}\end{subfigure}\ghs
      \begin{subfigure}{\fourfigwid}\includegraphics[width=\textwidth,height=\imgheiclustering]{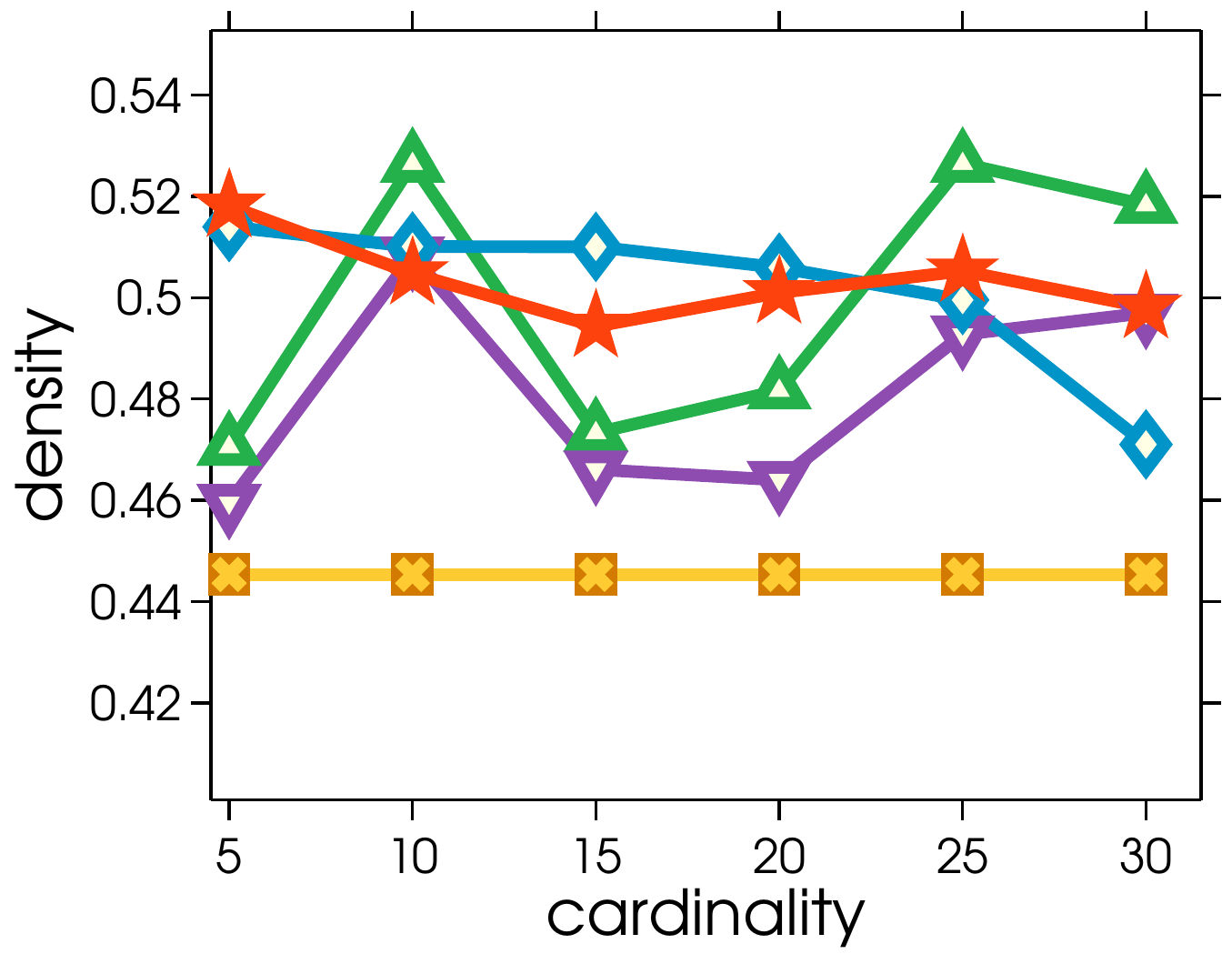}\caption{polbooks}\end{subfigure}\ghs
      \begin{subfigure}{\fourfigwid}\includegraphics[width=\textwidth,height=\imgheiclustering]{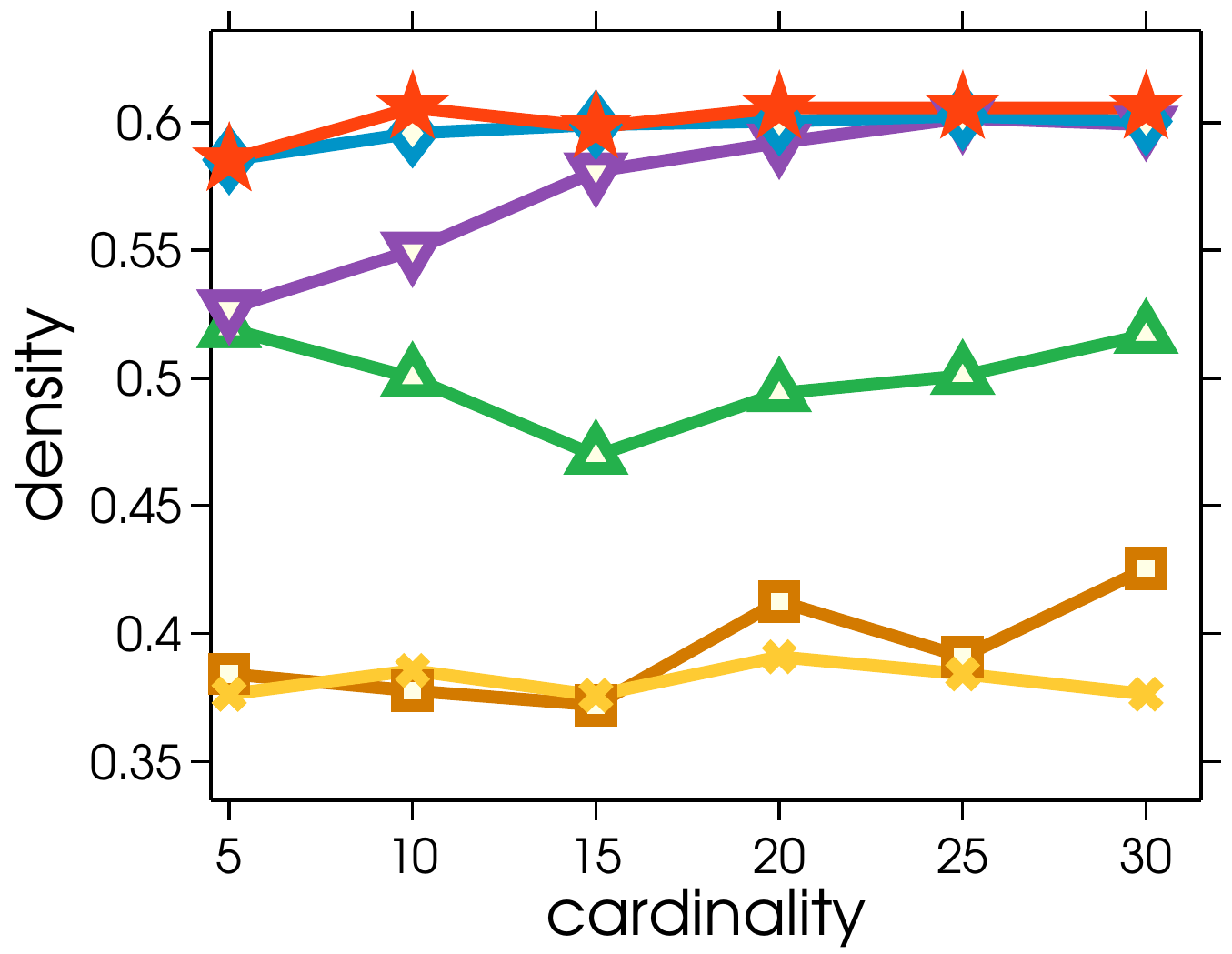}\caption{afootball}\end{subfigure}\ghs
      \begin{subfigure}{\fourfigwid}\includegraphics[width=\textwidth,height=\imgheiclustering]{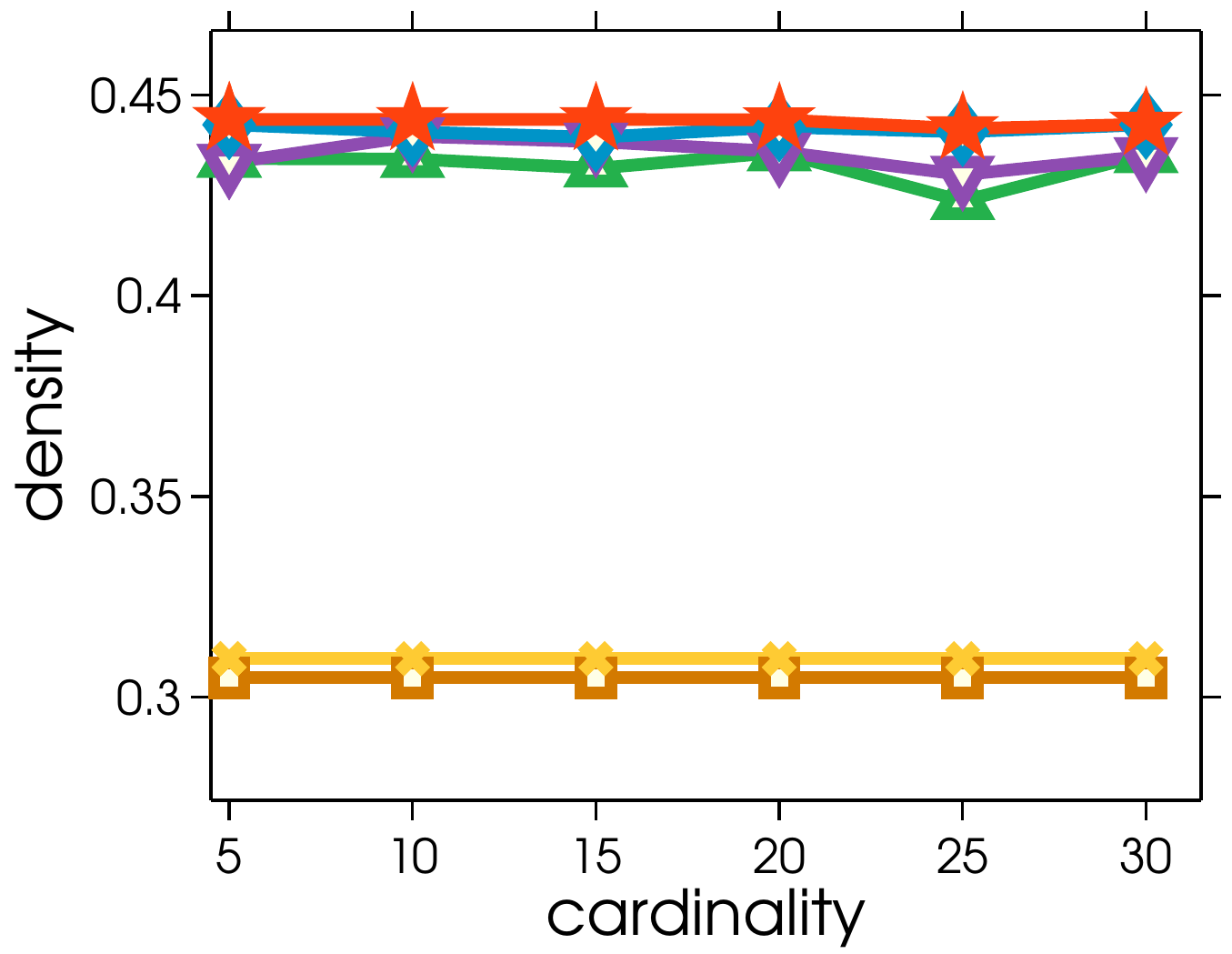}\caption{jazz}\end{subfigure}
\caption{Experimental results for modularity clustering.}
\label{fig:modularity}
\end{figure*}

\subsection{Modularity Clustering} \label{sect:mod}

Modularity was first introduced in \cite{newman2004finding} as a performance measure for the quality of community structure found by a clustering algorithm. Given a modularity matrix $\bbb{Q}\in \mathbb{R}^{n\times n}$ with $\bbb{Q}_{ij}=\bbb{W}_{ij} - deg(\bbb{v}_i)deg(\bbb{v}_j)/(2m)$, modularity clustering can be formulated as the following optimization problem \cite{newman2006modularity,brandes2008modularity,ChanY11}:
\beq \label{eq:modularity}
\min_{\bbb{X} \in \{-1,+1\}^{n\times n}}~\frac{1}{8m}tr(\bbb{X}^T\bbb{Q}\bbb{X}),~s.t.~\bbb{X}\bbb{1}=(2-k)\bbb{1}
\eeq
\noi Observing that $\bbb{Y} = (\bbb{X}+\bbb{1})/2\in \mathbb{R}^{n\times n}$, we obtain $\bbb{Y}\in \{0,1\}^{n\times n}$. Combining the linear constraint $\bbb{X}\bbb{1}=(2-k)\bbb{1}$, we have $\bbb{Y}\bbb{1}=\bbb{1}$ and  $\bbb{X}\bbb{X}^T$$=$$(2\bbb{Y}-\bbb{1})(2\bbb{Y}-\bbb{1})^T = 4 \bbb{Y} \bbb{Y}^T-3\bbb{1}\bbb{1}^T$. Based on these analyses, one can rewrite (\ref{eq:modularity}) as the following equivalent optimization problem:
\beq \label{eq:modularity2}
\min_{\bbb{Y} \in \{0,1\}^{n\times n}}~\frac{1}{8m}tr(\bbb{Y}^T\bbb{Q}\bbb{Y})+constant,~s.t.~\bbb{Y}\bbb{1}=\bbb{1}.\nn
\eeq
\bbb{Compared Methods.} We compare against 4 methods on 8 network data sets (See Table \ref{tab:clustering}). (i) Iterative rounding algorithm (IRA) \cite{ChanY11} in every stage solves a convex quadratic program and picks a fixed number of the vertices with largest values to assign to the cluster. However, such heuristic algorithm does not have any convergence guarantee. We use the code provided by the authors\footnote{\url{http://www.cse.ust.hk/~dyyeung/paper/publist.html}} and set the parameter $\rho=0.5$ in this method. (ii) LP relaxation solves a capped simplex problem. (iii) L2box-ADMM solves a spherical constrained optimization problem. (iv) Iterative Hard Thresholding (IHT) considers setting the current solution to the indicator vector of its top-k entries while decreasing the objective function. Due to its suboptimal performance in our previous experiment, LP relaxation is used as its initialization.

\bbb{Experimental Results.} Several observations can be drawn from Figure \ref{fig:modularity}. (i) IHT does not necessarily improve upon the LP relaxation method. (ii) IRA consistently outperforms IHT and LP in all the experiments. (iii) L2box-ADMM gives comparable result to IRA. (iv) Our proposed methods generally outperform IRA and L2box-ADM in the experiments.

\begin{figure*}[!t]
\captionsetup[subfigure]{justification=centering}
    \centering
      \begin{subfigure}{\fourfigwid}\includegraphics[width=\textwidth,height=\imgheimrf]{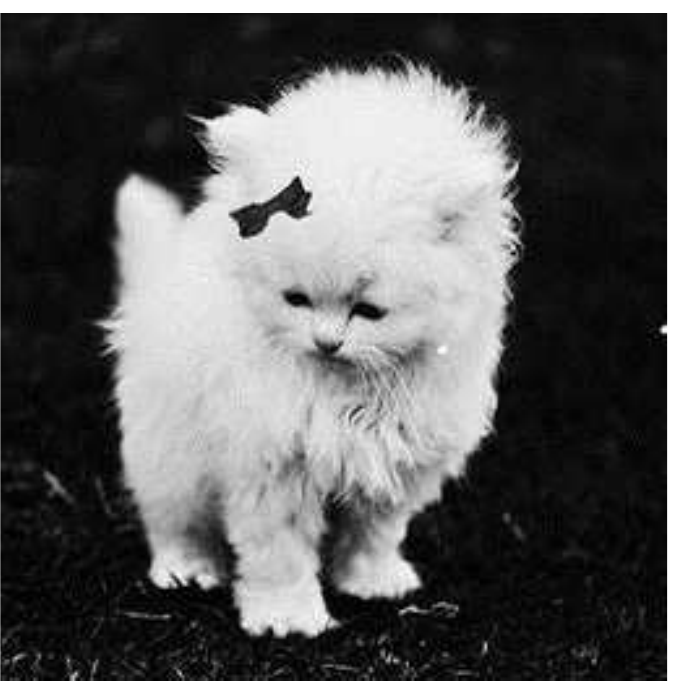}\caption{cat image\\ \bigskip}\end{subfigure}\ghs
      \begin{subfigure}{\fourfigwid}\includegraphics[width=\textwidth,height=\imgheimrf]{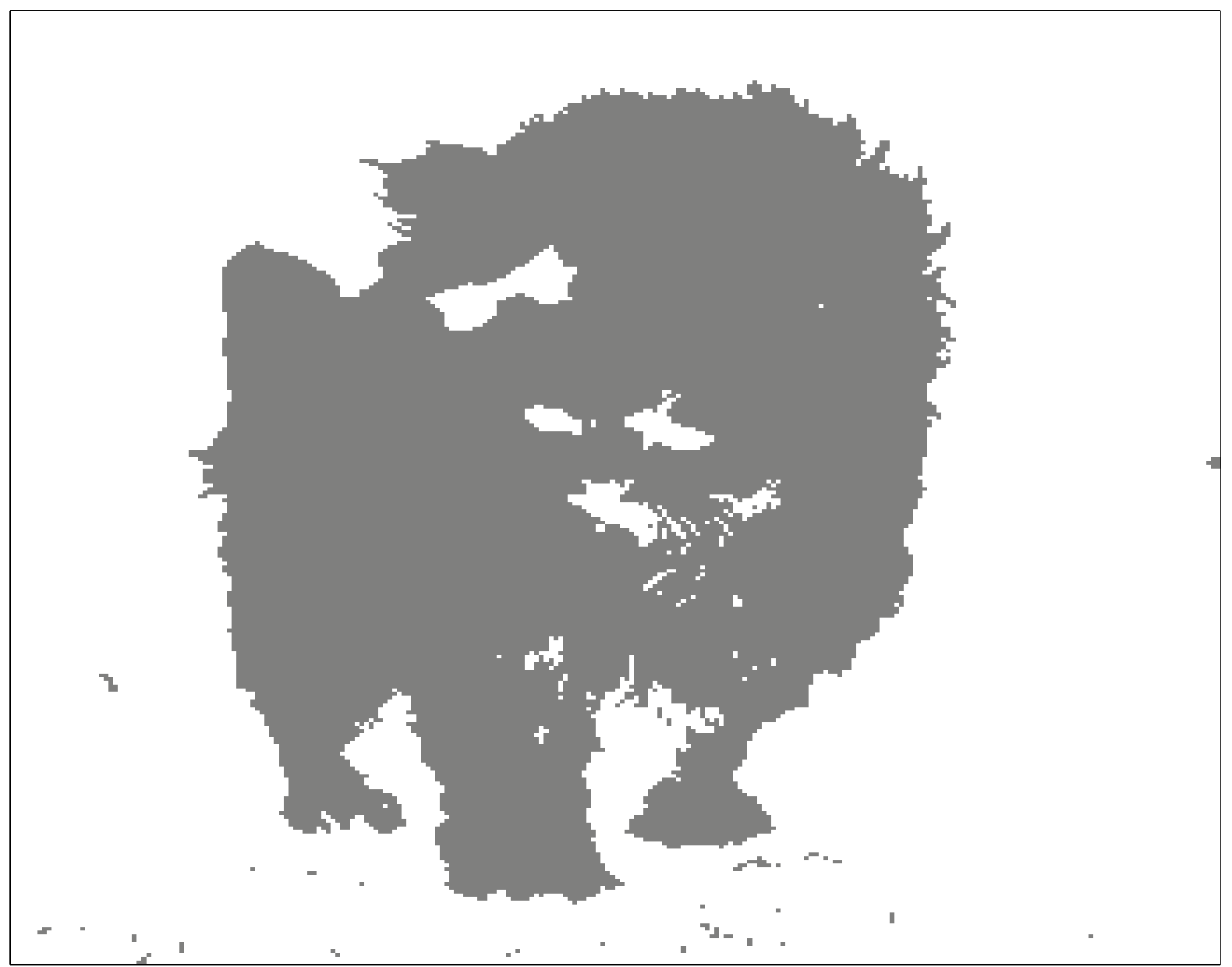}\caption{unary image\\$f =143599.88$}\end{subfigure}\ghs
      \begin{subfigure}{\fourfigwid}\includegraphics[width=\textwidth,height=\imgheimrf]{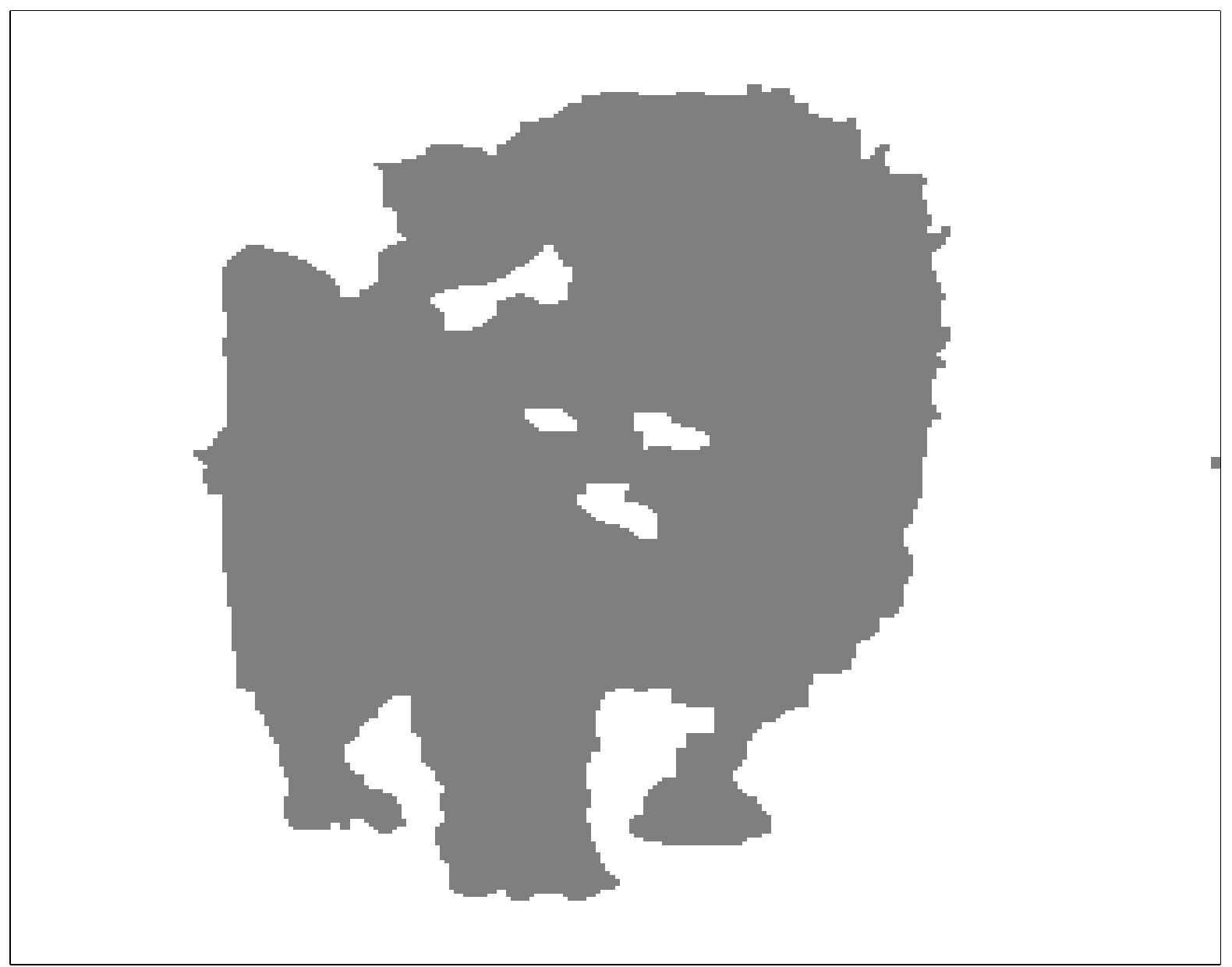}\caption{Graph Cut\\$f =143086.58$}\end{subfigure}\ghs
      \begin{subfigure}{\fourfigwid}\includegraphics[width=\textwidth,height=\imgheimrf]{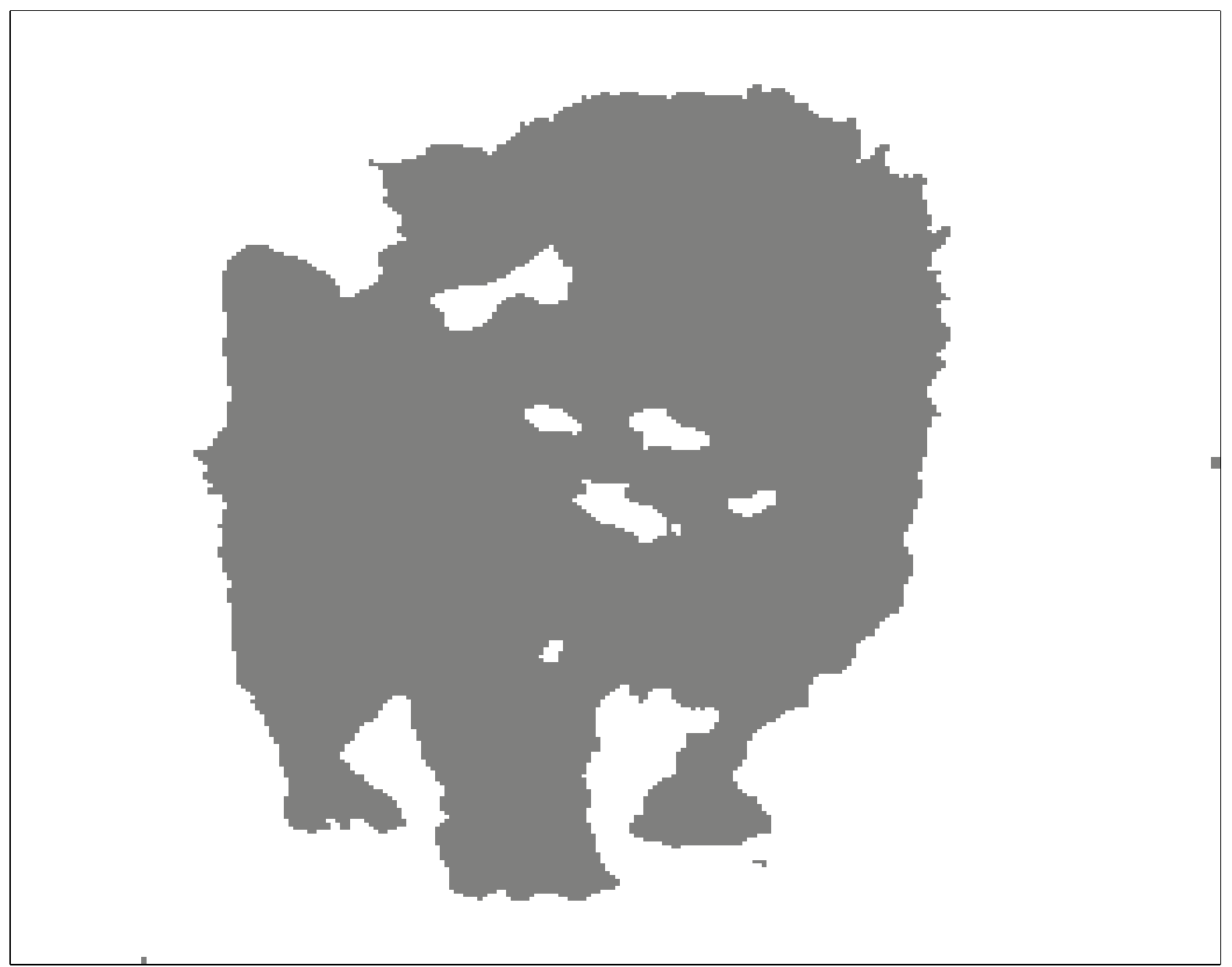}\caption{LP\\$f =143145.88$}\end{subfigure}

      \begin{subfigure}{\fourfigwid}\includegraphics[width=\textwidth,height=\imgheimrf]{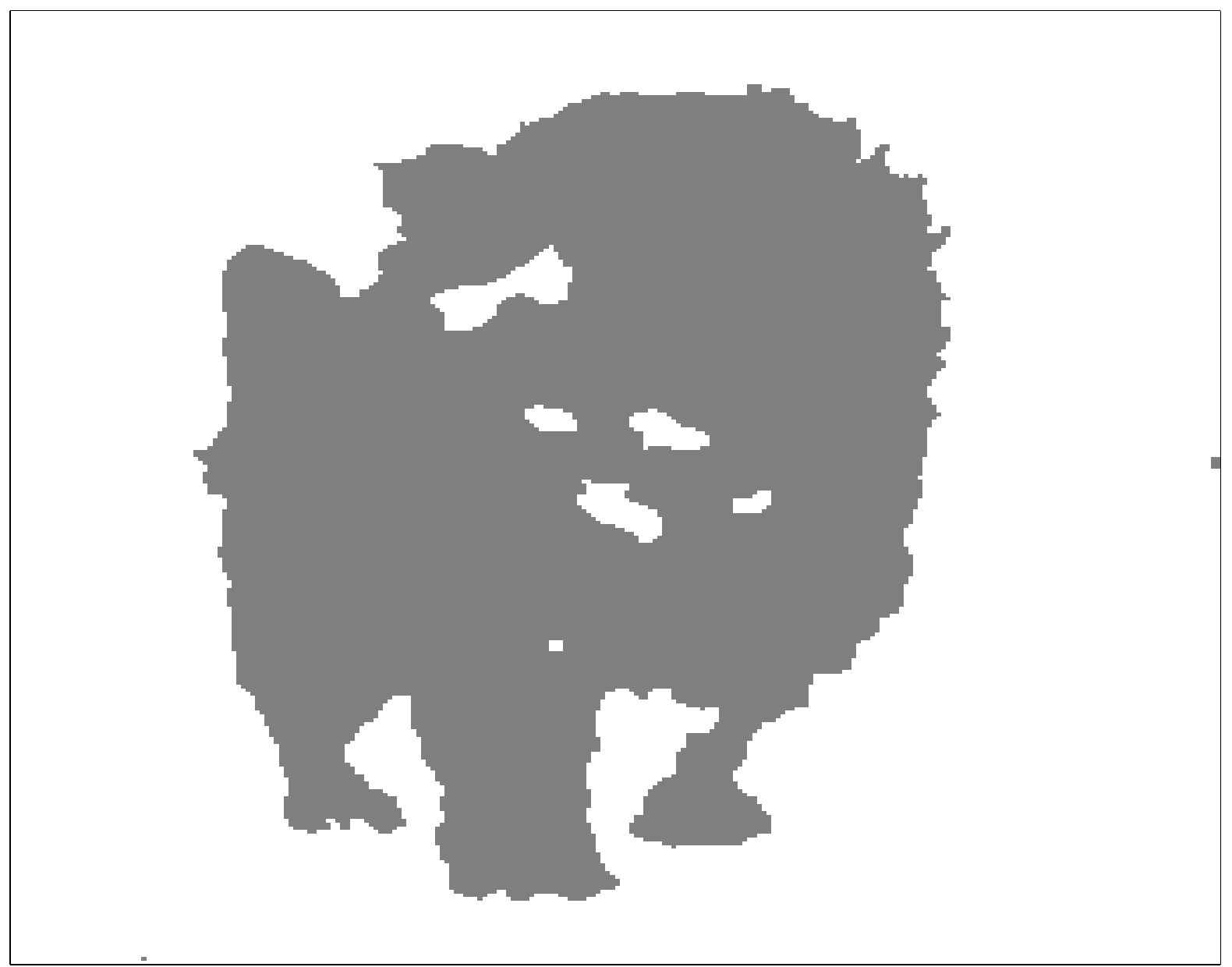}\caption{L2box-ADMM \\$f =143109.38$}\end{subfigure}\ghs
      \begin{subfigure}{\fourfigwid}\includegraphics[width=\textwidth,height=\imgheimrf]{cat_diadm-eps-converted-to.pdf}\caption{QPM-$\ell_0$ norm\\$f =143162.46$}\end{subfigure}\ghs
      \begin{subfigure}{\fourfigwid}\includegraphics[width=\textwidth,height=\imgheimrf]{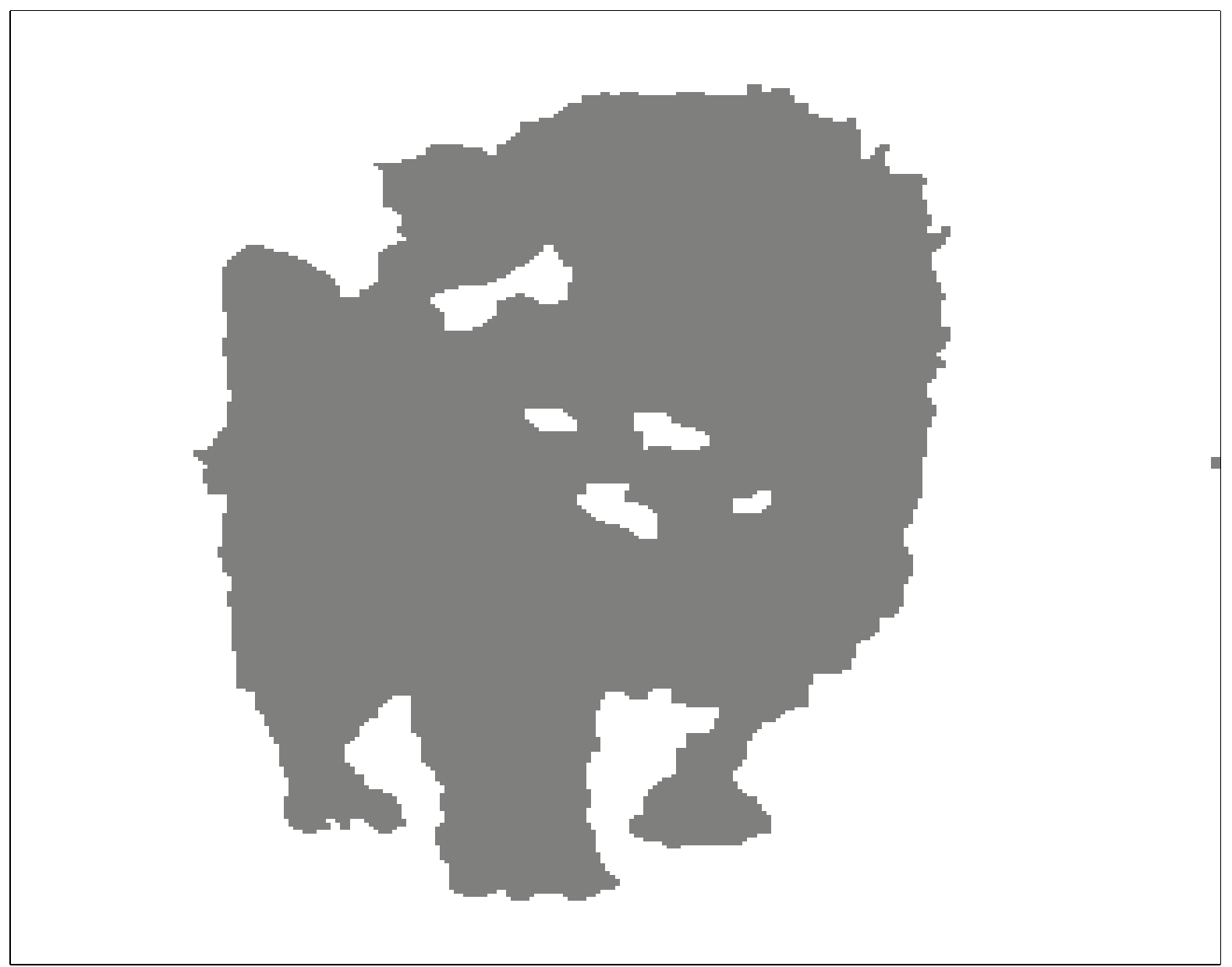}\caption{MPEC-EPM\\$f =143091.91$}\end{subfigure}\ghs
      \begin{subfigure}{\fourfigwid}\includegraphics[width=\textwidth,height=\imgheimrf]{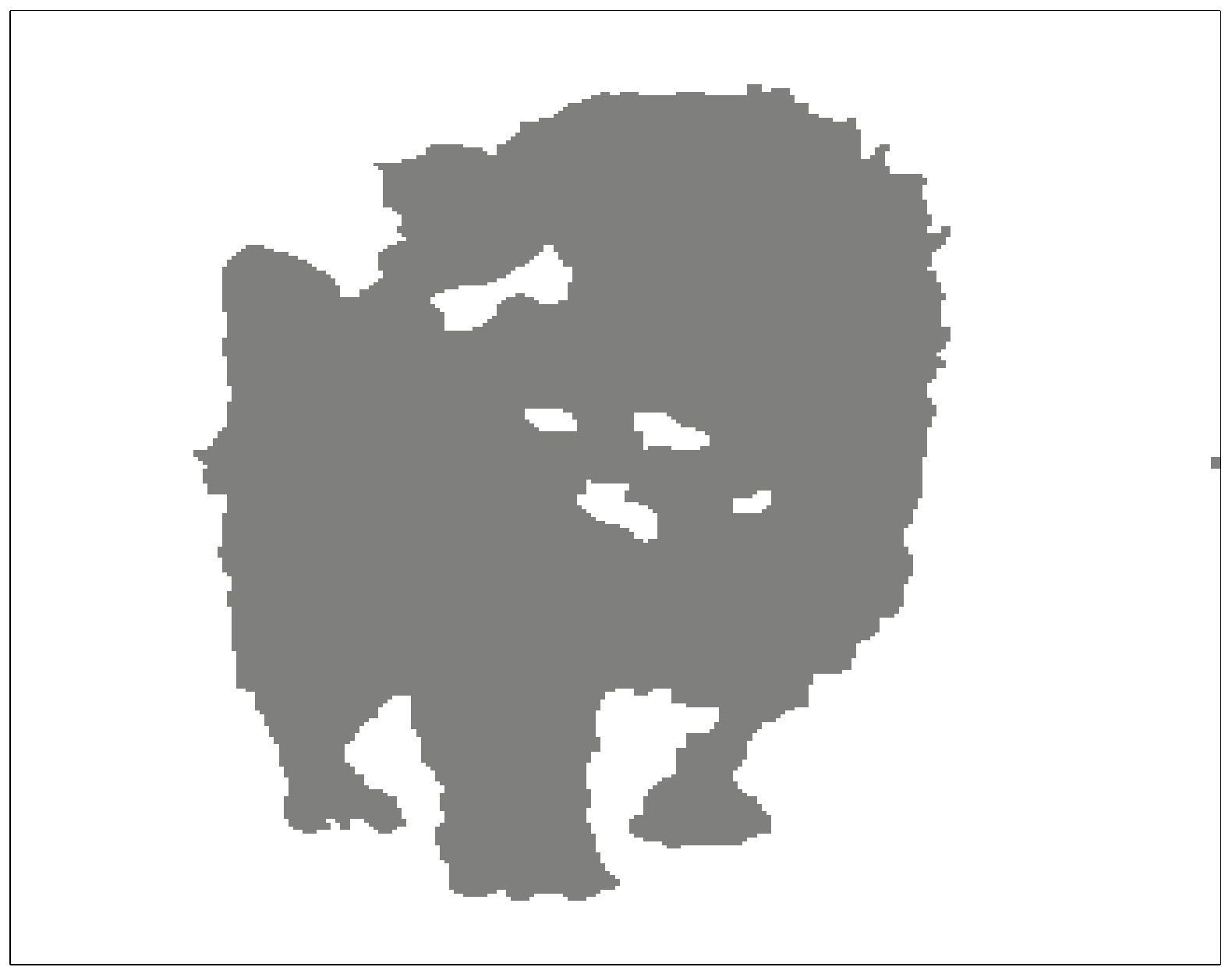}\caption{MPEC-ADM\\$f =143093.83$}\end{subfigure}

      \begin{subfigure}{\fourfigwid}\includegraphics[width=\textwidth,height=\imgheimrf]{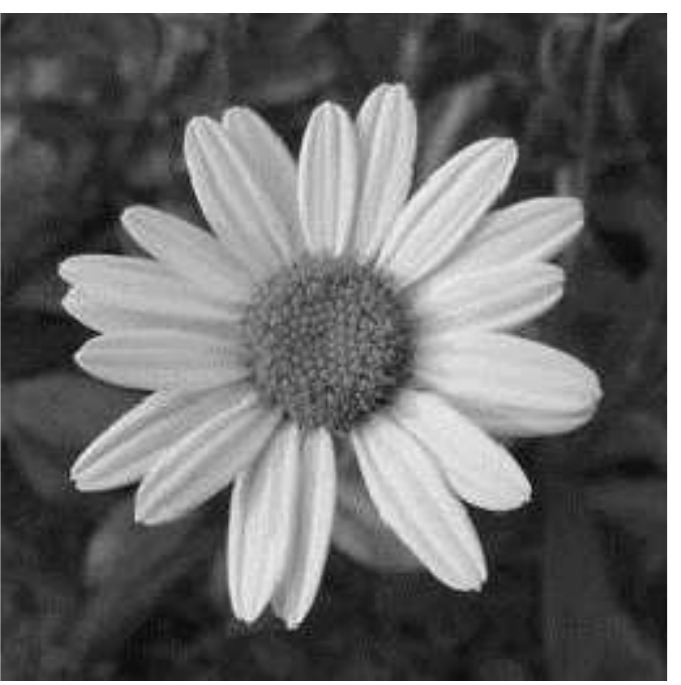}\caption{flower image\\ \bigskip }\end{subfigure}\ghs
      \begin{subfigure}{\fourfigwid}\includegraphics[width=\textwidth,height=\imgheimrf]{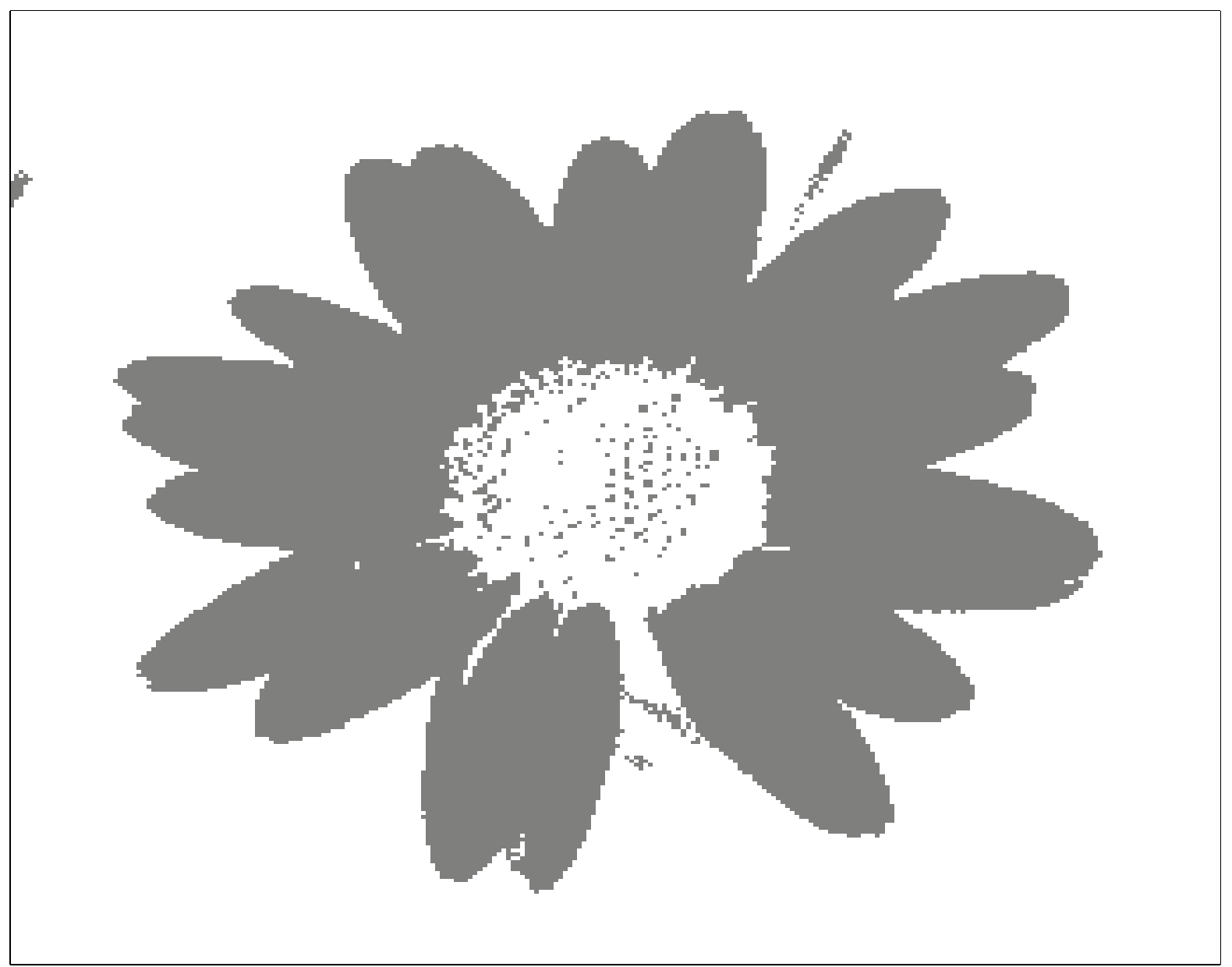}\caption{unary image\\$f = 37152.66$}\end{subfigure}\ghs
      \begin{subfigure}{\fourfigwid}\includegraphics[width=\textwidth,height=\imgheimrf]{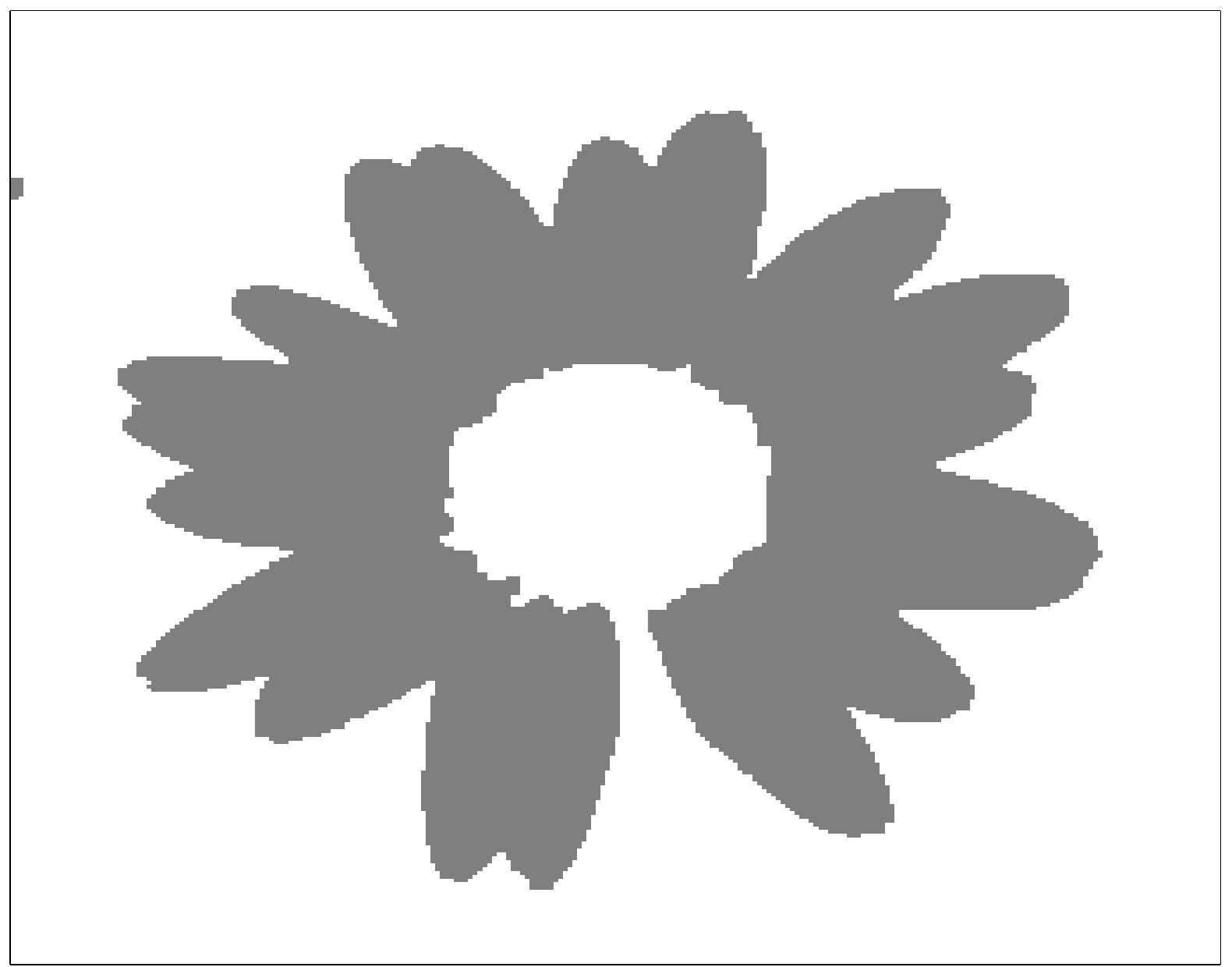}\caption{Graph Cut\\$f = 36616.50$}\end{subfigure}\ghs
      \begin{subfigure}{\fourfigwid}\includegraphics[width=\textwidth,height=\imgheimrf]{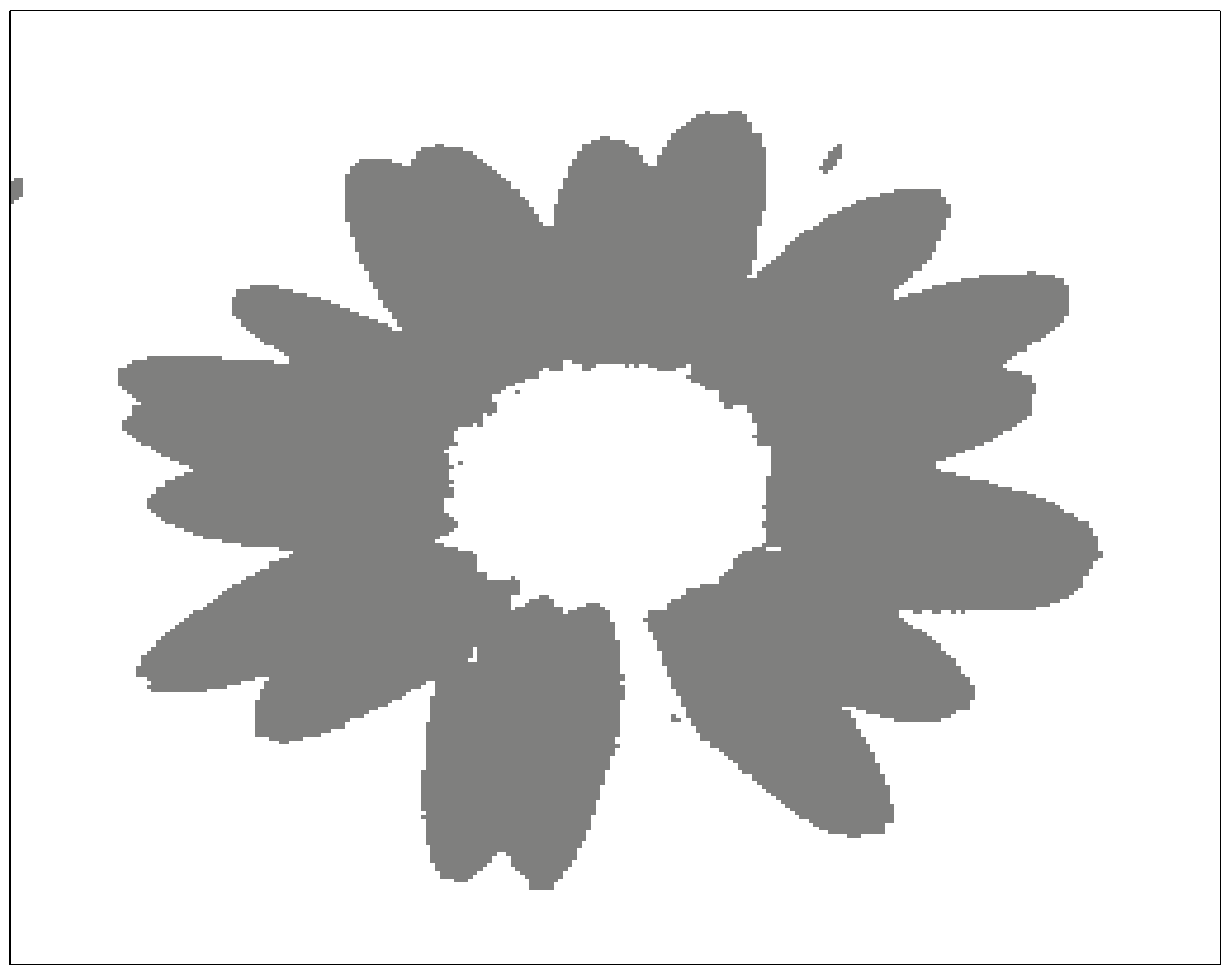}\caption{LP\\$f = 36653.40$}\end{subfigure}

      \begin{subfigure}{\fourfigwid}\includegraphics[width=\textwidth,height=\imgheimrf]{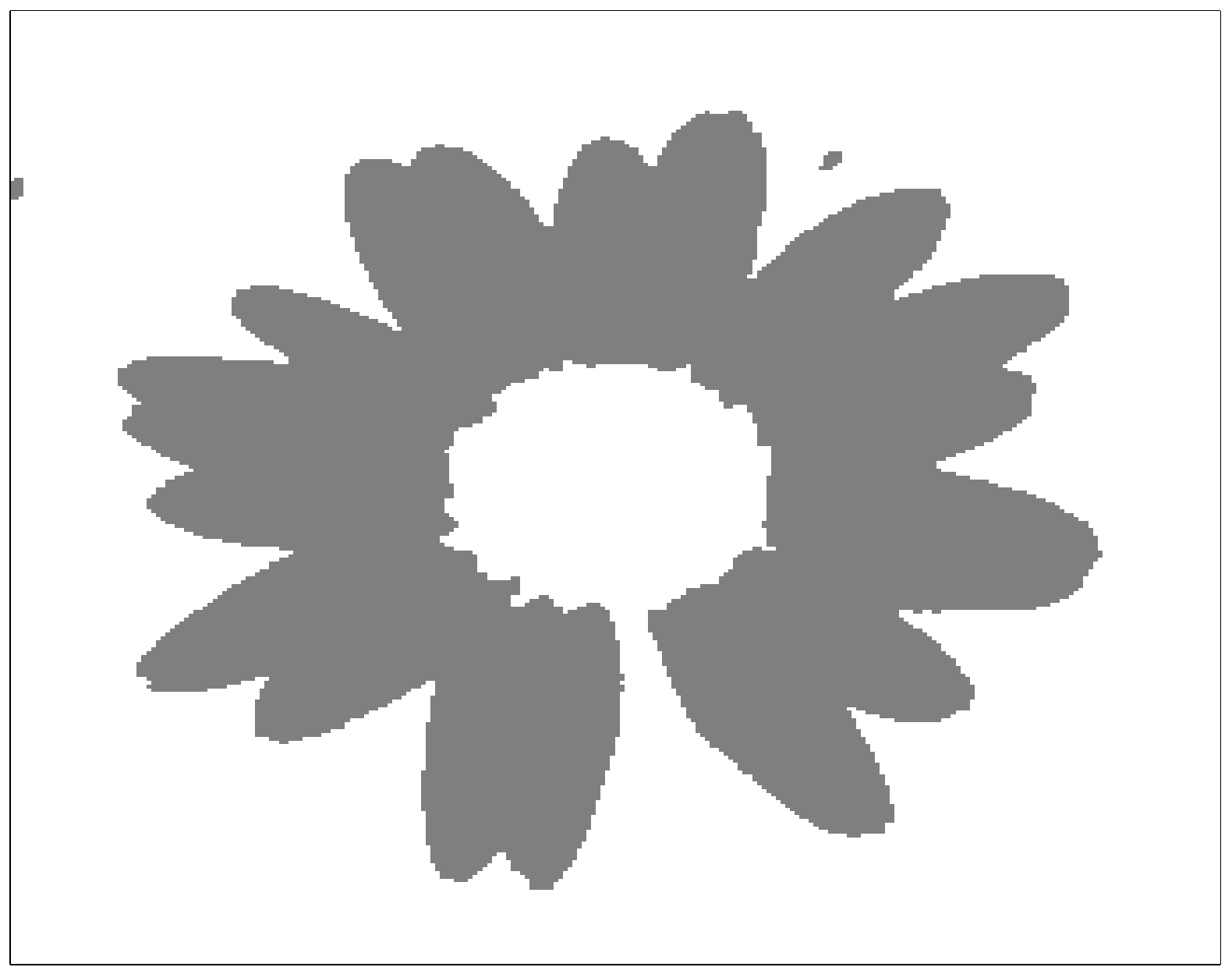}\caption{L2box-ADMM\\$f = 36626.17$}\end{subfigure}\ghs
      \begin{subfigure}{\fourfigwid}\includegraphics[width=\textwidth,height=\imgheimrf]{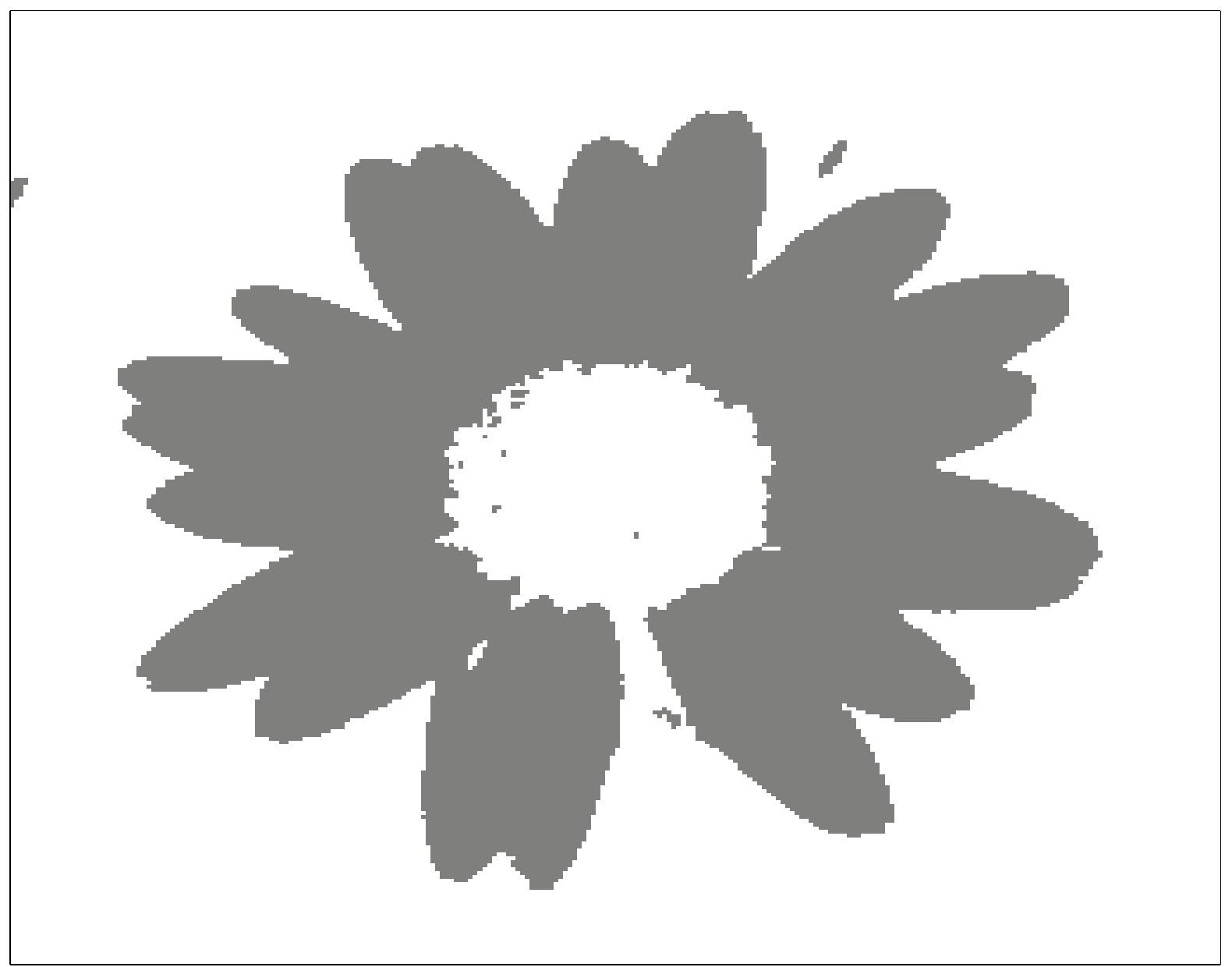}\caption{L0-QPM\\$f = 36689.11$}\end{subfigure}\ghs
      \begin{subfigure}{\fourfigwid}\includegraphics[width=\textwidth,height=\imgheimrf]{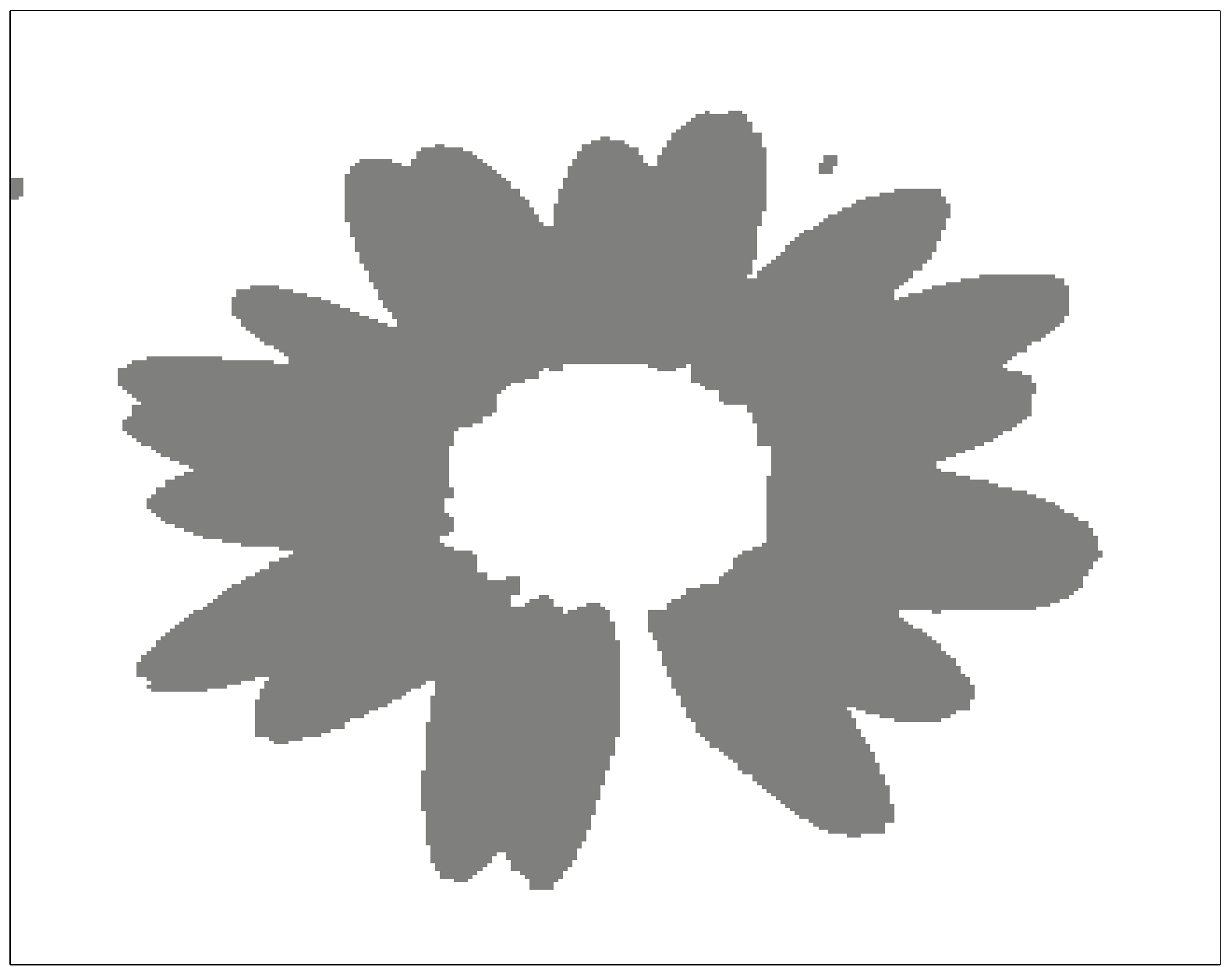}\caption{MPEC-EPM\\$f = 36619.97$}\end{subfigure}\ghs
      \begin{subfigure}{\fourfigwid}\includegraphics[width=\textwidth,height=\imgheimrf]{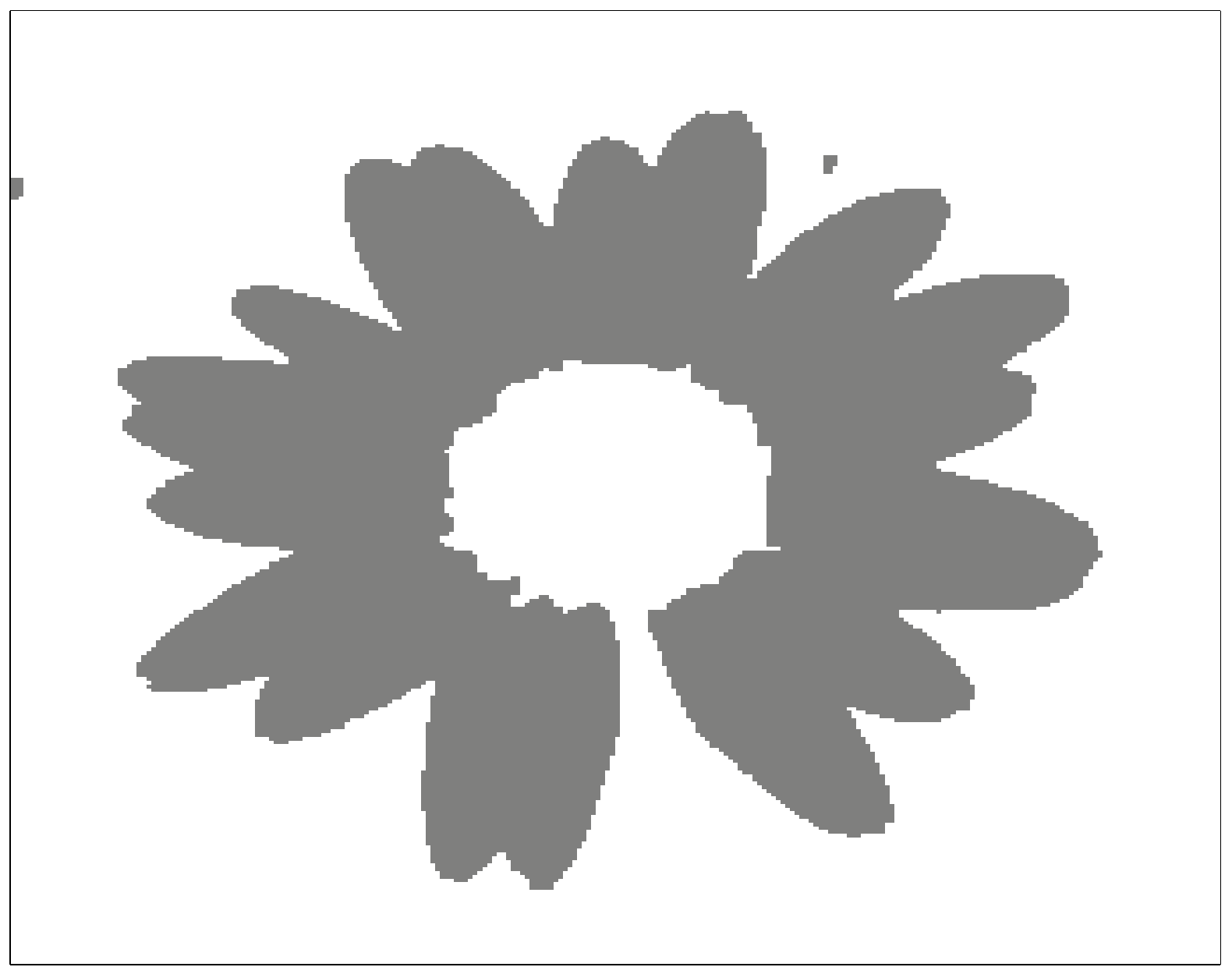}\caption{MPEC-ADM\\$f= 36621.70$}\end{subfigure}
\caption{Markov random fields on `cat' image and `flower' images.}
\label{fig:mrf}
\end{figure*}

\begin{table}[!h]
\small\tiny \footnotesize  \scriptsize
\caption{CPU time (in seconds) comparisons.}
\center
\begin{tabular}{|c|c|c|c|c|}
  \hline
 Graph & LP & L2box-ADM & MPEC-EPM& MPEC-ADM \\
  \hline
wordassociation & 1 & 7 & 2 & 13   \\
 enron & 2 & 40 & 29  & 85 \\
 uk-2007-05 & 6 & 75& 65& 77\\
cnr-2000 & 16& 210 & 209 & 245 \\
 dblp-2010 & 15& 234& 282& 253 \\
 in-2004 & 79 & 834& 1023& 1301\\
 amazon-2008 & 49 & 501& 586& 846 \\
 dblp-2011 & 59&554&621&1007 \\
\hline
\end{tabular}
\label{tab:cpu}
\end{table}
\subsection{Markov Random Fields} \label{sect:mrf}
The Markov Random Field (MRF) optimization \cite{Boykov2001,cour2007solving,HuangCG14} is widely used in many labeling applications, including image restoration, edge detection, and image segmentation. Generally speaking, it involves solving the following problem:
\beq
\min_{\bbb{x}\in\{0,1\}^{n}}~\frac{1}{2}\bbb{x}^T\bbb{L}\bbb{x} + \bbb{x}^T\bbb{b}\nn
\eeq
\noi where $\bbb{b}\in\mathbb{R}^n$ is determined by the unary term defined for the graph and $\bbb{L}\in\mathbb{R}^{n\times n}$ is the Laplacian matrix, which is based on the binary term relating pairs of graph nodes together. The quadratic term is usually considered a smoothness prior on the node labels.

\bbb{Compared Methods.} We perform image segmentation on the `cat' image and `flower' images. (i) Graph cut method \cite{Boykov2001} is included in our experiments. This method is known to achieve the global optimal solution for this specific class of binary problem. (ii) LP relaxation solves a box constrained quadratic programming problem. (iii) L2box-ADMM solves the $\ell_2$ box non-separable reformulation directly using classical alternating direction method of multipliers (ADMM) \cite{WuG16a}. (iv) L0-QPM norm solves the semi-continuous $\ell_0$ norm reformulation of the binary optimization problem by quadratic penalty method \cite{LuZ13,yuan2016l0mpec}.


\bbb{Experimental Results.} Figure \ref{fig:mrf} demonstrates a qualitative result for image segmentation. MPEC-EPM and MPEC-ADM produce solutions that are very close to the globally optimal one. Moreover, both our methods achieve lower objectives than the other compared methods.

\subsection{Convergence Curve and Computational Efficiency}
This subsection demonstrates the convergence curve and computational efficiency of the proposed algorithms. We only report the results on dense subgraph discovery.

\begin{figure*}
\captionsetup[subfigure]{justification=centering}
    \centering
      \begin{subfigure}{0.5\textwidth}\includegraphics[width=\textwidth,height=13pt]{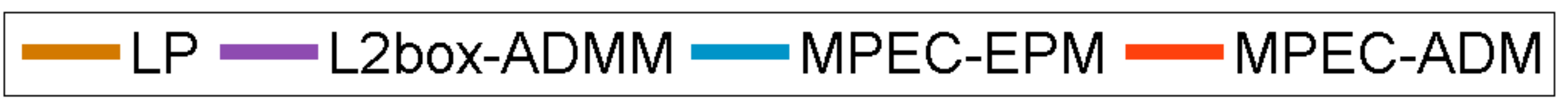}\end{subfigure}

      \vspace{3pt}

      \begin{subfigure}{\fourfigwid}\includegraphics[width=\textwidth,height=\imgheiconv]{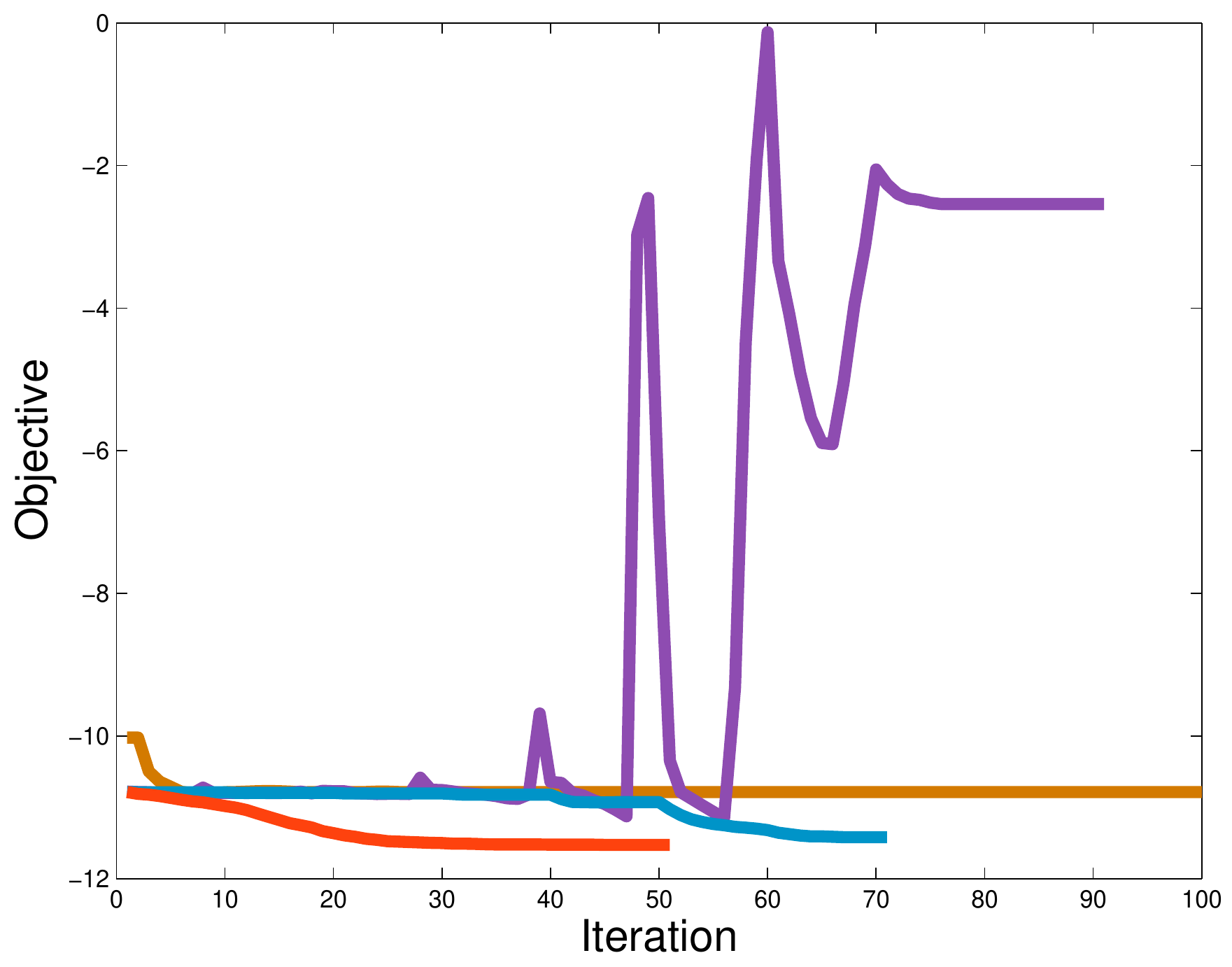}\caption{karate}\end{subfigure}\ghs
      \begin{subfigure}{\fourfigwid}\includegraphics[width=\textwidth,height=\imgheiconv]{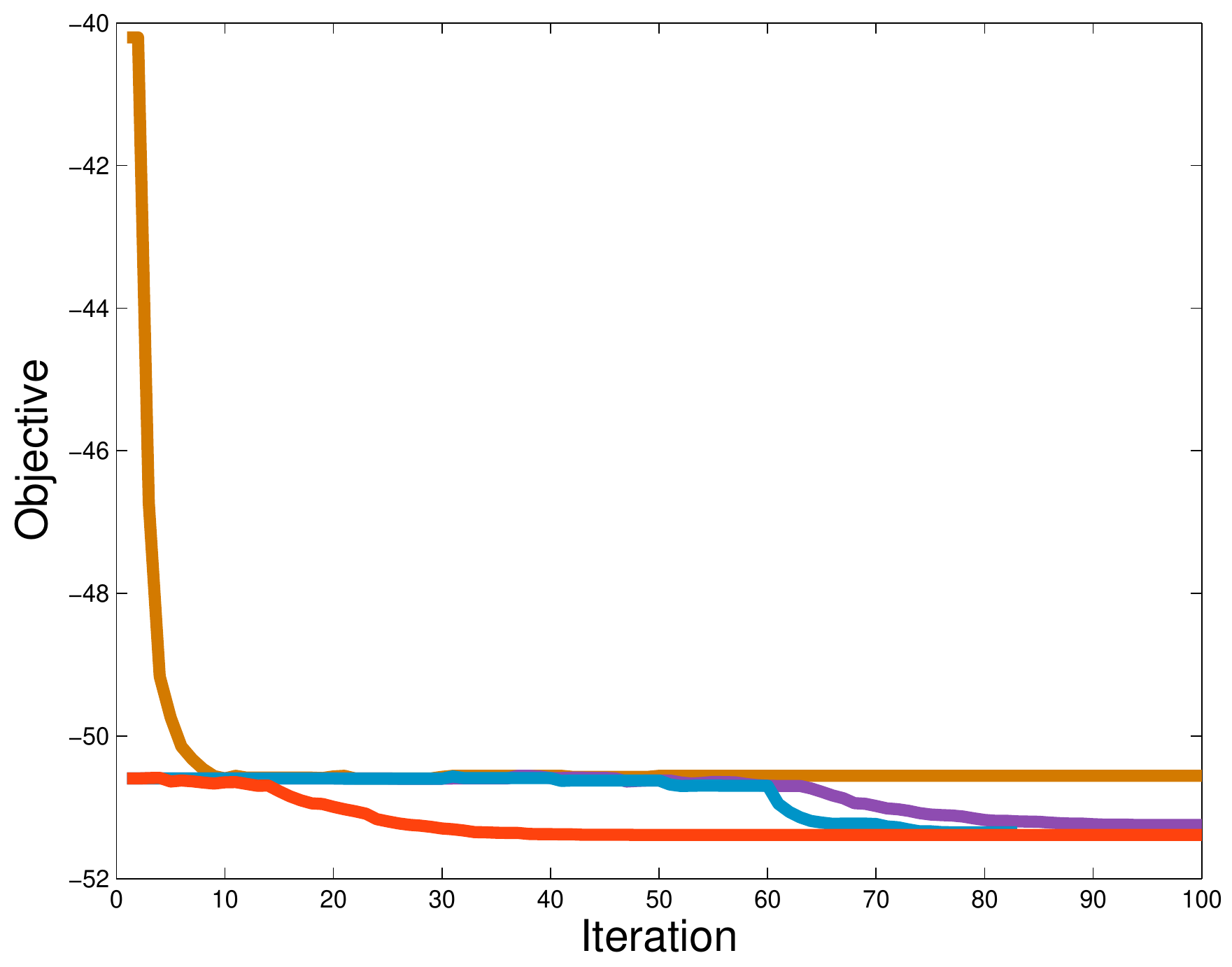}\caption{collab}\end{subfigure}\ghs
      \begin{subfigure}{\fourfigwid}\includegraphics[width=\textwidth,height=\imgheiconv]{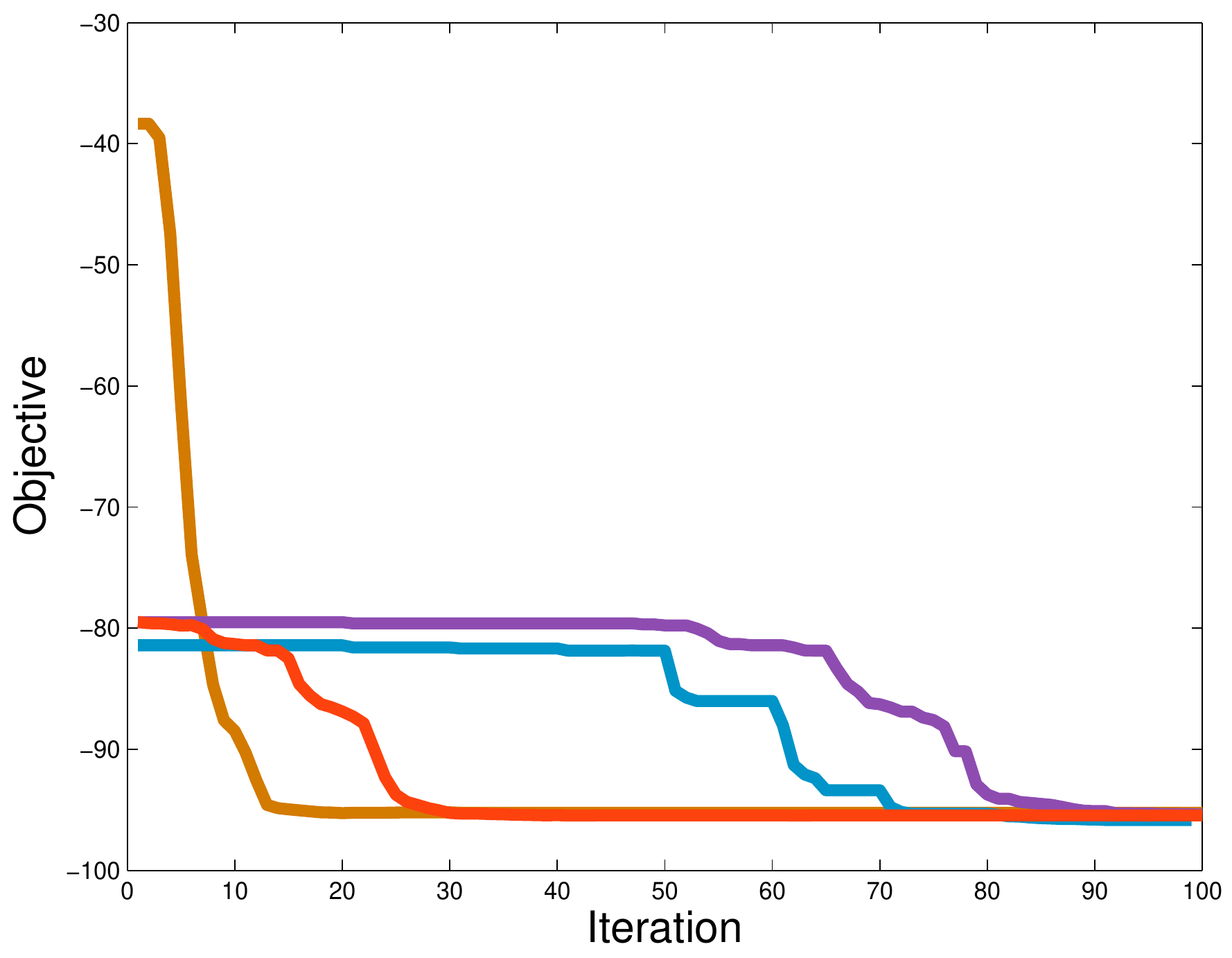}\caption{dolphins}\end{subfigure}\ghs
      \begin{subfigure}{\fourfigwid}\includegraphics[width=\textwidth,height=\imgheiconv]{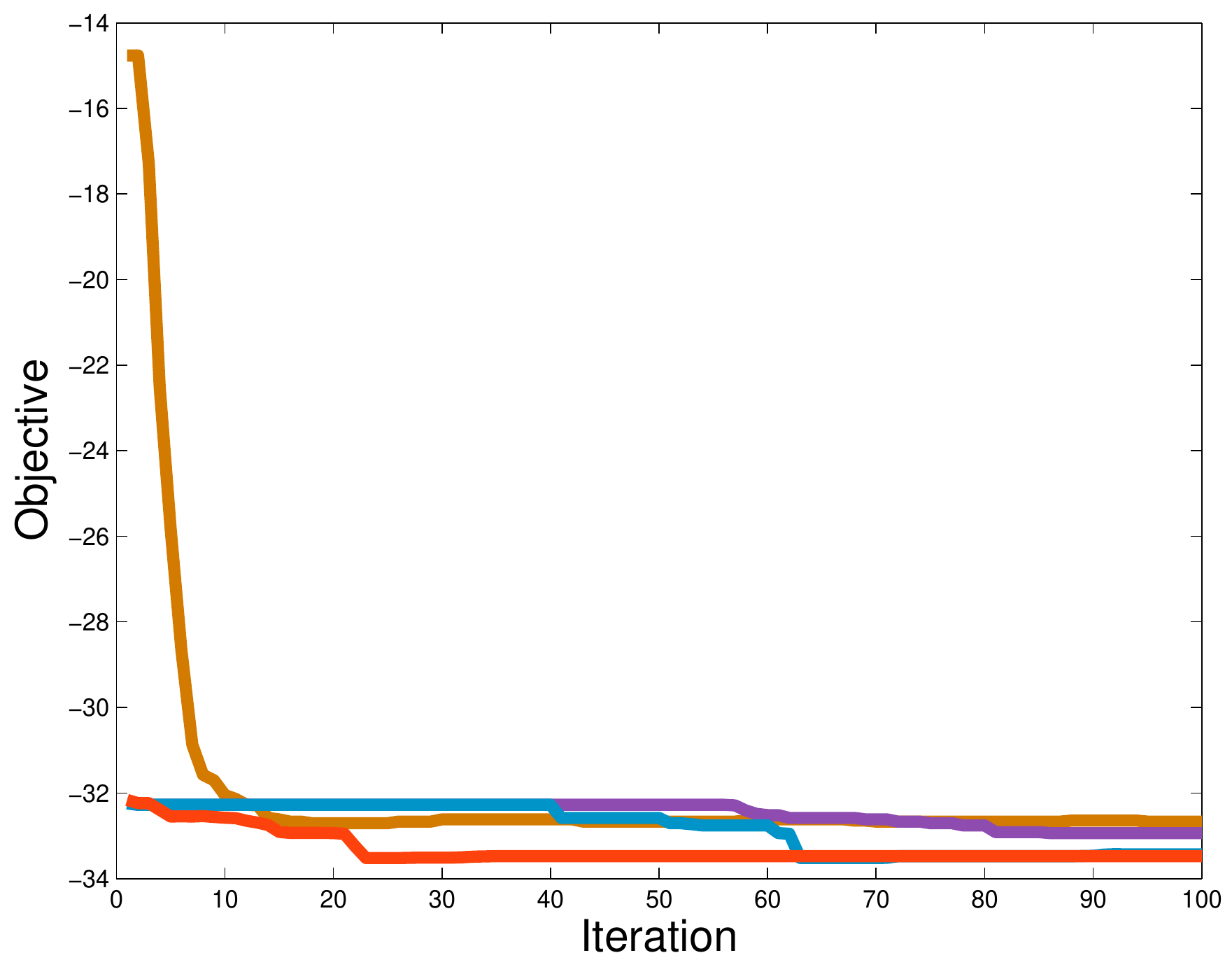}\caption{email}\end{subfigure}

      \begin{subfigure}{\fourfigwid}\includegraphics[width=\textwidth,height=\imgheiconv]{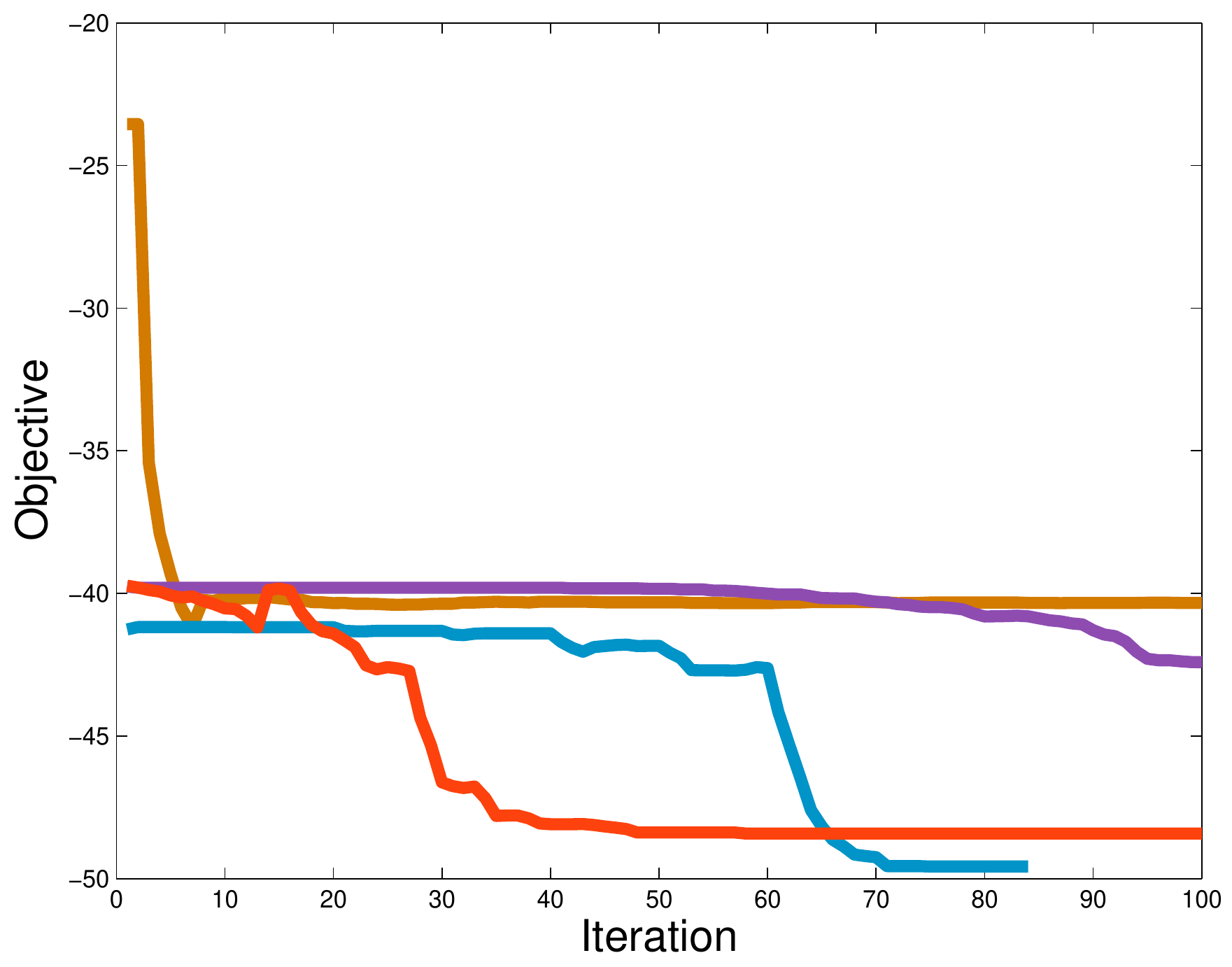}\caption{lesmis}\end{subfigure}\ghs
      \begin{subfigure}{\fourfigwid}\includegraphics[width=\textwidth,height=\imgheiconv]{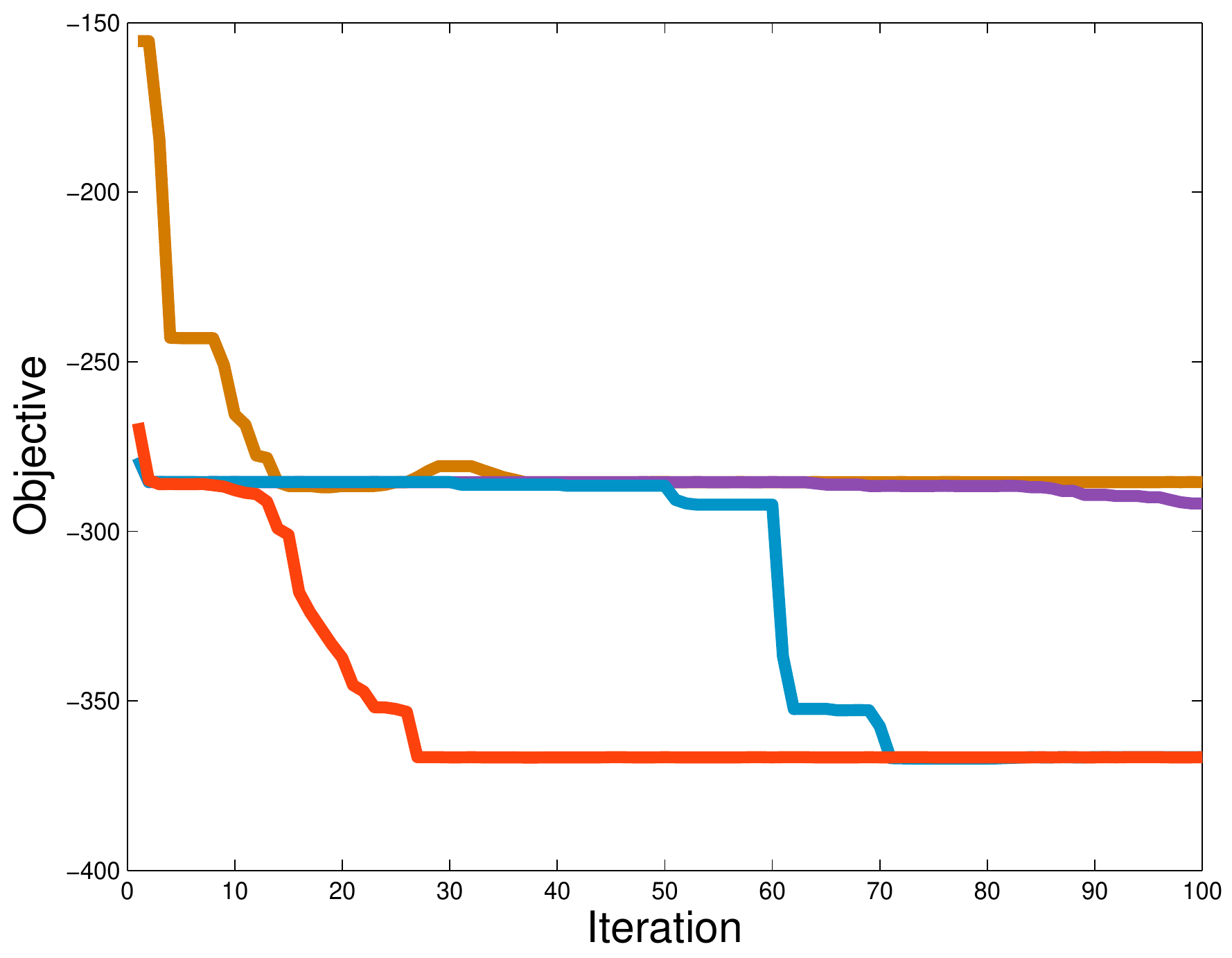}\caption{polbooks}\end{subfigure}\ghs
      \begin{subfigure}{\fourfigwid}\includegraphics[width=\textwidth,height=\imgheiconv]{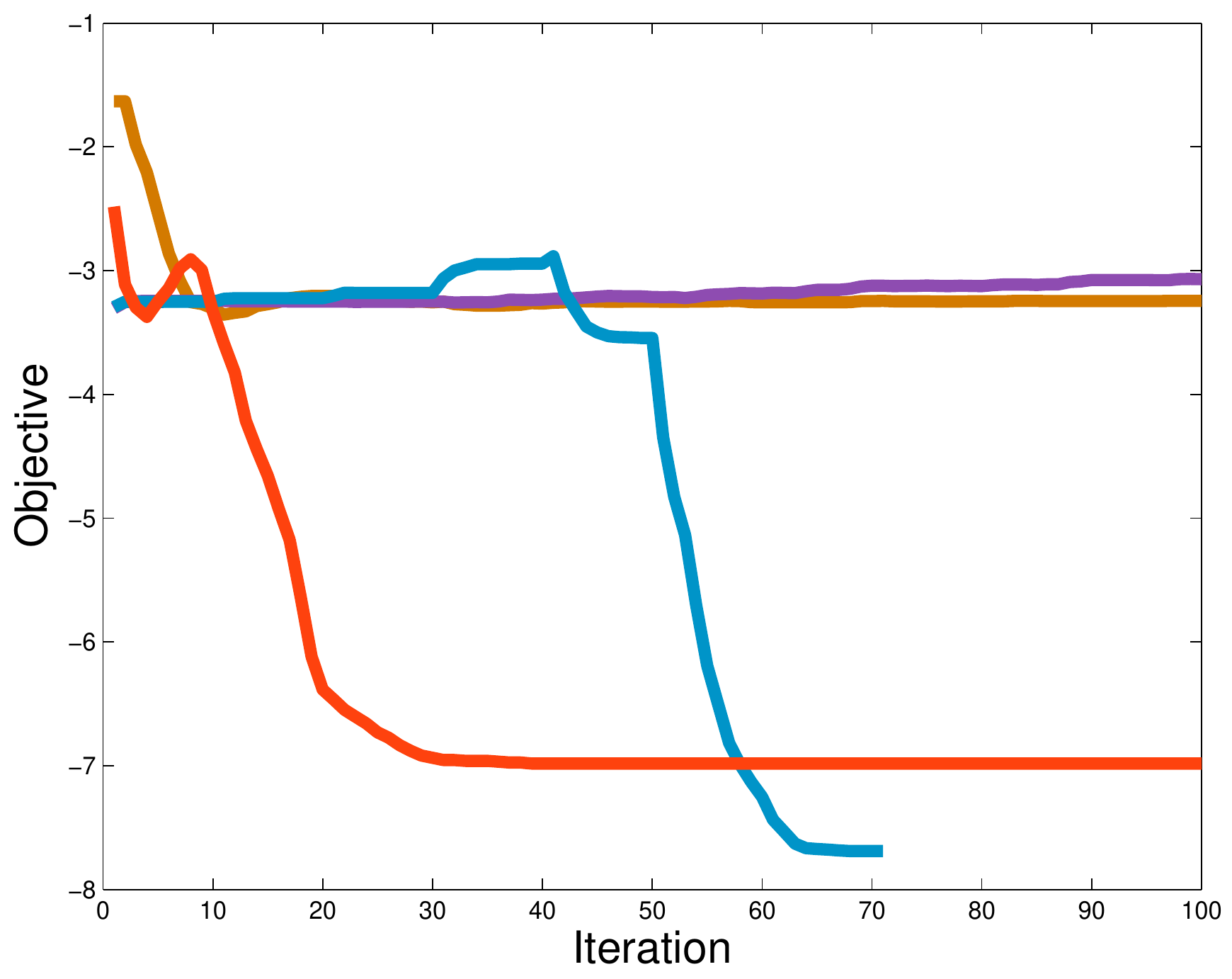}\caption{afootball}\end{subfigure}\ghs
      \begin{subfigure}{\fourfigwid}\includegraphics[width=\textwidth,height=\imgheiconv]{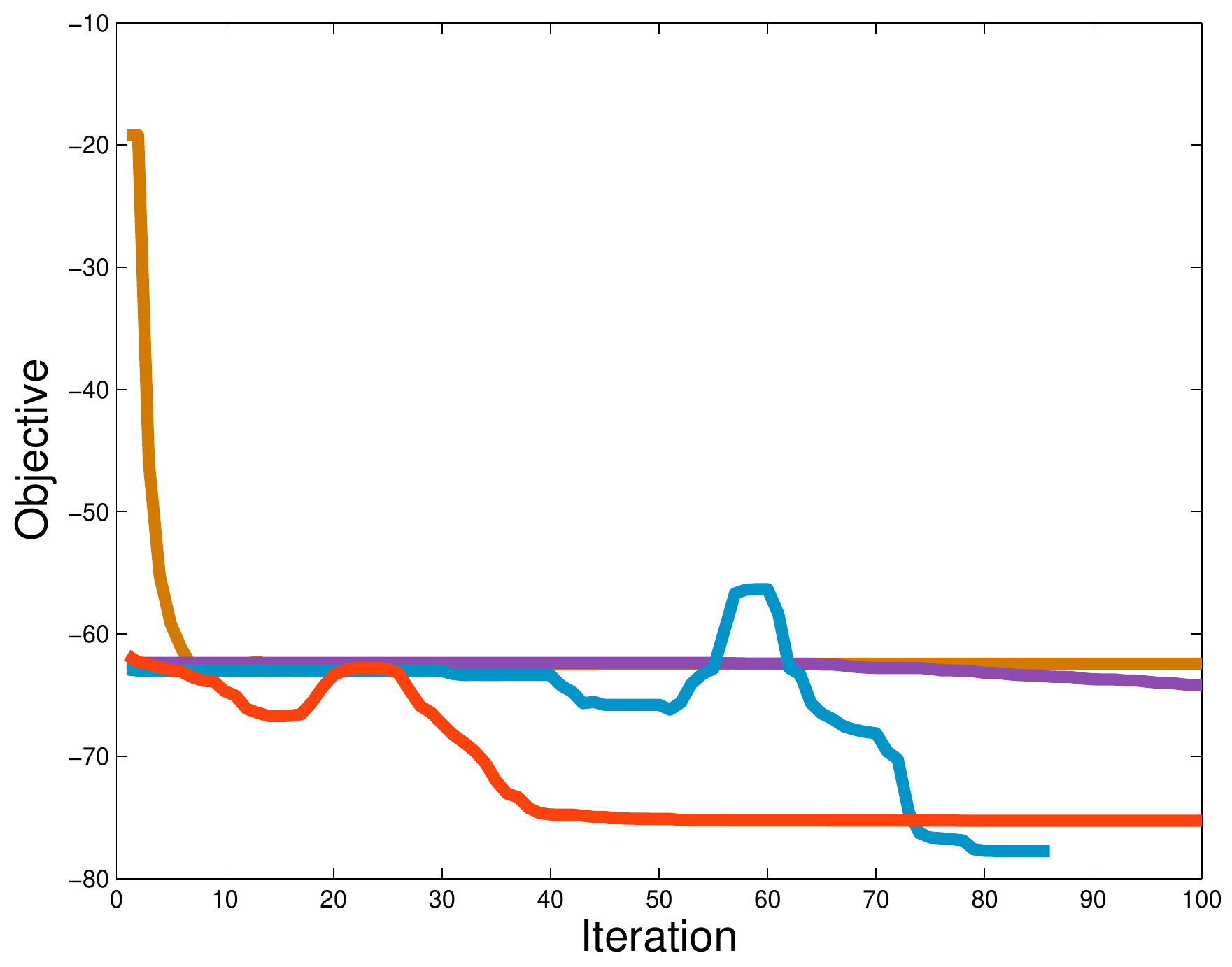}\caption{jazz}\end{subfigure}

\caption{Convergence curve for dense subgraph discovery on different datasets with $k=1000$.}
\label{fig:convergence:1}
\end{figure*}

\begin{figure*}
\captionsetup[subfigure]{justification=centering}
    \centering
      \begin{subfigure}{0.5\textwidth}\includegraphics[width=\textwidth,height=12pt]{fig_subgraph_convergence_legend-eps-converted-to.pdf}\end{subfigure}

      \vspace{3pt}

      \begin{subfigure}{\fourfigwid}\includegraphics[width=\textwidth,height=\imgheiconv]{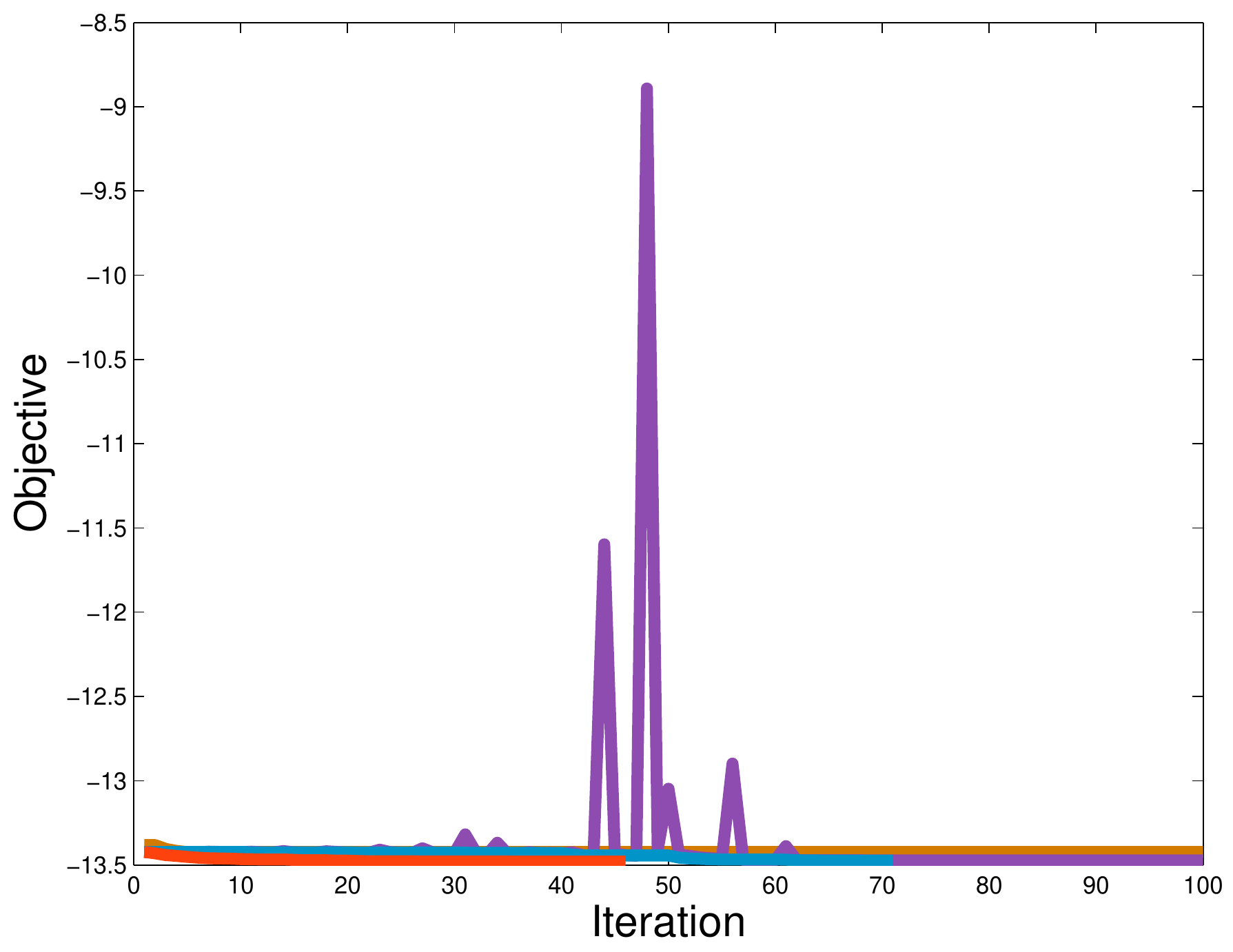}\caption{karate}\end{subfigure}\ghs
      \begin{subfigure}{\fourfigwid}\includegraphics[width=\textwidth,height=\imgheiconv]{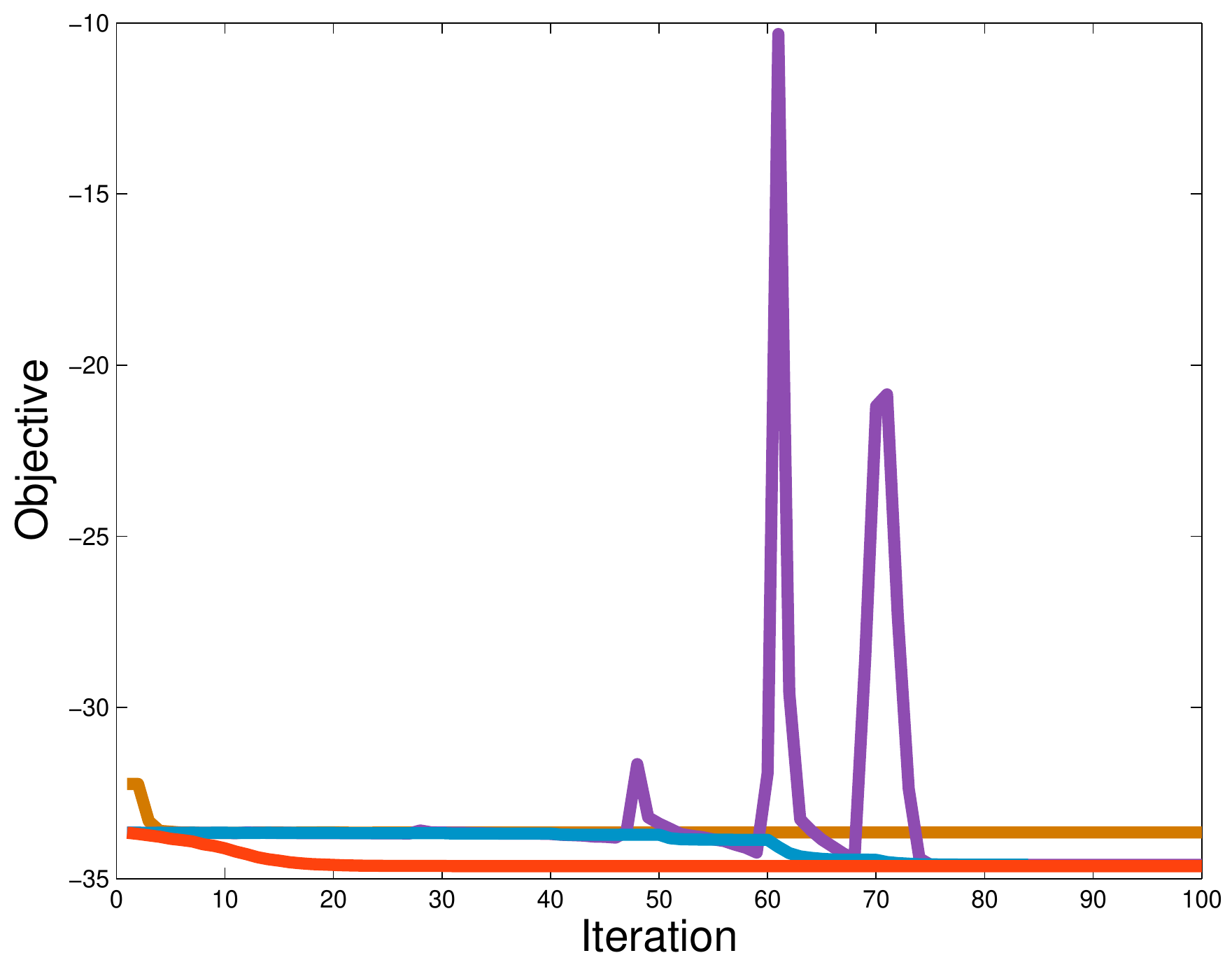}\caption{collab}\end{subfigure}\ghs
      \begin{subfigure}{\fourfigwid}\includegraphics[width=\textwidth,height=\imgheiconv]{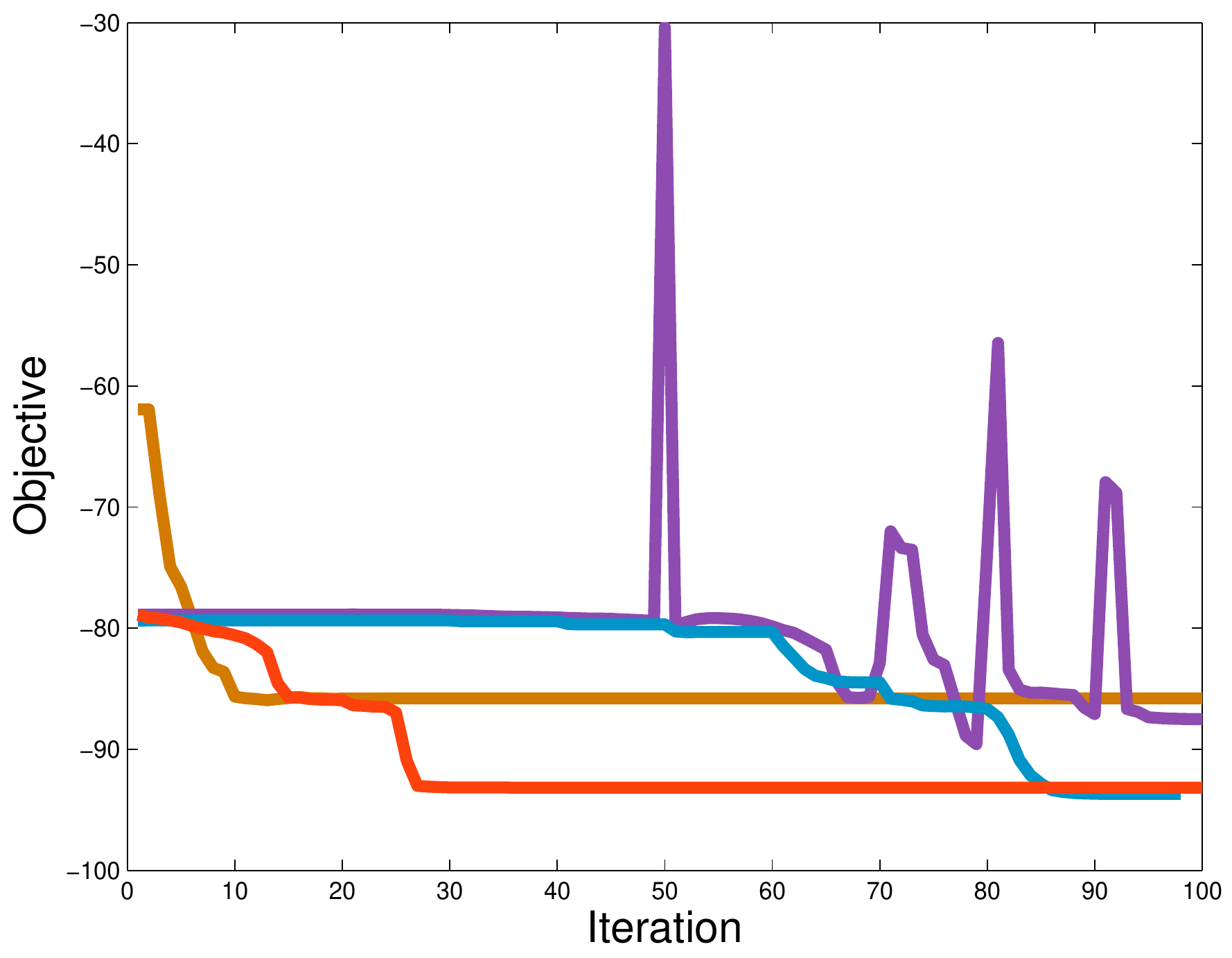}\caption{dolphins}\end{subfigure}\ghs
      \begin{subfigure}{\fourfigwid}\includegraphics[width=\textwidth,height=\imgheiconv]{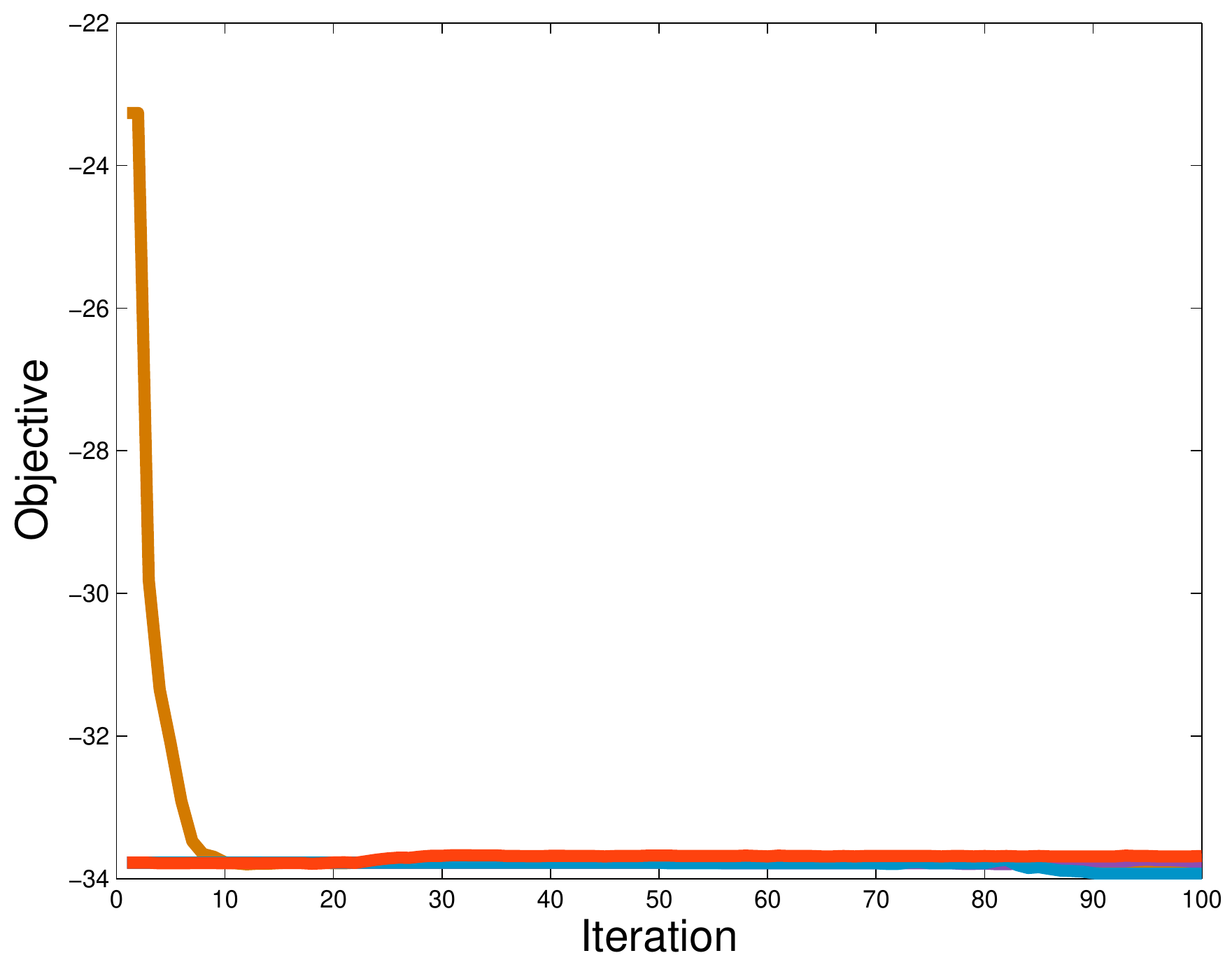}\caption{email}\end{subfigure}

      \begin{subfigure}{\fourfigwid}\includegraphics[width=\textwidth,height=\imgheiconv]{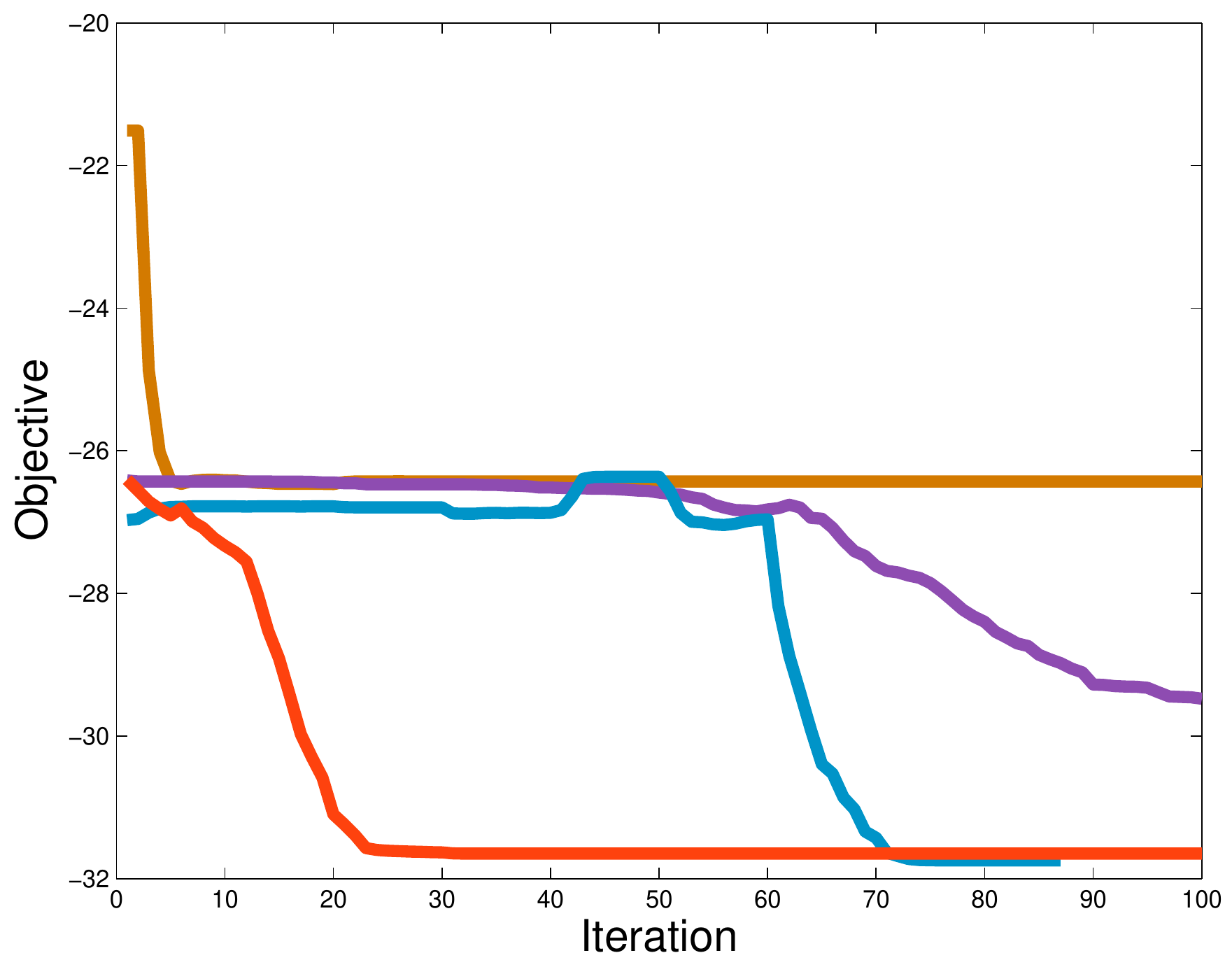}\caption{lesmis}\end{subfigure}\ghs
      \begin{subfigure}{\fourfigwid}\includegraphics[width=\textwidth,height=\imgheiconv]{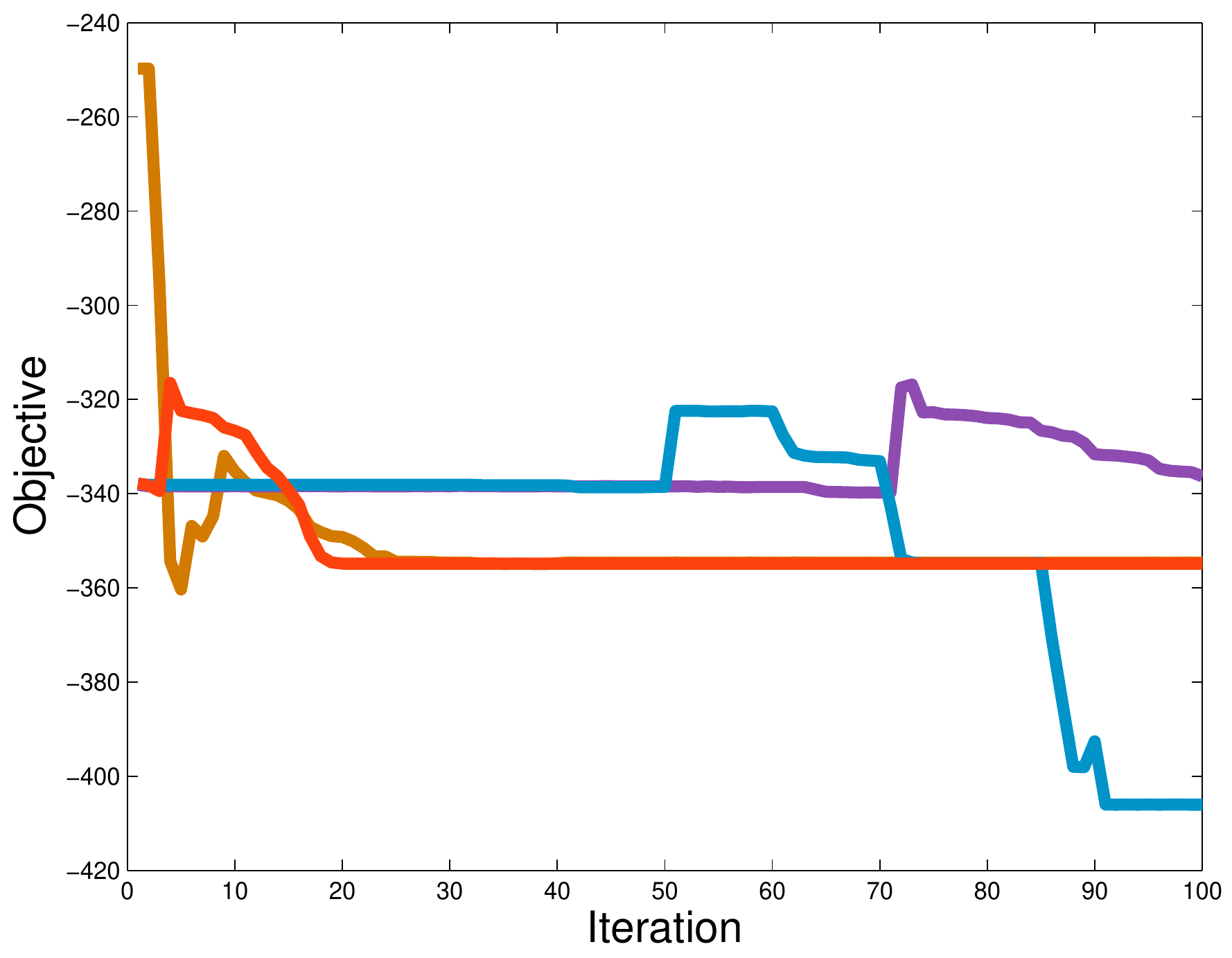}\caption{polbooks}\end{subfigure}\ghs
      \begin{subfigure}{\fourfigwid}\includegraphics[width=\textwidth,height=\imgheiconv]{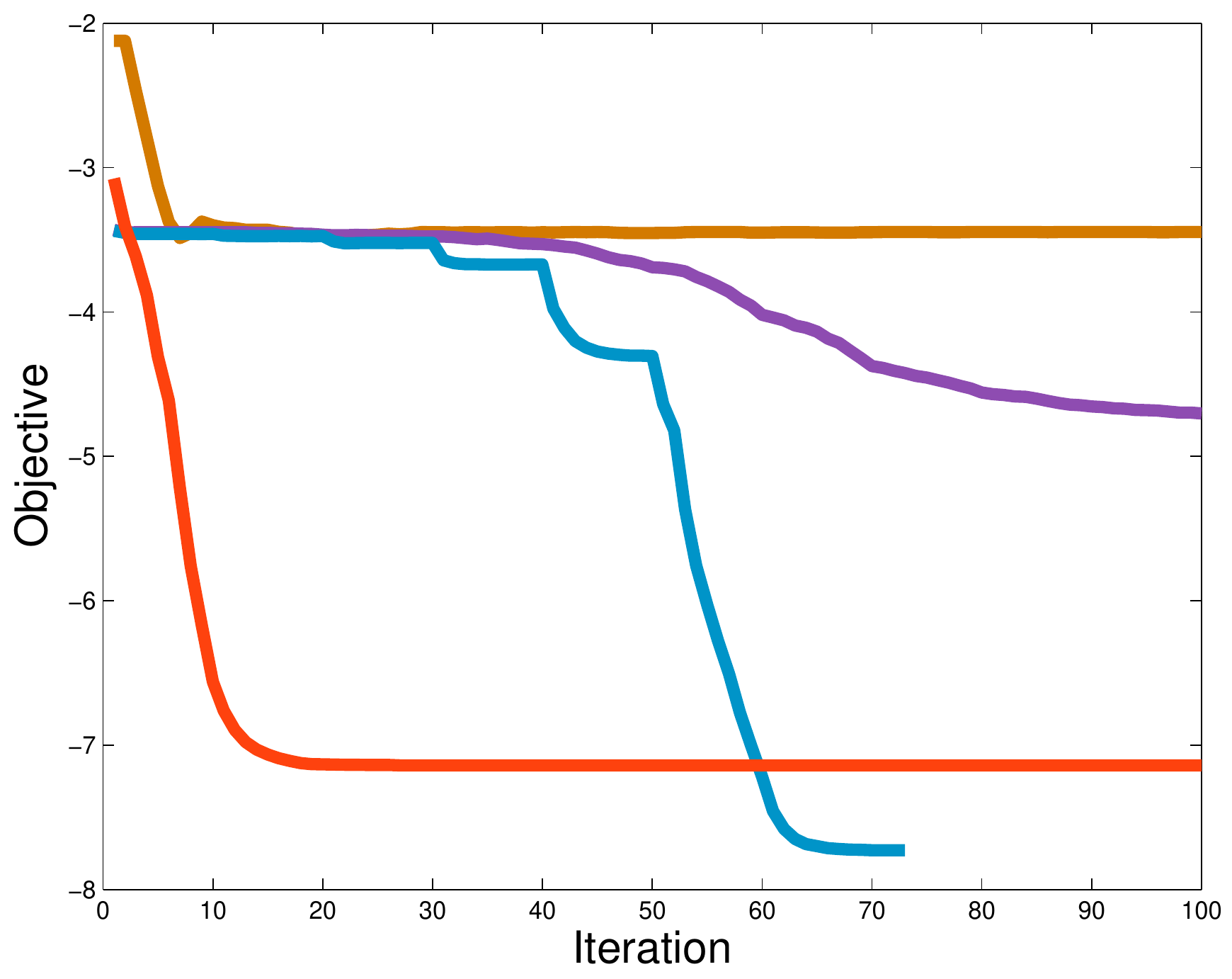}\caption{afootball}\end{subfigure}\ghs
      \begin{subfigure}{\fourfigwid}\includegraphics[width=\textwidth,height=\imgheiconv]{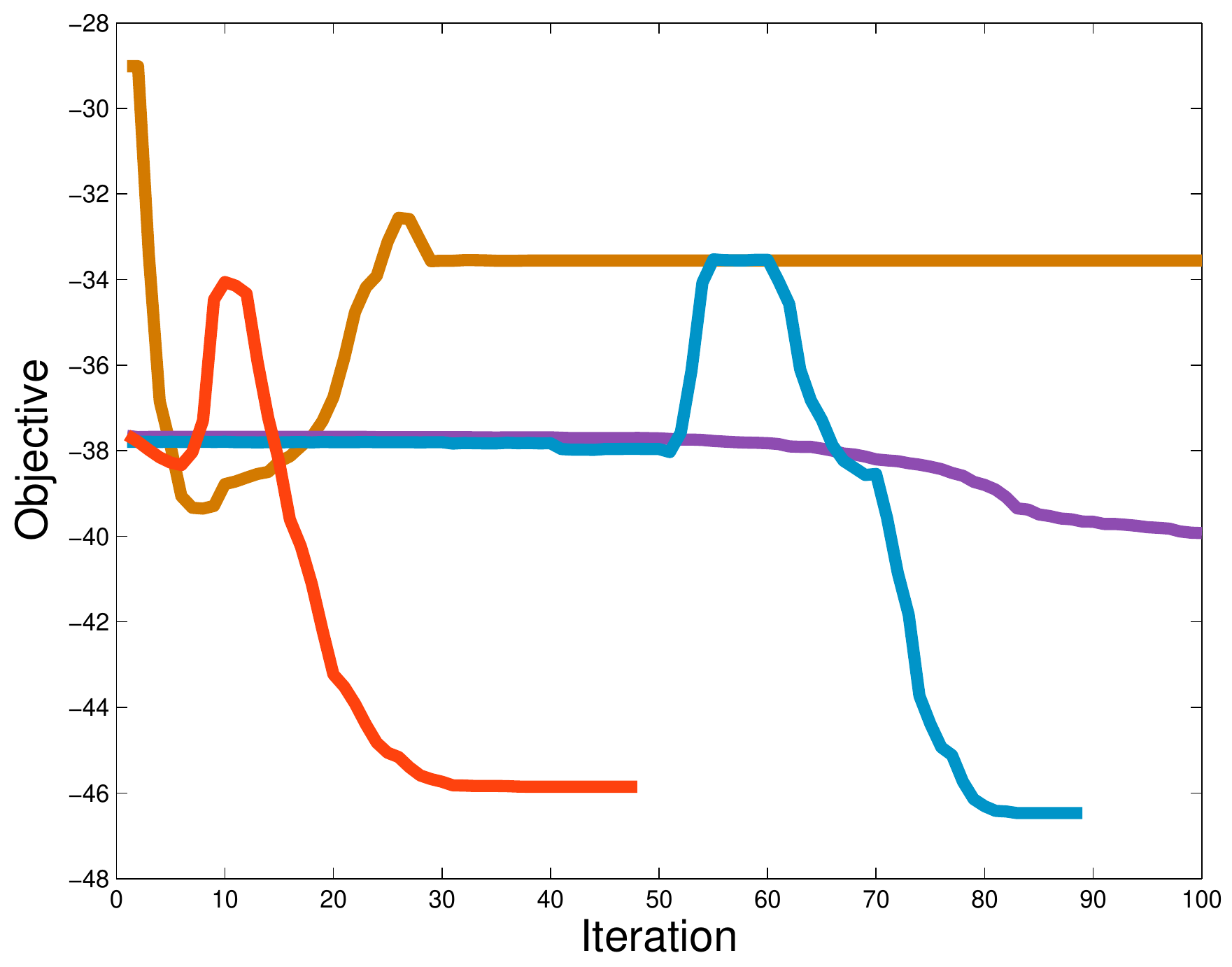}\caption{jazz}\end{subfigure}

\caption{Convergence curve for dense subgraph discovery on different datasets with $k=4000$.}
\label{fig:convergence:2}
\end{figure*}

\noi \bbb{Convergence Curve:} We demonstrate the convergence curve of the methods \{LP,~L2box-ADMM,~MPEC-EPM,~MPEC-ADM\} for dense subgraph discovery on different data sets. As can be seen in Figure \ref{fig:convergence:1} and Figure \ref{fig:convergence:2}, the proposed MPEC-based methods converges within 100 iterations. Moreover, we observe that the objective values generally decrease monotonically, and we attribute this to the greedy property of the penalty method for MPEC-EPM and monotone property of the dual variable $\rho$ update for MPEC-ADM.

\noi \bbb{Computational Efficiency:} We provide some running time comparisons for the methods \{LP,~L2box-ADMM,~MPEC-EPM,~MPEC-ADM\} on different data sets with different $k\in\{100, 1000, 2000, 3000, 4000, 5000\}$. As can be seen in Table \ref{tab:cpu}, even for the data set such as `dblp-2011' that contains about one million nodes and 7 million edges, all the methods can terminate with in 15 minutes. Moreover, the runtime efficiency of our methods are several times slower than LP and comparable with and L2box-ADMM. This is expected, since (i) our methods MPEC-EPM and MPEC-ADM need to call the LP procedure multiple times, and (ii) all the methods are all alternating methods and have the same computational complexity.

\subsection{Some Implementation Details }
This subsection presents some implementation details of MPEC-EPM and MPEC-ADM, including method of solving the $\bbb{x}$-subproblems and parameters setting of the algorithms.

The convex $\bbb{x}$-subproblems can be solved using Nesterov's projective gradient methods. For the applications of constrained image segmentation and Markov random fields, the projection step involved in the $\bbb{x}$-subproblems is easy since it only contains box constraints\footnote{It solves: $\min_{\bbb{0}\leq \bbb{x}\leq \bbb{1}}\|\bbb{a}-\bbb{x}\|_2^2$, where $\bbb{a}$ is given.}; for the applications of dense subgraph discovery, modularity clustering and graph bisection, the projection step can be hard since it contains an additional linear constraint besides box constraints (also known as capped simplex constraint)\footnote{It solves: $\min_{\bbb{0}\leq \bbb{x}\leq \bbb{1},~\bbb{x}^T\bbb{1}=k}\|\bbb{a}-\bbb{x}\|_2^2$, where $\bbb{a}$ and $k$ are given.}. Fortunately, this projection step can be solved by a break point search algorithm \cite{helgason1980polynomially} exactly in $n\log(n)$ time. In our experiments, we use the Matlab implementation provided in the appendix of \cite{yuan2016l0mpec} .

The following parameters are used in our algorithms. For all methods, we use proximal gradient descent algorithm to solve the inner $\bbb{x}$ subproblems. We stop the proximal gradient descent procedure when a relative change is smaller than $\epsilon=10^{-5}$, i.e. ${\|\bbb{x}^{k+1}-\bbb{x}^k\|}/{\|\bbb{x}^k\|}\leq \epsilon$, where $k$ is the iteration counter for the $\bbb{x}$-subproblem. In addition, we set $\rho^0=0.01,~T=10,~\sigma=\sqrt{10}$ for MPEC-EPM and set $\alpha^0=0.001,~T=10,~\sigma=\sqrt{10}$ for MPEC-ADM. Finally, for L2box-ADM, we update the penalty parameter by a factor of $\sqrt{10}$ in every $T=10$ iterations with initial value set to $0.1$ \footnote{We tune this parameter in the range \{0.01,~0.1,~1\} and find the value $0.1$ generally gives comparable results.}.

\section{Conclusions and Future Work}\label{sect:conc}

This paper presents a new class of continuous MPEC-based optimization methods to solve general binary optimization problems. Although the optimization problem is non-convex, we design two methods (exact penalty and alternating direction) to solve the equivalent problem. We also shed some theoretical lights to the equivalent formulations and optimization algorithms. Experimental results on binary problems demonstrate that our methods generally outperform existing solutions in terms of solution quality.

Our future work focuses on several directions. \textbf{(i)} We will investigate the optimality qualification of our multi-stage convex relaxation method for some specific objective functions, e.g., as is done in \cite{Goemans95,zhang2010analysis,CandesLS15,jain2013low}. \textbf{(ii)} We are also interested in extending the proposed algorithms to solve orthogonality and spherical optimization problems \cite{wen2013feasible,Chen2016} in computer vision and machine learning.


\section*{Acknowledgments}
\noi We would like to thank Prof. Shaohua Pan and Dr. Li Shen (South China University of Technology) for their valuable discussions on the MPEC techniques. Research reported in this publication was supported by competitive research funding from King Abdullah University of Science and Technology (KAUST).


\bibliographystyle{ieee}
\bibliography{my}
\setcounter{lemma}{0}
\setcounter{theorem}{0}

\begin{IEEEbiography}[{\includegraphics[width=1in,height=1.25in]{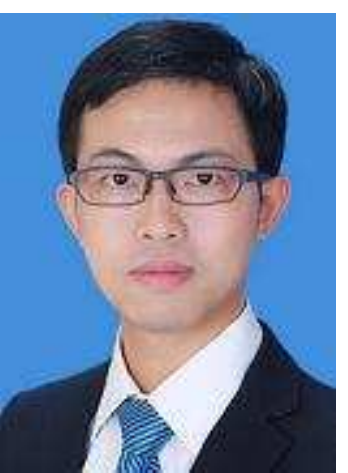}}]{Ganzhao Yuan} was born in Guangdong, China. He received his Ph.D. in School of Computer Science and Engineering, South China University of Technology (SCUT) in 2013. He is currently a research associate professor at School of Data and Computer Science in Sun Yat-sen University (SYSU). His research interests primarily center around large-scale nonlinear optimization and its applications in computer vision and machine learning. He has published papers in ICML, SIGKDD, AAAI, CVPR, VLDB, and ACM Transactions on Database System (TODS).
\end{IEEEbiography}

\begin{IEEEbiography}[{\includegraphics[width=1in,height=1.25in]{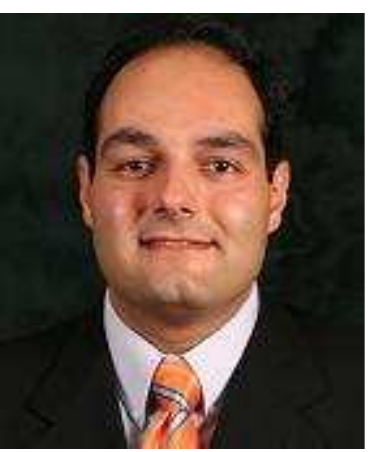}}]{Bernard Ghanem} was born in Betroumine, Lebanon. He received his Ph.D. in Electrical and Computer Engineering from the University of Illinois at Urbana-Champaign (UIUC) in 2010. He is currently an assistant professor at King Abdullah University of Science and Technology (KAUST), where he leads the Image and Video Understanding Lab (IVUL). His research interests focus on designing, implementing, and analyzing approaches to address computer vision problems (e.g. object tracking and action recognition/detection in video), especially at large-scale.
\end{IEEEbiography}


\clearpage



\section*{Appendix A: Proofs of Equivalent Reformulations}\label{app:equ}


This section presents the proofs of the other equivalent separable/non-separable reformulations for binary constraint which are claimed in Table 1.

\begin{lemma} \label{lemma:binary:mpec2}
\bbb{$\ell_{\infty}$ box non-separable MPEC}. We define
\beq
\Pi \triangleq \{ (\bbb{x},\bbb{v})~|~\bbb{x}^T\bbb{v}=n,~\|\bbb{x}\|_{\infty}\leq 1,~\|\bbb{v}\|_{\infty}\leq 1 \}. \label{eq:Pi}\nn
\eeq
\noi Assume that $(\bbb{x},\bbb{v}) \in \Pi$, then we have $\bbb{x}\in \{-1,+1\}^n$, $\bbb{v}\in \{-1,+1\}^n$, and $\bbb{x}=\bbb{v}$.

\begin{proof}
The proof of this lemma is similar to the proof of $\ell_2$ box non-separable MPEC. Firstly, we prove that $\bbb{x} \in \{-1,+1\}^n$. We have: $n = \bbb{x}^T\bbb{v} \leq  \|\bbb{x}\|_{1} \cdot \|\bbb{v}\|_{\infty} \leq \|\bbb{x}\|_{1}$. Thus, we obtain $\|\bbb{x}\|_1\geq n$. We define $\bbb{z}=|\bbb{x}|$. Combining $\|\bbb{x}\|_{\infty}\leq 1$, we have the following constraint sets for $\bbb{z}$: $\sum_{i}~ \bbb{z}_i \geq n,~\bbb{0}\leq \bbb{z}\leq \bbb{1}$. Therefore, we have $\bbb{z}=\bbb{1}$ and it holds that $\bbb{x}\in\{-1,+1\}^n$. Secondly, using the same methodology, we can prove that $\bbb{v}\in\{-1,+1\}^n$. Finally, we have $\bbb{x}=\bbb{v}$ since $\bbb{x}\in\{-1,+1\}^n$,~$\bbb{v}\in\{-1,+1\}^n$ and $\la \bbb{x},\bbb{v}\ra=n$.
\end{proof}
\end{lemma}

\begin{lemma}
\bbb{$\ell_{\infty}$ box separable MPEC.} We define
\beq \label{eq:Psi}
\Psi \triangleq \{ (\bbb{x},\bbb{v})~|~\bbb{x} \odot \bbb{v}=\bbb{1},~\|\bbb{x}\|_{\infty}\leq 1,~\|\bbb{v}\|_{\infty}\leq 1\}.\nn
\eeq
\noi Assume that $(\bbb{x},\bbb{v}) \in \Psi$, then we have $\bbb{x}\in \{-1,+1\}^n$, $\bbb{v}\in \{-1,+1\}^n$, and $\bbb{x}=\bbb{v}$.
\begin{proof}
We observe that all the constraints in $\Psi$ can be decomposed into $n$ independent components. We now focus on $i$th component. (i) Assuming that $\bbb{x}_i\geq0$, we have $0\leq \bbb{x}_i\leq 1$ and $0\leq {1}/{\bbb{x}_i}\leq 1$, we have $\bbb{x}_i=1$ and $\bbb{v}_i = {1}/{\bbb{x}_i}=1$. (ii) Assuming that $\bbb{x}_i\leq0$, we have $-1 \leq \bbb{x}_i\leq 0$ and $-1\leq {1}/{\bbb{x}_i}\leq 0$, we have $\bbb{x}_i=-1$ and $\bbb{v}_i = {1}/{\bbb{x}_i}=-1$. This finishes the proof.

\end{proof}
\end{lemma}

\begin{lemma}
\bbb{$\ell_2$ box separable MPEC.} We define
\beq\label{eq:Upsilon}
\Upsilon \triangleq \{ (\bbb{x},\bbb{v})~|~ \bbb{x} \odot \bbb{v}=1,~\|\bbb{x}\|_{\infty}\leq 1,~\|\bbb{v}\|_2^2\leq n \} \nn
\eeq
\noi Assume that $(\bbb{x},\bbb{v}) \in \Upsilon$, then we have $\bbb{x}\in \{-1,+1\}^n$, $\bbb{v}\in \{-1,+1\}^n$, and $\bbb{x}=\bbb{v}$.
\begin{proof}
We notice that $\bbb{v}=\frac{1}{\bbb{x}}$. We define $\bbb{z}=|\bbb{v}|$. Combining $-\bbb{1}\leq \bbb{x} \leq \bbb{1}$, we obtain $\bbb{z}\geq \bbb{1}$. We have the following constraint sets for $\bbb{z}$: $\bbb{z}\geq \bbb{1},~\bbb{z}^T\bbb{z}\leq n$. Therefore, we have $\bbb{z}=\bbb{1}$. Finally, we achieve $\bbb{v} \in \{-1,+1\}^n$ and $\bbb{x}=\frac{1}{\bbb{v}} = \bbb{v}$.
\end{proof}
\end{lemma}

\begin{lemma}
\bbb{$\ell_2$ box non-separable reformulation.} The following equivalence holds:
$$\{-1+1\}^n \Leftrightarrow \{\bbb{x}~|~-\bbb{1}\leq \bbb{x}\leq \bbb{1},~\|\bbb{x}\|_2^2 = n\}$$
\begin{proof}
First, it holds that: $n=\bbb{x}^T\bbb{x}\leq \|\bbb{x}\|_1\|\bbb{x}\|_{\infty}\leq \|\bbb{x}\|_1$. Therefore, we have $-\bbb{1}\leq \bbb{x}\leq \bbb{1},~\|\bbb{x}\|_1 \geq n$. Note that these constraint sets are symmetric. Letting $\bbb{z}=|\bbb{x}|$, we obtain:~$\bbb{0}\leq \bbb{z} \leq \bbb{1},~\sum_{i}\bbb{z}_i\geq n$. Thus, we have $\bbb{z}=\bbb{1}$ and it holds that $\bbb{x}\in\{-1,+1\}^n$.
\end{proof}
\end{lemma}


%
%
%
%
%


\section*{Appendix B: Solving the Rank-One Subproblem}\label{app:sub}
This subsection describes a nearly closed-form solution for solving the rank-one subproblem which is involved in our alternating direction method. For general purpose, we consider the following optimization problem:
\beq \label{eq:iden}
\min_{\|\bbb{x}\|_2 \leq \beta}~\frac{1}{2}\bbb{x}^T(\gamma \bbb{I}+\bbb{b}\bbb{b}^T)\bbb{x} + \la \bbb{x},\bbb{c} \ra
\eeq
\noi where $\gamma, \beta \in \mathbb{R},~\bbb{b},\bbb{c} \in \mathbb{R}^n$ are given. We assume $\beta\neq0$ or $\gamma \neq0$. Clearly, (\ref{eq:iden}) is equivalent to the following minimax optimization problem:
\beq \label{eq:iden:minimax}
\min_{\bbb{x}}\max_{\theta\geq 0}~\tfrac{1}{2}\bbb{x}^T(\gamma \bbb{I}+\bbb{b}\bbb{b}^T)\bbb{x} +  \la \bbb{x},\bbb{c} \ra +\tfrac{\theta}{2}(\|\bbb{x}\|^2_2 - \beta^2)
\eeq
\noi Setting the gradient respect of $\bbb{x}$ to zero, we obtain:
\beq \label{eq:recover}
\bbb{x} = -(\gamma \bbb{I}+\theta \bbb{I} +\bbb{b}\bbb{b}^T)^{-1}\bbb{c}
\eeq
\noi Putting this equality into (\ref{eq:iden:minimax}), we have the following minimization with respect to $\theta$:
\beq \label{eq:theta1}
 \max_{\theta\geq 0} -\frac{1}{2}\bbb{c}^T (\gamma \bbb{I}+\theta \bbb{I}+\bbb{b}\bbb{b}^T)^{-1}\bbb{c} - \frac{1}{2}\theta \beta^2
\eeq
\noi Using the well-known Sherman-Morrison inverse formula \footnote{$(\eta \bbb{I}+\bbb{b}\bbb{b}^T)^{-1}  = \frac{1}{\eta}\bbb{I}  - \frac{1}{\eta^2+\eta \bbb{b}^T \bbb{b}} \bbb{b}\bbb{b}^T,~\forall \eta>0$}, (\ref{eq:theta1}) reduces to the following optimization problem:
\beq
\max_{\theta\geq 0} ~\frac{1}{2} \left( \frac{t^2}{ (\gamma+\theta)(\gamma+\theta+r)} -\frac{s}{\gamma+\theta}  - s\theta \beta^2 \right)\nn
\eeq
\noi where $r=\bbb{b}^T\bbb{b},~s = \bbb{c}^T\bbb{c},~t = \bbb{c}^T\bbb{b}$. The optimal solution $\theta^*$ can be found using a simple one-dimensional bisection line search procedure. After that, the optimal solution $\bbb{x}^*$ for the original optimization problem in (\ref{eq:iden}) can be recovered using (\ref{eq:recover}).

\end{document}